\documentclass[11pt]{article}
\usepackage{tikz}
\usetikzlibrary{automata}
\usetikzlibrary{shapes.geometric}
\usetikzlibrary{calc}
\usetikzlibrary{decorations.pathmorphing}
\usetikzlibrary{decorations.markings}
\usepackage{mathtools}
\usepackage{graphicx}
\usepackage{amsmath,amsfonts,latexsym,amscd,amssymb,theorem}
\usepackage{graphicx}
\usepackage[ansinew]{inputenc}
\usepackage{epsfig}
\usepackage{amsmath}
\usepackage{amssymb}
\usepackage{afterpage}
\usepackage{selinput}
\def\a{\alpha}

\newcommand{\ba}{\begin{eqnarray}}
\newcommand{\ea}{\end{eqnarray}}

\newtheorem{thm}{Theorem}[section]

\newtheorem{theorem}[thm]{Theorem}

\newtheorem{definition}[thm]{Definition}
\newtheorem{lemma}[thm]{Lemma}

\newtheorem{corollary}[thm]{Corollary}
\newtheorem{remark}[thm]{Remark}

\makeatletter
\newcommand*{\rom}[1]{\expandafter\@slowromancap\romannumeral #1@}
\makeatother

\begin{document}
\author{Ayman El Zein$^{1,2}$}
\footnotetext[1]{Lebanese University, KALMA Laboratory, Beirut, Lebanon}
\footnotetext[2]{University Savoie Mont Blanc, LAMA Laboratory, Chamb\'{e}ry, France}
\title{Oriented Hamiltonian Cycles in Tournaments: a Proof of Rosenfeld's Conjecture}
\maketitle 
\begin{abstract}
Rosenfeld in 1974 conjectured that there is an integer $N>8$ such that every tournament of order $n\geq N$ contains every non-directed cycle of order $n$. We prove that, with exactly $35$ exceptions, every tournament of order $n\geq 3$ contains each non-directed cycle of order $m$, $3\leq m\leq n$.
\end{abstract}
\section{Introduction}
It was early proved by Camion \cite{C} that a tournament has a Hamiltonian directed cycle if and only if it is strong. Concerning non-directed cycles, Thomassen \cite{Thomassen}, in $1973$, proved that for any $n\geq 50$, any tournament of order $n$ contains a Hamiltonian antidirected cycle. Rosenfeld \cite{R'}, in $1974$, after improving Thomassen's result to $n\geq 28$, conjectured that there is an integer $N\geq 9$ such that any tournament of order $n\geq N$ contains any non-directed Hamiltonian cycle. This conjecture was studied in different approach and many interesting result were successively established. At the beginning, Gr\"{u}nbaum \cite{G'} proved the existence of cycles with block of length $n-1$, then, in $1983$, the existence of cycles of two blocks was proved by Benhocine and Wojda \cite{BW}. In the same year, Petrovi\'{c} \cite{P} improved Thomassen's and Rosenfeld's bound by proving the same result for $n\geq 16$. The first general proof of Rosenfeld's conjecture was established by Thomason \cite{T} for tournaments of order $n\geq 2^{128}$. Thomason was contented only to indicate that it should be true for tournaments of order at least $9$. Then, Havet \cite{H} improved this bound by proving the conjecture for tournaments of order $n\geq 68$ and for reducible tournaments of order $n\geq 9$. In this paper, we completely settle Rosenfeld's cycles conjecture: any tournament contains each Hamiltonian non-directed cycle with $30$ exceptions, all of order less than $9$ , that will be described explicitly together with the non-allowed cycles. Note that, Rosenfeld proposed the same conjecture for paths, which was proved by Havet and Thomass\'{e} \cite{HT}. The existence of Hamiltonian paths provides simply the existence of any non-Hamiltonian path, but this is not applied for cycles. Our argument can be used to extend Rosenfeld's conjecture by showing that any tournament of order $n$ contains any non-directed cycle of order $m<n$ with exactly $5$ exceptions, where $m=n-1$.\\
A \textit{tournament} is an orientation of a complete graph. If $X$ is a subset of vertices of $T$, the subtournament induced by $X$ is denoted by $T(X)$. For abbreviation, we write $x\rightarrow y$ whenever $(x,y)\in E(T)$. Similarly, if $X$ and $Y$ are two subtournaments of $T$, we write $X\rightarrow Y$ if $x\rightarrow y$ for every pair $(x,y)\in X\times Y$. A vertex $x$ is said to be minimal in $T$ if $d^-(x)=0$.\\
 Let $P=x_1...x_n$ be a path, $x_1$ is called the \textit{origin} of $P$ while $x_n$ is called its \textit{end}. $^*P$ (resp. $P^*$) denotes the path $x_2...x_n$ (resp. $x_1...x_{n-1}$), the path $x_n...x_1$ is denoted by $P^{-1}$. $P$ is said to be an \textit{outpath}(resp. \textit{inpath}) if $x_1\rightarrow x_2$ (resp. $x_1\leftarrow x_2$). $P$ is said to be a \textit{directed outpath} (resp. \textit{directed inpath}) if $x_i\rightarrow x_{i+1}$ (resp. $x_i\leftarrow x_{i+1}$) for all $i$, $1\leq i\leq n-1$. A \textit{block} of $P$ is a maximal directed subpath of $P$. In general, an oriented path is formed by successive blocks. The $i^{th}$ block of $P$ is denoted by $B_i(P)$ while its length is denoted by $b_i(P)$. The \textit{type} of an outpath $P$ is the sequence $+(b_1(P),b_2(P),...,b_s(P))$ while that of an inpath $Q$ is $-(b_1(Q),b_2(Q),...,b_t(Q))$. Occasionally, we write $P=\pm(b_1(P),...,b_s(P))$ instead of $P$, if $P$ is of type $\pm(b_1(P),...,b_s(P))$. When two paths $P$ and $Q$ are of same type, we write $P\equiv Q$. Given a path $P$ and a tournament $T$, we define the set $Or(P,T)=\{v\in T;v$ is an origin of a path $Q$ in $T$ such that $Q\equiv P\}$.\\
 Let $C=x_0...x_{n-1}x_0$ be a cycle. In a similar way, we define the blocks $B_i(C)$ of $C$ together with their lengths $b_i(C)$. The type of $C$ is the sequence $(b_1(C),...,b_s(C))$, we write $C=(b_1(C),...,b_s(C))$. When two cycles $C$ and $C'$ are of same type, we write $C\equiv C'$. A rotation of a cycle $C=a_0a_1...a_{n-1}a_0$ is an enumeration $a'_0a'_1...a'_{n-1}a'_0$ of the vertices of $C$ such that $a'_0=a_{n-1}$ and $a'_1=a_0$, after a rotation we re-write $C=a_0a_1...a_{n-1}a_0$, where $a_i=a'_i$ for all $i\in \{0,...,n-1\}$.\\
 The dual of a digraph $D$ is the digraph $\overline{D}$ such that $(y,x)\in E(\overline{D})$ whenever $(x,y)\in E(D)$.\\
 A tournament $T$ is said to be strong if $T$ contains an $xy-$directed outpath for every pair $(x,y)\in T$. Otherwise, $T$ is said to be \textit{reducible}.\\
 Let $X$ be a set of vertices of $T$. The \textit{outsection} generated by $X$, denoted by $S^+(X)$, is the set of vertices $y$ such that $T$ contains an $xy-$directed outpath for some $x\in X$. $s^+(X)$ denotes the cardinal of $S^+(X)$. Likewise, the \textit{insection} generated by $X$, denoted by $S^-(X)$, is the set of vertices $y$ such that $T$ contains an $xy-$directed inpath for some $x\in X$. $s^-(X)$ denotes the cardinal of $S^-(X)$. If $x_1,...x_s$ are $s$ vertices in $T$, we write $S^{\pm}(x_1,...,x_s)$ instead of $S^{\pm}(\{x_1,...,x_s\})$ and $s^{\pm}(x_1,...,x_s)$ instead of $s^{\pm}(\{x_1,...,x_s\})$. A vertex $x$ of $T$ is said to be an \textit{outgenerator} (resp. \textit{ingenerator}) of $T$ if $S^+(x)=V(T)$ (resp. $S^-(x)=V(T)$).\\
 In this paper, we say that a path $P$ (resp. a cycle $C$) is contained in a tournament $T$ if there exists a path $Q$ (a cycle $C'$) lying in $T$ such that $Q\equiv P$ (resp. $C'\equiv C$). If $X$ is set of vertices of $T$, $O_X$ (resp. $I_X$) designate an arbitrary Hamiltonian directed outpath (resp. inpath)of $T(X)$. Let $x$ be a vertex in $T$ and $Y$ a subtournament of $T$, the set $N^+_Y(x)$ (resp. $N^-_Y(x)$) denotes the set of outneighbors (resp. inneighbors) of $x$ in $Y$. The cardinal of $N^+_Y(x)$ (resp. $N^-_Y(x)$) is denoted by $d^+_Y(x)$ (resp. $d^-_Y(x)$).
\section{Preliminary study}
In order to prove Rosenfeld's paths conjecture, Havet and Thomass\'{e} \cite{HT} introduced the concept of exceptions as follows:\\
An exception is a pair $(T;P)$, where $T$ is a tournament illustrated in figures \ref{figure 1}, \ref{figure 2}, \ref{figure 3} and \ref{figure 4} and $P$ is an outpath with $|T|=|P|$ such that there exists two vertices $x$ and $y$ of $T$ neither of which is an origin of $P$ but $s^+(x,y)\geq b_1(P)+1$. The dual of an exception is also an exception.\\
\begin{figure}
\begin{tikzpicture}
\tikzset{enclosed/.style={draw,circle,inner sep=2pt,minimum size=4pt,fill=black}}
\tikzset{->-/.style={decoration={
            markings,
            mark=at position #1 with
            {\arrow{>}}},postaction={decorate}}}
\node[enclosed,label={left,yshift=.2cm:1}](1)at(0,0){};
\node[enclosed,label={left,yshift=.2cm:2}](2)at(1.5,2.6){};
\node[enclosed,label={right,yshift=.2cm:3}](3)at(3,0){};
\draw[black,->-=.5] (1)--(2);
\draw[black,->-=.5] (2)--(3);
\draw[black,->-=.5] (3)--(1);

\node[minimum size=.1pt,label={left:3A}](3A)at(2.15,-.5){};

\node[ellipse,minimum width=4cm,minimum height=1.3cm,draw](4Aa)at(7.6,3){};
\node[minimum size=.1pt](4Aa1)at(6.7,2.52){};
\node[minimum size=.1pt](4Aa2)at(8.5,2.52){};
\node[enclosed,label={left,yshift=.2cm:1}](4A1)at(6.4,3){};
\node[enclosed,label={right,yshift=.2cm:2}](4A2)at(8.8,3){};
\draw[black,->-=.5] (4A1)--(4A2);
\node[enclosed,label={right,yshift=.2cm:3}](4A3)at(9.1,0){};
\node[enclosed,label={left,yshift=.2cm:4}](4A4)at(6.1,0){};
\draw[black,->-=.5] (4A3)--(4A4);
\draw[black,->-=.5] (4A4)--(4Aa1);
\draw[black,->-=.5] (4Aa2)--(4A3);

\node[minimum size=.1pt,label={left:4A}](4A)at(8.25,-.5){};

\node[enclosed,label={left,yshift=.2cm:1}](5A1)at(11.58,2.15){};
\node[enclosed,label={right,yshift=.2cm:2}](5A2)at(13.5,3.8){};
\node[enclosed,label={right,yshift=.2cm:3}](5A3)at(15.42,2.15){};
\node[enclosed,label={right,yshift=.2cm:4}](5A4)at(14.7,0){};\node[enclosed,label={left,yshift=.2cm:5}](5A5)at(12.3,0){};

\draw[black,->-=.5] (5A1)--(5A2);
\draw[black,->-=.5] (5A1)--(5A3);
\draw[black,->-=.5] (5A2)--(5A3);
\draw[black,->-=.5] (5A2)--(5A4);
\draw[black,->-=.5] (5A3)--(5A4);
\draw[black,->-=.5] (5A3)--(5A5);
\draw[black,->-=.5] (5A4)--(5A5);
\draw[black,->-=.5] (5A4)--(5A1);
\draw[black,->-=.5] (5A5)--(5A1);
\draw[black,->-=.5] (5A5)--(5A2);

\node[minimum size=.1pt,label={left:5A}](5A)at(14.15,-.5){};

\node[ellipse,minimum width=4cm,minimum height=4cm,draw](5Ba)at(1.5,-3.2){};
\node[minimum size=.1pt](5Ba1)at(0.6,-4.8){};
\node[minimum size=.1pt](5Ba2)at(2.4,-4.8){};
\node[enclosed,label={left,yshift=.2cm:1}](5B1)at(0.3,-4){};
\node[enclosed,label={right,yshift=.2cm:2}](5B2)at(1.5,-1.93){};\node[enclosed,label={right,yshift=.2cm:3}](5B3)at(2.7,-4){};
\draw[black,->-=.5] (5B1)--(5B2);
\draw[black,->-=.5] (5B2)--(5B3);
\draw[black,->-=.5] (5B1)--(5B3);
\node[enclosed,label={right,yshift=.2cm:4}](5B4)at(3,-7){};
\node[enclosed,label={left,yshift=.2cm:5}](5B5)at(0,-7){};
\draw[black,->-=.5] (5B4)--(5B5);
\draw[black,->-=.5] (5B5)--(5Ba1);
\draw[black,->-=.5] (5Ba2)--(5B4);

\node[minimum size=.1pt,label={left:5B}](5B)at(2.15,-7.5){};

\node[ellipse,minimum width=4cm,minimum height=4cm,draw](5Ca)at(7.6,-3.2){};
\node[minimum size=.1pt](5Ca1)at(6.7,-4.8){};
\node[minimum size=.1pt](5Ca2)at(8.5,-4.8){};
\node[enclosed,label={left,yshift=.2cm:1}](5C1)at(6.4,-4){};
\node[enclosed,label={right,yshift=.2cm:2}](5C2)at(7.6,-1.93){};\node[enclosed,label={right,yshift=.2cm:3}](5C3)at(8.8,-4){};
\draw[black,->-=.5] (5C1)--(5C2);
\draw[black,->-=.5] (5C2)--(5C3);
\draw[black,->-=.5] (5C3)--(5C1);
\node[enclosed,label={right,yshift=.2cm:4}](5C4)at(9.1,-7){};
\node[enclosed,label={left,yshift=.2cm:5}](5C5)at(6.1,-7){};
\draw[black,->-=.5] (5C4)--(5C5);
\draw[black,->-=.5] (5C5)--(5Ca1);
\draw[black,->-=.5] (5Ca2)--(5C4);

\node[minimum size=.1pt,label={left:5C}](5C)at(8.25,-7.5){};

\node[ellipse,minimum width=4cm,minimum height=1.3cm,draw](5D)at(13.5,-4){};
\node[minimum size=.1pt](5Da1)at(12.6,-4.48){};
\node[minimum size=.1pt](5Da2)at(14.4,-4.48){};
\node[enclosed,label={left,yshift=.2cm:1}](5D1)at(12.3,-4){};
\node[enclosed,label={right,yshift=.2cm:2}](5D2)at(14.7,-4){};
\draw[black,->-=.5] (5D1)--(5D2);
\node[enclosed,label={right,yshift=.2cm:3}](5D3)at(15,-7){};
\node[enclosed,label={left,yshift=.2cm:4}](5D4)at(12,-7){};
\draw[black,->-=.5] (5D3)--(5D4);
\draw[black,->-=.5] (5D4)--(5Da1);
\draw[black,->-=.5] (5Da2)--(5D3);
\node[enclosed,label={left,yshift=.01cm:5}](5D5)at(13.5,-5.5){};
\draw[black,->-=.5] (5D1)--(5D5);
\draw[black,->-=.5] (5D3)--(5D5);
\draw[black,->-=.5] (5D5)--(5D2);
\draw[black,->-=.5] (5D5)--(5D4);

\node[minimum size=.1pt,label={left:5D}](5D)at(14.15,-7.5){};

\node[ellipse,minimum width=4cm,minimum height=1.3cm,draw](5Ea)at(2.4,-9.5){};
\node[minimum size=.1pt](5Ea1)at(1.5,-9.98){};
\node[minimum size=.1pt](5Ea2)at(3.3,-9.98){};
\node[enclosed,label={left,yshift=.2cm:1}](5E1)at(1.2,-9.5){};
\node[enclosed,label={right,yshift=.2cm:2}](5E2)at(3.6,-9.5){};
\draw[black,->-=.5] (5E1)--(5E2);
\node[enclosed,label={right,yshift=.2cm:3}](5E3)at(4.4,-12.5){};
\node[ellipse,minimum width=4cm,minimum height=1.3cm,draw](5Eb)at(1.25,-12.5){};
\node[enclosed,label={left,yshift=.2cm:4}](5E4)at(0.05,-12.5){};
\node[enclosed,label={right,yshift=.2cm:5}](5E5)at(2.45,-12.5){};
\draw[black,->-=.5] (5E3)--(5Eb);
\draw[black,->-=.5] (5Ea)--(5E3);
\draw[black,->-=.5] (5E4)--(5E5);
\draw[black,->-=.5] (5Eb)--(5Ea);

\node[minimum size=.1pt,label={left:5E}](5E)at(3.1,-13.65){};

\node[ellipse,minimum width=5cm,minimum height=3.2cm,draw](6Aa)at(8.1,-9.95){};
\node[minimum size=.1pt](6Aa1)at(6.5,-11.07){};
\node[minimum size=.1pt](6Aa2)at(9.7,-11.07){};
\node[ellipse,minimum width=2.2cm,minimum height=1.3cm,draw](6Ab)at(9.1,-10.5){};
\node[enclosed,label={left,yshift=.2cm:1}](6A1)at(6.9,-10.5){};
\node[enclosed,label={right,yshift=.2cm:2}](6A2)at(8.1,-9){};
\node[enclosed,label={left,yshift=.01cm:3}](6A3)at(8.55,-10.5){};
\node[enclosed,label={right,yshift=.01cm:4}](6A4)at(9.65,-10.5){};\node[enclosed,label={right,yshift=.2cm:5}](6A5)at(9.7,-12.5){};\node[enclosed,label={left,yshift=.2cm:6}](6A6)at(6.5,-12.5){};

\draw[black,->-=.5] (6A1)--(6A2);
\draw[black,->-=.5] (6A2)--(6Ab);
\draw[black,->-=.5] (6Ab)--(6A1);
\draw[black,->-=.5] (6A4)--(6A3);
\draw[black,->-=.5] (6Aa2)--(6A5);
\draw[black,->-=.5] (6A5)--(6A6);
\draw[black,->-=.5] (6A6)--(6Aa1);

\node[minimum size=.1pt,label={left:6A}](6A)at(8.75,-13.65){};

\node[enclosed,label={left,yshift=.2cm:6}](6B6)at(12.58,-10.35){};

\node[ellipse,minimum width=2.5cm,minimum height=1.2cm,draw](6Ba)at(14.5,-8.8){};

\node[enclosed,label={right,yshift=.3cm:1}](6B1)at(13.7,-8.8){};
\node[enclosed,label={left,yshift=.3cm:2}](6B2)at(15.3,-8.8){};
\node[enclosed,label={right,yshift=.2cm:3}](6B3)at(16.42,-10.35){};
\node[enclosed,label={right,yshift=.2cm:4}](6B4)at(15.7,-12.5){};\node[enclosed,label={left,yshift=.2cm:5}](6B5)at(13.3,-12.5){};
\draw[black,->-=.5] (6B1)--(6B2);
\draw[black,->-=.5] (6B6)--(6Ba);
\draw[black,->-=.5] (6B6)--(6B3);
\draw[black,->-=.5] (6Ba)--(6B5);
\draw[black,->-=.5] (6Ba)--(6B4);
\draw[black,->-=.5] (6B3)--(6B4);
\draw[black,->-=.5] (6B3)--(6Ba);
\draw[black,->-=.5] (6B5)--(6B4);
\draw[black,->-=.5] (6B4)--(6B6);
\draw[black,->-=.5] (6B5)--(6B6);
\draw[black,->-=.5] (6B5)--(6B3);

\node[minimum size=.1pt,label={left:6B}](6B)at(15.15,-13.65){};

\end{tikzpicture}
\caption{Finite exceptions -1-}
\label{figure 1}
\end{figure}
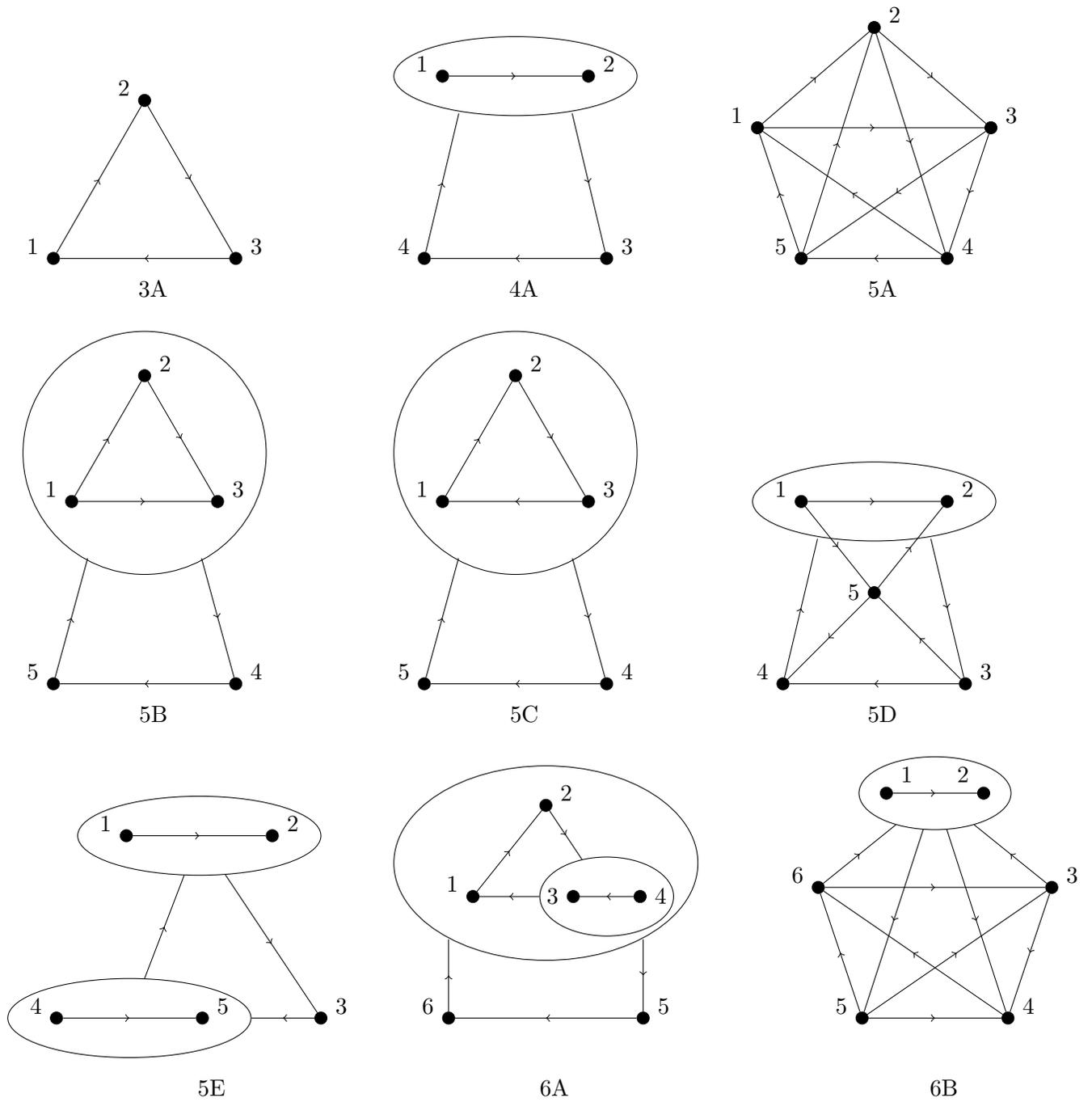
\begin{figure}
\begin{tikzpicture}
\tikzset{enclosed/.style={draw,circle,inner sep=2pt,minimum size=4pt,fill=black}}
\tikzset{->-/.style={decoration={
            markings,
            mark=at position #1 with
            {\arrow{>}}},postaction={decorate}}}
            
\node[ellipse,minimum width=3cm,minimum height=3cm,draw](a)at(0,0){};
\node[enclosed,label={left,yshift=.2cm:1}](1)at(-0.9,-0.7){};
\node[enclosed,label={right,yshift=.2cm:2}](2)at(0,0.9){};\node[enclosed,label={right,yshift=.2cm:3}](3)at(0.9,-0.7){};
\draw[black,->-=.5] (1)--(2);
\draw[black,->-=.5] (2)--(3);
\draw[black,->-=.5] (3)--(1);
\node[enclosed,label={right,yshift=.2cm:6}](6)at(1.5,-3.3){};
\node[ellipse,minimum width=3cm,minimum height=1.3cm,draw](b)at(-1.5,-3.3){};
\node[enclosed,label={right,yshift=.3cm:4}](4)at(-2.5,-3.3){};
\node[enclosed,label={left,yshift=.3cm:5}](5)at(-0.5,-3.3){};
\draw[black,->-=.5] (4)--(5);
\draw[black,->-=.5] (6)--(b);
\draw[black,->-=.5] (b)--(a);
\draw[black,->-=.5] (a)--(6);

\node[minimum size=.1pt,label={left:6C}](6C)at(.65,-4){};

\node[ellipse,minimum width=3cm,minimum height=1.3cm,draw](6Da)at(5,0.3){};
\node[ellipse,minimum width=1.3cm,minimum height=3cm,draw](6Db)at(8,-1.8){};
\node[ellipse,minimum width=3cm,minimum height=1.3cm,draw](6Dc)at(5,-3.3){};
\node[enclosed,label={right,yshift=.3cm:1}](6D1)at(4,0.3){};
\node[enclosed,label={left,yshift=.3cm:2}](6D2)at(6,0.3){};
\node[enclosed,label={left,yshift=-.2cm:3}](6D3)at(8,-0.8){};
\node[enclosed,label={left,yshift=.2cm:4}](6D4)at(8,-2.8){};
\node[enclosed,label={left,yshift=.3cm:5}](6D5)at(6,-3.3){};
\node[enclosed,label={right,yshift=.3cm:6}](6D6)at(4,-3.3){};

\draw[black,->-=.5] (6D1)--(6D2);
\draw[black,->-=.5] (6D3)--(6D4);
\draw[black,->-=.5] (6D5)--(6D6);
\draw[black,->-=.5] (6Da)--(6Db);
\draw[black,->-=.5] (6Db)--(6Dc);
\draw[black,->-=.5] (6Dc)--(6Da);

\node[minimum size=.1pt,label={left:6D}](6D)at(6.65,-4){};

\node[enclosed,label={left,yshift=.2cm:6}](6E6)at(10.08,-1.15){};

\node[ellipse,minimum width=2.5cm,minimum height=1.2cm,draw](6Ea)at(12,0.4){};

\node[enclosed,label={right,yshift=.3cm:1}](6E1)at(11.2,0.4){};
\node[enclosed,label={left,yshift=.3cm:2}](6E2)at(12.8,0.4){};

\node[enclosed,label={right,yshift=.2cm:3}](6E3)at(13.92,-1.15){};
\node[enclosed,label={right,yshift=.2cm:4}](6E4)at(13.2,-3.3){};\node[enclosed,label={left,yshift=.2cm:5}](6E5)at(10.8,-3.3){};
\draw[black,->-=.5] (6E1)--(6E2);
\draw[black,->-=.5] (6E6)--(6Ea);
\draw[black,->-=.5] (6E6)--(6E3);
\draw[black,->-=.5] (6Ea)--(6E3);
\draw[black,->-=.5] (6Ea)--(6E4);
\draw[black,->-=.5] (6E3)--(6E4);
\draw[black,->-=.5] (6E3)--(6E5);
\draw[black,->-=.5] (6E4)--(6E5);
\draw[black,->-=.5] (6E4)--(6E6);
\draw[black,->-=.5] (6E5)--(6E6);
\draw[black,->-=.5] (6E5)--(6Ea);

\node[minimum size=.1pt,label={left:6E}](6E)at(12.65,-4){};

\node[enclosed,label={right,yshift=.2cm:1}](6F1)at(1.2,-5.5){};
\node[enclosed,label={right,yshift=.2cm:2}](6F2)at(1.2,-9.9){};
\node[enclosed,label={left,yshift=.2cm:3}](6F3)at(-2.5,-7.7){};
\node[enclosed,label={left,yshift=.2cm:4}](6F4)at(-1.2,-9.9){};
\node[enclosed,label={left,yshift=.2cm:5}](6F5)at(-1.2,-5.5){};\node[enclosed,label={right,yshift=.2cm:6}](6F6)at(2.5,-7.7){};

\draw[black,->-=.4] (6F1)--(6F2);
\draw[black,->-=.5] (6F1)--(6F5);
\draw[black,->-=.2] (6F2)--(6F3);
\draw[black,->-=.5] (6F2)--(6F6);
\draw[black,->-=.2] (6F3)--(6F1);
\draw[black,->-=.5] (6F3)--(6F4);
\draw[black,->-=.4] (6F4)--(6F1);
\draw[black,->-=.5] (6F4)--(6F2);
\draw[black,->-=.2] (6F4)--(6F6);
\draw[black,->-=.4] (6F5)--(6F2);
\draw[black,->-=.5] (6F5)--(6F3);
\draw[black,->-=.4] (6F5)--(6F4);
\draw[black,->-=.5] (6F6)--(6F1);
\draw[black,->-=.4] (6F6)--(6F3);
\draw[black,->-=.2] (6F6)--(6F5);

\node[minimum size=.1pt,label={left:6F}](6F)at(.65,-10.6){};

\node[enclosed,label={right,yshift=.2cm:2}](6G2)at(7.6,-5.5){};
\node[enclosed,label={right,yshift=.2cm:4}](6G4)at(7.6,-9.9){};
\node[enclosed,label={left,yshift=.2cm:6}](6G6)at(3.9,-7.7){};
\node[enclosed,label={left,yshift=.2cm:5}](6G5)at(5.2,-9.9){};
\node[enclosed,label={left,yshift=.2cm:1}](6G1)at(5.2,-5.5){};\node[enclosed,label={right,yshift=.2cm:3}](6G3)at(8.9,-7.7){};

\draw[black,->-=.4] (6G1)--(6G4);
\draw[black,->-=.4] (6G1)--(6G5);
\draw[black,->-=.5] (6G2)--(6G1);
\draw[black,->-=.4] (6G2)--(6G5);
\draw[black,->-=.2] (6G3)--(6G1);
\draw[black,->-=.5] (6G3)--(6G2);
\draw[black,->-=.4] (6G3)--(6G6);
\draw[black,->-=.4] (6G4)--(6G2);
\draw[black,->-=.5] (6G4)--(6G3);
\draw[black,->-=.2] (6G5)--(6G3);
\draw[black,->-=.5] (6G5)--(6G4);
\draw[black,->-=.5] (6G5)--(6G6);
\draw[black,->-=.5] (6G6)--(6G1);
\draw[black,->-=.2] (6G6)--(6G2);
\draw[black,->-=.2] (6G6)--(6G4);

\node[minimum size=.1pt,label={left:6G}](6G)at(7.05,-10.6){};

\node[ellipse,minimum width=3cm,minimum height=3cm,draw](6Ha)at(12,-6.6){};
\node[enclosed,label={left,yshift=.2cm:1}](6H1)at(11.1,-7.3){};
\node[enclosed,label={right,yshift=.2cm:2}](6H2)at(12,-5.7){};\node[enclosed,label={right,yshift=.2cm:3}](6H3)at(12.9,-7.3){};
\draw[black,->-=.5] (6H1)--(6H2);
\draw[black,->-=.5] (6H2)--(6H3);
\draw[black,->-=.5] (6H3)--(6H1);
\node[enclosed,label={right,yshift=.2cm:6}](6H6)at(13.5,-9.9){};
\node[ellipse,minimum width=3cm,minimum height=1.3cm,draw](6Hb)at(10.5,-9.9){};
\node[enclosed,label={right,yshift=.3cm:4}](6H4)at(9.5,-9.9){};
\node[enclosed,label={left,yshift=.3cm:5}](6H5)at(11.5,-9.9){};
\draw[black,->-=.5] (6H5)--(6H4);
\draw[black,->-=.5] (6Hb)--(6H6);
\draw[black,->-=.5] (6Ha)--(6Hb);
\draw[black,->-=.5] (6H6)--(6Ha);

\node[minimum size=.1pt,label={left:6H}](6H)at(12.5,-10.6){};

\node[ellipse,minimum width=4cm,minimum height=1.3cm,draw](6Ia)at(0.2,-13){};
\node[minimum size=.1pt](6Ia1)at(-0.7,-13.48){};
\node[minimum size=.1pt](6Ia2)at(1.1,-13.48){};
\node[enclosed,label={left,yshift=.1cm:5}](6I5)at(-1,-13){};
\node[enclosed,label={right,yshift=.1cm:4}](6I4)at(1.4,-13){};
\draw[black,->-=.5] (6I5)--(6I4);
\node[enclosed,label={right,yshift=.2cm:3}](6I3)at(2.2,-16){};

\node[ellipse,minimum width=4cm,minimum height=1.3cm,draw](6Ib)at(-1.3,-16){};
\node[enclosed,label={left,yshift=.1cm:2}](6I2)at(-2.5,-16){};
\node[enclosed,label={right,yshift=.1cm:1}](6I1)at(-0.1,-16){};
\draw[black,->-=.5] (6I2)--(6I1);
\draw[black,->-=.5] (6I3)--(6Ib);
\draw[black,->-=.5] (6Ib)--(6Ia1);
\draw[black,->-=.5] (6Ia2)--(6I3);
\node[enclosed,label={right,yshift=.1cm:6}](6I6)at(0.2,-14.8){};
\draw[black,->-=.5] (6I5)--(6I6);
\draw[black,->-=.5] (6I6)--(6I4);
\draw[black,->-=.5] (6I3)--(6I6);
\draw[black,->-=.7] (6I1)--(6I6);
\draw[black,->-=.7] (6I6)--(6I2);

\node[minimum size=.1pt,label={left:6I}](6I)at(0.5,-16.7){};

\node[ellipse,minimum width=3cm,minimum height=3cm,draw](6Ja)at(6,-12.7){};
\node[enclosed,label={left,yshift=.2cm:2}](6J2)at(5.1,-13.4){};
\node[enclosed,label={right,yshift=.2cm:1}](6J1)at(6,-11.8){};\node[enclosed,label={right,yshift=.2cm:6}](6J6)at(6.9,-13.4){};
\draw[black,->-=.5] (6J2)--(6J1);
\draw[black,->-=.5] (6J2)--(6J6);
\draw[black,->-=.5] (6J1)--(6J6);
\node[enclosed,label={right,yshift=.2cm:3}](6J3)at(7.5,-16){};
\node[ellipse,minimum width=3cm,minimum height=1.3cm,draw](6Jb)at(4.5,-16){};
\node[enclosed,label={right,yshift=.3cm:4}](6J4)at(3.5,-16){};
\node[enclosed,label={left,yshift=.3cm:5}](6J5)at(5.5,-16){};
\draw[black,->-=.5] (6J5)--(6J4);
\draw[black,->-=.5] (6Jb)--(6J3);
\draw[black,->-=.5] (6Ja)--(6Jb);
\draw[black,->-=.5] (6J3)--(6Ja);

\node[minimum size=.1pt,label={left:6J}](6J)at(6.2,-16.7){};

\node[enclosed,label={left,yshift=.2cm:6}](6K6)at(10.08,-13.85){};

\node[ellipse,minimum width=2.5cm,minimum height=1.2cm,draw](6Ka)at(12,-12.3){};

\node[enclosed,label={right,yshift=.2cm:1}](6K1)at(11.08,-12.3){};
\node[enclosed,label={left,yshift=.2cm:2}](6K2)at(12.88,-12.3){};

\node[enclosed,label={right,yshift=.2cm:3}](6K3)at(13.92,-13.85){};
\node[enclosed,label={right,yshift=.2cm:4}](6K4)at(13.2,-16){};\node[enclosed,label={left,yshift=.2cm:5}](6K5)at(10.8,-16){};
\draw[black,->-=.5] (6K1)--(6K2);
\draw[black,->-=.5] (6K6)--(6Ka);
\draw[black,->-=.5] (6K6)--(6K3);
\draw[black,->-=.5] (6K3)--(6Ka);
\draw[black,->-=.5] (6Ka)--(6K4);
\draw[black,->-=.5] (6K4)--(6K3);
\draw[black,->-=.5] (6K5)--(6K3);
\draw[black,->-=.5] (6K5)--(6K4);
\draw[black,->-=.5] (6K6)--(6K4);
\draw[black,->-=.5] (6K5)--(6K6);
\draw[black,->-=.5] (6Ka)--(6K5);

\node[minimum size=.1pt,label={left:6K}](6K)at(12.65,-16.7){};

\end{tikzpicture}
\caption{Finite exceptions -2-}
\label{figure 2}
\end{figure}
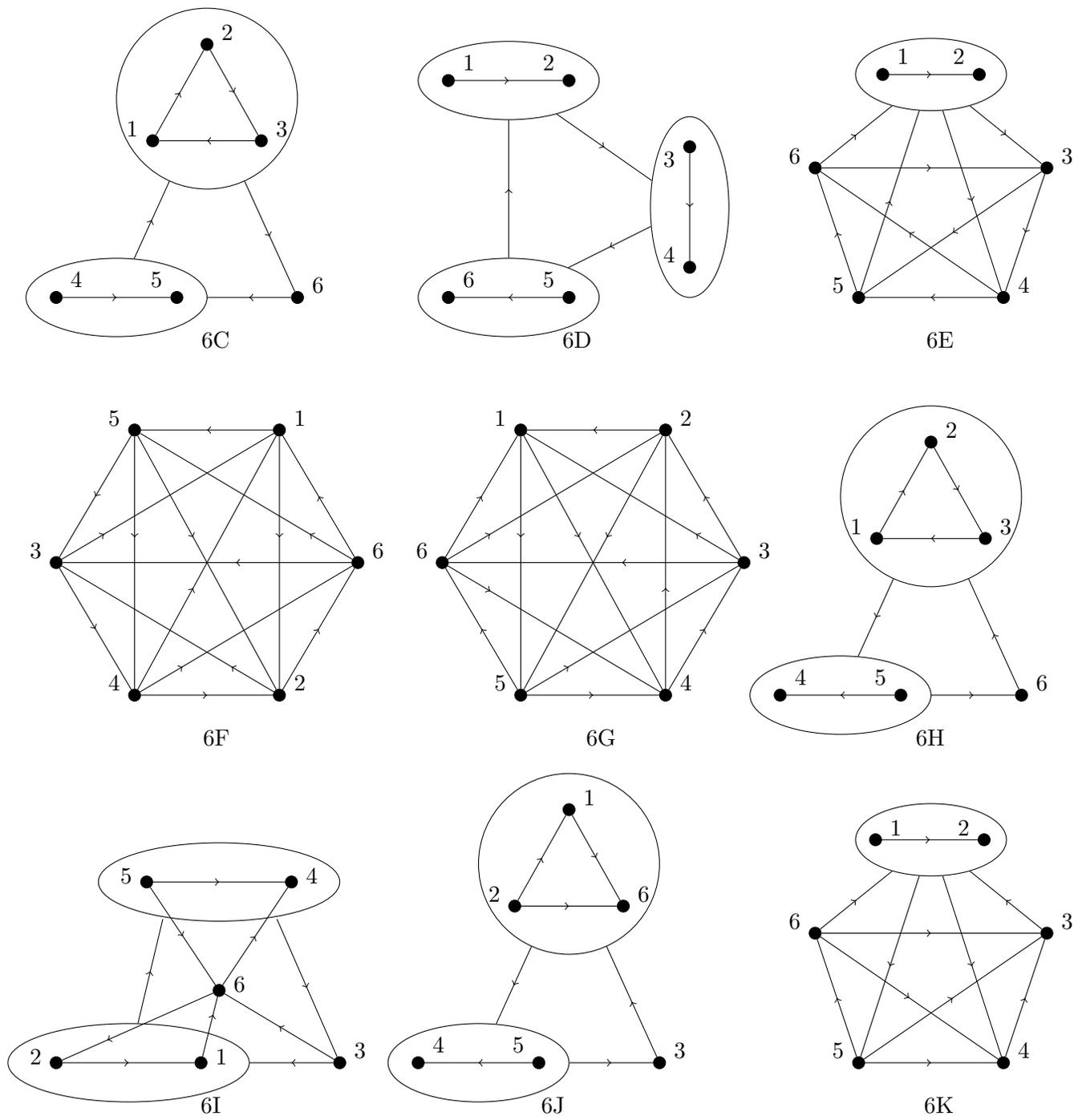
\begin{figure}
\begin{tikzpicture}
\tikzset{enclosed/.style={draw,circle,inner sep=2pt,minimum size=4pt,fill=black}}
\tikzset{->-/.style={decoration={
            markings,
            mark=at position #1 with
            {\arrow{>}}},postaction={decorate}}}
            
\node[ellipse,minimum width=3cm,minimum height=3cm,draw](a)at(0,0){};
\node[enclosed,label={left,yshift=.2cm:1}](6L1)at(-0.9,-0.7){};
\node[enclosed,label={right,yshift=.2cm:2}](6L2)at(0,0.9){};\node[enclosed,label={right,yshift=.2cm:3}](6L3)at(0.9,-0.7){};
\draw[black,->-=.5] (6L1)--(6L2);
\draw[black,->-=.5] (6L1)--(6L3);
\draw[black,->-=.5] (6L2)--(6L3);
\node[enclosed,label={right,yshift=.2cm:6}](6L6)at(1.5,-3.3){};
\node[ellipse,minimum width=3cm,minimum height=1.3cm,draw](b)at(-1.5,-3.3){};
\node[enclosed,label={right,yshift=.3cm:5}](6L5)at(-2.5,-3.3){};
\node[enclosed,label={left,yshift=.3cm:4}](6L4)at(-0.5,-3.3){};
\draw[black,->-=.5] (6L4)--(6L5);
\draw[black,->-=.5] (6L6)--(b);
\draw[black,->-=.5] (b)--(a);
\draw[black,->-=.5] (a)--(6L6);

\node[minimum size=.1pt,label={left:6L}](6L)at(.65,-4){};

\node[enclosed,label={right,yshift=.2cm:1}](1)at(6.2,1){};
\node[enclosed,label={right,yshift=.2cm:2}](2)at(7.6,-.7){};
\node[enclosed,label={right,yshift=.2cm:3}](3)at(7,-2.5){};
\node[enclosed,label={right,yshift=-.2cm:4}](4)at(5.3,-3.3){};
\node[enclosed,label={left,yshift=.2cm:5}](5)at(3.6,-2.5){};
\node[enclosed,label={left,yshift=.2cm:6}](6)at(3,-.7){};
\node[enclosed,label={left,yshift=.2cm:7}](7)at(4.4,1){};

\draw[black,->-=.5] (1)--(2);
\draw[black,->-=.4] (1)--(3);
\draw[black,->-=.65] (1)--(5);
\draw[black,->-=.5] (2)--(3);
\draw[black,->-=.3] (2)--(4);
\draw[black,->-=.7] (2)--(6);
\draw[black,->-=.5] (3)--(4);
\draw[black,->-=.4] (3)--(5);
\draw[black,->-=.6] (3)--(7);
\draw[black,->-=.5] (4)--(5);
\draw[black,->-=.4] (4)--(6);
\draw[black,->-=.7] (4)--(1);
\draw[black,->-=.5] (5)--(6);
\draw[black,->-=.35] (5)--(7);
\draw[black,->-=.7] (5)--(2);
\draw[black,->-=.5] (6)--(7);
\draw[black,->-=.2] (6)--(1);
\draw[black,->-=.75] (6)--(3);
\draw[black,->-=.5] (7)--(1);
\draw[black,->-=.2] (7)--(2);
\draw[black,->-=.3] (7)--(4);

\node[minimum size=.1pt,label={left:7A}](7A)at(5.95,-4){};

\node[ellipse,minimum width=3cm,minimum height=3cm,draw](7Ba)at(11.4,0){};
\node[enclosed,label={left,yshift=.2cm:3}](7B3)at(10.6,-0.6){};
\node[enclosed,label={right,yshift=.2cm:1}](7B1)at(11.4,0.8){};\node[enclosed,label={right,yshift=.2cm:2}](7B2)at(12.2,-0.6){};
\draw[black,->-=.5] (7B1)--(7B2);
\draw[black,->-=.5] (7B3)--(7B1);
\draw[black,->-=.5] (7B3)--(7B2);
\node[enclosed,label={right,yshift=.2cm:7}](7B7)at(12.9,-3.3){};
\node[ellipse,minimum width=3cm,minimum height=3cm,draw](7Bb)at(9.9,-3.3){};
\node[enclosed,label={left,yshift=.2cm:6}](7B6)at(9.1,-3.9){};
\node[enclosed,label={right,yshift=.2cm:4}](7B4)at(9.9,-2.5){};
\node[enclosed,label={right,yshift=.2cm:5}](7B5)at(10.7,-3.9){};
\draw[black,->-=.5] (7B4)--(7B5);
\draw[black,->-=.5] (7B5)--(7B6);
\draw[black,->-=.5] (7B6)--(7B4);
\draw[black,->-=.5] (7Ba)--(7B7);
\draw[black,->-=.5] (7B7)--(7Bb);
\draw[black,->-=.5] (7Bb)--(7Ba);

\node[minimum size=.1pt,label={left:7B}](7B)at(12.35,-4){};

\node[ellipse,minimum width=3cm,minimum height=3cm,draw](7Ca)at(2.8,-6.8){};
\node[enclosed,label={left,yshift=.2cm:3}](7C3)at(2,-7.4){};
\node[enclosed,label={right,yshift=.2cm:1}](7C1)at(2.8,-6){};\node[enclosed,label={right,yshift=.2cm:2}](7C2)at(3.6,-7.4){};
\draw[black,->-=.5] (7C1)--(7C2);
\draw[black,->-=.5] (7C3)--(7C1);
\draw[black,->-=.5] (7C2)--(7C3);
\node[enclosed,label={right,yshift=.2cm:7}](7C7)at(4.3,-10.1){};
\node[ellipse,minimum width=3cm,minimum height=3cm,draw](7Cb)at(1.3,-10.1){};
\node[enclosed,label={left,yshift=.2cm:6}](7C6)at(0.5,-10.7){};
\node[enclosed,label={right,yshift=.2cm:4}](7C4)at(1.3,-9.3){};
\node[enclosed,label={right,yshift=.2cm:5}](7C5)at(2.1,-10.7){};
\draw[black,->-=.5] (7C4)--(7C5);
\draw[black,->-=.5] (7C5)--(7C6);
\draw[black,->-=.5] (7C6)--(7C4);
\draw[black,->-=.5] (7Ca)--(7C7);
\draw[black,->-=.5] (7C7)--(7Cb);
\draw[black,->-=.5] (7Cb)--(7Ca);

\node[minimum size=.1pt,label={left:7C}](7C)at(3.5,-11.2){};

\node[ellipse,minimum width=3cm,minimum height=1.5cm,draw](7Da)at(6.5,-7.1){};
\node[ellipse,minimum width=3cm,minimum height=3cm,draw](7Db)at(10.5,-8.7){};
\node[ellipse,minimum width=3cm,minimum height=1.5cm,draw](7Dc)at(6.5,-10.1){};
\node[enclosed,label={right,yshift=.3cm:1}](7D1)at(5.5,-7.1){};
\node[enclosed,label={left,yshift=.3cm:2}](7D2)at(7.5,-7.1){};
\node[enclosed,label={left,yshift=.2cm:3}](7D3)at(10.5,-7.8){};
\node[enclosed,label={left,yshift=.2cm:5}](7D5)at(9.8,-9.1){};
\node[enclosed,label={right,yshift=.2cm:4}](7D4)at(11.2,-9.1){};
\node[enclosed,label={left,yshift=.3cm:6}](7D6)at(7.5,-10.1){};
\node[enclosed,label={right,yshift=.3cm:7}](7D7)at(5.5,-10.1){};

\draw[black,->-=.5] (7D1)--(7D2);
\draw[black,->-=.5] (7D3)--(7D4);
\draw[black,->-=.5] (7D4)--(7D5);
\draw[black,->-=.5] (7D5)--(7D3);
\draw[black,->-=.5] (7D6)--(7D7);
\draw[black,->-=.5] (7Da)--(7Db);
\draw[black,->-=.5] (7Db)--(7Dc);
\draw[black,->-=.5] (7Dc)--(7Da);

\node[minimum size=.1pt,label={left:7D}](7D)at(8.5,-11){};

\node[ellipse,minimum width=5cm,minimum height=4cm,draw](7Ea)at(3.5,-14.5){};
\node[ellipse,minimum width=3cm,minimum height=3cm,draw](7Eb)at(2.7,-14.5){};
\node[ellipse,minimum width=3cm,minimum height=1.5cm,draw](7Ec)at(-1.5,-16){};
\node[enclosed,label={right,yshift=.3cm:1}](7E1)at(-1.5,-13){};
\node[enclosed,label={right,yshift=.3cm:2}](7E2)at(5.3,-14.5){};
\node[enclosed,label={left,yshift=.2cm:3}](7E3)at(2.5,-13.7){};
\node[enclosed,label={left,yshift=.2cm:5}](7E5)at(1.8,-15){};
\node[enclosed,label={right,yshift=.2cm:4}](7E4)at(3.2,-15){};
\node[enclosed,label={left,yshift=.3cm:6}](7E6)at(-0.5,-16){};
\node[enclosed,label={right,yshift=.3cm:7}](7E7)at(-2.5,-16){};

\draw[black,->-=.5] (7E1)--(7Ea);
\draw[black,->-=.5] (7E3)--(7E4);
\draw[black,->-=.5] (7E4)--(7E5);
\draw[black,->-=.5] (7E5)--(7E3);
\draw[black,->-=.5] (7E6)--(7E7);
\draw[black,->-=.5] (7Eb)--(7E2);
\draw[black,->-=.5] (7Ea)--(7Ec);
\draw[black,->-=.5] (7Ec)--(7E1);

\node[minimum size=.1pt,label={left:7E}](7E)at(.65,-17){};

\node[ellipse,minimum width=3cm,minimum height=1.5cm,draw](7Fa)at(8,-13){};
\node[ellipse,minimum width=3cm,minimum height=3cm,draw](7Fb)at(12,-14.5){};
\node[ellipse,minimum width=3cm,minimum height=1.5cm,draw](7Fc)at(8,-16){};
\node[enclosed,label={right,yshift=.3cm:5}](7F5)at(7,-13){};
\node[enclosed,label={left,yshift=.3cm:4}](7F4)at(9,-13){};
\node[enclosed,label={left,yshift=.2cm:1}](7F1)at(12,-13.7){};
\node[enclosed,label={left,yshift=.2cm:2}](7F2)at(11.3,-15){};
\node[enclosed,label={right,yshift=.2cm:3}](7F3)at(12.7,-15){};
\node[enclosed,label={left,yshift=.3cm:6}](7F6)at(9,-16){};
\node[enclosed,label={right,yshift=.3cm:7}](7F7)at(7,-16){};

\draw[black,->-=.5] (7F1)--(7F3);
\draw[black,->-=.5] (7F3)--(7F2);
\draw[black,->-=.5] (7F2)--(7F1);
\draw[black,->-=.5] (7F4)--(7F5);
\draw[black,->-=.5] (7F6)--(7F7);
\draw[black,->-=.5] (7Fa)--(7Fb);
\draw[black,->-=.5] (7Fb)--(7Fc);
\draw[black,->-=.5] (7Fc)--(7Fa);
\draw[black,->-=.5] (7F5)--(7F7);

\node[minimum size=.1pt,label={left:7F}](7F)at(10,-17){};

\end{tikzpicture}
\caption{Finite exceptions -3-}
\label{figure 3}
\end{figure}
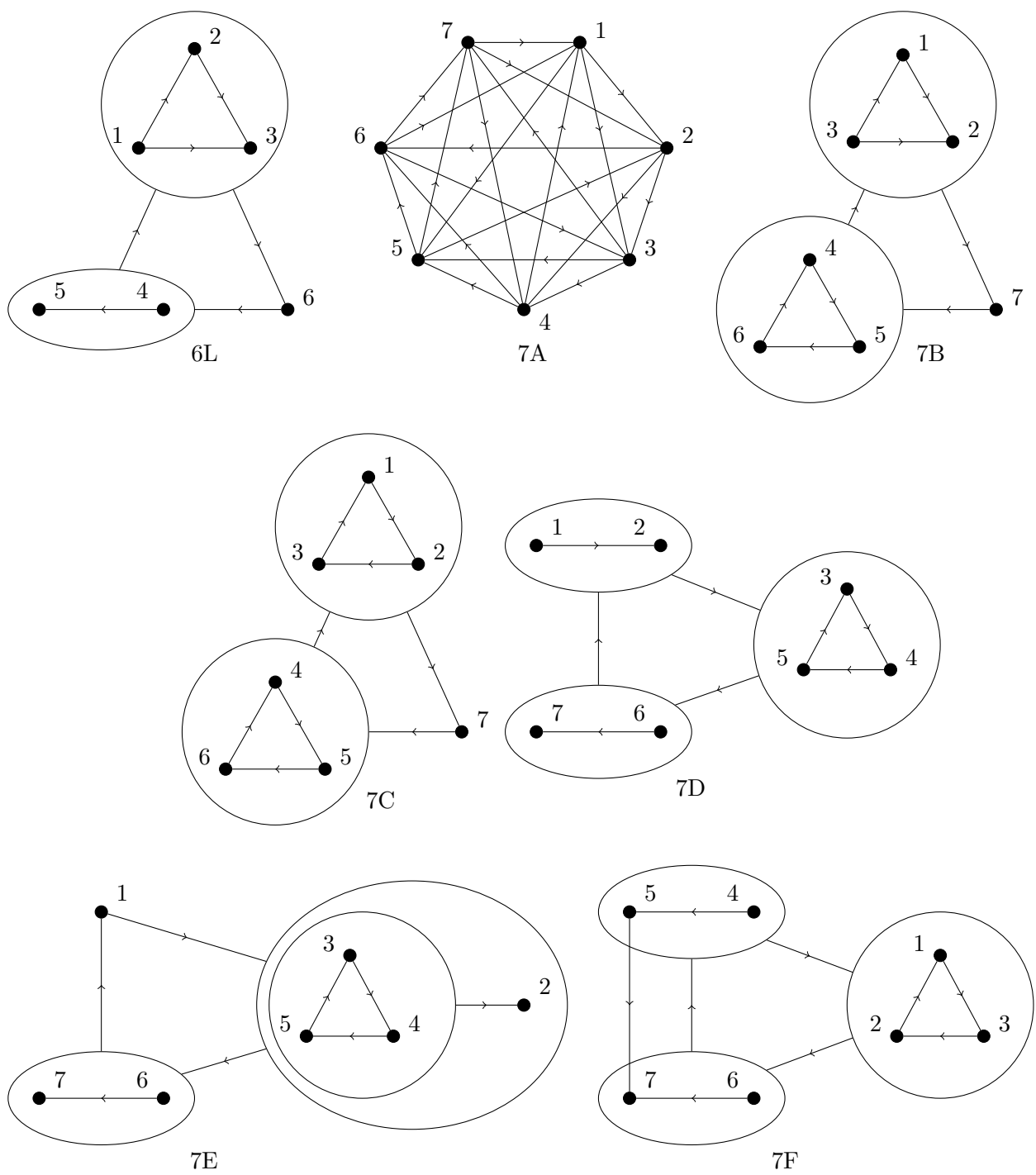
\begin{figure}
\begin{tikzpicture}
\tikzset{enclosed/.style={draw,circle,inner sep=2pt,minimum size=4pt,fill=black}}
\tikzset{->-/.style={decoration={
            markings,
            mark=at position #1 with
            {\arrow{>}}},postaction={decorate}}}
            
\node[ellipse,minimum width=5cm,minimum height=4cm,draw](a)at(5,-1.5){};
\node[ellipse,minimum width=1.3cm,minimum height=3cm,draw](b)at(6,-1.5){};
\node[ellipse,minimum width=3cm,minimum height=1.5cm,draw](c)at(0,-3){};
\node[enclosed,label={right,yshift=.3cm:1}](1)at(6,-2.3){};
\node[enclosed,label={left,yshift=.2cm:2}](2)at(4,-2.3){};
\node[enclosed,label={left,yshift=.2cm:3}](3)at(4,-0.7){};
\node[enclosed,label={left,yshift=.2cm:4}](4)at(1,-3){};
\node[enclosed,label={right,yshift=.2cm:5}](5)at(-1,-3){};
\node[enclosed,label={left,yshift=.2cm:6}](6)at(0,0){};
\node[enclosed,label={right,yshift=.2cm:7}](7)at(6,-0.7){};
\node[minimum size=.1pt](b1)at(5.63,-0.7){};
\node[minimum size=.1pt](b2)at(5.63,-2.3){};

\draw[black,->-=.5] (2)--(3);
\draw[black,->-=.5] (3)--(b1);
\draw[black,->-=.5] (b2)--(2);
\draw[black,->-=.5] (7)--(1);
\draw[black,->-=.5] (a)--(6);
\draw[black,->-=.5] (6)--(c);
\draw[black,->-=.5] (c)--(a);
\draw[black,->-=.5] (4)--(5);

\node[minimum size=.1pt,label={left:7G}](7G)at(2.65,-4){};

\node[enclosed,label={left,yshift=.2cm:6}](7H6)at(9.2,-.85){};

\node[ellipse,minimum width=2.5cm,minimum height=1.2cm,draw](7Ha)at(11.12,.7){};

\node[enclosed,label={right,yshift=.3cm:1}](7H1)at(10.2,.7){};
\node[enclosed,label={left,yshift=.3cm:2}](7H2)at(12,.7){};

\node[enclosed,label={right,yshift=.2cm:3}](7H3)at(13.04,-.85){};

\node[ellipse,minimum width=2.5cm,minimum height=1.2cm,draw](7Hb)at(12.32,-3){};

\node[enclosed,label={left,yshift=.3cm:4}](7H4)at(13.22,-3){};
\node[enclosed,label={right,yshift=.3cm:7}](7H7)at(11.42,-3){};

\node[enclosed,label={left,yshift=.2cm:5}](7H5)at(9.92,-3){};

\draw[black,->-=.5] (7H1)--(7H2);
\draw[black,->-=.5] (7H6)--(7Ha);
\draw[black,->-=.5] (7H6)--(7H3);
\draw[black,->-=.5] (7H3)--(7Ha);
\draw[black,->-=.5] (7Ha)--(7Hb);
\draw[black,->-=.5] (7Hb)--(7H3);
\draw[black,->-=.5] (7H5)--(7H3);
\draw[black,->-=.5] (7H5)--(7Hb);
\draw[black,->-=.6] (7H6)--(7Hb);
\draw[black,->-=.5] (7H5)--(7H6);
\draw[black,->-=.5] (7Ha)--(7H5);
\draw[black,->-=.5] (7H7)--(7H4);

\node[minimum size=.1pt,label={left:7H}](7H)at(11.72,-4){};

\node[ellipse,minimum width=3cm,minimum height=3cm,draw](7Ia)at(2,-6){};
\node[enclosed,label={left,yshift=.2cm:5}](7I5)at(1.2,-6.6){};
\node[enclosed,label={right,yshift=.2cm:4}](7I4)at(2,-5.2){};\node[enclosed,label={right,yshift=.2cm:7}](7I7)at(2.8,-6.6){};
\node[enclosed,label={right,yshift=.2cm:6}](7I6)at(3.5,-9.3){};
\node[ellipse,minimum width=3cm,minimum height=3cm,draw](7Ib)at(.5,-9.3){};
\node[enclosed,label={left,yshift=.2cm:3}](7I3)at(-.3,-9.9){};
\node[enclosed,label={right,yshift=.2cm:2}](7I2)at(0.5,-8.3){};
\node[enclosed,label={right,yshift=.2cm:1}](7I1)at(1.3,-9.9){};
\draw[black,->-=.5] (7I4)--(7I7);
\draw[black,->-=.5] (7I5)--(7I4);
\draw[black,->-=.5] (7I7)--(7I5);
\draw[black,->-=.5] (7Ia)--(7I6);
\draw[black,->-=.5] (7I6)--(7Ib);
\draw[black,->-=.5] (7Ib)--(7Ia);
\draw[black,->-=.5] (7I3)--(7I2);
\draw[black,->-=.5] (7I2)--(7I1);
\draw[black,->-=.5] (7I3)--(7I1);

\node[minimum size=.1pt,label={left:7I}](7I)at(2.65,-10.5){};

\node[ellipse,minimum width=4cm,minimum height=6cm,draw](7Ja)at(10.5,-7.3){};
\node[ellipse,minimum width=3cm,minimum height=3cm,draw](7Jb)at(10.5,-6.3){};
\node[enclosed,label={left,yshift=.2cm:1}](7J1)at(6.2,-7.9){};
\node[enclosed,label={left,yshift=.2cm:2}](7J2)at(6.2,-6.7){};\node[enclosed,label={right,yshift=.2cm:3}](7J3)at(7.4,-7.2){};
\node[enclosed,label={left,yshift=.2cm:4}](7J4)at(10,-5.7){};
\node[enclosed,label={right,yshift=.2cm:5}](7J5)at(11,-6.2){};
\node[enclosed,label={left,yshift=-.2cm:6}](7J6)at(10,-6.7){};
\node[enclosed,label={right,yshift=.2cm:7}](7J7)at(10.5,-8.9){};
\node[minimum size=.1pt](7Ja1)at(8.9,-5.9){};
\node[minimum size=.1pt](7Ja2)at(8.9,-8.9){};

\draw[black,->-=.5] (7J2)--(7J1);
\draw[black,->-=.5] (7J1)--(7J3);
\draw[black,->-=.5] (7J3)--(7J2);
\draw[black,->-=.5] (7J4)--(7J5);
\draw[black,->-=.5] (7J5)--(7J6);
\draw[black,->-=.5] (7J6)--(7J4);
\draw[black,->-=.5] (7Jb)--(7J7);
\draw[black,->-=.5] (7Ja1)--(7J2);
\draw[black,->-=.5] (7Ja2)--(7J1);
\draw[black,->-=.5] (7J3)--(7Jb);
\draw[black,->-=.5] (7J7)--(7J3);

\node[minimum size=.1pt,label={left:7J}](7J)at(9.2,-10.5){};

\node[enclosed,label={right,yshift=.2cm:6}](8A6)at(2.4,-11.5){};
\node[enclosed,label={right,yshift=.2cm:7}](8A7)at(3.8,-13.2){};
\node[enclosed,label={right,yshift=.2cm:8}](8A8)at(3.2,-15){};
\node[enclosed,label={right,yshift=.2cm:1}](8A1)at(.7,-16.45){};
\node[enclosed,label={left,yshift=.2cm:3}](8A3)at(-.2,-15){};
\node[enclosed,label={left,yshift=.2cm:4}](8A4)at(-.8,-13.2){};
\node[enclosed,label={left,yshift=.2cm:5}](8A5)at(.6,-11.5){};
\node[enclosed,label={left,yshift=.2cm:2}](8A2)at(2.3,-16.45){};
\node[minimum size=.1pt](8Aa)at(1.5,-15.9){};
\node[ellipse,minimum width=2.5cm,minimum height=1.3cm,draw](8Ab)at(1.5,-16.45){};

\draw[black,->-=.5] (8A1)--(8A2);
\draw[black,->-=.5] (8Aa)--(8A3);
\draw[black,->-=.6] (8Aa)--(8A4);
\draw[black,->-=.5] (8Aa)--(8A6);
\draw[black,->-=.5] (8A3)--(8A4);
\draw[black,->-=.6] (8A3)--(8A5);
\draw[black,->-=.25] (8A3)--(8A7);
\draw[black,->-=.5] (8A4)--(8A5);
\draw[black,->-=.7] (8A4)--(8A6);
\draw[black,->-=.43] (8A4)--(8A8);
\draw[black,->-=.5] (8A5)--(8A6);
\draw[black,->-=.8] (8A5)--(8A7);
\draw[black,->-=.54] (8A5)--(8Aa);
\draw[black,->-=.5] (8A6)--(8A7);
\draw[black,->-=.7] (8A6)--(8A8);
\draw[black,->-=.42] (8A6)--(8A3);
\draw[black,->-=.5] (8A7)--(8A8);
\draw[black,->-=.2] (8A7)--(8Aa);
\draw[black,->-=.5] (8A7)--(8A4);
\draw[black,->-=.5] (8A8)--(8Aa);
\draw[black,->-=.4] (8A8)--(8A3);
\draw[black,->-=.25] (8A8)--(8A5);

\node[minimum size=.1pt,label={left:8A}](8A)at(2.3,-17.5){};

\node[ellipse,minimum width=4cm,minimum height=6cm,draw](8Ba)at(10.5,-13.9){};
\node[ellipse,minimum width=3cm,minimum height=3cm,draw](8Bb)at(10.5,-12.9){};
\node[enclosed,label={left,yshift=.2cm:1}](8B1)at(6.2,-14.5){};

\node[ellipse,minimum width=2.5cm,minimum height=1cm,draw](8Bc)at(6.2,-13.3){};
\node[enclosed,label={right,yshift=.01cm:2}](8B2)at(6.8,-13.3){};
\node[enclosed,label={left,yshift=.01cm:8}](8B8)at(5.6,-13.3){};

\node[enclosed,label={right,yshift=.2cm:3}](8B3)at(7.4,-13.9){};
\node[enclosed,label={left,yshift=.2cm:4}](8B4)at(10,-12.3){};
\node[enclosed,label={right,yshift=.2cm:5}](8B5)at(11,-12.8){};
\node[enclosed,label={left,yshift=-.2cm:6}](8B6)at(10,-13.3){};
\node[enclosed,label={right,yshift=.2cm:7}](8B7)at(10.5,-15.5){};
\node[minimum size=.1pt](8Ba1)at(8.9,-12.5){};
\node[minimum size=.1pt](8Ba2)at(8.9,-15.5){};

\draw[black,->-=.5] (8Bc)--(8B1);
\draw[black,->-=.5] (8B1)--(8B3);
\draw[black,->-=.5] (8B3)--(8Bc);
\draw[black,->-=.5] (8B4)--(8B5);
\draw[black,->-=.5] (8B5)--(8B6);
\draw[black,->-=.5] (8B6)--(8B4);
\draw[black,->-=.5] (8Bb)--(8B7);
\draw[black,->-=.5] (8Ba1)--(8Bc);
\draw[black,->-=.5] (8Ba2)--(8B1);
\draw[black,->-=.5] (8B3)--(8Bb);
\draw[black,->-=.5] (8B7)--(8B3);
\draw[black,->-=.5] (8B8)--(8B2);

\node[minimum size=.1pt,label={left:8B}](8B)at(9.2,-17.5){};

\end{tikzpicture}
\caption{Finite exceptions -4-}
\label{figure 4}
\end{figure}
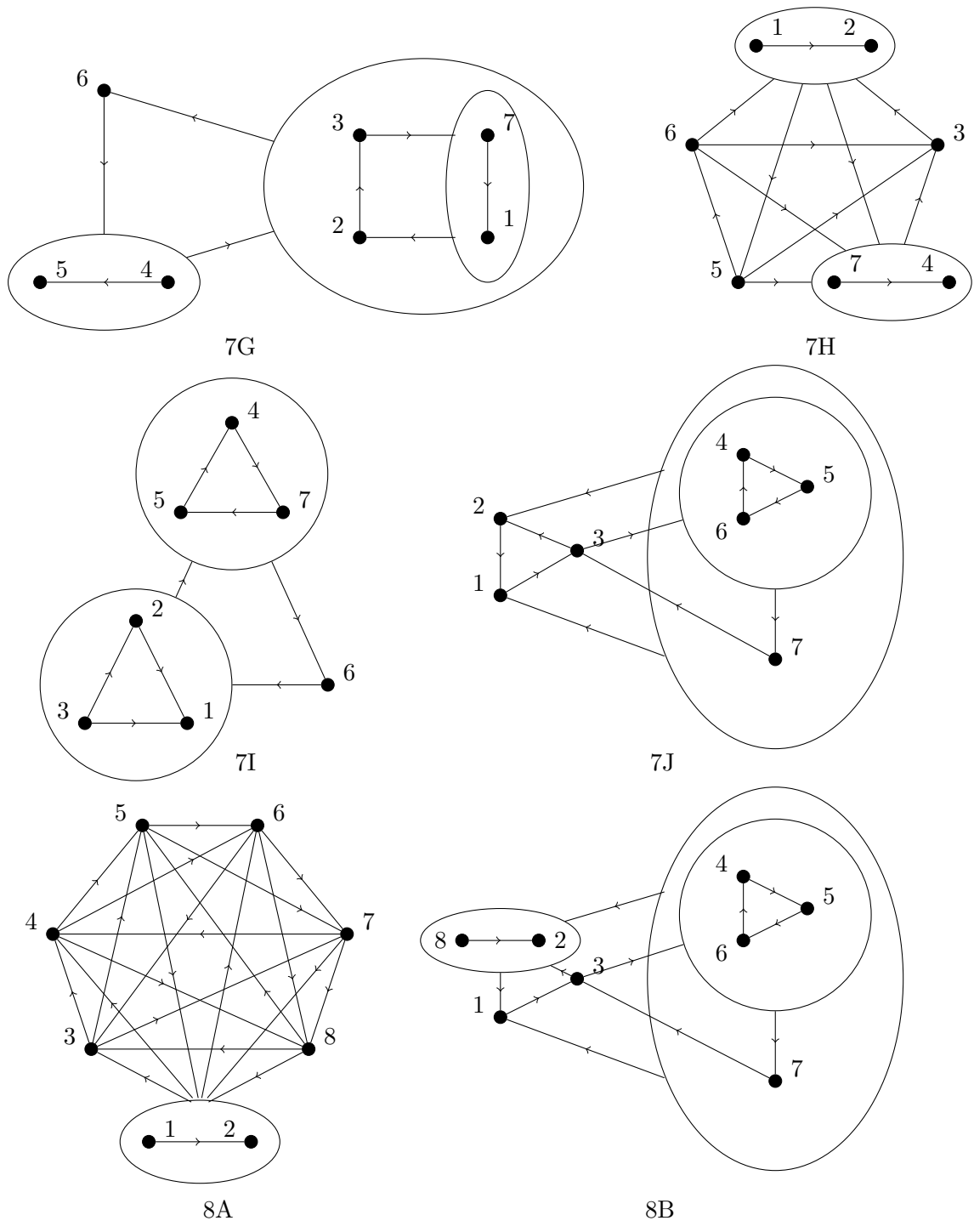
They established two categories of exceptions: the finite exceptions and the infinite families of exceptions. The notation of an exception is the following: $[T;P;S;P_1,...,P_k]$ where $T$ is a tournament illustrated in figures \ref{figure 1}, \ref{figure 2}, \ref{figure 3}, \ref{figure 4}, \ref{figure 5} and \ref{figure 6}, $P$ is an outpath, $S$ is the set of vertices of $T$ which are not origin of $P$ and the paths $P_1,...,P_k$ are the paths of $T$ whose origins are precisely the vertices of $V(T)\setminus S$.\\
Exc 0: $[3A;(1,1);\{1,2,3\}]$
\\Exc 1: $[4A;(1,1,1);\{1,2,3\};4213]$
\\Exc 2: $[4A;(1,2);\{3,4\};1324;2314]$
\\Exc 3: $[4A;(2,1);\{1,2,\};3421;4132]$
\\Exc 4: $[5A;(1,1,1,1);\{1,2,3,4,5\}]$
\\Exc 5: $[5B;(2,1,1);\{1,2,3\};45213;51423]$
\\Exc 6: $[5C;(1,1,2);\{4,5\};12534;23514;31524]$
\\Exc 7: $[5C;(2,1,1);\{1,2,3,4\};51432]$
\\Exc 8: $[5D;(1,1,1,1);\{2,5\};12543;35124;42153]$
\\Exc 9: $[5E;(1,1,1,1);\{2,4,53\};12453;35421]$
\\Exc 10: $[5E;(1,2,1);\{3,5\};12435;23145;45312]$
\\Exc 11: $[5E;(2,2);\{1,2\};34215;42315;52314]$
\\Exc 12: $[5E;(1,1,2);\{1,2\};35412;41523;51423]$
\\Exc 13: $[6A;(3,1,1);\{3,4\};156324;256143;562341;612345]$
\\Exc 14: $[6B;(2,1,1,1);\{3,4\};154326;254316;562143;612345]$
\\Exc 15: $[6C;(1,1,2,1);\{1,2,3,6\};435261;534261]$
\\Exc 16: $[6C;(1,2,1,1);\{4,5,6\};163425;263415;362415]$
\\Exc 17: $[6D;(2,1,1,1);\{2,4,6\};124365;346521;562143]$
\\Exc 18: $[6D;(1,2,2);\{2,4,6\};126345;341562;564123]$
\\Exc 19: $[6D;(1,1,1,2);\{2,4,6\};126543;341265;563421]$
\\Exc 20: $[6E;(1,1,1,1,1);\{1,2\};341256;465213;516324;621435]$
\\Exc 21: $[6E;(2,1,1,1);\{1,2\};346521;452136;562143;634125]$
\\Exc 21: $[6E;(2,1,1,1);\{1,2\};346521;452136;562143;634125]$
\\Exc 22: $[6F;(1,1,1,1,1);\{1,2,3\};421563;532641;613452]$
\\Exc 23: $[6G;(1,1,1,1,1);\{4,6\};145632;216453;326415;546132]$
\\Exc 24: $[6H;(1,1,1,1,1);\{1,2,3,4\};543162;613425]$
\\Exc 25: $[6H;(1,1,1,2);\{4,5\};142536;243516;341526;613452]$
\\Exc 26: $[6H;(1,1,3);\{4,5,6\};145623;245631;345612]$
\\Exc 27: $[6H;(1,3,1);\{4,6\};126534;236514;316524;543261]$
\\Exc 28: $[6H;(2,1,2);\{4,5\};124563;234561;314562;614235]$
\\Exc 29: $[6I;(1,1,1,1,1);\{4,6\};145632;213654;365421;546231]$
\\Exc 30: $[6J;(1,1,1,1,1);\{4,6\};162453;261453;312465;542631]$
\\Exc 31: $[6K;(1,2,2);\{3,4\};146532;246531;541632;634125]$
\\Exc 32: $[6L;(1,2,1,1);\{5,6\};163425;263415;361425;456132]$
\\Exc 33: $[7A;(1,1,1,1,1,1);\{1,2,3,4,5,6,7\}]$
\\Exc 34: $[7B;(1,1,2,1,1);\{1,2,3\};4576132;5674132;6475132;7541263]$
\\Exc 35: $[7B;(2,1,3);\{1,2,3\};4315627;5316427;6314527;7435612]$
\\Exc 36: $[7B;(2,3,1);\{1,2\};3125476;4567132;5647132;6457132;7421356]$
\\Exc 37: $[7C;(1,1,1,1,1,1);\{4,5,6\};1243567;2341567;3142567;7541632]$
\\Exc 38: $[7C;(1,1, 2,1,1);\{1,2,3\};4156327;5164327;6145327;7541263]$
\\Exc 39: $[7C;(2,1,3);\{1,2,3\};4315627;5316427;6314527;7435612]$
\\Exc 40: $[7D;(1,1,1,2,1);\{1,2\};3412756;4512736;5312746;6215437;7215436]$
\\Exc 41: $[7D;(1,1,1,3);\{6,7\};1546327;2546317;3745216;4753216;5734216]$
\\Exc 42: $[7D;(2,2,1,1);\{6,7\};1342675;2341675;3465127;4563127;5364127]$
\\Exc 43: $[7E;(1,1,2,1,1);\{2,7\};1236745;3214756;4213756;5213746;6734215]$
\\Exc 44: $[7F;(1,1,1,3);\{6,7\};1732546;2713546;3721546;4127635;5127634]$
\\Exc 45: $[7G;(2,1,2,1);\{1,7\};2654317;3654721;4367125;5367124;6517234]$
\\Exc 46: $[7H;(2,2,2);\{4,7\};1746532;2746531;3126574;5321674;6247531]$
\\Exc 47: $[7I;(1,1,2,1,1);\{4,5,7\};1456237;2456137;3456127;6135427]$
\\Exc 48: $[7J;(1,1,2,1,1);\{1,2\};3245167;4235167;5234167;6234157;7234156]$
\\Exc 49: $[8A;(1,1,1,1,1,1,1);\{1,2\};35461278;46527183;56487213;67341285;74358216;\\85347216]$
\\Exc 50: $[8A;(2,1,1,1,1,1);\{1,2\};34652718;46752138;56734128;68214375;78216453;\\83412576]$
\\Exc 51: $[8B;(2,1,2,1,1);\{2,8\};13245867;32145867;42156873;52164873;62145873;\\73245861]$.\\

The infinite families of exceptions are denoted by $E_i(n)=(F_i(n);P)$, where $F_i(n)$ is the tournament on $n$ vertices illustrated in figures \ref{figure 5} and \ref{figure 6}. For each $E_i(n)$, we define the set $S$ of vertices of $F_i(n)$ which are not origin of $P$ together with the conditions on the tournament. And finally, we give the paths $P$ with origin $x\notin S$. Below the list of infinite families of exceptions:
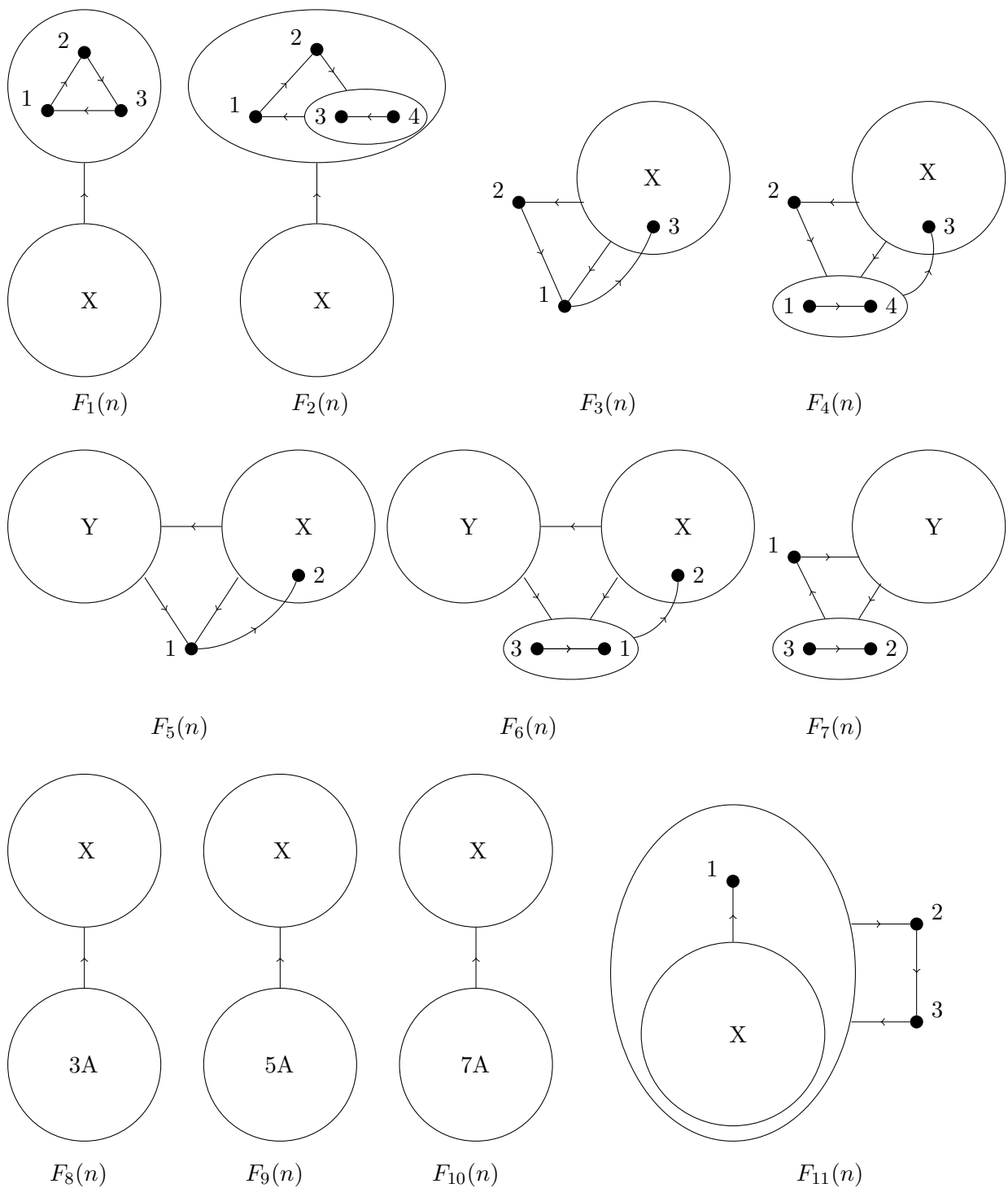
\begin{figure}
\begin{tikzpicture}
\tikzset{enclosed/.style={draw,circle,inner sep=2pt,minimum size=4pt,fill=black}}
\tikzset{->-/.style={decoration={
            markings,
            mark=at position #1 with
            {\arrow{>}}},postaction={decorate}}}
            
\node[ellipse,minimum width=2.5cm,minimum height=2.5cm,draw](a)at(0,0){};
\node[enclosed,label={left,yshift=.2cm:1}](1)at(-0.6,-0.4){};
\node[enclosed,label={left,yshift=.2cm:2}](2)at(0,0.55){};
\node[enclosed,label={right,yshift=.2cm:3}](3)at(0.6,-0.4){};
\node[ellipse,minimum width=2.5cm,minimum height=2.5cm,draw](b)at(0,-3.5){};
\node[minimum size=.1pt,label={left:X}](X)at(0.5,-3.5){};
\draw[black,->-=.5] (1)--(2);
\draw[black,->-=.5] (2)--(3);
\draw[black,->-=.5] (3)--(1);
\draw[black,->-=.5] (b)--(a);

\node[minimum size=.1pt,label={left:$F_1(n)$}](F_1)at(1,-5.2){};

\node[ellipse,minimum width=4.2cm,minimum height=2.5cm,draw](2Fa)at(3.8,0){};
\node[ellipse,minimum width=2cm,minimum height=0.9cm,draw](2Fc)at(4.6,-0.5){};
\node[enclosed,label={left,yshift=.2cm:1}](2F1)at(2.8,-0.5){};
\node[enclosed,label={left,yshift=.2cm:2}](2F2)at(3.8,0.6){};
\node[enclosed,label={left,yshift=.01cm:3}](2F3)at(4.2,-0.5){};
\node[enclosed,label={right,yshift=.01cm:4}](2F4)at(5.05,-0.5){};
\node[ellipse,minimum width=2.5cm,minimum height=2.5cm,draw](2Fb)at(3.8,-3.5){};
\node[minimum size=.1pt,label={left:X}](2FX)at(4.3,-3.5){};
\draw[black,->-=.5] (2F1)--(2F2);
\draw[black,->-=.5] (2F2)--(2Fc);
\draw[black,->-=.5] (2Fc)--(2F1);
\draw[black,->-=.5] (2F4)--(2F3);
\draw[black,->-=.5] (2Fb)--(2Fa);

\node[minimum size=.1pt,label={left:$F_2(n)$}](F_2)at(4.6,-5.2){};

\node[ellipse,minimum width=2.5cm,minimum height=2.5cm,draw](3Fa)at(9.3,-1.5){};
\node[minimum size=.1pt,label={left:X}](3FX)at(9.7,-1.5){};
\node[enclosed,label={left,yshift=.2cm:1}](3F1)at(7.85,-3.6){};
\node[enclosed,label={left,yshift=.2cm:2}](3F2)at(7.1,-1.9){};
\node[enclosed,label={right,yshift=.01cm:3}](3F3)at(9.3,-2.3){};
\node[minimum size=.1pt](3Fa1)at(8.3,-1.9){};
\node[minimum size=.1pt](3Fa2)at(8.7,-2.4){};

\draw[black,->-=.5] (3F2)--(3F1);
\draw[black,->-=.5] (3Fa2)--(3F1);
\draw[black,->-=.5] (3Fa1)--(3F2);
\draw[black,->-=.5] (3F1) to [out=0,in=-110,looseness=0.8] (3F3);

\node[minimum size=.1pt,label={left:$F_3(n)$}](F_3)at(9.3,-5.2){};

\node[ellipse,minimum width=2.5cm,minimum height=2.5cm,draw](4Fa)at(13.8,-1.5){};
\node[minimum size=.1pt,label={left:X}](4FX)at(14.2,-1.4){};
\node[enclosed,label={left,yshift=.01cm:1}](4F1)at(11.85,-3.6){};
\node[enclosed,label={left,yshift=.2cm:2}](4F2)at(11.6,-1.9){};
\node[enclosed,label={right,yshift=.01cm:3}](4F3)at(13.8,-2.3){};
\node[minimum size=.1pt](4Fa1)at(12.8,-1.9){};
\node[minimum size=.1pt](4Fa2)at(13.2,-2.4){};
\node[ellipse,minimum width=2.2cm,minimum height=1cm,draw](4Fb)at(12.35,-3.6){};
\node[enclosed,label={right,yshift=.01cm:4}](4F4)at(12.85,-3.6){};

\draw[black,->-=.5] (4F2)--(4Fb);
\draw[black,->-=.5] (4Fa2)--(4Fb);
\draw[black,->-=.5] (4Fa1)--(4F2);
\draw[black,->-=.5] (4Fb) to [out=10,in=-70,looseness=0.9] (4F3);
\draw[black,->-=.5] (4F1)--(4F4);

\node[minimum size=.1pt,label={left:$F_4(n)$}](F_4)at(13,-5.2){};

\node[ellipse,minimum width=2.5cm,minimum height=2.5cm,draw](5Fa)at(3.5,-7.2){};
\node[minimum size=.1pt,label={left:X}](5FX)at(4,-7.2){};
\node[enclosed,label={left,yshift=.01cm:1}](5F1)at(1.75,-9.2){};
\node[ellipse,minimum width=2.5cm,minimum height=2.5cm,draw](5Fb)at(0,-7.2){};
\node[minimum size=.1pt,label={left:Y}](5FY)at(.5,-7.2){};
\node[enclosed,label={right,yshift=.01cm:2}](5F2)at(3.5,-8){};
\node[minimum size=.1pt](5Fa1)at(2.6,-7.9){};
\node[minimum size=.1pt](5Fb1)at(.9,-7.9){};

\draw[black,->-=.5] (5Fa)--(5Fb);
\draw[black,->-=.5] (5Fa1)--(5F1);
\draw[black,->-=.5] (5Fb1)--(5F1);
\draw[black,->-=.5] (5F1) to [out=0,in=-110,looseness=0.8] (5F2);

\node[minimum size=.1pt,label={left:$F_5(n)$}](F_5)at(2.3,-10.5){};

\node[ellipse,minimum width=2.5cm,minimum height=2.5cm,draw](a)at(9.7,-7.2){};
\node[minimum size=.1pt,label={left:X}](X)at(10.2,-7.2){};
\node[enclosed,label={left,yshift=.01cm:3}](3)at(7.4,-9.2){};
\node[enclosed,label={right,yshift=.01cm:2}](2)at(9.7,-8){};
\node[minimum size=.1pt](a1)at(8.8,-7.9){};
\node[minimum size=.1pt](c1)at(7.1,-7.9){};
\node[ellipse,minimum width=2.2cm,minimum height=1cm,draw](b)at(7.95,-9.2){};
\node[enclosed,label={right,yshift=.01cm:1}](1)at(8.5,-9.2){};
\node[ellipse,minimum width=2.5cm,minimum height=2.5cm,draw](c)at(6.2,-7.2){};
\node[minimum size=.1pt,label={left:Y}](Y)at(6.7,-7.2){};

\draw[black,->-=.5] (3)--(1);
\draw[black,->-=.5] (a1)--(b);
\draw[black,->-=.5] (c1)--(b);
\draw[black,->-=.5] (b) to [out=10,in=-90,looseness=0.9] (2);
\draw[black,->-=.5] (3)--(1);
\draw[black,->-=.5] (a)--(c);

\node[minimum size=.1pt,label={left:$F_6(n)$}](F_6)at(8,-10.5){};

\node[ellipse,minimum width=2.5cm,minimum height=2.5cm,draw](a)at(13.8,-7.2){};
\node[minimum size=.1pt,label={left:Y}](Y)at(14.3,-7.2){};
\node[enclosed,label={left,yshift=.01cm:3}](3)at(11.85,-9.2){};
\node[enclosed,label={left,yshift=.2cm:1}](1)at(11.6,-7.7){};
\node[minimum size=.1pt](a1)at(12.8,-7.7){};
\node[minimum size=.1pt](a2)at(13.1,-8){};
\node[ellipse,minimum width=2.2cm,minimum height=1cm,draw](b)at(12.35,-9.2){};
\node[enclosed,label={right,yshift=.01cm:2}](2)at(12.85,-9.2){};

\draw[black,->-=.5] (1)--(a1);
\draw[black,->-=.5] (a2)--(b);
\draw[black,->-=.5] (b)--(1);
\draw[black,->-=.5] (3)--(2);

\node[minimum size=.1pt,label={left:$F_7(n)$}](F_7)at(13,-10.5){};

\node[ellipse,minimum width=2.5cm,minimum height=2.5cm,draw](a)at(0,-12.5){};
\node[ellipse,minimum width=2.5cm,minimum height=2.5cm,draw](b)at(0,-16){};
\node[minimum size=.1pt,label={left:3A}](3A)at(0.5,-16){};
\node[minimum size=.1pt,label={left:X}](X)at(0.45,-12.5){};

\draw[black,->-=.5] (b)--(a);

\node[minimum size=.1pt,label={left:$F_8(n)$}](F_8)at(0.65,-17.8){};

\node[ellipse,minimum width=2.5cm,minimum height=2.5cm,draw](a)at(3.2,-12.5){};
\node[ellipse,minimum width=2.5cm,minimum height=2.5cm,draw](b)at(3.2,-16){};
\node[minimum size=.1pt,label={left:5A}](5A)at(3.7,-16){};
\node[minimum size=.1pt,label={left:X}](X)at(3.65,-12.5){};

\draw[black,->-=.5] (b)--(a);

\node[minimum size=.1pt,label={left:$F_9(n)$}](F_9)at(3.85,-17.8){};

\node[ellipse,minimum width=2.5cm,minimum height=2.5cm,draw](a)at(6.4,-12.5){};
\node[ellipse,minimum width=2.5cm,minimum height=2.5cm,draw](b)at(6.4,-16){};
\node[minimum size=.1pt,label={left:7A}](7A)at(6.9,-16){};
\node[minimum size=.1pt,label={left:X}](X)at(6.85,-12.5){};

\draw[black,->-=.5] (b)--(a);

\node[minimum size=.1pt,label={left:$F_{10}(n)$}](F_10)at(7.05,-17.8){};

\node[ellipse,minimum width=4cm,minimum height=5.5cm,draw](a)at(10.6,-14.5){};
\node[enclosed,label={left,yshift=.2cm:1}](1)at(10.6,-13){};

\node[ellipse,minimum width=3cm,minimum height=3cm,draw](b)at(10.6,-15.5){};
\node[minimum size=.1pt,label={left:X}](X)at(11.1,-15.5){};

\node[minimum size=.1pt](a1)at(12.4,-13.7){};
\node[minimum size=.1pt](a2)at(12.4,-15.3){};

\node[enclosed,label={right,yshift=.2cm:2}](2)at(13.6,-13.7){};
\node[enclosed,label={right,yshift=.2cm:3}](3)at(13.6,-15.3){};

\draw[black,->-=.5] (b)--(1);
\draw[black,->-=.5] (2)--(3);
\draw[black,->-=.5] (3)--(a2);
\draw[black,->-=.5] (a1)--(2);

\node[minimum size=.1pt,label={left:$F_{11}(n)$}](F_11)at(13,-17.8){};

\end{tikzpicture}
\caption{The infinite families of exceptions -1-}
\label{figure 5}
\end{figure}
\begin{figure}
\begin{tikzpicture}
\tikzset{enclosed/.style={draw,circle,inner sep=2pt,minimum size=4pt,fill=black}}
\tikzset{->-/.style={decoration={
            markings,
            mark=at position #1 with
            {\arrow{>}}},postaction={decorate}}}
            
\node[ellipse,minimum width=4cm,minimum height=5.5cm,draw](a)at(0,0){};
\node[ellipse,minimum width=2.2cm,minimum height=1cm,draw](c)at(0,1.5){};
\node[enclosed,label={right,yshift=.01cm:1}](1)at(0.5,1.5){};
\node[enclosed,label={left,yshift=.01cm:4}](4)at(-0.5,1.5){};

\node[ellipse,minimum width=3cm,minimum height=3cm,draw](b)at(0,-1){};
\node[minimum size=.1pt,label={left:X}](X)at(0.5,-1){};

\node[minimum size=.1pt](a1)at(1.8,0.8){};
\node[minimum size=.1pt](a2)at(1.8,-0.8){};

\node[enclosed,label={right,yshift=.2cm:2}](2)at(3,0.8){};
\node[enclosed,label={right,yshift=.2cm:3}](3)at(3,-0.8){};

\draw[black,->-=.5] (b)--(c);
\draw[black,->-=.5] (4)--(1);
\draw[black,->-=.5] (2)--(3);
\draw[black,->-=.5] (3)--(a2);
\draw[black,->-=.5] (a1)--(2);

\node[minimum size=.1pt,label={left:$F_{12}(n)$}](F_12)at(0.65,-3.5){};

\node[ellipse,minimum width=4cm,minimum height=5.5cm,draw](a)at(8,0){};
\node[enclosed,label={left,yshift=.2cm:1}](1)at(8,1.3){};

\node[ellipse,minimum width=3cm,minimum height=3cm,draw](b)at(8,-1){};
\node[minimum size=.1pt,label={left:X}](X)at(8.5,-1){};

\node[minimum size=.1pt](a1)at(9.7,1.2){};
\node[minimum size=.1pt](a2)at(9.8,-0.8){};

\node[enclosed,label={right,yshift=.2cm:2}](2)at(12,1.2){};

\node[ellipse,minimum width=3cm,minimum height=3cm,draw](c)at(12,-0.8){};
\node[enclosed,label={right,yshift=.2cm:3}](3)at(12.7,-1.3){};
\node[enclosed,label={left,yshift=.2cm:4}](4)at(11.3,-1.3){};
\node[enclosed,label={right,yshift=.2cm:5}](5)at(12,0){};

\draw[black,->-=.5] (b)--(1);
\draw[black,->-=.5] (2)--(c);
\draw[black,->-=.5] (c)--(a2);
\draw[black,->-=.5] (a1)--(2);
\draw[black,->-=.5] (3)--(4);
\draw[black,->-=.5] (4)--(5);
\draw[black,->-=.5] (5)--(3);

\node[minimum size=.1pt,label={left:$F_{13}(n)$}](F_13)at(8.65,-3.5){};

\node[ellipse,minimum width=4cm,minimum height=5.5cm,draw](a)at(5,-7){};
\node[ellipse,minimum width=2.2cm,minimum height=1cm,draw](c)at(5,-5.5){};
\node[enclosed,label={right,yshift=.01cm:1}](1)at(5.5,-5.5){};
\node[enclosed,label={left,yshift=.01cm:6}](6)at(4.5,-5.5){};

\node[ellipse,minimum width=3cm,minimum height=3cm,draw](b)at(5,-8){};
\node[minimum size=.1pt,label={left:X}](X)at(5.5,-8){};

\node[minimum size=.1pt](a1)at(6.7,-5.8){};
\node[minimum size=.1pt](a2)at(6.8,-7.8){};

\node[enclosed,label={right,yshift=.2cm:2}](2)at(9,-5.8){};

\node[ellipse,minimum width=3cm,minimum height=3cm,draw](d)at(9,-7.8){};
\node[enclosed,label={right,yshift=.2cm:3}](3)at(9.7,-8.3){};
\node[enclosed,label={left,yshift=.2cm:4}](4)at(8.3,-8.3){};
\node[enclosed,label={right,yshift=.2cm:5}](5)at(9,-7){};

\draw[black,->-=.5] (b)--(c);
\draw[black,->-=.5] (6)--(1);
\draw[black,->-=.5] (2)--(d);
\draw[black,->-=.5] (d)--(a2);
\draw[black,->-=.5] (a1)--(2);
\draw[black,->-=.5] (3)--(4);
\draw[black,->-=.5] (4)--(5);
\draw[black,->-=.5] (5)--(3);

\node[minimum size=.1pt,label={left:$F_{14}(n)$}](7I)at(5.65,-10.5){};

\end{tikzpicture}
\caption{The infinite families of exceptions -2-}
\label{figure 6}
\end{figure}
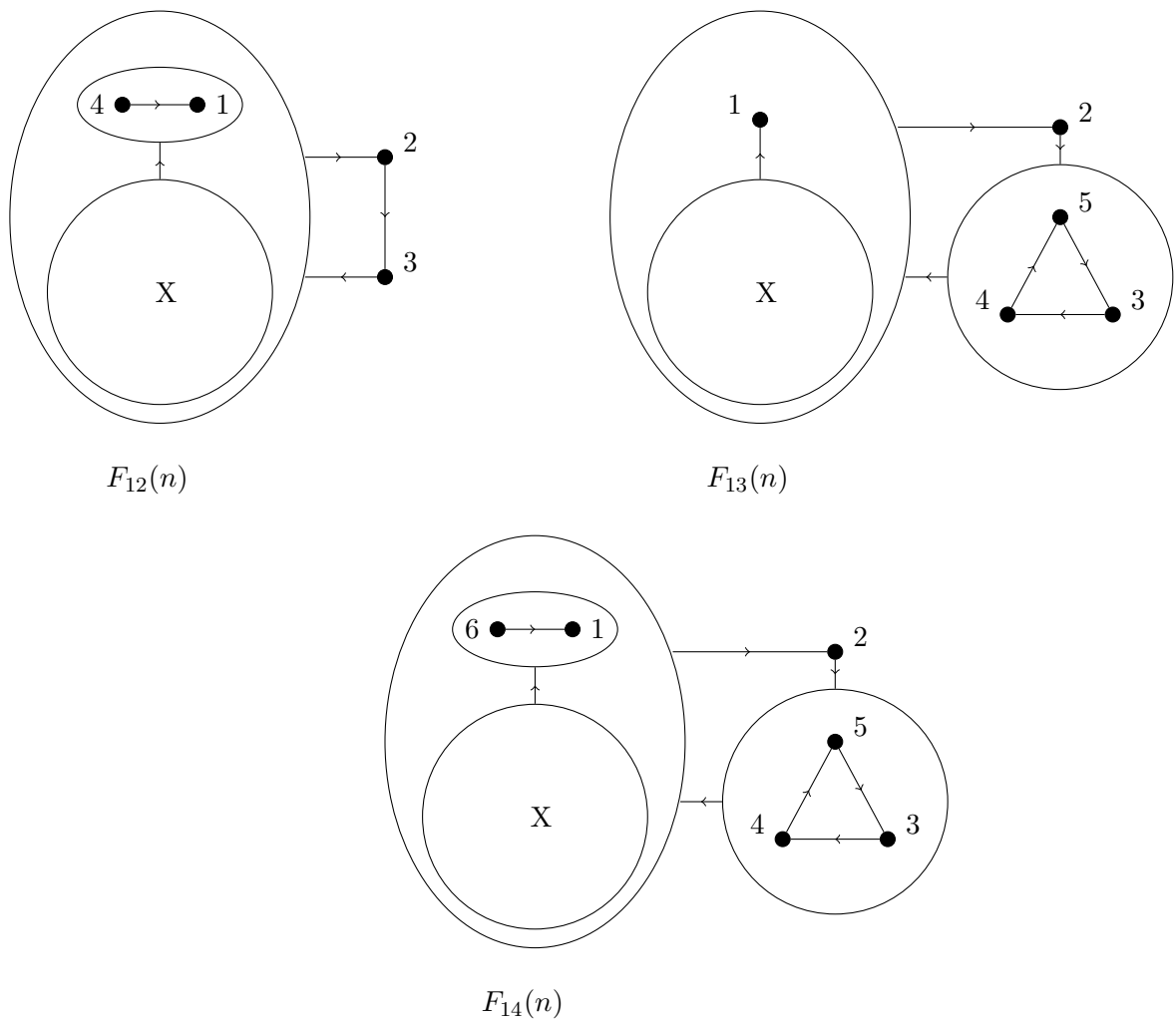
Exception $E_1(n)=(F_1(n),(1,n-2))$; $S=\{1,2,3\}$. Conditions: $|X|\geq 1$. Paths: for any $u\in X,\,P=u132I_{X-u}$.
\\Exception $E_2(n)=(F_2(n),(2,n-3))$; $S=\{3,4\}$. Conditions: $|X|\geq 1$. Paths: $P=1234I_X$, $P=2314I_X$ and for any $u\in X,\,P=u132I_{X-u}$.
\\Exception $E_3(n)=(F_3(n),(1,n-2))$; $S=\{1,3\}$. Conditions: $N^+(3)\neq\{2\}$ and 3 is an ingenerator of $T(X)$ (in particular there exists a Hamiltonian directed inpath $u3R_1$ of $T(X)$ and another $3vR_2$ of $T(X)$). Paths: $P=21u3R_1$ and for any $y\in X\setminus \{3\}$, $P=y132I_{X-y}$.
\\Exception $E_4(n)=(F_4(n),(2,n-3))$; $S=\{1,4\}$. Conditions: $N^+(3)\neq\{2\}$ and 3 is an ingenerator of $T(X)$ (in particular there exists a Hamiltonian directed inpath $u3R_1$ of $T(X)$ and another $3vR_2$ of $T(X)$). Paths: $P=214u3R_1$, $P=3241vR_2$ and for any $y\in X\setminus \{3\}$, $P=u132I_{X-u}$.
\\Exception $E_5(n)=(F_5(n),(1,n-2))$; $S=\{1,2\}$. Conditions: $n\geq 5$, $|Y|\geq 2$ and 2 is an ingenerator of $T(X)$. Paths: for any $x\in X\setminus \{2\}$, $P=x1I_YI_{X-x}$, and for any $y\in Y$, $P=y1I_{Y-y}I_X$.
\\Exception $E_6(n)=(F_6(n),(2,n-3))$; $S=\{1,3\}$. Conditions: $|Y|\geq 2$ and 2 is an ingenerator of $T(X)$. Paths: for any $x\in X$, $P=xz13I_{Y-z}I_{X-x}$ (for a given $z\in Y$) and for any $y\in Y$, $P=y31I_{Y-y}I_X$.
\\Exception $E_7(n)=(F_7(n),(1,1,n-3))$; $S=\{2,3\}$. Conditions: $T(Y)$ is not $3-$cycle and $|Y|\geq 3$. Paths: since $T(Y)$ is not a $3-$cycle, there is a path $Q=-(1,n-5)$ in $T(Y)$ (this is clear if $T(Y)$ is not reducible; if $T(Y)$ is strong, there exists a vertex $w\in Y$ such that $T(Y)-w$ is also strong, such a vertex is certainly an origin of $Q$ since it has an inneighbor that is an origin of a Hamiltonian directed path of $T(Y)-w$), we then have $P=1Q23$ and for any $y\in Y$, $P=y231O_{Y-y}$.
\\Exception $E_8(n)=(F_8(n),(n-4,1,1,1))$; $S=X$. Conditions: $3A$ is the $3-$cycle, its set of vertices is $\{1,2,3\}$, we furthermore need $|X|\geq 2$. Paths: $P=1O_{X-u}2u3$, $P=2O_{X-u}1u3$ and $P=3O_{X-u}2u1$ for a given $u\in X$.
\\Exception $E'_8(n)=(F_8(n),(n-4,2,1))$; $S=X$. Conditions: $3A$ is the $3-$cycle, its set of vertices is $\{1,2,3\}$, we furthermore need $|X|\geq 2$. Paths: $P=1O_{X-u}32u$, $P=2O_{X-u}13u$ and $P=3O_{X-u}21u$ for a given $u\in X$.
\\Exception $E_9(n)=(F_9(n),(n-6,1,1,1,1,1))$; $S=X$. Conditions: $5A$ is the $2-$regular tournament, its set of vertices is $\{1,2,3,4,5\}$, we furthermore need $|X|\geq 2$. Paths: $P=1O_{X-u}2u453$, $P=2O_{X-u}3u514$, $P=3O_{X-u}4u125$, $P=4O_{X-u}5u231$ and $P=5O_{X-u}1u342$ for a given $u\in X$.
\\Exception $E'_9(n)=(F_9(n),(n-6,2,1,1,1))$; $S=X$. Conditions: $5A$ is the $2-$regular tournament, its set of vertices is $\{1,2,3,4,5\}$, we furthermore need $|X|\geq 2$. Paths: $P=1O_{X-u}32u45$, $P=2O_{X-u}43u51$, $P=3O_{X-u}54u12$, $P=4O_{X-u}15u23$ and $P=5O_{X-u}21u34$ for a given $u\in X$.
\\Exception $E_{10}(n)=(F_{10}(n),(n-8,1,1,1,1,1,1,1))$; $S=X$. Conditions: $7A$ is the Paley tournament, its set of vertices is $\{1,2,3,4,5,6,7\}$, we furthermore need $|X|\geq 2$. Paths: $P=1O_{X-u}2u45376$, $P=2O_{X-u}3u56417$, $P=3O_{X-u}4u67521$, $P=4O_{X-u}5u71632$, $P=5O_{X-u}6u12743$, $P=6O_{X-u}7u23154$ and $P=7O_{X-u}1u34265$ for a given $u\in X$.
\\Exception $E'_{10}(n)=(F_{10}(n),(n-8,2,1,1,1,1,1))$; $S=X$. Conditions: $7A$ is the Paley tournament, its set of vertices is $\{1,2,3,4,5,6,7\}$, we furthermore need $|X|\geq 2$. Paths: $P=1O_{X-u}32u4657$, $P=2O_{X-u}43u5761$, $P=3O_{X-u}54u6172$, $P=4O_{X-u}65u7213$, $P=5O_{X-u}76u1324$, $P=6O_{X-u}17u2435$ and $P=7O_{X-u}21u3546$ for a given $u\in X$.
\\Exception $E_{11}(n)=(F_{11}(n),(1,1,n-3))$; $S=\{1,2\}$. Conditions: $|X|\geq 2$. Paths: $P=31O_X2$ and for any $u\in X,\,P=u1O_{X-u}23$.
\\Exception $E_{12}(n)=(F_12(n),(2,1,n-4))$; $S=\{1,4\}$. Conditions: $|X|\geq 2$. Paths: $P=231O_X4$ and for any $u\in X,\,P=u41O_{X-u}23$.
\\Exception $E_{13}(n)=(F_{13}(n),(1,1,n-3))$; $S=\{1,2\}$. Conditions: $|X|\geq 2$. Paths: $P=3425O_X1$, $P=4523O_X1$, $P=5324O_X1$ and for all $u,v\in X,\,P=u1v2345O_{X\setminus \{u,v\}}$.
\\Exception $E_{14}(n)=(F_{14}(n),(2,1,n-4))$; $S=\{1,6\}$. Conditions: $|X|\geq 2$. Paths: for every vertex $u\in X$, $P=24u53O_{X-u}61$, $P=3u16245O_{X-u}$, $P=4u16253O+{X-u}$, $P=5u16234O_{X-u}$ and $P=u61345O_{X-u}2$.\\
They conclude the following:
\begin{theorem}\label{2.1}
Let $T$ be a tournament of order $n$, $P$ an outpath of order $n$ and $x,y$ two distinct vertices of $T$ such that $s^+(x,y)\geq b_1(P)+1$. Then at least one of the following holds:\begin{itemize}
\item[(i)] $x$ or $y$ is an origin of $P$ in $T$.
\item[(ii)] $(T;P)$ is an exception.
\end{itemize}
\end{theorem}
The following remark is due to Bou Hanna \cite{Bou Hanna}:
\begin{remark} \label{2.2}
If $T$ is a tournament of order $4$, $6$ or $8$, there is at most $2$ vertices $x$ such that $T-x$ may be $3A$, $5A$ or $7A$.
\end{remark}
For the sake of our proof, we need to establish the following lemmas:
\begin{lemma} \label{2.3}
Let $P$ be a non-directed outpath of order $n$ and $T$ a tournament of order $n$. If $P$ is not directed and $(T;P)$ is not a Gr\"{u}nbaum's exception then $T$ contains $P$ of origin of indgree at least $1$ unless $(T;P)=(4B;(1,2))$.
\end{lemma}
\begin{proof}
This is clear if $\delta^-(T)\geq 1$. So we can assume that there exists $x\in T$ of indegree zero. We can divide $T-x$ into two sets of vertices $V_1$ and $V_2$ such that $T(V_1)$ and $T(V_2)$ is not $3A$, $5A$ nor $7A$ and $|V_1|=b_1+b_2-1$ unless $(T-x;(b_1,b_2-1)=(3A;(1,1))$. Let $a_1...a_t=(b_1,b_2-1)$ in $T(V_1)$ and $a_{t+2}...a_n=(b_3-1,...,b_s)$ in $T(V_2)$. We have $a_1...a_txa_{t+2}...a_n\equiv P$ and $d^+(a_1)\geq 1$. If $(T-x;(b_1,b_2-1))=(3A;(1,1))$ then $T=4B$ and $P=(1,2)$.    $\square$
\end{proof}
\begin{lemma} \label{2.4}
Let $P$ be a non-directed outpath of order $n\geq 3$ and $T$ a tournament of order $n$. If $(T;P)$ is not the Gr\"{u}nbaum's exceptions then $T$ contains $P$ of origin of outdegree at least $2$ unless $(T;P)=(4B;(1,1,1))$ or $(T;P)=(4B;(2,1))$. In these cases, $1$, $2$ and $3$ are origins of $P$. Moreover, if $P$ is directed and $T$ is distinct from $3A$, there is an origin of $P$ of outdegree at least $2$. Else, every vertex of $3A$ is an origin of $P$.
\end{lemma}
\begin{proof}
Set $X=\{v\in T;d^+(v)\geq 2\}$, we have $|X|\geq n-3$ and $s^+(X)=n$, since otherwise, there is a vertex $x\in T-X$ such that $x\rightarrow X$, a contradiction. If $(T;P)$ is an exception, one can check the result. So we may suppose that $(T;P)$ is not an exception, if $|X|\geq 2$, the result follows by theorem \ref{2.1}. Otherwise, $|X|=1$, if $n=3$ then $T$ is transitive while if $n=4$ then $T=4B$, the result follows.   $\square$
\end{proof}
\begin{lemma} \label{2.5}
Let $P$ be an outpath of order $n$ having at least three blocks and $T$ a tournament of order $n$. If $(T;P)$ is not a Gr\"{u}nbaum's exception, then $T$ contains $a_1...a_n\equiv P$ such that $d^+(a_2)\geq 1$ unless $(T;P)=(\overline{4B};(1,1,1))$ \end{lemma}
\begin{proof}
This is trivial if $\delta^-(T)\geq 1$. So we can assume that $T$ contains a maximal vertex $x$. Set $P=(b_1,...,b_s)$. We can divide $T-x$ into two sets of vertices $V_1$ and $V_2$ such that $|V_1|=b_1+b_2+b_3$ and $T(V_1)$ and $T(V_2)$ are not $3A$, $5A$ nor $7A$, unless $(T(V_1);(b_1,b_2,b_3-1))=(3A;(1,1))$; in this case $(T;P)=(\overline{4B};(1,1,1))$. Let $a_1...a_t=(b_1,...,b_3-1)$ in $T(V_1)$ and $a_{t+2}...a_n=-(b_4-1,...b_s)$. We have $a_1...a_txa_{t+2}...a_n\equiv P$ and $d^+(a_2)\ 1$.  $\square$
\end{proof}
\begin{lemma} \label{2.6}
Let $T$ be a tournament of order $n=2k\geq 4$. $T$ contains an antidirected Hamiltonian outpath $x_1...x_n$ such that $d^-(x_1)\geq 2$ and $d^+(x_n)\geq 1$ if and only if $T\notin \{4A,\overline{4B}\}$.
\end{lemma}
\begin{proof}
Set $X=\{v\in T;d^-(v)\geq 2\}$, we have $|X|\geq n-3$. If $(T;P)$ is an exception one can check the result. Otherwise, by theorem \ref{2.1}, $T$ contains an antidirected Hamiltonian outpath $x_1...x_n$ such that $x_1\in X$. If $d^+(x_n)=0$, the path $x_1x_nx_{n-1}...x_2$ gives the result.   $\square$
\end{proof}
\begin{remark}\label{2.7}
Let $T$ be a tournament on $n$ vertices and $P$ be a path of order $n$ with at least three blocks, such that $b_1(P)=1$. If $T$ has a minimal vertex then $T$ has at least $|T|-2$ origins of $P$.
\end{remark}
\begin{proof}
This is a corollary of Theorem \ref{2.1}.  $\square$
\end{proof}
\begin{lemma}\label{2.8}
Let $(T;P)$ be an exception, $Or(P;T)\cap Or(P;\overline{T})\neq \emptyset$ if and only if $(T;P)$ is not one of the exceptions $0,1,4,7,18,19,22,33,E_1(4)$ and $E_1(6)$ (with $X=3A$).
\end{lemma}
Inspired by the concept of exceptions, which are crucial in finding paths in tournaments, the following concept will be useful in finding cycles.

\begin{definition}
Let $T$ be a tournament of order $n$, $P$ a non-directed outpath of order $n$ and $x$ be a vertex of $T$ such that $d^+(x)\geq 2$. The couple $(T-x;^*P)$ is said to be a biexception with respect to $x$ if $x$ is not an origin of $P$ in $T$, $(T-x;^*P)$ is an exception and there exists a vertex $y\in T$ such that $x$ is the end of $P$ whenever $y$ is its origin in $T$.\\
The notation of a biexception is Exc $(k,i)$, where $(T-x;^*P)$ is the biexception $k$, $i=0$ when $^*P$ is an outpath and $i=1$ when $^*P$ is an inpath. The set of vertices of $T-x$ which are not origins of $P$ in $T$ when $x$ is not its end is denoted by $S$. The biexceptions are the following:\\
Exc $(0,0)(1)$: $P=(2,1)$, $N^+(x)=\{1,2,3\}$, $S=\{1,2,3\}$.\\
Exc $(0,0)(2)$: $P=(2,1)$, $N^+(x)=\{1,2\}$, $S=\{1,2\}$.\\
Exc $(1,0)(1)$: $P=(2,1,1)$, $N^+(x)=\{1,3\}$, $S=\{1,2,3\}$.\\
Exc $(1,0)(2)$: $P=(2,1,1)$, $N^+(x)=\{2,3\}$, $S=\{1,2\}$.\\
Exc $(1,0)(3)$: $P=(2,1,1)$, $N^+(x)=\{1,2,3\}$, $S=\{1,2\}$.\\
Exc $(2,0)$: $P=(2,2)$, $N^+(x)=\{3,4\}$, $S=\{3\}$.\\
Exc $(4,0)(1)$: $P=(2,1,1,1)$, $N^+(x)=\{1,2\}$, $S=\{5\}$.\\
Exc $(4,0)(2)$: $P=(2,1,1,1)$, $N^+(x)=\{1,2,3\}$, $S=\{1,2\}$.\\
Exc $(6,0)$: $P=(2,1,2)$, $N^+(x)=\{4,5\}$, $S=\{4\}$.\\
Exc $(7,0)(1)$: $P=(3,1,1)$, $N^+(x)=\{1,2,4\}$, $S=\{1\}$.\\
Exc $(7,0)(2)$: $P=(3,1,1)$, $N^+(x)=\{1,4\}$, $S=\{3\}$.\\
Exc $(9,0)(1)$: $P=(2,1,1,1)$, $N^+(x)=\{2,4,5\}$, $S=\{2\}$.\\
Exc $(9,0)(2)$: $P=(2,1,1,1)$, $N^+(x)=\{4,5\}$, $S=\{2,5\}$.\\
Exc $(15,0)(1)$: $P=(2,1,2,1)$, $N^+(x)=\{1,2,6\}$, $S=\{1\}$.\\
Exc $(15,0)(2)$: $P=(2,1,2,1)$, $N^+(x)=\{1,6\}$, $S=\{3\}$.\\
Exc $(16,0)(1)$: $P=(2,2,1,1)$, $N^+(x)=\{4,5\}$, $S=\{6\}$.\\
Exc $(16,0)(2)$: $P=(2,2,1,1)$, $N^+(x)=\{4,5,6\}$, $S=\{6\}$.\\
Exc $(26,0)(1)$: $P=(2,1,3)$, $N^+(x)=\{4,5,6\}$, $S=\{4,5\}$.\\
Exc $(26,0)(2)$: $P=(2,1,3)$, $N^+(x)=\{4,6\}$, $S=\{4,5\}$.\\
Exc $(26,0)(3)$: $P=(2,1,3)$, $N^+(x)=\{5,6\}$, $S=\{4,5\}$.\\
Exc $(27,0)$: $P=(2,3,1)$, $N^+(x)=\{4,6\}$, $S=\{4\}$.\\
Exc $(31,0)$: $P=(2,2,2)$, $N^+(x)=\{3,4\}$, $S=\{4\}$.\\
Exc $(33,0)(1)$: $P=(2,1,1,1,1,1)$, $N^+(x)=\{1,2,3,5\}$, $S=\{1\}$.\\
Exc $(33,0)(2)$: $P=(2,1,1,1,1,1)$, $N^-(x)=\{1,2,4,7\}$, $S=\{7\}$.\\
Exc $(48,0)$: $P=(2,1,2,1,1)$, $N^+(x)=\{1,2\}$, $S=\{2\}$.\\
Exc $(E_1(n-1),0)$: $P=(2,n-3)$, $N^+(x)=\{1,2\}$, $S=\{1\}$.\\
Exc $(E_3(n-1),0)$: $P=(2,n-3)$, $N^+(x)=\{1,3\}$, $S=\{1\}$.\\
Exc $(E_5(n-1),0)$: $P=(2,n-3)$, $N^+(x)=\{1,2\}$, $S=\{1\}$.\\
Exc $(E_8(n-1),0)$: $P=(n-4,1,1,1)$, $N^+(x)\subseteq X$, $S=X$.\\
Exc $(E'_8(n-1),0)$: $P=(n-4,2,1)$, $N^+(x)\subseteq X$, $S=X$.\\
Exc $(E_9(n-1),0)$: $P=(n-6,1,1,1,1,1)$, $N^+(x)\subseteq X$, $S=X$.\\
Exc $(E'_9(n-1),0)$: $P=(n-6,2,1,1,1)$, $N^+(x)\subseteq X$, $S=X$.\\
Exc $(E_{10}(n-1),0)$: $P=(n-8,1,1,1,1,1,1,1)$, $N^+(x)\subseteq X$, $S=X$.\\
Exc $(E'_{10}(n-1),0)$: $P=(n-8,2,1,1,1,1,1)$, $N^+(x)\subseteq X$, $S=X$.\\
Exc $(E_{11}(n-1),0)$: $P=(2,1,n-4)$, $N^+(x)=\{1,2\}$, $S=\{1\}$.\\
Exc $(E_{13}(n-1),0)$: $P=(2,1,n-4)$, $N^+(x)=\{1,2\}$, $S=\{1\}$.\\
Exc $(0,1)$: $P=-(1,1,1)$, $N^-(x)=\{1,2\}$, $S=\{1,2\}$.\\
Exc $(1,1)(1)$: $P=-(1,1,1,1)$, $N^-(x)=\{1,2\}$, $S=\{1,3\}$.\\
Exc $(1,1)(2)$: $P=-(1,1,1,1)$, $N^-(x)=\{1,3\}$, $S=\{1\}$.\\
Exc $(1,1)(3)$: $P=-(1,1,1,1)$, $N^-(x)=\{2,3\}$, $S=\{1,2,3,4\}$.\\
Exc $(1,1)(4)$: $P=-(1,1,1,1)$, $N^-(x)=\{1,2,3\}$, $S=\{1,3\}$.\\
Exc $(2,1)$: $P=-(1,1,2)$, $N^-(x)=\{3,4\}$, $S=\{4\}$.\\
Exc $(3,1)$: $P=-(1,2,1)$, $N^-(x)=\{1,2\}$, $S=\{1,4\}$.\\
Exc $(4,1)(1)$: $P=-(1,1,1,1,1)$, $N^-(x)=\{1,2\}$, $S=\{3\}$.\\
Exc $(4,1)(2)$: $P=-(1,1,1,1,1)$, $N^-(x)=\{1,2,3\}$, $S=\{3\}$.\\
Exc $(4,1)(3)$: $P=-(1,1,1,1,1)$, $N^-(x)=\{1,2,4\}$, $S=\{4\}$.\\
Exc $(5,1)$: $P=-(1,2,1,1)$, $N^-(x)=\{1,2,3\}$, $S=\{5\}$.\\
Exc $(6,1)$: $P=-(1,1,1,2)$, $N^-(x)=\{4,5\}$, $S=\{4,5\}$.\\
Exc $(7,1)(1)$: $P=-(1,2,1,1)$, $N^-(x)=\{1,2,3,4\}$, $S=\{4,5\}$.\\
Exc $(7,1)(2)$: $P=-(1,2,1,1)$, $N^-(x)=\{1,2,3\}$, $S=\{4,5\}$.\\
Exc $(7,1)(3)$: $P=-(1,2,1,1)$, $N^-(x)=\{1,2,4\}$, $S=\{1,2\}$.\\
Exc $(8,1)$: $P=-(1,1,1,1,1)$, $N^-(x)=\{2,5\}$, $S=\{1,5\}$.\\
Exc $(9,1)(1)$: $P=-(1,1,1,1,1)$, $N^-(x)=\{2,4,5\}$, $S=\{1,2,4\}$.\\
Exc $(9,1)(2)$: $P=-(1,1,1,1,1)$, $N^-(x)=\{2,5\}$, $S=\{1\}$.\\
Exc $(9,1)(3)$: $P=-(1,1,1,1,1)$, $N^-(x)=\{2,4\}$, $S=\{4\}$.\\
Exc $(9,1)(4)$: $P=-(1,1,1,1,1)$, $N^-(x)=\{4,5\}$, $S=\{4\}$.\\
Exc $(10,1)$: $P=-(1,1,2,1)$, $N^-(x)=\{3,5\}$, $S=\{3,4,5\}$.\\
Exc $(11,1)$: $P=-(1,2,2)$, $N^-(x)=\{1,2\}$, $S=\{1,4\}$.\\
Exc $(12,1)$: $P=-(1,1,1,2)$, $N^-(x)=\{1,2\}$, $S=\{1,4\}$.\\
Exc $(15,1)(1)$: $P=-(1,1,1,2,1)$, $N^-(x)=\{1,2,3,6\}$, $S=\{6\}$.\\
Exc $(15,1)(2)$: $P=-(1,1,1,2,1)$, $N^-(x)=\{1,2,3\}$, $S=\{6\}$.\\
Exc $(16,1)(1)$: $P=-(1,1,2,1,1)$, $N^-(x)=\{4,5,6\}$, $S=\{4,5\}$.\\
Exc $(16,1)(2)$: $P=-(1,1,2,1,1)$, $N^-(x)=\{4,5\}$, $S=\{4\}$.\\
Exc $(16,1)(3)$: $P=-(1,1,2,1,1)$, $N^-(x)=\{4,6\}$, $S=\{4,5\}$.\\
Exc $(16,1)(4)$: $P=-(1,1,2,1,1)$, $N^-(x)=\{5,6\}$, $S=\{4,5\}$.\\
Exc $(22,1)$: $P=-(1,1,1,1,1,1)$, $N^-(x)=\{1,2,3\}$, $S=\{1,2,3,4,5,6\}$.\\
Exc $(24,1)$: $P=-(1,1,1,1,1,1)$, $N^-(x)=\{1,2,3,4\}$, $S=\{4,5\}$.\\
Exc $(26,1)(1)$: $P=-(1,1,1,3)$, $N^-(x)=\{4,5,6\}$, $S=\{6\}$.\\
Exc $(26,1)(2)$: $P=-(1,1,1,3)$, $N^-(x)=\{4,5\}$, $S=\{6\}$.\\
Exc $(26,1)(3)$: $P=-(1,1,1,3)$, $N^-(x)=\{4,6\}$, $S=\{6\}$.\\
Exc $(32,1)$: $P=-(1,1,2,1,1)$, $N^-(x)=\{5,6\}$, $S=\{4,5\}$.\\
Exc $(33,1)(1)$: $P=-(1,1,1,1,1,1,1)$, $N^-(x)=\{1,2,5,7\}$, $S=\{2\}$.\\
Exc $(33,1)(2)$: $P=-(1,1,1,1,1,1,1)$, $N^-(x)=\{1,2,3,6\}$, $S=\{3\}$.\\
Exc $(E_1(n-1),1))$: $P=-(1,1,n-3)$, $N^-(x)=\{1,2,3\}$, $S=X$.\\
Exc $(E_2(5),1)$: $P=-(1,2,2)$, $N^-(x)=\{3,4\}$, $S=\{5\}$.\\
Exc $(E_5(5),1)(1)$: $P=-(1,1,3)$, $N^-(x)=\{1,2\}$ and $T(Y)=3A$, $S=\{1,2\}$.\\
Exc $(E_5(n-1),1)(2)$: $P=-(1,1,n-3)$, $N^-(x)\{1,2\}$, $T(Y)\neq 3A$ and $X=\{2\}$, $S=\{2\}$.\\
Exc $(E_6(5),1)$: $P=-(1,2,2)$, $N^-(x)=\{1,3\}$ and $Y=\{u\rightarrow v\}$, $S=\{3,u\}$.\\
Exc $(E_7(6),1)$: $P=-(1,1,1,3)$, $N^-(x)=\{2,3\}$ and $T(Y)=3A$, $S=\{1\}$.\\
Exc $(E_8(5),1)(1)$: $P=-(1,1,1,1,1)$, $N^-(x)=X$, $S=X$.\\
Exc $(E_8(6),1)(2)$: $P=-(1,2,1,1,1)$, $T(X)$ is a transitive triangle of minimal vertex $u$ and $N^-(x)=X$, $S=\{u\}$.
\end{definition}
\begin{lemma}\label{2.10}
Let $T$ be a tournament on $n$ vertices, $P$ an outpath of order $n$ and $x$ a vertex of $T$ of outdegree at least two. If $(T-x;^*P)$ is a finite exception, $^*P$ is an outpath and $x$ is not an origin of $P$ in $T$ then any vertex of $T$ distinct from $x$ is an origin of $P$ in $T$ with end distinct from $x$ unless $(T;P)$ is a biexception with respect to $x$.
\end{lemma}
\begin{proof}
We enumerate the finite exceptions and extend them in all possible ways by a vertex $x$ of outdegree at least two. To shorten the proof we use the following notation: Exc $k$: $P$ $[N^+(x);S;P_1;...;P_t]$, where $P_i$ are the paths of origin in $T\setminus (\{x\}\cup S)$ and of end distinct from $x$.\\
Exc $0$: $(2,1)\,[\{1,2,3\};\{1,2,3\}]\,[\{1,2\};\{1,2\};3x21]$ The other possible outneighbors of $x$ are equivalent\\
Exc $1$: $(2,1,1)\,[\{1,2\};\emptyset;13x42;23x41;3x142;4x213]\,[\{1,3\};\{1,2,3\};4x312]$\\$[\{2,3\};\{1,2\};34x12;413x2]\,[\{1,2,3\};\{1,2\};341x2;4x312]$\\
Exc $2$: $(2,2)\,[\{3,4\};\{3\};413x2;134x2;234x1]$\\
Exc $3$: $(3,1)\,[\{1,2\};\emptyset;13x24;23x14;3x124;4x132]$\\
Exc $4$: $(2,1,1,1)\,[\{1,2\};\{5\};1243x5;23x514;3x2514;4125x3]$\\$[\{1,3\};\emptyset;1342x5;2x1453;3514x2;4x3125;513x42]$\\$[\{1,2,3\};\{1,2\};34x521;4x3125;523x14]$\\$[\{1,2,4\};\emptyset;134x25;234x15;341x25;451x32;512x43]$\\$[\{1,2,3,4\};\emptyset;134x25;234x15;341x25;451x32;512x43]$\\$[\{1,2,3,4,5\};\emptyset;134x25;245x31;351x42;412x53;523x14]$\\ The other possible outneighbors of $x$ are equivalent.\\
Exc $5$: $(3,1,1)\,[\{1,2,3\};\emptyset;1452x3;2451x3;3451x2;4512x3;5x1423]$\\$[\{1,2\};\emptyset;14x352;24x153;3451x2;45x213;5x1423]$\\$[\{1,3\};\emptyset;14x352;24x153;34x152;4x1352;5x1423]$\\$[\{2,3\};\emptyset;14x352;24x351;34x251;4x2351;5x2413]$\\
Exc $6$: $(2,1,2)\,[\{4,5\};\{4\};1243x5;2341x5;3142x5;51x234]$\\
Exc $7$: $(3,1,1)\,[\{1,2,3,4\};\emptyset;1453x2;2451x3;3451x2;4512x3;5x1423]$\\$[\{1,2,3\};\emptyset;1453x2;2451x3;3451x2;4512x3;5x1423]$\\$[\{1,2,4\};\{1\};245x31;3451x2;4531x2;5x1423]$\\$[\{1,2\};\emptyset;14x253;24x153;34x251;4x2351;5x1423]$\\$[\{1,4\};\{3\};145x23;245x31;452x31;5x1423]$\\ The other possible outneighbors of $x$ are equivalent.\\
Exc $8$: $(2,1,1,1)\,[\{2,5\};\emptyset;13x425;2351x4;342x51;4x5132;5124x3]$\\
Exc $9$: $(2,1,1,1)\,[\{2,4,5\};\{2\};135x24;345x21;412x53;512x43]$\\$[\{4,5\};\{2,5\};12x354;3x5421;4521x3]$\\
Exc $10$: $(2,2,1)\,[\{3,5\};\emptyset;1x3245;2x3145;345x12;4125x3;523x41]$\\
Exc $11$: $(3,2)\,[\{1,2\};\emptyset;13x254;23x154;3x1254;4x1253;5x1243]$\\
Exc $12$: $(2,1,2)\,[\{1,2\};\emptyset;13x452;23x451;3x1452;4x1523;5x1423]$\\
Exc $13$: $(4,1,1)\,[\{3,4\};\emptyset;15623x4;245x361;346x152;46321x5;5x43162;6x43512]$\\
Exc $14$: $(3,1,1,1)\,[\{3,4\};\emptyset;14635x2;24635x1;346x152;46321x5;5x34126;6x34152]$\\
Exc $15$: $(2,1,2,1)\,[\{1,2,3,6\};\emptyset;123x654;231x654;312x654;416x235;516x234;6x14532]$\\$[\{1,2,6\};\{1\};23x6415;312x654;4x61235;5x61234;6435x21]$\\$[\{1,2,3\};\emptyset;12345x6;23145x6;31245x6;4x15263;5x14263;6x14532]$\\$[\{1,2\};\emptyset;12345x6;23145x6;31245x6;4x15263;5x14263;6x14532]$\\$[\{1,6\};\{3\}16542x3;2x13654;4x15362;5x14362;6524x13]$\\ The other possible outneighbors of $x$ are equivalent.\\
Exc $16$: $(2,2,1,1)\,[\{4,5\};\{6\};12x3456;23x1456;31x2456;45132x6;51246x3]\,[\{4,6\};\emptyset;$\\$126x534;236x514;316x524;45x2163;5x61432;65x2413]$\\$[\{5,6\};\emptyset;165x432;265x413;365x421;4x61532;516x234;64x1532]$\\$[\{4,5,6\};\{6\};1234x56;2314x56;3124x56;4125x63;5124x63\}]$\\
Exc $17$: $(3,1,1,1)\,[\{2,4,6\};\emptyset;1234x65;24563x1;3456x21;46125x3;5612x43;62341x5]$\\$[\{2,4\};\emptyset;12465x3;24563x1;34621x5;45631x2;56243x1;6123x54]$\\ The other possible outneighbors of $x$ are equivalent.\\
Exc $18$: $(2,2,2)\,[\{2,4,6\};\emptyset;13x5462;234x561;35x1624;456x123;51x3246;612x345]$\\$[\{2,4\};\emptyset;13x5462;234x561;35x1624;45631x2;51x3246;612x345]$\\ The other possible outneighbors of $x$ are equivalent.\\
Exc $19$: $(2,1,1,2)\,[\{2,4,6\};\emptyset;124x653;23x5641;346x215;45x1263;562x431;61x3425]$\\$[\{2,4\};\emptyset;1243x65;23x5641;36214x5;45x1263;562x431;6x41253]$\\ The other possible outneighbors of $x$ are equivalent.\\
Exc $20$: $(2,1,1,1,1)\,[\{1,2\};\emptyset;13x4652;23x4651;34x5261;45x6213;56x3412;63x5214]$\\
Exc $21$: $(3,1,1,1)\,[\{1,2\};\emptyset;13526x4;23516x4;345x621;463x521;5x14236;634x521]$\\
Exc $22$: $(2,1,1,1,1)\,[\{1,2\};\emptyset;1236x54;2314x65;3124x65;4125x63;5264x31;6314x52]$\\$[\{1,2,3\};\emptyset;1236x54;2314x65;3125x64;4125x63;5316x42;6124x53]$\\ The other possible outneighbors are equivalent.\\
Exc $23$: $(2,1,1,1,1)\,[\{4,6\};\emptyset;1423x56;2145x36;3612x54;4213x56;5623x146145x32]$\\
Exc $24$: $(2,1,1,1,1)\,[\{1,2,3,4\};\emptyset;1245x63;2345x61;3145x62;461x325;5x41263;6x41523]$\\$[\{1,2,4\};\emptyset;1245x63;2345x61;354x162;46x3512;5x41263;6x41523]$\\$[\{1,2,3\};\emptyset;1245x63;2345x61;3145x62;461x325;561x324;6x13425]$\\$[\{1,2\};\emptyset;1245x63;2345x61;3145x62;462x135;562x134;6x13425]$\\$[\{1,4\};\emptyset;1245x63;2345x61;3145x62;46x3512;5613x24;6x13425]$\\
Exc $25$: $(2,1,1,2)\,[\{4,5\};\emptyset;1243x65;2341x65;3142x65;46135x2;56134x2;6x41532]$\\
Exc $26$: $(2,1,3)\,[\{4,5,6\};\{4,5\};12x3456;23x1456;31x2456;61x2345]$\\$[\{4,5\};\emptyset;12x3456;23x1456;31x2456;46x1235;56x1234;61x2345]$\\$[\{4,6\};\{4,5\};12x3456;23x1456;31x2456;61x2345]$\\$[\{5,6\};\{4,5\};12x3456;23x1456;31x2456;61x2345]$\\
Exc $27$: $(2,3,1)\,[\{4,6\};\{4\};1x64523;2x64531;3x64512;5x43261;614x523]$\\
Exc $28$: $(3,1,2)\,[\{4,5\};\emptyset;1234x56;2314x56;3124x56;46x5123;56x4123;6x54123]$\\
Exc $29$: $(2,1,1,1,1)\,[\{4,6\};\emptyset;154x632;254x631;3645x21;4312x56;5642x31;6213x54]$\\
Exc $30$: $(2,1,1,1,1)\,[\{4,6\};\emptyset;1642x53;2641x53;3645x21;43x1625;536x421;6435x21]$\\
Exc $31$: $(2,2,2)\,[\{3,4\};\{4\};1x34256;2x34156;314x256;543x612;643x125]$\\
Exc $32$: $(2,2,1,1)\,[\{5,6\};\emptyset;1x62435;2x61435;3x61425;416x235;536x421;645x132]$\\
Exc $33$: $(2,1,1,1,1,1)$. The Paley tournament is arc transitive, so without loss of generality we may suppose that $\{1,2\}\subseteq N^+(x)$. We have $341x2576\equiv 461x2543\equiv 561x2743\equiv 63541x27\equiv P$. Now let us prove that $2$ is an origin of $P$: if $x\rightarrow 5$ then $245x1673\equiv P$, so we assume that $5\rightarrow x$. If $x\rightarrow 6$ then $241x6573$. Thus $2$ is an origin of $P$ in any case. In a similar way, it is easy to prove that $1$ is not an origin of $P$ if and only if $N^+(x)=\{1,2,3,5\}$ and that $7$ is not an origin of $P$ if and only if $N^+(x)=\{1,2,4,7\}$.\\
Exc $34$: $(2,1,2,1,1)\,[\{1,2,3\};\emptyset;17x45362;27x45163;37x45162;4536x217;5634x217;$\\$6435x217;7x145263][\{1,2\};\emptyset;17x45362;27x45163;37x45162;4536x217;5634x217;$\\$6435x217;7x145263][\{1,3\};\emptyset;17x45362;27x45163;37x45163;4526x137;5624x137;$\\$6425x137;7x145263][\{2,3\};\emptyset;17x45362;27x45163;37x45162;4536x217;5634x217;$\\$6435x217;7x243156]$\\
Exc $35$: $(3,1,3)\,[\{1,2,3\};\emptyset;175x6432;275x6431;375x6412;4x123756;5x123764;$\\$6x123745;7x124563][\{1,2\};\emptyset;175x6432;275x6431;375x6412;4x123756;5x123764;$\\$6x123745;7x124563]$\\$[\{1,3\};\emptyset;175x6432;275x6431;375x6412;4x123756;5x123764;$\\$6x123745;7x124563][\{2,3\};\emptyset;175x6432;275x6431;375x6412;45x16327;56x14327;$\\$64x15327;7x324561]$\\
Exc $36$: $(3,3,1)\,[\{1,2\};\emptyset;17x24653;27x14653;37x14652;4x172356;5x172364;6x172345;$\\$7x124653]$\\
Exc $37$: $(2,1,1,1,1,1)\,[\{4,5,6\};\emptyset;174x6532;274x6513;374x6521;451632x7;561432x7;$\\$641532x7;7x461532][\{4,5\};\emptyset;175x4632;275x4613;375x4621;451632x7;561432x7;$\\$641532x7;7x461532]$ The other possible outneighbors of $x$ are equivalent.\\
Exc $38$: $(2,1,2,1,1)\,[\{1,2,3\};\emptyset;17x45362;27x45163;37x45162;45x61723;56x41723;$\\$64x51723;741x2356][\{1,2\};\emptyset;17x45362;27x45163;37x45162;45x61723;56x41723;$\\$64x51723;741x2356]$ The other possible outneighbors are equivalent.\\
Exc $39$: $(3,1,3)\,[\{1,2,3\};\emptyset;175x6423;275x6431;375x6412;4563x127;5643x127;$\\$6453x127;745x6123][\{1,3\};\emptyset;175x6423;275x6431;375x6412;4563x127;5643x127;$\\$6453x127;745x6123]$ The other possible outneighbors are equivalent.\\
Exc $40$: $(2,1,1,2,1)\,[\{1,2\};\emptyset;13x45267;23x45167;34512x67;45312x67;5341x67;$\\$6x215437;7x215436]$\\
Exc $41$: $(2,1,1,3)\,[\{6,7\};\emptyset;1x734265;2x734165;3x745216;4x753216;5x734216;$\\$61x34275;71x34265]$\\
Exc $42$: $(3,2,1,1)\,[\{6,7\};\emptyset;1x673254;2x673154;3x674251;4x675231;5x673241;$\\$613754x2;716354x2]$\\
Exc $43$: $(2,1,2,1,1)\,[\{2,7\};\emptyset;137x2645;26x13745;3x761245;4x761253;5x761234]$\\
Exc $44$: $(2,1,1,3)\,[\{6,7\};\emptyset;136x7254;216x7354;326x7154;457x6312;516x7234;$\\$64127x35;74126x35]$\\
Exc $45$: $(3,1,2,1)\,[\{1,7\};\emptyset;165x2374;2657x134;3657x124;4537x162;5x374162;$\\$6537x124;765x2314]$\\
Exc $46$: $(3,2,2)\,[\{4,7\};\emptyset;1x746532;2x746531;3x742156;431265x7;5x741632;$\\$6x742153;731265x4]$\\
Exc $47$: $(2,1,2,1,1)\,[\{4,5\};\emptyset;1x456237;2x456137;3x456127;46x32517;56x32417;$\\$6x427531;76x32514]$\\$[\{4,5,7\};\emptyset;1x456237;2x456137;3x456127;46x32517;56x32417;$\\$6x427531;76x32514]$ The other possible outneighbors are equivalent.\\
Exc $48$: $(2,1,2,1,1)\,[\{1,2\};\{2\};13x45267;32145x67;4x235167;5x234167;$\\$6x234157;7x234156]$\\
Exc $49$: $(2,1,1,1,1,1,1)\,[\{1,2\};\emptyset;134265x78;234165x78;351x27846;
452x18367;5x1784623;$\\$671x28534;781x26453;851x24376]$\\
Exc $50$: $(3,1,1,1,1,1)\,[\{1,2\};\emptyset;1345827x6;2345817x6;3x1256487;4571x2836;5712836x4;$\\$68127x453;7812564x3;8512743x6]$\\
Exc $51$: $(3,1,2,1,1)\,[\{2,8\};\emptyset;13475x862;21x834756;3x8245167;4x8235761;5x8234761;$\\$6x8234751;738245x61;81x234756]$  $\square$
\end{proof}
\begin{lemma}\label{2.11}
Let $T$ be a tournament on $n$ vertices, $P$ an inpath of order $n$ and $x$ a vertex of $T$ of indegree at least two. If $(T-x;^*P)$ is a finite exception, $^*P$ is an outpath and $x$ is not an origin of $P$ in $T$ then any vertex of $T$ distinct from $x$ is an origin of $P$ in $T$ with end distinct from $x$ unless $(T;P)$ is a biexception with respect to $x$.
\end{lemma}
\begin{proof}
We use the same notations of the previous lemma:\\
Exc $0$: $-(1,1,1)\,[\{1,2,3\};\emptyset;13x2;21x3;32x1]$\\$[\{1,2\};\{1,2\};32x1]$\\ The other possible inneighbors of $x$ are equivalent.\\
Exc $1$: $-(1,1,1,1)\,[\{1,2\};\{1,3\};213x4;4x312]$\\$[\{1,3\};\{1\};21x34;312x4;4x213]$\\$[\{2,3\};\{1,2,3,4\}]$\\$[\{1,2,3\};\{1,3\};21x34;43x12]$\\
Exc $2$: $-(1,1,2)\,[\{3,4\};\{4\};142x3;241x3;312x4]$\\
Exc $3$: $-(1,2,1)\,[\{1,2\};\{1,4\};21x43;3x421]$\\
Exc $4$: $-(1,1,1,1,1)\,[\{1,2\};\{3\};145x32;213x54;4x3125;5342x1]$\\$[\{1,3\};\emptyset;152x43;21x354;312x54;435x21;534x21]$\\$[\{1,2,3\};\{3\};1453x2;21x354;43x125;5342x1]$\\$[\{1,2,4\};\{4\};145x32;21x354;32x415]$\\$[\{1,2,3,4\};\emptyset;1453x2;2514x3;32x415;43x125;54x231]$\\$[\{1,2,3,4,5\};\emptyset;1453x2;2514x3;3125x4;4231x5;5342x1]$\\
Exc $5$: $-(1,2,1,1)\,[\{1,2,3\};\{5\};152x34;251x34;351x24;4x5312]$\\$[\{1,2\};\emptyset;1523x4;2513x4;3x5214;4x5213;5x3412]$\\$[\{1,3\};\emptyset;1534x2;2134x5;3512x4;4x5312;5x2413]$\\$[\{2,3\};\emptyset;1x5324;2534x1;3524x1;4x5312;5x1423]$\\
Exc $6$: $-(1,1,1,2)\,[\{4,5\};\{4,5\};1532x4;2513x4;3521x4]$\\
Exc $7$: $-(1,2,1,1)\,[\{1,2,3,4\};\{4,5\};152x34;251x34;351x24]$\\$[\{1,2,3\};\{4,5\};152x34;251x34;351x24]$\\$[\{1,2,4\};\{1,2\};3x5214;4123x5;5x3412]$\\$[\{1,2\};\emptyset;1523x4;2145x3;351x24;4x3152;5x3412]$\\$[\{1,4\};\emptyset;1345x2;2x5134;3x5214;4123x5;5x2431]$\\ The other possible inneighbors of $x$ are equivalent.\\
Exc $8$: $-(1,1,1,1,1)\,[\{2,5\};\{1,5\};254x31;3x1425;4x3125]$\\
Exc $9$: $-(1,1,1,1,1)\,[\{2,4,5\};\{1,2,4\};3x1425;54x231]$\\$[\{4,5\};\{4\};1x2453;2x1453;3x1425;542x31]$\\$[\{2,4\};\{4\};145x32;245x31;3x1425;5x3124]$\\$[\{2,5\};\{1\};251x43;3x1425;4x3125;541x32]$\\
Exc $10$: $-(1,1,2,1)\,[\{3,5\};\{3,4,5\};1x2435;2x1435]$\\
Exc $11$: $-(1,2,2)\,[\{1,2\};\{1,4\};21x543;3x4215;5413x2]$\\
Exc $12$: $-(1,1,1,2)\,[\{1,2\};\{1,4\};213x45;3x5412;5421x3]$\\
Exc $13$: $-(1,3,1,1)\,[\{3,4\};\emptyset;1x62543;2x61543;34x6215;463x215;5x62341;6x12543]$\\
Exc $14$: $-(1,2,1,1,1)\,[\{3,4\};\emptyset;1x63452;2x63451;3645x21;4635x21;5x61432;6x15432]$\\
Exc $15$: $-(1,1,1,2,1)\,[\{1,2,3,6\};\{6\};1362x54;2163x54;3261x54;46x1235;56x1234]$\\$[\{1,2,6\};\emptyset;145x362;245x361;3261x54;46x1235;56x1234;623x415]$\\$[\{1,2,3\};\{6\};13x2654;21x3654;32x1654;4x61235;5x61234]$\\$[\{1,2\};\emptyset;1362x54;245x361;3261x54;4x61235;5x61234;6x34521]$\\$[\{1,6\};\emptyset;145x362;245x361;3261x54;46x1235;56x1234;623x415]$\\ The other possible inneighbors of $x$ are equivalent.\\
Exc $16$: $-(1,1,2,1,1)\,[\{4,5,6\};\{4,5\};14x6235;24x6315;34x6125;612x435]$\\$[\{4,5\};\{4\};1x62435;2x61435;3x62415;5413x62;6x15432]$\\$[\{4,6\};\{4,5\};14x6235;24x6315;34x6125;612x435]$\\$[\{5,6\};\{4,5\};15x6234;25x6314;35x6124;612x534]$\\
Exc $17$: $-(1,2,1,1,1)\,[\{2,4,6\};\emptyset;156x243;21354x6;312x465;43516x2;534x621;65132x4]$\\$[\{2,4\};\emptyset;1623x54;2613x54;312x465;4235x16;5461x32;65132x4]$\\ The other possible inneighbors of $x$ are equivalent.\\
Exc $18$: $-(1,1,2,2)\,[\{2,4,6\};\emptyset;1x53246;213x456;3x15462;435x612;5x32614;651x234]$\\$[\{2,4\};\emptyset;1x53246;213x456;3x15462;4231x56;5x32614;651x234]$\\ The other possible inneighbors of $x$ are equivalent.\\
Exc $19$: $-(1,1,1,1,2)\,[\{2,4,6\};\emptyset;16x4532;26x4531;32x6154;42x6153;54x2316;64x2315]$\\$[\{2,4\};\emptyset;1x63542;213x654;3x15642;423x165;54x2316;64x2315]$\\ The other possible inneighbors of $x$ are equivalent.\\
Exc $20$: $-(1,1,1,1,1,1)\,[\{1,2\};\emptyset;156x423;256x413;324x651;413x652;546x312;645x312]$\\
Exc $21$: $-(1,2,1,1,1)\,[\{1,2\};\emptyset;1623x54;2613x54;3124x65;4x51632;5x63142;6514x32]$\\
Exc $22$: $-(1,1,1,1,1,1)\,[\{1,2,3\};\{1,2,3,4,5,6\}]$\\$[\{1,2\};\emptyset;14253x6;254x361;3x64215;4x36512;51246x3;6x35241]$\\ The other possible inneighbors of $x$ are equivalent.\\
Exc $23$: $-(1,1,1,1,1,1)\,[\{4,6\};\emptyset;16453x2;2x35461;35461x2;45632x1;51462x3;65342x1]$\\
Exc $24$: $-(1,1,1,1,1,1)\,[\{1,2,3,4\};\{4,5\};16x3425;26x1435;36x1425;6541x23]$\\$[\{1,2,3\};\emptyset;13425x6;21435x6;32415x6;4x51263;5x41263;6541x23]$\\$[\{1,2,4\};\emptyset;13425x6;21435x6;32415x6;41523x6;5x36214;6541x23]$\\$[\{1,2\};\emptyset;13425x6;21435x6;32415x6;41523x6;5x36214;6541x23]$\\$[\{1,4\};\emptyset;13425x6;21435x6;32415x6;41523x6;5x36214;64x1523]$\\ The other possible inneighbors of $x$ are equivalent.\\
Exc $25$: $-(1,1,1,1,2)\,[\{4,5\};\emptyset;1x65432;2x65413;3x65421;412x653;512x643;6x13452]$\\
Exc $26$: $-(1,1,1,3)\,[\{4,5,6\};\{6\};1x26354;2x36154;3x16254;45x6123;512x346]$\\$[\{4,5\};\{6\};1x26354;2x36154;3x16254;45x6123;512x346]$\\$[\{4,6\};\{6\};1x26354;2x36154;3x16254;412x356;512x346]$\\$[\{5,6\};\emptyset;1x26354;2x36154;3x16254;45x6123;512x346;654x123]$\\
Exc $27$: $-(1,1,3,1)\,[\{4,6\};\emptyset;162x435;261x435;361x425;4152x63;5142x63;65432x1]$\\
Exc $28$: $-(1,2,1,2)\,[\{4,5\};\emptyset;1x63254;2x61354;3x62154;4156x32;5146x32;6x14235]$\\
Exc $29$: $-(1,1,1,1,1,1)\,[\{4,6\};\emptyset;1x35462;2x13654;35462x1;45632x1;5x41632;631x524]$\\
Exc $30$: $-(1,1,1,1,1,1)\,[\{4,6\};\emptyset;132x564;231x564;34x6521;453x126;5x13624;621x354]$\\
Exc $31$: $-(1,1,2,2)\,[\{3,4\};\emptyset;132x654;231x654;35461x2;45361x2;5x21463;6x21453]$\\
Exc $32$: $-(1,1,2,1,1)\,[\{5,6\};\{4,5\};15x6234;25x6134;35x6124;613x524]$\\
Exc $33$: $-(1,1,1,1,1,1,1)$. The Paley tournament is arc transitive, so without loss of generality we may suppose that $\{1,2\}\subseteq N^-(x)$. We have $42x13675\equiv 51x24376\equiv 62x15347\equiv 7351x264\equiv P$. Now if $3\rightarrow x$ then $146573x2\equiv P$, so we may assume that $x\rightarrow 3$. If $6\rightarrow x$ then $145376x2\equiv P$ and else $17423x65\equiv P$. So $1$ is always an origin of $P$. Now if $3\rightarrow x$ then $21x37564\equiv P$ so we may suppose that $x\rightarrow 3$. If $x\rightarrow 7$ then $213x7564\equiv P$, so $7\rightarrow x$. If $x\rightarrow 5$ then $21763x54\equiv P$, so $5\rightarrow x$. If $4\rightarrow x$ then $213675x4\equiv P$, so $x\rightarrow 4$. If $6\rightarrow x$ then $274351x6\equiv P$, so $x\rightarrow 6$, in this case $2$ is not origin of $P$. Now if $x\rightarrow 3$ then $3x645127\equiv P$, so $3\rightarrow x$. If $7\rightarrow x$ then $315462x7\equiv P$, so $x\rightarrow 7$. If $x\rightarrow 6$ then $31257x64\equiv P$, so $6\rightarrow x$. If $4\rightarrow x$ then $367521x4\equiv P$, so $x\rightarrow 4$. If $5\rightarrow x$ then $36x52714\equiv P$, so $x\rightarrow 5$. In this case $3$ is not an origin of $P$.\\
Exc $34$: $-(1,1,1,2,1,1)\,[\{1,2,3\};\emptyset;137x6254;237x6145;3621x574;47651x32;57461x32;$\\$67541x32;7124x635][\{1,2\};\emptyset;137x6245;237x6145;3621x574;4x651732;5x461732;$\\$6x541732;7124x635][\{1,3\};\emptyset;137x6254;237x6145;3621x574;47651x32;57461x32;$\\$67541x32;7124x635][\{2,3\};\emptyset;137x6254;237x6145;3621x574;4x651732;5x461732;$\\$6x541732;7124x635]$\\
Exc $35$: $-(1,2,1,3)\,[\{1,2,3\};\emptyset;1437x562;2417x563;3417x562;461x3275;541x3276;$\\$651x3274;7x415632]$ The other possible inneighbors of $x$ are equivalent.\\
Exc $36$: $-(1,2,3,1)\,[\{1,2\};\emptyset;1437x256;2437x156;3417x256;4617x253;5417x263;$\\$6517x243;7x316542]$\\
Exc $37$: $-(1,1,1,1,1,1,1)\,[\{4,5,6\};\emptyset;1x243567;2x341567;3x142567;46153x72;54163x72;$\\$65143x72;712435x6][\{4,5\};\emptyset;1x243567;2x341567;3x142567;46153x72;54163x72;$\\$65143x72;712634x6]$ The other possible inneighbors of $x$ are equivalent.\\
Exc $38$: $-(1,1,1,2,1,1)\,[\{1,2,3\};\emptyset;137x6245;217x6345;327x6145;4x651723;5x461723;$\\$6x541723;7124x653][\{1,2,\};\emptyset;137x6245;217x6345;327x6145;4x651723;5x461723;$\\$6x541723;7124x653]$ The other possible inneighbors of $x$ are equivalent.\\
Exc $39$: $-(1,2,1,3)\,[\{1,2,3\};\emptyset;1437x562;2417x563;3417x562;4x536127;5x634127;$\\$6x435127;7x435612][\{1,2\};\emptyset;1437x562;2417x563;3417x562;4x536127;5x634127;$\\$6x435127;7x435612]$ The other possible inneighbors of $x$ are equivalent.\\
Exc $40$: $-(1,1,1,1,2,1)\,[\{1,2\};\emptyset;167x4325;267x4315;3x746512;4x736512;5x735412;$\\$634x7512;734x6512]$\\
Exc $41$: $-(1,1,1,1,3)\,[\{6,7\};\emptyset;1x267543;2x167543;3x167542;4x167352;5x167432;$\\$634x5217;734x5216]$\\
Exc $42$: $-(1,2,2,1,1)\,[\{6,7\};\emptyset;1x342675;2x341675;3x126475;4x126375;5x126374;$\\$63721x54;73621x54]$\\
Exc $43$: $-(1,1,1,2,1,1)\,[\{2,7\};\emptyset;1x632745;213x4657;314x5627;415x3627;514x3627;$\\$632x1547;761x3245]$\\
Exc $44$: $-(1,1,1,1,3)\,[\{6,7\};\emptyset;1x237546;2x317546;3x127546;4x567123;5x467123]$\\
Exc $45$: $-(1,2,1,2,1)\,[\{1,7\};\emptyset;1426x573;2x317654;3265x417;4652x317;5642x317;$\\$6x314527;7426x513]$\\
Exc $46$: $-(1,2,2,2)\,[\{4,7\};\emptyset;1x247563;2x147563;3x147562;4562x731;5x671324;$\\$6x571324;7562x431]$\\
Exc $47$: $-(1,1,1,2,1,1)\,[\{4,5\};\emptyset;1x647532;2x647531;3x647521;425763x1;527463x1;$\\$6x134725;724563x1]$\\$[\{4,5,7\};\emptyset;1x647532;2x647531;3x647521;425763x1;527463x1;$\\$6x134725;724563x1]$ The other possible inneighbors of $x$ are equivalent.\\
Exc $48$: $-(1,1,1,2,1,1)\,[\{1,2\};\emptyset;147x3256;234x5167;3x762156;435x6271;534x6271;$\\$634x5271;7x435261]$\\
Exc $49$: $-(1,1,1,1,1,1,1,1)\,[\{1,2\};\emptyset;1853426x7;2853416x7;3x4125768;4x3125768;5x3127846;$\\$6x3127845;7x3125648;8x4125763]$\\
Exc $50$: $-(1,2,1,1,1,1,1)\,[\{1,2\};\emptyset;182367x54;281367x54;3x7216458;4x8213675;$\\$5x8214376;6x7214358;7x8216453;8x7216453]$\\
Exc $51$: $-(1,2,1,2,1,1)\,[\{2,8\};\emptyset;1x4835267;2738451x6;3x4856172;4x3852167;5x3842167;\\6x3842157;7x3842156;8732451x6]$
\end{proof}
\begin{lemma}\label{2.12}
Let $T$ be a tournament on $n$ vertices, $P$ an outpath of order $n$ and $x$ a vertex of $T$ of outdegree at least two. If $(T-x;^*P)$ belongs to one of the infinite families of exceptions, $^*P$ is an outpath and $x$ is not an origin of $P$ in $T$ then any vertex of $T$ distinct from $x$ is an origin of $P$ in $T$ with end distinct from $x$ unless $(T;P)$ is a biexception with respect to $x$.
\end{lemma}
\begin{proof}
We use the same notations of the previous lemmas and we denote by $a$ and $b$ a random vertices of $X$ and $Y$ respectively:\\
Exc $E_1(n-1)$: $(2,n-3)\,[\{1,2,3\};\emptyset;123xI_X;231xI_X;l312xI_X;ax321I_{X-a}]$\\$[\{1,2\};\{1\};231xI_X;312xI_X;ax132I_{X-a}]$\\ The other possible outneighbors of $x$ are equivalent.\\
Exc $E_2(n-1)$: $(3,n-4)\,[\{3,4\};\emptyset;12x34I_X;2x314I_X;3124xI_X;4123xI_X;ax4321I_{X-a}]$\\
Exc $E_3(n-1)$: $(2,n-3)\,[\{1,3\};\{1\};2x31vR_2;321xvR_2;a21xI_{X-a}]$\\
Exc $E_4(n-1)$: $(3,n-4)\,[\{1,4\};\emptyset;13x42vR_2;2x14u3R_1;3x142vR_2;43x12vR_2;ax142I_{X-a}]$\\
Exc $E_5(n-1)$: $(2,n-3)\,[\{1,2\};\{1\};ab1xI_{Y-b}I_{X-a};bx1I_{Y-b}I_X]$\\
Exc $E_6(n-1)$: $(3,n-4)\,[\{1,3\};\emptyset;12x3I_YI_{X-2};32x1I_YI_{X-2};ax3I_YI_{X-a};bx31I_{Y-b}I_X]$\\
Exc $E_7(n-1)$: $(2,1,n-4)\,[\{2,3\};\emptyset;1x2O_Y3;21xO_Y3;31xO_Y2;bx2O_{Y-b}31]$\\
Exc $E_8(n-1)$: $(n-4,1,1,1)\,[N^+(x)\subseteq X;X;1O_X2x3;2O_X1x3;3O_X1x2]$\\
Exc $E'_8(n-1)$: $(n-4,2,1)\,[N^+(x)\subseteq X;X;1O_X32x;2O_X13x;3O_X21x]$\\
Exc $E_9(n-1)$: $(n-6,1,1,1,1,1)\,[N^+(x)\subseteq X;X;1O_X2x453;2O_X3x514;3O_X4x125;\\4O_X5x231;5O_X1x342]$\\
Exc $E'_9(n-1)$: $(n-6,2,1,1,1)\,[N^+(x)\subseteq X;X;1O_X32x45;2O_X43x51;3O_X54x12;\\4O_X15x23;5O_X21x34]$\\
Exc $E_{10}(n-1)$: $(n-8,1,1,1,1,1,1,1)\,[N^+(x)\subseteq X;X;1O_X2x45376;2O_X3x56417;\\3O_X4x67521;4O_X5x71632;5O_X6x12743;6O_X7x23154;7O_X1x34265]$\\
Exc $E'_{10}(n-1)$: $(n-8,2,1,1,1,1,1)\,[N^+(x)\subseteq X;X;1O_X32x4657;2O_X43x5761;\\3O_X54x6172;4O_X65x7213;5O_X76x1342;6O_X17x2435;7O_X21x3546]$\\
Exc $E_{11}(n-1)$: $(2,1,n-4)\,[\{1,2\};\{1\};23xO_X1;3x1O_X2;ax1O_{X-a}23]$\\
Exc $E_{12}(n-1)$: $(3,1,n-5)\,[\{1,4\};\emptyset;12x43O_X;2x413O_X;3x12O_X4;42x13O_X;a2x13O_{X-a}4]$\\
Exc $E_{13}(n-1)$: $(2,1,n-4)\,[\{1,2\};\{1\};2x1345O_X;3x145O_X2;4x153O_X2;5x134O_X2;$\\$a2x345O_{X-a}1]$\\
Exc $E_{14}(n-1)$: $(3,1,n-5)\,[\{1,6\};\emptyset;12x6345O_X;2x61345O_X;3x6145O_X2;4x6153O_X2;\\5x6134O_X2;a12x345O_{X-a}6]$  $\square$
\end{proof}
\begin{lemma}\label{2.13}
Let $T$ be a tournament on $n$ vertices, $P$ an inpath of order $n$ and $x$ a vertex of $T$ of indegree at least two. If $(T-x;^*P)$ belongs to one of the infinite families of exceptions, $^*P$ is an outpath and $x$ is not an origin of $P$ in $T$ then any vertex of $T$ distinct from $x$ is an origin of $P$ in $T$ with end distinct from $x$ unless $(T;P)$ is a biexception with respect to $x$.
\end{lemma}
\begin{proof}
We use the same notations of the previous lemma:\\
Exc $E_1(n-1)$: $-(1,1,n-3)\,[\{1,2,3\};X;13x2I_X;21x3I_X;32x1I_X]$\\$[\{1,2\};\emptyset;1a3x2I_{X-a};2a3x1I_{X-a};32x1I_X;ax321I_{X-a}]$\\
Exc $E_2(n-1)$: $-(1,2,n-4)\,[\{3,4\};\emptyset;1x234I_X;2a31x4I_{X-a};34x21I_X;4a12x3I_{X-a};$\\$axa'3421I_{X\setminus \{a,a'\}}]$ If $|X|\geq 2$ and if $|X|=1$, $a$ is not an origin of $P$.\\
Exc $E_3(n-1)$: $-(1,1,n-3)\,[\{1,3\};\emptyset;1a2x3R_1;2r1u3(R_1-r);3v2x1R_2;ax231I_{X\setminus \{a,3\}}]$ Where $r$ is the end of $R_1$.\\
Exc $E_4(n-1)$: $-(1,2,n-4)\,[\{1,4\};\emptyset;1v32x4R_2;2r14u3(R_1-r);3x241I_{X-3};$\\$4v32x1R_2;ax241I_{X-a}]$ Where $r$ is the end of $R_1$.\\
Exc $E_5(n-1)$: $-(1,1,n-3)$. Three cases may arise:\begin{itemize}
\item If $X=\{2\}$ and $T(Y)$ is not $3A$, let $Q=(1|Y|-2)$ in $Y$, $[\{1,2\};\{2\};1Qx2;bxI_{Y-b}21]$.
\item If $X=\{2\}$ and $T(Y)=3A$, $[\{1,2\};\{1,2\};bxI_{Y-b}21]$.
\item If $|X|\geq 2$, since $2$ is an ingenerator of $X$, let $2Q$ be a Hamiltonian directed inpath of $X$ with end $v$ and let $u$ be an inneighbor of $2$ in $X$.\\$[\{1,2\};\emptyset;1vI_Yx2Q^*;2uI_Yx1I_{X\setminus \{2,u\}};axI_Y21I_{X\setminus \{a,2\}};bxI_{Y-b}21I_{X-2}]$.
\end{itemize}
Exc $E_6(n-1)$: $-(1,2,n-4)$. Three cases may arise:\begin{itemize}
\item If $X=\{2\}$ and $|Y|\geq 3$, let $u\rightarrow v\rightarrow w$ be three vertices of $Y$, \\$[\{1,3\};\emptyset;1uvw2x3I_{Y\setminus \{u,v,w\}};2xb13I_{Y-b};3uvw2x1I_{Y\setminus \{u,v,w\}};bxb'13I_{Y\setminus \{b,b'\}}2]$.
\item If $X=\{2\}$ and $Y\{u\rightarrow v\}$, $[\{1,3\};\{3,u\};13xvu2;2xu13v;vu12x3]$.
\item If $|X|\geq 2$, let $u$ be an inneighbor of $2$ in $X$, \\$[\{1,3\};\emptyset;1u2I_Yx3I_{X\setminus \{2,u\}};3u2I_Yx1I_{X\setminus \{2,u\}};axb13I_{Y-b}I_{X-a};bxu213I_{Y-b}I_{X\setminus \{2,u\}}]$.
\end{itemize}
Exc $E_7(n-1)$: $-(1,1,1,n-4)$. Let $u\rightarrow v$ be an arc of $T(Y)$, two cases may arise:\begin{itemize}
\item If $T(Y)$ is not $3A$, let $Q=-(1,|Y|-2)$ be a Hamiltonian path of $Y$, \\$[\{2,3\};\emptyset;1xQ32;2uvxO_{Y\setminus \{u,v\}}31;3uvxO_{X\setminus \{u,v\}}21;axa'1I_{Y\setminus \{a,a'\}}32]$.
\item If $T(Y)$ is not $3A$, $S=\{1\}$.
\end{itemize}
Exc $E_8(n-1)$: $-(1,n-5,1,1,1)$. We have $N^-(x)\subseteq X$. Let $u\rightarrow v$ be an arc of $T(X)$. Two cases may arise:\begin{itemize}
\item If $|X|\geq 4$, let $w\in N^-(x)$.\\$[N^-(x);\emptyset;1x2O_{X\setminus \{u,v\}}3vu;2x3O_{X\setminus \{u,v\}}1vu;3x1O_{X\setminus \{u,v\}}2vu;a1wx3O_{X\setminus \{a,u,v,w\}}2vu]$\\(The path $O_{X\setminus \{a,u,v,w\}}$ may be empty).
\item If $|X|=3$, say $X=\{u,v,w\}$. $1,2$ and $3$ are origins as above. If $u$ is an outneighbor of $x$, we have $ux1v2w3\equiv P$. If $u$ has an inneighbor $v$ that is an inneighbor of $x$, we have $uvx13w2$. So the sole vertex that is not an origin of $P$ is the inneighbor of $x$ that is minimal in $X$.
\item If $|X|=2$. $[\{u,v\};\{u,v\};1x32vu;2x13vu;3x21vu]$
\end{itemize}
Exc $E'_8(n-1)$: $-(1,n-5,2,1)$. Two cases may arise:\begin{itemize}
\item If $|X|\geq 3$, let $Q=(n-7,1)$ be a path in $T(X)$ and $u=X-V(Q)$.\\$[N^-(x);\emptyset;1x2Q3u;2x3Q1u;3x1Q2u;a1O_{X-a}2x3]$
\item If $|X|=2$, say $X=\{u,v\}$\\$[\{u,v\};\{1,2,3\};u1v2x3;v1u2x3]$
\end{itemize}
Exc $E_9(n-1)$: $-(1,n-7,1,1,1,1,1)$. Let $u,v\in X$.\\$[N^-(x)\subseteq X;\emptyset;1x3O_{X\setminus \{u,v\}}2u4v5;a1O_{X-a}23x54]$. The other vertices of $5A$ are equivalent.\\
Exc $E'_9(n-1)$: $-(1,n-7,2,1,1,1)$. Let $u,v\in X$.\\$[N^-(x)\subseteq X;\emptyset;1x4O_{X\setminus \{u,v\}}32u5v;a1O_{X-a}2x435]$. The other vertices of $5A$ are equivalent.\\
Exc $E_{10}(n-1)$: $-(1,n-9,1,1,1,1,1,1,1)$. Let $u,v\in X$.\\$[N^-(x)\subseteq X;\emptyset;1x3O_{X\setminus \{u,v\}}2u4v675;a1O_{X-a}23x7564]$. The other vertices of $7A$ are equivalent.\\
Exc $E'_{10}(n-1)$: $-(1,n-9,2,1,1,1,1,1)$. Let $u,v\in X$.\\$[N^-(x)\subseteq X;\emptyset;1x4O_{X\setminus \{u,v\}}32u5v67;a1O_{X-a}2x73546]$. The other vertices of $7A$ are equivalent.\\
Exc $E_{11}(n-1)$: $-(1,1,1,n-4)$. Let $u\rightarrow v$ be an arc of $T(X)$, \\$[\{1,2\};\emptyset;1xvu23O_{X\setminus \{u,v\}};2uvx3O_{X\setminus \{u,v\}}1;3x1O_X2;a31xO_{X-a}2]$\\
Exc $E_{12}(n-1)$: $-(1,2,1,n-5)\,[\{1,4\};\emptyset;1a23xO_{X-a}4;2x34O_X1;3xa2O_{X-a}41;\\4a23xO_{X-a}1;axa'1423O_{X\setminus \{a,a'\}}]$\\
Exc $E_{13}(n-1)$: $-(1,1,1,n-4)$. Let $u\rightarrow v$ be an arc of $T(X)$, \\$[\{1,2\};\emptyset;1342x5O_X;2uvx345O_{X\setminus \{u,v\}}1;3x425O_X1;4x523O_X1;5x324O_X1;a34x5O_{X-a}12]$\\
Exc $E_{14}(n-1)$: $-(1,2,1,n-5)\,[\{1,6\};\emptyset;13455xO_X62;2134x5O_X6;3x254O_X61;4x235O_X61;\\5x243O_X61;6345x)_X12;ax3425O_{X-2}61]$  $\square$
\end{proof}
\\From now on, we call the lemmas \ref{2.10}, \ref{2.11}, \ref{2.12} and \ref{2.13} the \textit{building lemmas}.

\section{Main result}
Below are all the exceptions $(T;C)$ where $T$ does not contains $C$:
\\$A_1=(3A;(2,1))$, $A_2=(5A;(2,1,1,1))$, $A_3=(7A;(2,1,1,1,1,1))$, $A_4=(4C;(1,1,1,1))$, $A_5=(6M;(1,1,1,1,1,1))$, $A_6=(8C;(1,1,1,1,1,1,1,1))$, $A_7=(6N;(4,2))$,\\$A_8=(4B;(2,2))$, $A_9=(5C;(4,1))$, $A_{10}=(8D;(1,1,1,1,1,1,1,1))$, $A_{11}=(6O;(1,1,1,1,1,1))$, $A_{12}=(6P;(1,1,1,1,1,1))$, $A_{13}=(6F;(2,1,2,1))$, $A_{14}=(6D;(2,2,1,1))$,\\$A_{15}=(5A;(1,1,1,1))$, $A_{16}=(5F;(1,1,1,1))$, $A_{17}=(7K;(1,1,1,1,1,1))$,\\$A_{18}=(7A;(2,1,2,1))$.\\

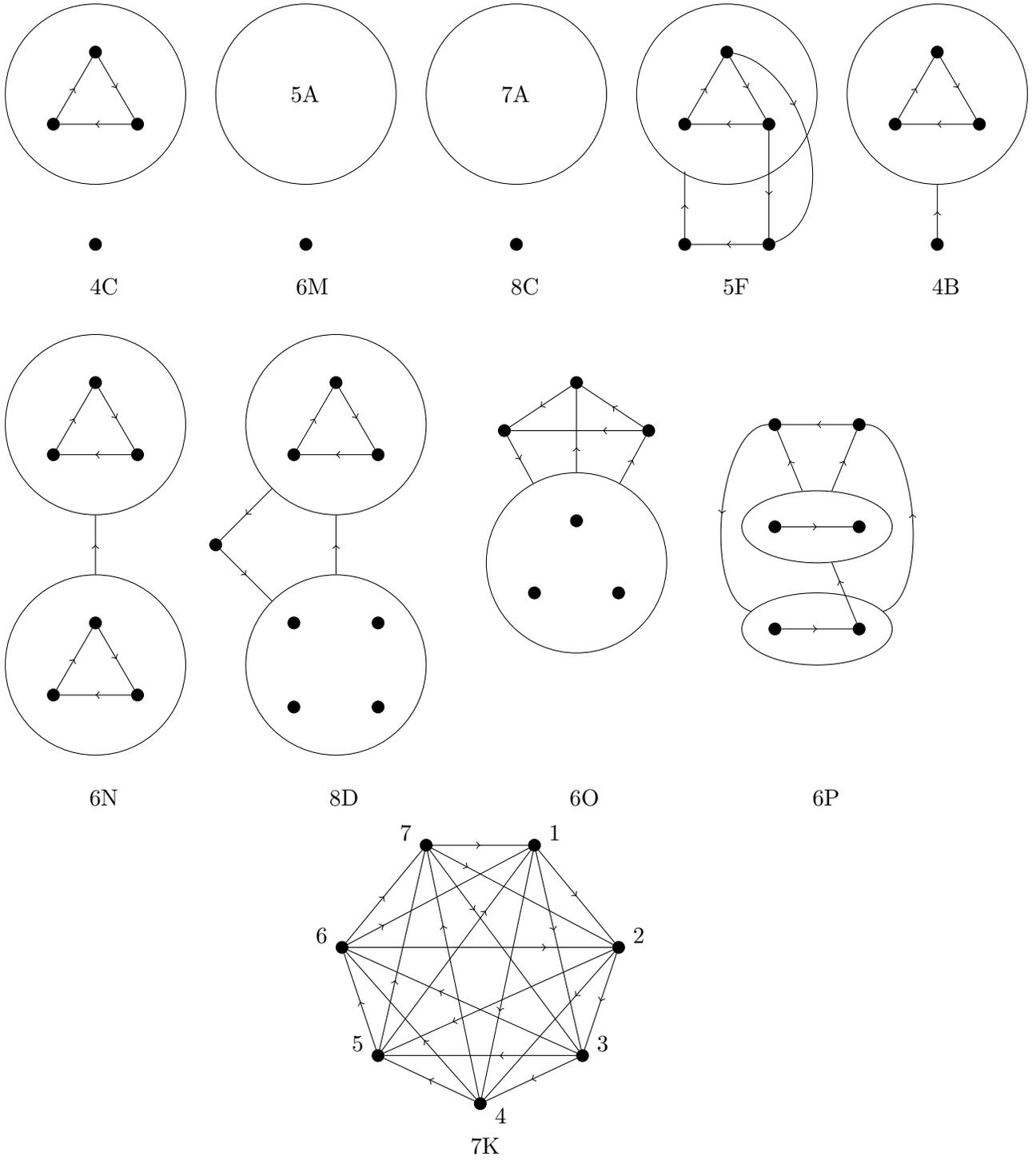
\begin{figure}
\begin{tikzpicture}
\tikzset{enclosed/.style={draw,circle,inner sep=2pt,minimum size=4pt,fill=black}}
\tikzset{->-/.style={decoration={
            markings,
            mark=at position #1 with
            {\arrow{>}}},postaction={decorate}}}
            
\node[ellipse,minimum width=3cm,minimum height=3cm,draw](a)at(0,0){};
\node[enclosed](1)at(-0.7,-0.5){};
\node[enclosed](2)at(0,0.7){};
\node[enclosed](3)at(0.7,-0.5){};
\node[enclosed](4)at(0,-2.5){};

\draw[black,->-=.5] (1)--(2);
\draw[black,->-=.5] (2)--(3);
\draw[black,->-=.5] (3)--(1);

\node[minimum size=.1pt,label={left:4C}](4C)at(0.65,-3.2){};

\node[ellipse,minimum width=3cm,minimum height=3cm,draw](a)at(3.5,0){};
\node[minimum size=.1pt,label={left:5A}](X)at(4,0){};
\node[enclosed](4)at(3.5,-2.5){};

\node[minimum size=.1pt,label={left:6M}](6M)at(4.15,-3.2){};

\node[ellipse,minimum width=3cm,minimum height=3cm,draw](a)at(7,0){};
\node[minimum size=.1pt,label={left:7A}](X)at(7.5,0){};
\node[enclosed](4)at(7,-2.5){};

\node[minimum size=.1pt,label={left:8C}](8C)at(7.65,-3.2){};

\node[ellipse,minimum width=3cm,minimum height=3cm,draw](a)at(0,-5.5){};
\node[ellipse,minimum width=3cm,minimum height=3cm,draw](b)at(0,-9.5){};

\node[enclosed](1)at(-.7,-6){};
\node[enclosed](2)at(0,-4.8){};
\node[enclosed](3)at(.7,-6){};

\draw[black,->-=.5] (1)--(2);
\draw[black,->-=.5] (2)--(3);
\draw[black,->-=.5] (3)--(1);

\node[enclosed](4)at(-.7,-10){};
\node[enclosed](5)at(0,-8.8){};
\node[enclosed](6)at(.7,-10){};

\draw[black,->-=.5] (4)--(5);
\draw[black,->-=.5] (5)--(6);
\draw[black,->-=.5] (6)--(4);

\draw[black,->-=.5] (b)--(a);

\node[minimum size=.1pt,label={left:6N}](6N)at(.65,-11.7){};

\node[ellipse,minimum width=3cm,minimum height=3cm,draw](a)at(14,0){};
\node[enclosed](1)at(13.3,-0.5){};
\node[enclosed](2)at(14,0.7){};
\node[enclosed](3)at(14.7,-0.5){};
\node[enclosed](4)at(14,-2.5){};

\draw[black,->-=.5] (1)--(2);
\draw[black,->-=.5] (2)--(3);
\draw[black,->-=.5] (3)--(1);
\draw[black,->-=.5] (4)--(a);

\node[minimum size=.1pt,label={left:4B}](4B)at(14.65,-3.2){};

\node[ellipse,minimum width=3cm,minimum height=3cm,draw](a)at(10.5,0){};
\node[enclosed](1)at(9.8,-.5){};
\node[enclosed](2)at(10.5,.7){};
\node[enclosed](3)at(11.2,-0.5){};
\node[enclosed](4)at(9.8,-2.5){};
\node[enclosed](5)at(11.2,-2.5){};
\node[minimum size=.1pt](a1)at(9.8,-1.15){};
\node[minimum size=.1pt](a2)at(11.2,-1.15){};

\draw[black,->-=.5] (1)--(2);
\draw[black,->-=.5] (2)--(3);
\draw[black,->-=.5] (3)--(1);
\draw[black,->-=.5] (4)--(a1);
\draw[black,->-=.5] (5)--(4);
\draw[black,->-=.6] (3)--(5);
\draw[black,->-=.35] (2) to [out=-10,in=20,looseness=1] (5);

\node[minimum size=.1pt,label={left:5F}](5F)at(11.15,-3.2){};

\node[ellipse,minimum width=3cm,minimum height=3cm,draw](a)at(4,-5.5){};
\node[ellipse,minimum width=3cm,minimum height=3cm,draw](b)at(4,-9.5){};

\node[enclosed](1)at(3.3,-6){};
\node[enclosed](2)at(4,-4.8){};
\node[enclosed](3)at(4.7,-6){};

\draw[black,->-=.5] (1)--(2);
\draw[black,->-=.5] (2)--(3);
\draw[black,->-=.5] (3)--(1);

\node[enclosed](4)at(3.3,-10.2){};
\node[enclosed](5)at(3.3,-8.8){};
\node[enclosed](6)at(4.7,-10.2){};
\node[enclosed](7)at(4.7,-8.8){};

\node[enclosed](8)at(2,-7.5){};

\draw[black,->-=.5] (a)--(8);
\draw[black,->-=.5] (8)--(b);
\draw[black,->-=.5] (b)--(a);

\node[minimum size=.1pt,label={left:8D}](7I)at(4.65,-11.7){};

\node[ellipse,minimum width=3cm,minimum height=3cm,draw](a)at(8,-7.8){};
\node[enclosed](1)at(7.3,-8.3){};
\node[enclosed](2)at(8,-7.1){};
\node[enclosed](3)at(8.7,-8.3){};

\node[enclosed](4)at(9.2,-5.6){};
\node[enclosed](5)at(8,-4.8){};
\node[enclosed](6)at(6.8,-5.6){};

\draw[black,->-=.5] (a)--(4);
\draw[black,->-=.3] (a)--(5);
\draw[black,->-=.5] (4)--(5);
\draw[black,->-=.3] (4)--(6);
\draw[black,->-=.5] (5)--(6);
\draw[black,->-=.5] (6)--(a);

\node[minimum size=.1pt,label={left:6O}](7I)at(8.65,-11.7){};

\node[ellipse,minimum width=2.5cm,minimum height=1.2cm,draw](b)at(12,-7.2){};
\node[ellipse,minimum width=2.5cm,minimum height=1.2cm,draw](c)at(12,-8.9){};

\node[enclosed](1)at(11.3,-5.5){};
\node[enclosed](2)at(12.7,-5.5){};

\node[enclosed](3)at(11.3,-7.2){};
\node[enclosed](4)at(12.7,-7.2){};

\node[enclosed](5)at(11.3,-8.9){};
\node[enclosed](6)at(12.7,-8.9){};

\draw[black,->-=.5] (2)--(1);
\draw[black,->-=.5] (3)--(4);
\draw[black,->-=.5] (5)--(6);
\draw[black,->-=.5] (b)--(1);
\draw[black,->-=.5] (b)--(2);
\draw[black,->-=.7] (6)--(b);

\draw[black,->-=.5] (c) to [out=15,in=0,looseness=0.7] (2);
\draw[black,->-=.5] (1) to [out=180,in=165,looseness=0.7] (c);

\node[minimum size=.1pt,label={left:6P}](7I)at(12.65,-11.7){};

\node[enclosed,label={right,yshift=.2cm:1}](1)at(7.3,-12.5){};
\node[enclosed,label={right,yshift=.2cm:2}](2)at(8.7,-14.2){};
\node[enclosed,label={right,yshift=.2cm:3}](3)at(8.1,-16){};
\node[enclosed,label={right,yshift=-.2cm:4}](4)at(6.4,-16.8){};
\node[enclosed,label={left,yshift=.2cm:5}](5)at(4.7,-16){};
\node[enclosed,label={left,yshift=.2cm:6}](6)at(4.1,-14.2){};
\node[enclosed,label={left,yshift=.2cm:7}](7)at(5.5,-12.5){};

\draw[black,->-=.5] (1)--(2);
\draw[black,->-=.4] (1)--(3);
\draw[black,->-=.65] (1)--(4);
\draw[black,->-=.5] (2)--(3);
\draw[black,->-=.3] (2)--(4);
\draw[black,->-=.7] (2)--(5);
\draw[black,->-=.5] (3)--(4);
\draw[black,->-=.4] (3)--(5);
\draw[black,->-=.6] (3)--(6);
\draw[black,->-=.5] (4)--(5);
\draw[black,->-=.4] (4)--(6);
\draw[black,->-=.7] (4)--(7);
\draw[black,->-=.5] (5)--(6);
\draw[black,->-=.35] (5)--(7);
\draw[black,->-=.7] (5)--(1);
\draw[black,->-=.5] (6)--(7);
\draw[black,->-=.2] (6)--(1);
\draw[black,->-=.75] (6)--(2);
\draw[black,->-=.5] (7)--(1);
\draw[black,->-=.2] (7)--(2);
\draw[black,->-=.3] (7)--(3);

\node[minimum size=.1pt,label={left:7K}](7I)at(7,-17.5){};

\end{tikzpicture}
\caption{Cycle's exceptions}
\label{figure 7}
\end{figure}
The exceptions $A_i$ for $1\leq i\leq 12$ are found by Havet \cite{H}. We show our main result in the other cases, for the Hamiltonian cycles, using the following theorem due to El-Sahili and Ghazo-Hanna \cite{ES}, by discussing the minimal indegree of the tournament.
\begin{theorem}\label{3.1}
Let $T$ be a tournament and $P$ a path such that $|T|=|P|$. The number of paths $P$ in $T$ is equal to that of $\overline{P}$ in $T$.
\end{theorem}
The following theorem concerns the tournaments with minimal indegree at most $1$ and the non-antidirected cycles.
\begin{theorem}\label{3.2}
Let $T$ be a tournament of order $n$ such that $\delta^-(T)\leq 1$ and $C$ a non-directed cycle of order $n$. If $C$ is not antidirected then $T$ contains $C$ if and only if $(T;C)$ is not one of the exceptions $A_i$, $i\in \{1,4,5,6,7,8,9,10,11,12\}$.
\end{theorem}
\begin{proof}
By theorem \ref{3.1}, we may suppose that $\delta^-(T)\leq \delta^+(T)$. Set $\delta=\delta^-(T)$. Let $v\in T$ be a vertex such that $d^-(v)=\delta$, $u$ be its inneighbor if it exists and $T_2=T(N^+(v))$. Let $B$ be a block of $C$ of maximal length and set $C=a_0a_1...a_{n-1}a_0$ such that $B$ is forward and it ends at $a_{\delta+2}$. If $\delta=0$, let $s=min \{i; a_i\rightarrow \{a_{i-1},a_{i+1}\}\}$. Clearly, the path $P=a_{i+1}a_{i+2}...a_{n-1}a_0...a_{i-1}$ is contained in $T_2$ if and only if $(T_2;P)\neq (3A;(1,1))$ and in this case $(T;C)=A_8$, we have $vPv\equiv C$. So we may suppose that $\delta=1$, we distinguish three cases:\\
\textbf{Case $1$}: If $|B|\geq 4$, rotate $C$ until $a_{n-1}\rightarrow a_{n-2}$ for the first time, then $b_1(a_0a_{n-1}a_{n-2}...a_3)=1$.\begin{itemize}
\item[(a)] If $d^-_{T_2}(u)\geq 2$ then, by theorem \ref{2.1}, by taking two inneighbors of $u$ in $T_2$, $u$ has an inneighbor $w$ origin of a path $P\equiv a_0a_{n-1}...a_3$ in $T_2$, and thus $wuvP^{-1}\equiv C$, or $(T_2;a_0a_{n-1}...a_3)$ is an exception, the possible exceptions are: Dual Exc $0,1,2,4,6,8,9$,
$10,12,15,16,18,19,20,22,23,$\\$24,25,26,27,29
,30,31,32,33,34,37,38,40,41,43,44,47,$ 
$48,49$. By checking each one, for Dual Exc $1,2,8,9,10,12,15,16,18,19,20,22,23,24,$
$25,26,27,29,30,31,32,34,37,38,40,41,$\\$43,44,47,48,49$ one can find an inneighbor of $u$ origin of a path $Q\equiv a_{n-1}a_{n-2}...a_2$ in $T_2$ and thus $uvQ^{-1}u\equiv C$. While for Dual Exc $0,4,6,33$, if $d^+_{T_2}(u)\geq 1$, $u$ has an outneighbor origin of a path $R\equiv a_{n-2}a_{n-3}...a_1$ and so $vR^{-1}uv$. Otherwise, we have:\\
Dual Exc $0$: $(T;C)$ is the exception $A_9$.\\
Dual Exc $4$: $543u21v5\equiv C$.\\
Dual Exc $6$: $15uv3241\equiv C$.\\
Dual Exc $33$: $765u4231v7\equiv C$.
\item[(b)] If $d^-_{T_2}(u)=1$, set $N^-(u)=\{w\}$. If $T_2-w$ contains a path $P\equiv a_2a_3...a_{n-2}$, then $wuPvw\equiv C$. Otherwise, $(T_2-w;a_2a_3...a_{n-2})$ is a Gr\"{u}nbaum's exception, so by theorem \ref{2.1} $u$ has an outneighbor $x$ origin of a path $Q\equiv a_0a_1...a_{n-3}$ in $T_2$ and so $Qvux\equiv C$, unless $T_2=w\rightarrow 3A$ and $a_0a_1...a_{n-3}$ is directed, which is the exception $A_7$.
\end{itemize}
\textbf{Case $2$}: If $|B|=3$, denote by $B'$ and $B"$ the blocks of $C$ containing the arc $a_0\rightarrow a_{n-1}$ and $a_4\rightarrow a_3$ respectively, without loss of generality, we can suppose that $|B'|\geq |B"|$, since otherwise, we can take $\overline{C}=a'_0a'_1...a'_{n-1}a'_0=\overline{a_3a_2...a_4a_3}$, we have two cases:\begin{itemize}
\item[(a)] If $d^-_{T_2}(u)\geq 2$, suppose that $s^+_{T_2}(N^-_{T_2}(u))\geq |B'|+1$, then $u$ has an inneighbor $w$ origin of a path $P\equiv a_0a_{n-1}...a_3$ and so $wuvP^{-1}\equiv C$ or $(T_2;a_0a_{n-1}...a_3)$ is an exception. By checking all the possible exceptions, $u$ has an inneighbor $w$ origin of a path $Q\equiv \overline{a_3a_4...a_{n-1}a_0}$ and thus $vuQv\equiv \overline{a_1a_2...a_1}=\overline{C}$ unless $(T_2;a_0a_{n-1}...a_3)$ is one of the exceptions: $1,20,22,23,24,28,29,30,34,$ $38,43,46,47,48,49$. If it is $28$, then $u614523vu\equiv C$. Otherwise, there exist $w,t\in N^-_{T_2}(u)$ and $y\in N^+_{T_2}(u)$ such that $y$ is an origin of a path $R\equiv a_2a_3...a_{n-3}$ in $T_2\setminus \{w,t\}$ and $w\rightarrow t$, thus $wuRvtw\equiv C$. So we may assume that $s^+_{T_2}(N^-_{T_2}(u))\leq |B'|$, $u$ has an outneighbor origin of a path $R\equiv a_{n-1}a_{n-2}...a_2$ and thus $uvR^{-1}u\equiv C$, unless $(T_2;a_{n-1}...a_2)=(4B;(2,1))$, in this case $v41u32v\equiv C$.
\item[(b)] If $d^-_{T_2}(u)=1$, set $N^-(u)=\{w\}$. Since $d^+_{T_2}(u)=|T_2|-1$, $u$ has an outneighbor origin of a path $R\equiv a_{n-1}a_{n-2}...a_2$ in $T_2$ and thus $uvR^{-1}u\equiv C$, unless $(T_2;a_{n-1}...a_2)$ is one of the exceptions $0,1,7$ and in these cases, by taking $u$ in the place of $v$, we obtain an inneighbor of $u$ having an indegree greater than $1$ and the result follows.
\end{itemize}
\textbf{Case $3$}: If $|B|=2$, if $u$ has an outneighbor $w$ origin of a path $P\equiv a_0a_{n-1}...a_3$ or $Q\equiv \overline{a_4a_5...a_1}$ then $wuvP^{-1}\equiv C$ or $uQvu\equiv \overline{C}$. Otherwise, suppose first that $d^+_{T_2}(u)\geq 2$, we have two cases:\begin{itemize}
\item[(a)] If $C=(2,1,2,1,...,2,1)$ then either $(T_2;a_0a_{n-1}...a_3)$ is the dual of the exception Exc $3$ or $45$, in these cases there exist $w,t\in N^+_{T_2}(u)$ and $y\in T_2\setminus \{w,t\}$ such that $w\rightarrow y$ and $t$ is an origin of a path $R\equiv a_2a_3...a_{n-3}$ in $T_2\setminus \{w,y\}$, we have $wuRvtw\equiv C$, or $s^-_{T_2}(N^+_{T_2}(u))=2$, let $y\in N^-_{T_2}(u)$ be an origin of a path $R\equiv a_1a_0a_{n-1}...a_6$ ($R$ may be an arc) in $T(N^-(u))$ and $\{w,t\}=N^+_{T_2}(u)$, we have $yuwvtR^{-1}\equiv a_1a_2...a_1=C$.
\item[(b)] If $a_0\rightarrow a_{n-1}$, then $(T_2;a_0a_{n-1}...a_3)$ is one of the exceptions $6,12,18,31$, let $w,t\in N^+_{T_2}(u)$ and $y\in N^-_{T_2}(u)$ such that $y$ is an origin of a path $R\equiv a_1a_0a_{n-1}...a_6$ in $T_2\setminus \{w,t\}$, we have $yuwtvR^{-1}\equiv a_1a_2...a_2=C$.
\item[(c)] If $a_{n-1}\rightarrow \{a_0,a_{n-2}\}$, without loss of generality we can suppose that $a_4\rightarrow \{a_3,a_5\}$, so $(T_2;a_0a_{n-1}...a_3)$ is one of the exceptions $0,4,8,9,33,37$. There exist $w\in N^-_{T_2}(u), t\in N^+_{T_2}(u)$ and $y\in T_2\setminus \{w,t\}$ such that $w\rightarrow y$ and $t$ is an origin of a path $R\equiv a_3a_4...a_{n-2}$ in $T_2\setminus \{w,y\}$; we have $ywuRvy\equiv C$.
\end{itemize}
So we may assume that $d^+_{T_2}(u)\leq 1$. If $|T|\geq 6$, there exists $w\in N^-_{T_2}(u)$ such that $T_2-w$ has at least two origins of a path $R\equiv a_4a_5...a_0$ and thus $vwuRv\equiv a_1a_2...a_1=C$, unless $T_2=4B$ and $R$ is directed, in this case $(T;C)=A_{13}$. And if $|T|\leq 5$, easily one can find $C$ unless $(T;C)=A_1$.    $\square$
\end{proof}
\\To prove the existence of non-antidirected cycles in tournaments of minimal indegree at least $2$, we establish four lemmas. Before stating them, we give an overview of the way to reach the result. For a given tournament $T$ with minimal indegree at least $2$, by theorem \ref{3.1}, we may suppose that $\delta^-(T)\leq \delta^+(T)$. Set $\delta=\delta^-(T)$. Let $v\in T$ be a vertex such that $d^-(v)=\delta$ and denote by $T_1=T(N^-(v))$ and $T_2=T(N^+(v))$. Let $B$ be a block of $C$ of maximal length and set $C=a_0a_1...a_{n-1}a_0$ such that $B$ is forward and it ends at $a_{\delta+2}$. We consider the four cases of orientations of the subpath $a_{n-1}a_0a_1$, each one discussed in a lemma. In each case, we rotate $C$ to minimize the length of the first block of $a_0a_{n-1}...a_{\delta+2}$ for the first time, keeping $a_{\delta}a_{\delta+1}a_{\delta+2}$ forward. Since the length of $B$ is greater than that of the first block of $a_0a_{n-1}...a_{\delta+2}$, this rotation allows us to obtain $b_1(a_0a_{n-1}...a_{\delta+2})\leq 2$. We consider two cases: $(T_1;a_1...a_{\delta})$ is not a Gr\"{u}nbaum's exception or it is an exception. Discussing the first one, if there exist $x_1...x_{\delta}\equiv a_1...a_{\delta}$ in $T_1$ and $x_0x_{n-1}...x_{\delta+2}\equiv a_0a_{n-1}...a_{\delta+2}$ in $T_2$ such that $x_0x_1\equiv a_0a_1$, then $x_0x_1...x_{\delta}vx_{\delta+2}...x_{n-1}x_0\equiv C$. Otherwise, let $x_1...x_{\delta}\equiv a_1...a_{\delta}$ be a path in $T_1$ such that $d^+_{T_2}(x_1)=max\,\{d^+_{T_2}(u);u\in Or(a_1...a_{\delta},T_1)\}$ if $a_1\rightarrow a_0$ and $d^-_{T_2}(x_1)=max\,\{d^-_{T_2}(u);u\in Or(a_1...a_{\delta},T_1)\}$ if $a_0\rightarrow a_1$. Set $X=N^+_{T_2}(x_1)$ if $a_1\rightarrow a_0$ and $X=N^-_{T_2}(x_1)$ otherwise. As $b_1(a_0a_{n-1}...a_{\delta+2})\leq 2$, by theorem \ref{2.1}, either $(T_2;a_0a_{n-1}...a_{\delta+2})$ is an exception or $s^+_{T_2}(X)\leq b_1(a_0a_{n-1}...a_{\delta+2})$ if $a_1\rightarrow a_0$ and $s^-_{T_2}(X)\leq b_1(a_0a_{n-1}...a_{\delta+2})$ if $a_0\rightarrow a_1$. In the first case, if $a_2\rightarrow a_1$ and $d^+_{T_2}(x_2)=0$, by lemma \ref{2.8} to find a path $P\equiv a_0a_{n-1}...a_{\delta+2}$ in $\overline{T_2}$ with origin $w$ such that $wx_1\equiv a_0a_1$, so $wx_1vx_3...x_{\delta}x_2P^{-1}\equiv C$ in $\overline{T}$ and then $T$ contains $C$, unless $(T_2;a_0a_{n-1}...a_{\delta+2})$ is one of the exceptions that are mentioned in lemma \ref{2.8}; these cases are treated one by one. Otherwise, we may apply lemma \ref{2.10} to find a path $y_1y_0....y_{\delta+2}\in a_1a_0...a_{\delta+2}$ in $T_2\cup \{x_1\}$ such that $y_1x_2\equiv a_1a_2$ and so $y_0y_1x_2...x_{\delta}vy_{\delta+2}...y_{n-1}y_0\equiv C$ unless $(T_2;a_0a_{n-1}...a_{\delta+2})$ is a biexception with respect to $x_1$ and this is a particular case that can be verified. And in the second case, the direction of the arcs between $x_1$ and $T_2$ and between $X$ and $T_2$ are known, so we can find $C$ unless $(T;C)$ is an exception, which is the desired result. On the other hand, if $(T_1;a_1...a_{\delta})$ is a Gr\"{u}nbaum's exception, according to each case, sometimes we add a vertex from $T_2$ to $T_1$ to obtain $T'_1$, $(T'_1;a_0...a_{\delta})$ is not a Gr\"{u}nbaum's exception, and sometimes $(T_1;a_0...a_{\delta-1})$ is not a Gr\"{u}nbaum's exception with $a_{\delta-1}a_{\delta}a_{\delta+1}$ is directed; in both cases we treat in a similar way as above.\\
It is clear that the four possible orientations of the subpath $a_{n-1}a_0a_1$ are: $a_{n-1}\rightarrow a_0\rightarrow a_1$, $a_{n-1}\leftarrow a_0\leftarrow a_1$, $a_{n-1}\leftarrow a_0\rightarrow a_1$ and $a_{n-1}\rightarrow a_0\leftarrow a_1$. So we start by the first one:
\begin{lemma}\label{3.3}
Let $T$ be a tournament of order $n$ with $\delta^-(T)\geq 2$ and $C$ a non-directed and non-antidirected cycle of order $n$. Let $B$ be a block of $C$ of maximal length and set $C=a_0a_1...a_{n-1}a_0$ such that $B$ is forward and it ends at $a_{\delta^-(T)+2}$. If $a_{n-1}\rightarrow a_0\rightarrow a_1$, then $T$ contains $C$ if and only if $(T;C)\neq A_{13}$.
\end{lemma}
\begin{proof}
By theorem \ref{3.1}, we may suppose that $\delta^-(T)\leq \delta^+(T)$. Set $\delta=\delta^-(T)$. Let $v\in T$ be a vertex such that $d^-(v)=\delta$ and denote by $T_1=T(N^-(v))$ and $T_2=T(N^+(v))$. Rotate $C$ until $a_{n-1}\rightarrow a_{n-2}$ for the first time. We are sure now that $a_{n-2}\leftarrow a_{n-1}\rightarrow a_0\rightarrow a_1$ and $a_{\delta}\rightarrow a_{\delta +1}\rightarrow a_{\delta +2}$. Suppose first that $(T_1;a_1...a_{\delta})$ is a Gr\"{u}nbaum's exception, let $a\in T_2$ such that $d^-_{T_2}(a)=\Delta^-(T_2)$. Set $T'_1=T_1\cup \{a\}$ and $T'_2=T_2-a$. We have two cases:\\
\textbf{Case $1$}: If $a_1\rightarrow a_2$. Suppose that $a_{n-3}\rightarrow a_{n-2}$. If $T'_1\neq 4B$ then $T'_1$ contains a path $x_0x_1...x_{\delta}\equiv a_0a_1...a_{\delta}$ such that $x_0,x_{\delta}\neq a$, if $x_0$ has an inneighbor in $T'_2$ origin of a path $P\equiv a_{n-1}a_{n-2}...a_{\delta+2}$ then $x_0...x_{\delta}vP^{-1}x_0\equiv C$. Otherwise, since $T_1$ is regular, $d^-_{T_1}(x_0)=\frac{|T_1|}{2}$, then, without loss of generality we may suppose that $d^-_{T'_2}(x_0)\geq 2$; in fact, if $\delta=5$ or $7$, it is clear, and if $\delta=3$, one can check that. By theorem \ref{2.1}, $(T'_2;a_{n-1}a_{n-2}...a_{\delta+2})$ is an exception. By the building lemmas, any vertex $w$ of $T'_2$ is origin of a path $Q\equiv a_0a_{n-1}...a_{\delta+2}$ in $T'_2\cup \{x_0\}$ with end distinct from $x_0$ and thus $wx_1...x_{\delta}vP^{-1}\equiv C$, unless $(T'_2;a_{n-1}a_{n-2}...a_{\delta+2})$ is one of the following biexceptions with respect to $x_0$: (From now on, the possible biexceptions are determined by the number of blocks of the path, the lengths of some blocks and the outdegree and indegree of $x_1$ in $T_2$)\\
Exc $(0,1)$: Since $d^+_{T_2}(x_2)=2$ then $x_1=a$. We have $x_03x_2x_3v1a2x_0\equiv C$.\\
Exc $(4,1)$: $d^-_{T'_2}(x_1)=1$ so $x_1\neq a$ and $d^+_{T'_2}(x_{\delta})\leq 3$ then $x_{\delta}$ is an origin of a path $P\equiv a_1a_0a_{n-1}...a_{\delta+2}$ in $(T'_2)\cup \{x_0,x_{\delta}\}$ with end distinct from $x_0$ so $x_{\delta}x_{\delta-1}...x_1vP^{-1}\equiv a_1a_2...a_0a_1=C$.\\
Exc $(8,1)$: Clearly $x_1\neq a$ so if $3\rightarrow x_{\delta}$ it is solved as above and if $x_{\delta}\rightarrow 3$ we have $x_0x_1...x_{\delta}341v25x_0\equiv C$.\\
Exc $(9,1)(1)$: If $d^+_{T'_2}(x_{\delta})\geq 1$ then it is clear that there is an outneighbor of $x_{\delta}$ origin of a path $P\equiv a_{\delta+1}...a_{n-1}$ in $(T'_2)\cup \{v\}$ with end and inneighbor of $x_0$, so $x_0...x_{\delta}Px_0\equiv C$. Else, $\delta=3$, if $x_1\neq a$, since $5$ is an origin of a path $P\equiv a_0a_{n-1}...a_{\delta+2}$ in $T'_2\cup \{x_0\}$ with end distinct from $x_0$ then $5x_3x_2x_1vP^{-1}\equiv C$, so $x_1=a$ and there is an inneighbor $u$ of $x_{\delta}$ origin of a path $P\equiv a_0a_{n-1}...a_{\delta+2}$ in $(T'_2)\cup \{x_2\}$ with end distinct from $x_2$, so $ ux_{\delta}ax_0vP^{-1}u\equiv C$.\\
Exc $(9,1)(2)$: $x_1\neq a$, so, as above, we may suppose that $V(T'_2)-S\subseteq N^+(x_{\delta}$ and then $x_0...x_{\delta}413v25x_0\equiv C$.\\
Exc $(9,1)(3)$ and $(4)$: $x_1\neq a$, so, as above, we may suppose that $V(T'_2)-S\subseteq N^+(x_{\delta}$ and then $x_0...x_{\delta}513v24x_0\equiv C$.\\
Exc $(15,1)$: $x_1\neq a$, so, as above, we may suppose that $V(T'_2)-S\subseteq N^+(x_{\delta}$ and then $x_0...x_{\delta}4532v61x_0\equiv C$.\\
Exc $(16,1)(1), (2)$ and $(3)$: $x_1\neq a$, so, as above, we may suppose that $V(T'_2)-S\subseteq N^+(x_{\delta}$ and then $x_0...x_{\delta}123v654x_0\equiv C$.\\
Exc $(16,1)(4)$: $x_1\neq a$, so, as above, we may suppose that $V(T'_2)-S\subseteq N^+(x_{\delta}$ and then $x_0...x_{\delta}123v456x_0\equiv C$.\\
Exc $(32,1)$: $x_1\neq a$, so, as above, we may suppose that $V(T'_2)-S\subseteq N^+(x_{\delta}$ and then $x_0...x_{\delta}136v425x_0\equiv C$.\\
Exc $(33,1)$: $x_1\neq a$, so, as above, we may suppose that $V(T'_2)-S\subseteq N^+(x_{\delta}$ and then $x_0...x_{\delta}463v7521x_0\equiv C$.\\
Exc $(E_8(5),1)$: $x_1\neq a$, so, as above, we may suppose that $V(T'_2)-S\subseteq N^+(x_{\delta}$ and then $x_0...x_{\delta}123vwux_0\equiv C$.\\
And if $T'_1=4B$ ($a_0...a_7=(2,1,3,1)$), if there exists $x\in T_1$ such that $d^-_{T'_2}(x)\geq 2$, let $y_1...y_4\equiv a_1...a_4$ be a path in $T'_1$ of origin $x$, if $y_1$ has an inneighbor $u$ origin of a path $P\equiv a_0a_{n-1}...a_6$ then $uy_1...y_4vP^{-1}\equiv C$, otherwise, $(T'_2;a_0a_{n-1}...a_6)$ is an exception. By the building lemmas, we may assume that $(T'_2;a_0a_{n-1}...a_6)$ is one of the following biexceptions with respect to $y_1$: (We can suppose that $N^-_{T'_2}(y_4)\subseteq S$, since otherwise, $y_4$ has an inneighbor origin of a path $P\equiv a_1a_0...a_6$ in $T'_2\cup \{y_1\}$ with end distinct from $y_1$, so $y_4y_3y_2vP^{-1}y_4\equiv a_2a_3...a_2=C$)\\
Dual Exc $(0,0)$: $d^-(3)=5$, a contradiction.\\
Dual Exc $(4,0)(1)$: We have $15y_2ay_1v342y_41\equiv C$.\\
Dual Exc $(4,0)(2)$: We have $43y_1ay_4v125y_24\equiv C$.\\
Dual Exc $(9,0)$: We have $54y_1ay_4v3y_2125\equiv C$.\\
Dual Exc $(15,0)$: We have $21y_1ay_4v5346y_22\equiv C$.\\
Dual Exc $(33,0)$: We have $21y_1ay_4v46573y_22\equiv C$.\\
Dual Exc $(E_1(4),0)$: We have $21y_1ay_4v43y_22\equiv C$.\\
Dual Exc $(E_8(5),0),(E_9(7),0),(E_{10}(9),0)$: Let $P$ be an antidirected outpath in $T'_2-u$, we have $uy_1...y_4Pvu\equiv C$.\\
So $d^-_{T'_2}(x)=1$ for all $x\in T_1$, let $u_1u_2u_3$ be a directed outpath in $T_1$, let $w_1\in N^+_{T_2}(u_3)$, we have $d^+_{T_2-w_1}(u_1)\geq |T_2-w_1|-2$ so $u_1$ has an outneighbor origin of a path $P\equiv a_{n-2}a_{n-3}...a_4$ in $T_2-w_1$ and thus $u_2u_3w_1vP^{-1}u_1u_2\equiv C$ unless $T_2-w_1$ is a transitive triangle, set $V(T_2)=\{a,w_1,w_2,w_3\}$ where $w_2\rightarrow w_3$. If $a\rightarrow w_3$ then $u_2u_3w_1vw_2w_3u_1au_2\equiv C$, so $w_3\rightarrow a$ and then we may suppose that $w_1\rightarrow a$. If $w_3\rightarrow u_1$ then $u_2u_3w_1vw_2w_3u_1au_2\equiv C$, so $u_1\rightarrow w_3$. Suppose that $w_2\rightarrow w_1$ then $u_2\rightarrow w_1$ since otherwise $u_3u_1w_3vw_2w_1u_2au_3\equiv C$, and thus $u_3\rightarrow w_3$ since otherwise $u_1u_2w_1vw_2w_3u_3au_1\equiv C$, but now $d^-(w)=4$ then $w_1\rightarrow u_1$ and so $u_2u_3w_3vw_2w_1u_1au_2\equiv C$. So we may assume that $w_1\rightarrow w_2$ and, without loss of generality, $T(w_1,w_2,w_3)=3A$. We remark that $u_i\rightarrow w_j$ implies that $u_{i+1}\rightarrow w_{j-1}$, so if $u_2\rightarrow w_3$ then $u_3vaw_1w_2u_1w_3u_2u_3\equiv C$ and if $u_2\rightarrow w_1$ then $u_3vaw_3w_1u_1w_2u_2u_3\equiv C$.\\
So we can assume that $a_{n-2}\rightarrow a_{n-3}$. Let $x_{n-1}x_0x_1...x_{\delta-2}\equiv a_{n-1}a_0a_1...a_{\delta-2}$ in $T_1$ such that $d^+_{T_2}(x_{n-1})=max \{d^+_{T_2}(x);x\in T_1\}$. If $x_{n-1}$ has an outneighbor origin of a path $P\equiv a_{\delta-1}...a_{n-3}$ in $T_2$ then $x_0x_1...x_{\delta-2}Pvx_{n-1}x_0\equiv C$. Otherwise, either $d^+_{T_2}(x_{\delta-2})\geq 2$ and $(T_2;a_{\delta-1}...a_{n-3})$ is a biexception with respect to $x_{\delta-2}$ or $d^+_{T_2}(x)=1$ for all $x\in T_1$. In the first case, $(T_2;a_{\delta-1}...a_{n-3})$ is one of the following biexceptions with respect to $x_{\delta-2}$:\\
Dual Exc $(E_1(5),1)$: We have $v1w23ux_{\delta-2}...x_0x_{n-1}v\equiv C$.\\
Dual Exc $(E_5(5),1)$: We have $13x_142x_{n-1}x_05v1\equiv C$.\\
Now, if $d^+_{T_2}(x)=1$ for all $x\in T_1$. Let $a\in T_2$ and $x_0x_1...x_{\delta}\equiv a_0a_1...a_{\delta}$ in $T_1\cup \{a\}$, we have $d^-_{T_2}(x_0)\geq |T_2-a|-1$ so there is an inneighbor of $x_0$ origin of a path $P\equiv a_{n-1}a_{n-2}...a_{\delta+2}$ in $T_2-a$ and thus $x_0...x_{\delta}vP^{-1}x_0\equiv C$.\\
\textbf{Case $2$}: If $a_2\rightarrow a_1$. Suppose that $a_{n-3}\rightarrow a_{n-2}$, $T'_1$ contains a path $x_0x_1...x_{\delta}\equiv a_0a_1...a_{\delta}$ such that $x_0,x_{\delta}\neq a$, if $x_0$ has an inneighbor in $T'_2$ origin of a path $P\equiv a_{n-1}a_{n-2}...a_{\delta+2}$ then $x_0...x_{\delta}vP^{-1}x_0\equiv C$. Otherwise, either $d^-_{T'_2}(x_0)\geq 2$ and $(T'_2;a_{n-1}a_{n-2}...a_{\delta+2})$ is a biexception with respect to $x_0$ or $T'_1=4B$ and $d^-_{T'_2}(x)=1$ for all $x\in T_1$. In the first case, $(T'_2;a_{n-1}a_{n-2}...a_{\delta+2})$ is one of the following biexceptions with respect to $x_0$:\\
Exc $(1,1)$: There is a vertex $b\in T_2$ such that $d^-_{T_2}(b)\geq 2$ and $T_2-b\neq 4A$, so $C$ can be found as above unless $T_1\cup \{b\}=4B$, in this case there is an origin of $P=-(1,1,1)$ in $T_1\cup \{b\}$ outneighbor of $3$ and so $v1423Pv\equiv C$.\\
Exc $(6,1)$: We have $x_{\delta}\rightarrow 1$ so $x_0...x_{\delta}142v35x_0\equiv C$.\\
Exc $(10,1)$: Since $T_2$ cannot be $6F$ nor $6H$, then if $T_1=5A$, let $y_1...y_5\equiv \overline{a_{n-1}a_{n-2}...a_{\delta+2}}$ in $T_1$, $d^-_{T_2}(y_1)\geq 3$ implies that there is an inneighbor $u$ of $y_1$ origin of $P\equiv a_0a_1...a_{\delta}$ in $T_2$ and thus $Pvy_5...y_1u\equiv \overline{C}$. So $T$ contains $C$. Otherwise, $T_1=3A$, if $x_{\delta}\rightarrow 1$ then $x_0...x_{\delta}12v453x_0\equiv C$, so we can assume that $x_1=a$ since otherwise, taking $x_{\delta}...x_1\equiv a_1...a_{\delta}$, $(T_2,a_{n-1}...a_{\delta+2})$ is not a biexception with respect to $x_0$. $T_2\setminus \{2\}\neq 5E$, so let $y_0...y_{\delta}\equiv a_0...a_{\delta}$ in $T_1\cup \{2\}$ such that $3,4$ or $5$ is an inneighbor of $y_1$, which is an origin of $P=(1,2,1)$ in $T_2\setminus \{2\}$, we have $y_0...y_{\delta}vP^{-1}y_0\equiv C$.\\
Exc $(12,1)$: We have $x_{\delta}\rightarrow 3$ so $x_0...x_{\delta}354v21x_0\equiv C$.\\
Exc $(22,1)$: There is $y_0...y_{\delta}\equiv a_{n-\delta -1}a_{n-\delta}...a_{n-1}$ in $T'_1$ having an inneighbor origin of $P\equiv a_{n-\delta -2}...a_1$ in $T'_2$, so $vP^{-1}y_0...y_{\delta}v\equiv C$.\\
Exc $(24,1)$: We have $x_{\delta}\rightarrow 1$ so $x_0...x_{\delta}1523v64x_0\equiv C$.\\
Exc $(26,1)$: Let $y_1...y_{\delta}\equiv \overline{a_3...a_{\delta+2}}$ in $T_1$, $(\overline{T_2};\overline{a_1a_0a_{n-1}...a_{\delta+4}})$ is not an exception so there is an inneighbor of $y_1$ origin of a path $P\equiv \overline{a_1a_0a_{n-1}...a_{\delta+4}}$ in $T_2$, we have $Py_1...y_{\delta+2}v\equiv \overline{C}$. Thus $T$ contains $C$.\\
Now, if $T'_1=4B$ and $d^-_{T'_2}(x)=1$ for all $x\in T_1$. If there exists $v'\in T$ such that $T(N^-(v'))$ is a transitive triangle, set $T"_1=T(N^-(v'))=\{u_1,u_2,u_3\}$ and $T"_2=T(N^+(v'))$ such that $u_1\rightarrow u_2\rightarrow u_3$, if $u_1$ has an inneighbor origin of a path $P\equiv a_{n-1}a_{n-2}...a_4$ in $T"_2$ then $u_1u_3u_2vP^{-1}u_1\equiv C$, otherwise, $(T_2;a_{n-1}a_{n-2}...a_4)$ is one of the biexceptions Exc $(4,1),(9,1),(15,1),(16,1)$ and $(33,1)$ with respect to $u_1$, but in these cases $(T_2;a_{n-4}a_{n-5}...a_1)$ is not an exception so there exists an origin $x$ of a path $P\equiv a_{n-4}a_{n-5}...a_1$ in $T_2$ such that $xu_1\equiv a_{n-4}a_{n-3}$, thus $vP^{-1}u_1u_3u_2v\equiv C$. So we may assume that for any vertex $x\in T$ of minimal indegree, $T(N^-(x))=3A$, thus the inneighbors $w_1$, $w_2$ and $w_3$ of $x_0$, $x_1$ and $x_3$ are pairwasly distinct, so $(T'_2;a_{n-1}a_{n-2}...a_5)$ is one of the exceptions $1,18,22,24,26,34,38$ and $47$, if it is Exc $18$ then $|\displaystyle \bigcup_{x\in T_1} N^-_{T_2}(x)|\geq 4$ thus, without loss of generality, we may suppose that there exists an inneighbor of $x_0$ origin of a path $P\equiv a_4a_3...a_{n-2}$ in $T_2$ and therefore $x_0x_1x_3vP^{-1}x_0\equiv a_5a_6...a_5=C$. Otherwise, there exists $x\in T'_2$ origin of a path $P\equiv a_{n-5}a_{n-6}...a_1$ such that $xx_3\equiv a_{n-5}a_{n-4}$ and thus $vP^{-1}x_3x_2x_1x_0v\equiv C$.\\
So we may assume that $a_{n-2}\rightarrow a_{n-3}$, set $T_1=\{u_1,u_2,u_3\}$, if $u_3$ has an inneighbor origin of a path $P\equiv a_2a_3...a_{n-3}$ in $T_2$ then $u_2u_3Pvu_1u_2\equiv C$, otherwise, without loss of generality, since there is no biexceptions such that $P=-(1,3,Q)$, we may suppose that $s^+_{T_2}(\displaystyle \bigcup_{i=1}^3 N^-_{T_2}(u_i))\leq 3$; if $|N^+_{T_2}(T_1)|\geq 2$, let $u,w\in N^+_{T_2}(T_1)$ such that $w$ is an origin of a path $P\equiv a_3a_4...a_{n-3}$ in $T_2-u$, we have $u_2uu_3Pvu_1u_2\equiv C$, otherwise, set $\{u\}=N^+_{T_2}(T_1)$, we have $u_1uu_2u_3vI_{T_2-u}u_1\equiv C$.
\\So we can assume that $(T_1;a_1...a_{\delta})$ is not a Gr\"{u}nbaum's exception, let $x_1...x_{\delta}\equiv a_1...a_{\delta}$ in $T_1$ such that $d^-_{T_2}(x_1)=max \{d^-_{T_2}(x); x\in Or(a_1...a_{\delta};T_1)\}$. If there is an inneighbor $u$ of $x_1$ origin of $P\equiv a_0a_{n-1}...a_{\delta+2}$ in $T_2$ then $ux_1...x_{\delta}vP^{1}\equiv C$. Otherwise, if $a_1\rightarrow a_2$, then, either $d^-_{T_2}(x_1)\geq 3$ and $(T_2,a_0a_{n-1}...a_{\delta+2})$ is an exception or $d^-_{T_2}(x_1)=2$. If $d^-_{T_2}(x_1)\geq 3$ and $(T_2,a_0a_{n-1}...a_{\delta+2})$ is an exception, by the building lemmas we may suppose that $(T_2;a_0a_{n-1}...a_{\delta+2})$ is one of the following biexceptions with respect to $x_1$:\\
Dual Exc $(0,0)(1)$: We have $x_{\delta}\rightarrow 3$ so $12x_1...x_{\delta}3v1\equiv C$.\\
Dual Exc $(1,0)(3)$: $4\rightarrow x_2$, a contradiction.\\
Dual Exc $(4,0)(2)$: If $d^+_{T_1}(x_1)=2$ then $d^+(x_1)=5$ and $T(N^+(x_1))\neq 5A$ since $\{x_2,v\}\rightarrow \{4,5\}$ that is solved unless $T(N^+(x_1))=\overline{5E}$ which is the next biexception; in any case $T$ contains $C$. Otherwise $T_1=3A$ and then $4$ is an inneighbor of  vertex in $T_1$, so $T$ contains $C$; or $|T_1|=2$ and then $T(N^-(4))\neq 5A$, so $T$ contains $C$.\\
Dual Exc $(9,0)(1),(15,0)(1),(16,0)(1)$ and $(33,0)$: We have $d^-(x_2)=\delta$ and $(T(N^+(x_2)),a_0a_{n-1}...a_{\delta+2})$ is not an exception since $T(N^+(x_2))$ has $v$ as a minimal vertex. Thus $T$ contains $C$.\\
Dual Exc $(26,0)(1)$: If $d^+(x_2)=6$ then $N^+(x_2)\neq \overline{6H}$ and if $d^+(x_2)=5$ then $T(N^+(6)\neq \overline{6H}$; thus $T$ contains $C$.\\
Now if $d^-_{T_2}(x_1)=2$, by lemma \ref{2.4}, $(T_1;a_1...a_{\delta})=(4B;(1,1,1)), (4B;(2,1)), (3A;(2))$ or $|T_1|=2$. If $|T_1|\geq 3$, by the building lemmas we may assume that $(T_2;a_0a_{n-1}...$\\$a_{\delta+2})$ is one of the following biexceptions with respect to $x_1$:\\
Dual Exc $(0,0)(1)$: $d^-(3)\geq 4$ is a contradiction.\\
Dual Exc $(1,0)$: $d^-(4)\geq 5$ is a contradiction.\\
In the other possible biexceptions, $N^-(x_1)\setminus S\neq \emptyset$ and so we can take $x_{\delta}$ as origin of $y_1...y_{\delta}\equiv a_1...a_{\delta}$, where $y_2=x_1$, the result follows.\\
And if $|T_1|=2$, $d^-(x_1)=2$ and $T(N^+(x_1))=\{v\}\rightarrow (N^+(x_1)\setminus \{v,x_2\})$ which implies that $(T_2;a_0a_{n-1}...a_{\delta+2})$ can be only one of the following biexceptions with respect to $x_1$:\\
Dual Exc $(31,0),(48,0)$: We have $d^-(x_2)=2$ and $T(N^+(x_2))$ is not one of the tournaments of these biexceptions.\\
Dual Exc $(E_{8,9,10}(|T_2|),0)$ and $(E'_{8,9,10}(|T_2|),0)$: There is an inneighbor of $x_1$ origin of $P\equiv a_{n-1}a_{n-2}...a+{\delta+1}$ in $T_2$ so $x_1x_2vP^{-1}x_1\equiv C$.
So, we may assume that $a_2\rightarrow a_1$. If $d^-_{T_2}(x_1)\geq 2$ and $d^-_{T_2}(x_2)=|T_2|$ then $(T_2;a_0a_{n-1}...a_{\delta+2})$ is an exception. If $(T_2;a_0a_{n-1}...a_{\delta+2})$ is not one of the dual of the exceptions $0,1,4,7,18,19,$\\$22,33,E_1(|T_2|)$ and $E_1(6)$ with $X=3A$, then by lemma \ref{2.8}, $N^+_{T_2}(x_1)$ contains an origin of $w$ of a path $P\equiv a_0a_{n-1}...a_{\delta+2}$ in $\overline{T_2}$ and so $wx_1vQx_2P^{-1}\equiv C$ in $\overline{T}$, where $Q\equiv a_3...a_{\delta}$ in $T_1\setminus \{x_1,x_2\}$. So, we may assume that $(T_2;a_0a_{n-1}...a_{\delta+2})$ is one of the following exceptions:\\
Dual Exc $0$: Set $N^-_{T_2}(x_1)=\{1,2\}$, we have $x_13v1x_22x_1\equiv C$.\\
Dual Exc $4,19$ or $33$: Any vertex of $T_2$ is an origin of a path $P\equiv a_1a_0...a_{\delta+2}$ in $T_2\cup \{x_2\}$, so if $T_1-x_2$ contains a path $y_2...y_{\delta}\equiv a_2...a_{\delta}$, we have $y_2...y_{\delta}vP^{-1}y_2\equiv a_2a_3...a_2=C$. Otherwise, $T_1-x_2=3A$, let $y_1y_2y_3$ be a directed outpath in $T_1-x_2$ and $u\in N^+_{T_2}(y_1)$, any vertex of $(T_2-u)\cup \{x_2\}$ is an origin of a path $P\equiv x_6x_7...x_1$ in $(T_2-u)\cup \{x_2\}$, thus $vuy_1y_2y_3Pv\equiv a_2a_3...a_2=C$.\\
Dual Exc $E_1(4)$: If $\delta=2$ then $x_14321vx_2x_1\equiv C$. If $\delta=3$ then $x_2x_1O_{T_2}vx_3x_2\equiv C$. If $\delta=4$, let $P\equiv a_{n-5}a_{n-6}...a_1$ be a path in $T_2$ and $Q=(1,2)$ be a path in $T_1$ of origin $x_2$, we have $QvP^{-1}x_2\equiv a_{n-4}a_{n-5}...a_{n-4}=C$.\\
Dual Exc $E_1(6)$: $\delta\geq 4$, if $\delta=4$, let $P=-(1,4)$ be a path in $T_2$, we have $x_1Pvx_4x_3x_2x_1\equiv C$. If $\delta=5$, let $Q=-(2,2)$ be a path in $T_1$, we have $vQO_{T_2}v\equiv a_{n-4}a_{n-3}...a_{n-4}=C$. If $\delta=6$, let $Q=-(1,2,1,1)$ be a path in $T_1$ of origin $x$ and $P=(2,3)$ be a path in $T_2$ of origin an outneighbor of $x$, we have $QvP^{-1}x\equiv a_{n-2}a_{n-1}...a_{n-2}=C$.\\
Otherwise, either $d^-_{T_2}(x_1)\geq 2$, $d^-_{T_2}(x_2)<|T_2|$ and $(T_2,a_0a_{n-1}...a_{\delta+2})$ is an exception or $d^-_{T_2}(x_1)=1$. In the first case, by the building lemmas we may assume that it is one of the following biexceptions with respect to $x_1$:\\
Dual Exc $(0,0)(2)$: If there is a vertex $x\in T$ such that $d^-(x)=3$ and $T(N^-(x))\neq 3A$, or $d^+(x)=3$ and $T(N^+(x))\neq 3A$, then $T$ contains $C$. Otherwise $(T;C)=(6Z;(2,1,2,1))$.\\
Dual Exc $(2,0)$: We have $N^-(1)\neq 4A$ so $\overline{T}$ contains $C$ and thus $T$ contains $C$.\\
Dual Exc $(4,0)$: If $x_{\delta}\rightarrow 1,5$ or $4$ (say $i$) then $1$ or $2$ is an origin of a path $P=-(1,1)$ in $T_2\setminus \{i-1,i\}$, we have $x_1...x_{\delta}i(i-1)vPx_1\equiv a_1...a_{n-1}a_0a_0=C$. So we can suppose that $1,4,5\in N^-(x_{\delta})$. If $a_{\delta-1}\rightarrow a_{\delta}$ then $1x_1...x_{\delta-1}v2x_{\delta}3541\equiv C$. Otherwise, if $x_{\delta-1}\leftarrow 1$ or $2$, say $1$, we have $2x_1...x_{\delta-1}13x_{\delta}45v2\equiv C$. Else, if $3\rightarrow x_{\delta-1}$ we have $x_1vx_2...x_{\delta-1}3514x_{\delta}2x_1\equiv C$. Assume now that $x_{\delta-1}\rightarrow 3$. Suppose first that $4\rightarrow x_{\delta-1}$, if $T(N^+(4))\neq 5A$ then $T$ contains $C$; else, we have $\{5,x_2\}\rightarrow x_{\delta}$ so $T(N^+(5))\neq 5A$. So we can assume that $x_{\delta-1}\rightarrow 4$, we have $N^+(x_{\delta-1})=\{1,2,3,4,v\}$ where $T(N^+(x_{\delta-1}))\neq 5A$; so $T$ contains $C$.\\
Dual Exc $(9,0)(1)$: If $x_{\delta}\rightarrow 1$ then $2x_1...x_{\delta}1453v2\equiv C$. Otherwise, if $d^+_{T_2}(x_{\delta})\geq 1$, let $u\in N^+_{T_2}(x_{\delta})$, there exists $x\in N^-_{T_2}(x_1)$ origin of a path $P=-(1,1)$ in $T_2\setminus \{u,w\}$ where $w\in N^+_{T_2}(u)$; we have $xx_1...x_{\delta}uwvP^{-1}\equiv C$. So we can suppose that $T_2\rightarrow x_{\delta}$. We have $N^-(x_2)\neq 5E$, so $\overline{T}$ contains $C$ and thus $T$ contains $C$.\\
Dual Exc $(9,0)(2)$: If $N^+(x_{\delta})\cap \{2,3,4,5\}\neq \emptyset$ then there exist $u\in N^+_{T_2}(x_{\delta})$ and $x\in N^-_{T_2}(x_1)$ such that $x$ is an origin of a path $P=-(1,1)$ in $T_2\setminus \{u,w\}$ where $w\in N^+_{T_2}(u)$; we have $xx_1...x_{\delta}uwvP^{-1}\equiv C$. So we can assume that $\{2,3,4,5\}\subseteq N^-(x_{\delta})$. If $a_{\delta-1}\rightarrow a_{\delta}$ then $4x_1...x_{\delta-1}v312x_{\delta}54\equiv C$. Otherwise, if $1$ or $2\rightarrow x_{\delta-1}$, say $1$, then $4x_1...x_{\delta-1}15x_{\delta}32v1\equiv C$; and if $4$ or $5\rightarrow x_{\delta-1}$, say $4$, then $5x_1...x_{\delta-1}43x_{\delta}21v5\equiv C$. So we may suppose that $N^+(x_{\delta-1})=\{1,2,4,5,v\}$ so $T(N^+(x_{\delta-1}))\neq \overleftarrow{5E}$ thus $T$ contains $C$.\\
Dual Exc $(15,0)$: If $x_{\delta}\rightarrow 6$ then $1x_1...x_{\delta}63v5421\equiv C$ and if $N^+(x_{\delta})\cap \{1,2,3\}\neq \emptyset$, say $1$, then $6x_1...x_{\delta}13v2456\equiv C$. So we may suppose that $\{1,2,3,6\}\subseteq N^-(x_{\delta})$, then $T(N^+(6))\neq \overline{6C}$ since otherwise $x_{\delta}\equiv 4$ or $5$ but $\{1,2,3\}\nsubseteq N^+(x_1)$ nor $N^-(x_1)$. Therefore $T$ contains $C$.\\
Dual Exc $(16,0)$: If $x_{\delta}\rightarrow 6$ then $4x_1...x_{\delta}62v5314\equiv C$. If $\{1,2,3\}\cap N^+(x_{\delta})\neq \emptyset$, say $1\in N^+(x_{\delta})$, then $4x_1...x_{\delta}13v6524\equiv C$. If $\{4,5\}\cap N^+(x_{\delta})\neq \emptyset$, say $4\in N^+(x_{\delta})$, then $5x_1...x_{\delta}46v3125\equiv C$. So we may assume that $N^-(x_{\delta})=T_2$, then $a_{\delta}\rightarrow a_{\delta-1}$. One can find $P\equiv a_{\delta}...a_{n-1}a_0$ in $T_2\cup \{x_{\delta},v\}$ with origin an inneighbor of $x_{\delta-1}$ and end an inneighbor of $x_1$; we have $x_1...x_{\delta-1}Px_1\equiv C$.\\
Dual Exc $(33,0)$: If $\{2,3,...,7\}\cap N^+(x_{\delta})\neq \emptyset$ then there exist $u\in N^+_{T_2}(x_{\delta})$ and $x\in N^-_{T_2}(x_1)$ such that $x$ is an origin of a path $P\equiv a_0a_{n-1}...a_{\delta+4}$ in $T_2\setminus \{u,w\}$ where $w\in N^+_{T_2}(u)$, then $xx_1...x_{\delta}uwvP^{-1}\equiv C$. So assume that $\{2,...,7\}\subseteq N^-(x_{\delta})$, if $a_{\delta-1}\rightarrow a_{\delta}$ then there exist $x\in N^-_{T_2}(x_1)$ origin of $P\equiv a_0a_{n-1}...a_{\delta+1}$ in $T_2\cup \{x_{\delta}\}$; we have $xx_1...x_{\delta-1}vP^{-1}\equiv C$. Otherwise, one can find $P\equiv a_{\delta}...a_{n-1}a_0$ in $T_2\cup \{x_{\delta},v\}$ with origin an inneighbor of $x_{\delta-1}$ and end an inneighbor of $x_1$; we have $x_1...x_{\delta-1}Px_1\equiv C$.\\
Dual Exc $(E_1(|T_2|),0)$: We have $T(N^+(2))\neq \overline{F_1}$ nor $\overline{F_5}$, so $T$ contains $C$.\\
Dual Exc $(E_3(|T_2|),0)$: If $x_{\delta}\rightarrow X$ then if $X\neq 3A$, let $Q=(1,|X|-2)$ be a Hamiltonian path of $X$, we have $1x_1...x_{\delta}Q2v1\equiv C$ and if $X=3A=\{3,4,5\}$ then $3x_1...x_{\delta}5421v3\equiv C$. So we may suppose that $X\nsubseteq N^+(x){\delta})$. If $x_{\delta}\rightarrow 2$ then $1x_1...x_{\delta}2I_Xv1\equiv C$. Assume now that $2\rightarrow x_{\delta}$, we have $T(N^+(2))\neq \overline{F_1}$, $\overline{F_3}$ nor $\overline{F_5}$, so $T$ contains $C$.\\
Dual Exc $(E_5(|T_2|),0)$: If $x_{\delta}\rightarrow 1$ then $2x_1...x_{\delta}1I_{T_2\setminus \{1,2\}}v2\equiv C$. So assume that $1\rightarrow x_{\delta}$. If $N^+(1)\neq \overline{F_1}$ nor $\overline{F_5}$ then $T$ contains $C$. Otherwise, $x_{\delta}\rightarrow x_1$, so $21x_2...x_{\delta}x_1I_{T_2\setminus \{1,2\}}v2\equiv C$.\\
Dual Exc $(E_8(|T_2|),0),(E_9(|T_2|),0)$ and $(E_{10}(|T_2|),0)$: Set $X=\{u\rightarrow w\}$. Without loss of generality, any vertex $x$ of $3A$, $5A$ and $7A$ respectively is an origin of a path $P\equiv a_{\delta+1}...a_{n-1}a_0$ in $T_2\cup \{v\}$; so if there is a such $x\in N^+(x_{\delta})$ then $x_1...x_{\delta}Px_1\equiv C$. Now we may assume that $3A$, $5A$ or $7A\rightarrow x_{\delta}$ respectively. If $a_{\delta-1}\rightarrow a_{\delta}$ then any vertex of $3A$, $5A$ or $7A$ is an origin of a path $P\equiv a_{\delta+1}...a_{n-1}a_0$ in $T_2\cup \{x_{\delta}\}$ with end $w$ so $wx_1...x_{\delta-1}vP\equiv C$. Otherwise, $x_{\delta-1}$ has an inneighbor $x$ in $T_2-u$. $x$ is an origin of a path $P\equiv a_{\delta}...a_{n-1}a_0$ in $T_2\cup \{x_{\delta},v\}$ with end $u$, so $ux_1...x_{\delta-1}P\equiv C$.\\
Now if $d^-_{T_2}(x_1)=1$, by lemma \ref{2.3}, $(T_1;a_1...a_{\delta})=(\overline{4B};-(1,2))$ or $T_1$ is transitive and $a_1...a_{\delta}$ is directed (it is clear that $\delta=2$ or $3$). Suppose first that $(T_1;a_1...a_{\delta})=(\overline{4B};-(1,2))$. Set $\overline{C}=u_0u_1...u_{n-1}u_0$ such that $u_0\equiv a_{\delta+1}$ and $u_1\equiv a_{\delta}$. We have $u_1...u_4=(2,1)$, so let $y_1...y_4=(2,1)$ in $T_1$. We have $d^-_{T_2}(y_1)\geq 3$, so $y_1$ has an inneighbor $u\in T_2$ origin of a path $P\equiv u_0u_{n-1}...a_{\delta+2}$ in $T_2$ and so $uy_1...y_4vP^{-1}\equiv C$, or $(T_2;u_0u_{n-1}...u_{\delta+2})$ is an exception. The building lemmas give the result unless $(T_2;u_0u_{n-1}...u_{\delta+2})$ is a biexception with respect to $y_1$. Since any vertex of $T_1$ distinct from $x_1$ can be an origin of a path $(2,1)$, we may assume that $|N^-(y_1)|\geq 3$ and $|S|\geq 3$, and since the path $u_1u_0...u_{\delta+2}$ has an even number of blocks with first block of length $2$, there is no biexception verifying these conditions.\\
Now if $T_1$ is a transitive triangle, if $C$ has a block of length $1$, let $s=min\,\{i; a_i$ \textit{is contained in a block of length one}$\}$. Without loss of generality we can suppose that $s=\delta+2$, indeed, by considering $C$ and $\overline{C}$, we may assume that $a_{s+1}\rightarrow a_s$, rotate $C$ until $a_{\delta+2}=a_s$. Thus $a_{\delta+3}\rightarrow a_{\delta+4}$, let $w\in N^+_{T_2}(x_2)$, there is an outneighbor of $x_1$ origin of a path $P\equiv a_1a_0a_{n-1}...a_{\delta+4}$ in $T_2-w$, and thus $x_1x_3x_2wvP^{-1}x_1\equiv a_2a_3...a_1a_2=C$ unless $|T_2|=3$, we may suppose that $T=7K$ and so $13654721=(2,2,2,1)\equiv C$. Assume now that $C=(2,2,...,2)$. If there exists $w\in N^+_{T_2}(x_2)$ such that $d^-_{T_2}(w)\geq 1$, let $t\in T_2$ be an inneighbor of $w$, and $P\equiv a_1a_0a_{n-1}...a_{\delta+5}$ be a path in $T_2\setminus \{w,t\}$ of origin an outneighbor of $x_1$ (this is possible since $d^+_{T_2\setminus \{w,t\}}(x_1)\geq |T_2\setminus \{w,t\}|-1$), we have $x_1x_3x_2wtvP^{-1}x_1\equiv a_2a_3...a_2=C$. Else, $N^+_T{_2}(x_2)=\{w\}$ and $d^+_{T_2}(w)=0$, we have $d^-_{T_2}(x_3)\geq 3$ so there is an inneighbor $t\in T_2$ of $x_3$ origin of a path $Q=(2,1)$ containing $w$. Let $P\equiv a_{\delta+4}...a_{n-2}$ be a path of origin an inneighbor of $x_2$ in $T_2-V(Q)$, we have $vQ^{-1}x_3x_1x_2Pv\equiv a_{n-1}a_0a_1...a_{n-1}=C$.\\
The last case is when $|T_1|=2$, let $u$ be the inneighbor of $x_1$ in $T_2$, if $d^-_{T_2}(u)\geq 1$, take $x_1$ in the place of $v$ and set $T'_1=T(N^-(x_1))$ and $T'_2=T(N^+(x_1))$. Let $y_1y_2\equiv a_1a_2$ in $T'_1$, and let $w\in N^-_{T'_2}(y_1)$. Since $T'_2=v\rightarrow (T'_2-v)$, $w$ is an origin of a path $P\equiv a_0a_{n-1}...a_{\delta+2}$ in $T'_2$, and thus $wy_1y_2x_1P^{-1}\equiv C$, unless $(T'_2\setminus \{w,v\};a_{n-2}a_{n-3}...a_{\delta+2})$ is a Gr\"{u}nbaum's exception.
If $d^-_{T'_2-v}(w)\geq 1$, any inneighbor $t$ of $w$ in $T'_2-v$ is an origin of a path $P\equiv a_{n-1}a_{n-2}...a_{\delta+2}$ in $T'_2-w$, so $wy_1y_2x_1P^{-1}w\equiv C$, else, any outneighbor $t$ of $w$ in $T'_2-v$ is an origin of a path $P\equiv a_{n-4}a_{n-5}...a_1$ in $T'_2-w$, so $x_1P^{-1}wy_1y_2x_1\equiv C$. So we may suppose that $d^-_{T_2}(u)=0$. We have $d^-_{T_2-u}(x_2)\geq 2$ and so $s^+_{T_2-u}(N^-_{T_2}(x_2))\geq 3$, since otherwise there is an inneighbor of $x_2$ of outdegree $1$. If $s=min\,\{4\leq i\leq n-2; d^+_C(a_i)=2\}$ exists, let $P\equiv a_{n-1}a_{n-2}...a_{s+1}$ be a path in $T_2-u$ of origin an inneighbor of $x_2$ and $Q\equiv a_3a_4...4_{s-1}$ be a Hamiltonian path of $T_2\setminus (V(P)\cup \{u\})$, we have $x_2x_1uQvP^{-1}x_2\equiv C$. Else, $C=(2,2,2,1)$, let $P\equiv a_1a_2...a_{n-3}$ be a Hamiltonian path of $(T_2-u)\cup \{v\}$ with origin and end distinct from $v$, we have $x_1Pux_2x_1\equiv C$.    $\square$
\end{proof}
\begin{lemma}\label{3.4}
Let $T$ be a tournament of order $n$ with $\delta^-(T)\geq 2$ and $C$ a non-directed and non-antidirected cycle of order $n$. Let $B$ be a block of $C$ of maximal length and set $C=a_0a_1...a_{n-1}a_0$ such that $B$ is forward and it ends at $a_{\delta^-(T)+2}$. If $a_{n-1}\leftarrow a_0\leftarrow a_1$, then $T$ contains $C$.
\end{lemma}
\begin{proof}
By theorem \ref{3.1}, we may suppose that $\delta^-(T)\leq \delta^+(T)$. Set $\delta=\delta^-(T)$. Let $v\in T$ be a vertex such that $d^-(v)=\delta$ and denote by $T_1=T(N^-(v))$ and $T_2=T(N^+(v))$. Rotate $C$ until $a_{n-2}\rightarrow a_{n-1}$ for the first time. We are sure now that $a_{n-2}\rightarrow a_{n-1}\leftarrow a_0\leftarrow a_1$ and $a_{\delta}\rightarrow a_{\delta +1}\rightarrow a_{\delta +2}$. Suppose first that $(T_1;a_1...a_{\delta})$ is a Gr\"{u}nbaum's exception, we have two cases:\\
\textbf{Case $1$}: If $a_2\rightarrow a_1$. For $x_i\in T_1$ let $P_i\equiv a_0a_1...a_{\delta-1}$ in $T_1$. If there is an outneighbor of an $x_i$ origin of a path $P\equiv a_{n-1}a_{n-2}...a_{\delta+1}$ in $T_2$ then $P_ivP^{-1}x_i\equiv C$. Otherwise, either there exists $x_i\in T_1$ such that $d^+_{T_2}(x_i)\geq 2$ or $T_1=3A$ and $d^+_{T_2}(x_i)=1$ for all $x_i\in T_1$. In the first case, denote by $B'$ the block of $C$ containing the arc $a_{n-2}\rightarrow a_{n-1}$ and let $Q\equiv a_2a_3...a_{\delta+1}$ be a path in $T_1$ of origin $x_i$, say $x_1$.\\
If $|B'|=1$,, $(T_2;a_{n-1}a_{n-2}...a_{\delta+1})$ is an exception, so by the building lemmas it is one of the following biexceptions with respect to $x_1$:\\
Dual Exc $(1,1)(1),(2)$ and $(4)$: $x_2$ has an outneighbor origin of a directed outpath $P$ of order $3$ in $T_2$ with end distinct from $1$, so $x_1x_3x_2Pv1x_1\equiv C$.\\
Dual Exc $(1,1)(3)$: It is not possible for all $x_i$ since otherwise $d^-(2)=5$ contradiction.\\
Dual Exc $(6,1)$: , we have $14Qv2351\equiv C$.\\
Dual Exc $(10,1)$: $x_2\rightarrow 2$ contradiction.\\
Dual Exc $(12,1)$: $4$ is an outneighbor of an $x_i$ contradiction.\\
Dual Exc $(22,1)$: We have $x_i\leftarrow 4,5,6$ for all $i$, so $315Q64v23\equiv C$.\\
Dual Exc $(24,1)$: We have $x_i\leftarrow 5,6$ for all $i$, so $315Q64v23\equiv C$.\\
If $|B'|=2$, suppose that $s^-_{T_2}(N^+_{T_2}(x_1))=2$ and set $N^+_{T_2}(x_1)=\{u\rightarrow w\}$. If $N^+(x_2)\cap (T_2\setminus \{u,w\})\neq \emptyset$, let $t\in N^+(x_2)\cap (T_2\setminus \{u,w\})$, we have $|T_2|\geq 4$ so there exists $P\equiv a_{n-2}a_{n-3}...a_{\delta+1}$ in $T_2-t$ of origin distinct from $u$ and $w$, thus $t^*P_1vP^{-1}x_1t\equiv C$. Otherwise, by induction on $i$, we may suppose that $T_1\leftarrow T_2\setminus \{u,w\}$. If $d^+_{T_1}(u)\geq 1$ or $d^+_{T_1}(w)\geq 1$, assume without loss of generality that it is $u$, then there exists a path $R\equiv a_{n-1}a_0...a_{\delta-1}$ in $T_1\cup \{u\}$ of origin $x$ and end distinct from $u$, let $S\equiv a_{n-2}a_{n-3}...a_{\delta+1}$ be a path in $T_2-u$ of origin distinct from $w$, we have $RvSx\equiv a_{n-1}a_0...a_{n-1}=C$. Otherwise, if $T_2\setminus \{u,w\}$ contains a path $R\equiv a_{n-3}a_{n-4}...a_{\delta+1}$ then $u^*P_1vR^{-1}x_1wu\equiv C$, else, let $t\in T_2\setminus \{u,w\}$, $R\equiv a_{n-3}a_{n-4}...a_{\delta+2}$ in $T_2\setminus \{u,w,t\}$ and $S\equiv a_1...a_{\delta}$ in $(T_1-x_1)\cup \{t\}$ with origin and end distinct from $t$, we have $uSvR^{-1}x_1wu\equiv C$. So assume that $s^-_{T_2}(N^+_{T_2}(x_2))\geq 3$ and $(T_2;a_{n-1}a_{n-2}...a_{\delta+1})$ is an exception. By the building lemmas, it is one of the following biexceptions with respect to $x_1$:\\
Dual Exc $(5,1)$ and $(7,1)$: We have $2vP^{-1}_114352\equiv C$.\\
If $|B'|=3$, if $s^-_{T_2}(N^+_{T_2}(x_1))\geq 4$ then, since there is no a such biexception, $T$ contains $C$. So we may suppose that $s^-_{T_2}(N^+_{T_2}(x_1))\leq 3$, if there exists $u\in N^+_{T_2}(x_1)\cap N^+_{T_2}(x_2)$, let $w\in N^+_{T_2}(x_1)-u$, $t\in N^-_{T_2}(x_1)$ such that $(T_2\setminus \{u,w,t\},a_{n-4}a_{n-5}...a_{\delta+1})$ is not a Gr\"{u}nbaum's exception and $R\equiv a_{n-4}...a_{\delta+1}$ in $T_2\setminus \{u,w,t\}$, we have $u^*P_1vR^{-1}x_1wtu\equiv C$. Otherwise, let $t\in N^+_{T_2}(x_2)$, if $N^+_{T_2}(x_1)=\{u,w\}$, then there exists $R\equiv a_{n-4}a_{n-5}...a_{\delta+1}$ in $T_2\setminus \{u,t\}$ of origin distinct from $w$ thus $t^*P_1vR^{-1}x_1ut\equiv C$, so we may assume that $s^-_{T_2}(N^+_{T_2}(x_1))=3$, then either $t$ is an origin of a path $R\equiv a_{n-1}a_{n-2}...a_{\delta+1}$ and thus$P_2vR^{-1}x_2\equiv C$ or $t$ is an origin of $R^{-1}$ and thus $vP^{-1}_2R^{-1}v\equiv a_1a_2...a_1=C$.\\
Assume now that $T_1=3A$ and $d^+(x_i)=1$ for all $x_i\in T_1$. If $|T_2|=3$ then if $T_2$ is transitive, $x_3$ has an outneighbor origin of an antidirected outpath $S$ in $T_2$ and thus $vx_1x_2x_3Sv\equiv a_1a_2...a_1=C$, and if $T_2$ is a circuit then $x_3x_2x_1vO_{T_2}x_3\equiv C$. So $|T_2|\geq 4$. If the path $R=a_2a_3...a_{n-2}$ has at least three blocks, one can find $R$ in $T_2\cup \{v\}$ with origin $x\in N^-(x_1)$ and end $y\in N^-(x_3)$, thus $x_2x_1Rx_3x_2\equiv C$. Otherwise, by taking $R'=a_2a_1...a_6$, we can suppose without loss of generality that $C=(3,3,1,1)$. If $N^+_{T_2}(x_i)=\{u\}$ for all $i$, then $T_2=4B$, let $V(T_2)=\{u_1,u_2,u_3\}$, we have $x_1u_1vuu_2x_3u_3x_2x_1\equiv C$. Otherwise, $T_2=\overline{4A}$ and $4\rightarrow T_1$. $x_1$ has an outneighbor in $T_2-4$ origin of an antidirected inpath $P$, thus $x_1x_34x_2vP^{-1}x_1\equiv C$.\\
\textbf{Case $2$}: If $a_1\rightarrow a_2$. For $x_i\in T_1$ set $P_i\equiv a_3a_4...a_{\delta+2}$ be a path in $T_1$ of origin $x_i$ and end $y_i$.  There exists a $P_i$ such that for a $u_i\in N^+_{T_2}(x_i)$ we have $N^-_{T_2-u_i}(y_i)\geq 2$. If there exists $Q\equiv a_{\delta+3}a_{\delta+4}...a_{n-1}a_0$ or $R\equiv a_{n-2}a_{n-3}...a_{\delta+1}$ in $T_2-u_i$ of origin and inneighbor of $y_i$, then $vuP_iQv\equiv a_1a_2...a_1=C$ or $vR^{-1}P^{-1}_iv\equiv a_{\delta}a_{\delta+1}...a_{\delta}=C$. Otherwise $(T_2-u_i;Q)$ is an exception, so by the building lemmas, it is one of the following biexceptions with respect to $y_i$: (in fact, let $z_i$ be the end of $P_i^*$, if $d^+_{T_2-u_i}(z_i)=0$ by taking the path $P_j$ of origin $y_i$, we have $y_j=z_i$ and since $d^-_{T_2-u_j}(z_i)\geq |T_2|-2$, $T$ contains $C$)\\
Dual Exc $(1,0)$: $(T_2,\overline{a_2a_1a_0a_{n-1}a_{n-2}})$ is not a biexception with respect to $y_i$ so $z_i$ is an origin of a path $P\equiv \overline{a_3a_2a_1a_0a_{n-1}a_{n-2}}$ in $T_2\cup \{y_i\}$ and thus $z_ix_ivP^{-1}z_i\equiv \overline{a_4a_5...a_4}=\overline{C}$. Therefor $T$ contains $C$.\\
Exc $(3,1)$: $(T_2;R)$ is the biexception Dual Exc $(1,0)$ with respect to $y_i$ so $T$ contains $C$.\\
Exc $(4,1),(8,1),(9,1),(33,1)$: $(T_2;\overline{a_2a_1a_0...a_{\delta+5}})$ is not a biexception with respect to $y_i$ so $z_i$ is an origin of a path $P\equiv \overline{a_4a_3a_2...a_{\delta+4}}$ in $T_2\cup \{y_i\}$, thus $\overline{^*(P_i^{-1})}u_ivP^{-1}z_i\equiv \overline{a_4a_5...a_4}=\overline{C}$, therefore $T$ contains $C$.\\
Exc $(E_8(5),1)$: $(T_2-1;Q)$ is not the biexception Exc $(E_8(5),1)$ with respect to $z_i$ so $T$ contains $C$.\\
Dual Exc $(E'_8(5),0),(E'_9(7),0)$ and $(E'_{10}(9),0)$: $(T_2-1;Q)$ is not one of the biexceptions Dual Exc $(E'_8(5),0)),(E'_9(7),0)$ and $(E'_{10}(9),0)$ with respect to $z_i$ so $T$ contains $C$.\\
Assume now that $(T_1;a_1...a_{\delta})$ is not a Gr\"{u}nbaum's exception, let $x_1...x_{\delta}\equiv a_1...a_{\delta}$ in $T_1$ such that $d^+_{T_2}(x_1)=max\,\{d^+_{T_2}(x);x\in Or(a_1...a_{\delta};T_1)\}$. If there is an outneighbor $u$ of $x_1$ origin of a path $P\equiv a_0a_{n-1}...a_{\delta+2}$ or an outneighbor $w$ of $x_{\delta}$ origin of a path $Q\equiv a_{\delta+1}a_{\delta+2}...a_{n-1}$ in $T_2$, then $ux_1...x_{\delta}vP^{-1}\equiv C$ or $vx_1...x_{\delta}Qv\equiv C$. Otherwise, suppose that $a_2\rightarrow a_1$. If $d^+_{T_2}(x_1)\geq 2$ and $d^+_{T_2}(x_2)=0$, then, by lemma \ref{2.8} $N^-_{T_2}(x_1)$ contains an origin $w$ of a path $P\equiv a_0a_{n-1}...a_{\delta+2}$ in $\overline{T_2}$ and so $wx_1vQx_2P^{-1}\equiv C$ where $Q\equiv a_3...a_{\delta}$ in $T_1\setminus \{x_1,x_2\}$ unless $(T_2;a_0a_{n-1}...a_{\delta+2})$ is one of the exceptions Exc $0,1,4,18,19,22,33,E_1(4),E_1(6)$, so if $T_1-x_2$ contains a path $y_2...y_{\delta}\equiv a_2...a_{\delta}$ then every vertex of $T_2$ is an origin of a path $R\equiv  a_{\delta+1}a_{\delta+2}...a_0$ in $T_2\cup \{x_2\}$ with end in $T_2$, thus $vy_2...y_{\delta}Rv\equiv a_1a_2...a_1=C$. Otherwise, $(T_1-x_2;a_2...a_{\delta})$ is a G\"{u}nbaum's exception, more precisely $T_1=4B$, then $a_6\rightarrow a_7$ and so the path $a_0a_{n-1}...a_{\delta+2}$ has an odd number of blocks, so $(T_2;a_0a_{n-1}...a_{\delta+2})$ is either one of the exceptions $4,33,E_1(4),E_1(6)$, let $y_2...y_5\equiv a_2...a_5$ be a path in $T_1$, $y_2$ has an outneighbor origin of a path $R\equiv a_1a_0...a_7$ in $T_2$ and thus $y_2...y_5vR^{-1}y_2\equiv a_2a_3...a_2=C$; or the exception Exc $19$, let $y_3...y_6\equiv a_3...a_6$ be a path in $T_1$, $y_3$ has an inneighbor origin of a path $R\equiv a_2a_1...a_8$ in $T_2$, thus $y_3...y_6vR^{-1}y_3\equiv a_3a_4...a_3=C$. Otherwise, if $d^+_{T_2}(x_1)\geq 2$ then $(T_2;a_0a_{n-1}...a_{\delta +2})$ is an exception, so by the building lemmas we may assume that it is one of the following biexceptions with respect to $x_1$:\\
Exc $(0,0)(2)$: $\overline{T}(N^-_{\overline{T_2}}(x_1))\neq 3A$ so $\overline{T}$ contains $C$ and so $T$ contains $C$.\\
Exc $(1,0)$: We have $\delta\geq 3$ so $x_{\delta}\rightarrow 4$ which is an origin of a path $Q$. A contradiction.\\
Exc $(2,0)$: $x_{\delta}\rightarrow 1$, a contradiction.\\
Exc $(4,0)$: $\delta=2$ so $21x_1v453x_22\equiv C$.\\
Exc $(6,0)$: We have $x_{\delta-1}\rightarrow x_{\delta}$ so $d^+_{T_2}(x_{\delta})\geq 1$ and $4\rightarrow (T_1-x_1)$ so $N^+_{T_2}(x_{\delta})\cap \{1,2,3,5\}\neq \emptyset$, a contradiction.\\
Exc $(9,0),(15,0),(16,0),(26,0)(27,0),(31,0),(33,0)$ and $(48,0)$: If $d^+_{T_2}(x_{\delta})\geq 1$, set $s=min\,\{i\geq \delta +2; d^+_{C}(a_s)=0\}$ and let $Q_1\equiv a_{\delta +1}a_{\delta +2}...a_{s-1}$ of origin $x$ and $Q_2\equiv a_0a_{n-1}...a_{s+1}$ of origin $a$ be two paths such that $T_2=V(Q_1)\cup V(Q_2)$, $x\in N^+(x_{\delta})$ and $a\in N^+(x_1)$, we have $ax_1...x_{\delta}Q_1vQ_2^{-1}\equiv C$. Otherwise, $\delta=2$, there exists $R\equiv a_1a_0a_{n-1}...a_{\delta+2}$ in $T_2\cup \{x_2\}$ of origin an outneighbor of $x_1$ and with end distinct from $x_2$, we have $x_1vP^{-1}x_1\equiv a_2a_3...a_2=C$.\\
Exc $(E_5(|T_2|),0)$: We have $T_2\rightarrow x_{\delta}$, so $d^+(x_2)< \delta$ if and only if $\delta=2$, then $\delta=2$ and thus the first block of the path $a_{\delta+2}a_{\delta+3}...a_0$ is at most one which is a contradiction.\\
Otherwise, by lemma \ref{2.4}, $(T_1;a_1...a_{\delta})=(\overline{4B};-(1,1,1)),(\overline{4B};-(2,1))$ or $(3A;$\\$-(2))$ or $|T_1|=2$: We have $d^+_{T_2}(x)=1$ for all $x\in Or(a_1...a_{\delta};T_1)$. So there exist $y_1...y_{\delta}\equiv \overline{a_1...a_{\delta}}$ in $T_1$ and $u\in N^-_{T_2}(y_1)$ such that $(T_2-u;\overline{a_{\delta+1}a_{\delta+2}...a_{n-2}})\neq (\overline{5C};-(2,1,1))$. Since $d^-_{T_2-u}(x_{\delta})\geq |T_2-u|-1$, then there is an inneighbor $w$ of $x_{\delta}$ origin of a path $P\equiv \overline{a_{\delta+1}a_{\delta+2}...a_{n-2}}$ in $T_2-u$; we have $ux_1...x_{\delta}Pvu\equiv \overline{C}$. Thus $T$ contains $C$.\\
Assume now that $a_1\rightarrow a_2$. Either $d^+_{T_2}(x_1)\geq 2$ and $(T_2;a_0a_{n-1}...a_{\delta+2})$ is an exception or $d^+_{T_2}(x_1)\leq 1$. So by the building lemmas, it is one of the following biexceptions with respect to $x_1$:\\
Exc $(1,0)$: $d^-(4)< \delta$ if and only if $\delta=2$. So $\delta=2$ and $T(N^+(4))\neq \overline{4A}$, then $\overline{T}$ contains $C$ and thus $T$ contains $C$.\\
Exc $(6,0),(27,0),(33,0)$: $(\overline{T}(N^+(x_2));a_0a_{n-1}...a_{\delta+2})$ is not a biexception so $\overline{T}$ contains $C$ and thus $T$ contains $C$.\\
Exc $(26,0)$: There is a block of length $3$ so $a_{\delta-1}\rightarrow a_{\delta}$. If $x_1\rightarrow 4$ then $\overline{T}(N^-(4))\neq 6H$ so $T$ contains $C$. Then we may assume that $N^+_{T_2}(x_1)=\{5,6\}$. Also if $x_2\rightarrow 5$ then $\overline{T}(N^-(5))\neq 6H$ and so $T$ contains $C$. We can suppose that $5\rightarrow x_2$. We have $\{1,2,3,5\}\rightarrow x_{\delta}$ since otherwise $x_{\delta}$ has an outneighbor that is an origin of a path $Q\equiv a_{\delta+1}a_{\delta+2}...a_{n-1}$ which is a contradiction. Now if $x_{\delta}\rightarrow 6$ we have $45x_2...x_{\delta-1}v126x_{\delta}3x_14\equiv C$ and if $6\rightarrow x_{\delta}$ we have $45x_2...x_{\delta-1}v6x_{\delta}321x_14\equiv C$.\\
Exc $(E_{13}(|T_2|),0)$: We have $\overline{T}(N^-(1))\neq F_{13}(|T_2|)$ so $T$ contains $C$.\\
Now if $d^+_{T_2}(x_1)\leq 1$. Set $\{u\}=N^+_{T_2}(x_1)$, if it exists, if there exist $w,t,y\in T_2$ such that $w\rightarrow \{x_2,t\}$, $\{t,y\}\rightarrow x_1$ and $y$ is an origin of a path $P\equiv a_{n-3}a_{n-4}...a_{\delta+2}$ then $twx_2...x_{\delta}vP^{-1}x_1t\equiv C$. Otherwise, suppose first that $N^+_{T_2}(x_1)=\{u\}$, if there exists $w\in N^-_{T_2}(x_2)$ such that $d^+_{T_2-u}(w)\geq 1$, let $w'\in N^+_{T_2-u}(w)$, then $(T_2\setminus \{w,w'\};a_{n-3}a_{n-4}...a_{\delta+2})$ is one of the following exceptions:\\
Dual Exc $0$: Set $u\equiv 1$, if $3\rightarrow w'$ then $x_1wx_2...x_{\delta}vw'321x_1\equiv C$, so $w'\rightarrow 3$. If $2\rightarrow w'$ then $u$ is an origin of a path $P\equiv a_0a_{n-1}...a_{\delta+2}$ in $T_2$, a contradiction. So $w'\rightarrow 2$  and thus $x_1wx_2...x_{\delta}3w'21x_1\equiv C$.\\
Exc $1$: $u\equiv 4$, if $w\rightarrow u$, let $\{y_2\rightarrow y_1\}=T_2\setminus \{u,w,1,2\}$ then $uwx_2...x_{\delta}vy_1y_2x_112u$\\$\equiv C$. Otherwise, $d^+_{T_2}(u)\geq 3$ then $T_2-u=5E$, we have $twx_2...x_{\delta}vu3x_112t\equiv C$.\\
Dual Exc $4$: Set $u\equiv 1$. If $w\rightarrow u$, let $\{y_2\rightarrow y_1\}=T_2\setminus \{w,u,2,4,5\}$, we have $uwx_2...x_{\delta}vy_1y_2x_1452u\equiv C$. Otherwise, since $T_2-u\neq 6D$ then $u$ is an origin of $P\equiv a_0a_{n-1}...a_{\delta+2}$ in $T_2$, a contradiction.\\
Exc $7$: $u\equiv 5$, so there exists an outneighbor of $w$ satisfying the conditions implying the result.\\
Dual Exc $33$: We have $d^+_{T_2}(u)\geq 3$, a contradiction.\\
Or $\{w\}\subseteq N^-_{T_2}(x_2)=\{u,w\}$ and $d^+_{T_2-u}(y)=0$ for all $y\in N^-_{T_2}(x_2)$. If $u\in N^-_{T_2}(x_2)$ then $d^+(u)<\delta^+(T)$ a contradiction, so $N^-_{T_2}(x_2)=\{w\}$. Take $x_2$ in the place of $v$ and set $T'_2=T(N^-(x_2))$ and $T'_2=T(N^+(x_2))$. Without loss of generality we may assume that $(T'_1;a_1...a_{\delta})$ is not a Gr\"{u}nbaum's exception, so let $y_1...y_{\delta}\equiv a_1...a_{\delta}$ be a path in $T'_1$. Obviously, $y_1$ or $y_{\delta}$ is an inneighbor of $v$, so if $(T'_2-v;a_{n-1}a_{n-2}...a_{\delta+2})$ is not a GR\"{u}nbaum's exception, let $Q\equiv a_{n-1}...a_{\delta+2}$, if $y_1\in N^-(v)$ then $vy_1...y_{\delta}x_2Q^{-1}v\equiv C$ and if $y_{\delta}\in N^-(v)$ then $x_2y_1...y_{\delta}vQ^{-1}x_2\equiv C$. So we may assume that $(T'_2-v;a_{n-1}a_{n-2}...a_{\delta+2})$ is a Gr\"{u}nbaum's exception, which implies $u$ is an origin of a path $P\equiv a_0a_{n-1}...a_{\delta+2}$ in $T_2$, a contradiction. Assume now that $d^+_{T_2}(x_1)=0$, then by lemma \ref{2.3} $(T_1;a_1...a_{\delta})=(4B;(1,2))$ or $a_1...a_{\delta}$ is directed. Suppose first that there exists $w\in N^-_{T_2}(x_2)$ such that $d^+_{T_2}(w)\geq 1$ and let $w'\in N^+_{T_2}(w)$, then $(T_2\setminus \{w,w'\};a_{n-2}a_{n-3}...a_{\delta+2})$ is a Gr\"{u}nbaum's exception. If $(T_1;a_1...a_{\delta})=(4B;(1,2))$, take $C=u_0u_1...u_{n-1}u_0=a_{\delta+1}a_{\delta}...a_{\delta+2}a_{\delta+1}$, $x_{\delta}...x_1\equiv u_1...u_{\delta}$ with $d^+_{T_2}(x_{\delta})\geq 1$, the result follows. Else, $a_1...a_{\delta}$ is directed, let $P\equiv \overline{a_{\delta+4}a_{\delta+5}...a_0a_1}$ be a path in $T_2-w$, we have $vx_{\delta}...x_2wx_1Pv\equiv \overline{a_2a_3...a_2}=\overline{C}$. So we may assume that $N^-_{T_2}(x_2)=\{w\}$ with $d^+_{T_2}(w)=0$, clearly $d^+_{T_2}(x_{\delta})\geq 2$, then by theorem \ref{2.1} $x_{\delta}$ has an outneighbor origin of a path $P\equiv a_{\delta+1}a_{\delta+2}...a_{n-1}$ in $T_2$ and so $vx_1...x_{\delta}Pv\equiv C$.    $\square$
\end{proof}

\begin{lemma}\label{3.5}
Let $T$ be a tournament of order $n$ with $\delta^-(T)\geq 2$ and $C$ a non-directed and non-antidirected cycle of order $n$. Let $B$ be a block of $C$ of maximal length and set $C=a_0a_1...a_{n-1}a_0$ such that $B$ is forward and it ends at $a_{\delta^-(T)+2}$. If $a_{n-1}\rightarrow a_0\leftarrow a_1$, then $T$ contains $C$ if and only if $(T;C)$ is not one of the exceptions $A_3$ and $A_{14}$.
\end{lemma}
\begin{proof}
Denote by $B'$ the block of $C$ containing the arc $a_{n-1}\rightarrow a_0$. We distinguish six cases:\\
\textbf{Case $1$}: If $|B|=|B'|\geq 4$. Rotate $C$ until $a_{\delta}\rightarrow a_{\delta-1}$. If $(T_1;a_1...a_{\delta})$ is not a Gr\"{u}nbaum's exception, there exists $x_1...x_{\delta}\equiv a_1...a_{\delta}$ such that $d^-_{T_2}(x_1)\geq 2$ and $s^-_{T_2}(N^-_{T_2}(x_1))\geq 3$, then, since $(T_2;a_0a_{n-1}...a_{\delta+2})$ is not an exception, by theorem \ref{2.1} there is an inneighbor $u$ of $x_1$ origin of a path $P\equiv a_0a_{n-1}...a_{\delta+2}$ in $T_2$ and thus $ux_1...x_{\delta}vP^{-1}\equiv C$. Assume now that $(T_1;a_1...a_{\delta})$ is a Gr\"{u}nbaum's exception, let $a\in T_2$ such that $d^-_{T_2}(a)=\Delta^-(T_2)$ and set $T'_1=T_1\cup \{a\}$ and $T'_2=T_2-a$. If there exists $x_0...x_{\delta}\equiv a_0...a_{\delta}$ in $T'_1$ such that $x_{\delta}\neq a$ then $d^-_{T'_2}(x_0)\geq 2$, so there is an inneighbor of $x_0$ origin of a path $P\equiv a_{n-1}a_{n-2}...a_{\delta+2}$ in $T'_2$ unless $(T'_2;a_{n-1}...a_{\delta+2})$ is one of the following biexceptions with respect to $x_0$:\\
Dual Exc $(6,0)$: $x_1\neq a$ so $\{x_1,x_{\delta}\}\rightarrow 5$, then by taking $5$ in the place of $v$ the result follows.\\
Dual Exc $(31,1)$: In a similar way take $6$ in the place of $v$.\\
Dual Exc $(E_{11}(5),0)$: In a similar way take $3$ in the place of $v$.\\
Otherwise, $T'_1=4B$, if there is an inneighbor $u$ of a vertex $y_i\in T_1$ origin of a path $P\equiv a_0a_{n-1}...a_{\delta+3}$ in $T'_2$ then $uy_iay_{i+1}y_{i+2}vP^{-1}\equiv C$. Else, if there is a vertex $y_i\in T_1$ such that $d^-_{T'_2}(y_i)=1$, we have $s^+_{T_2}(N^+_{T_2}(y_i))\geq |T_2|-1$ so there exists $u\in N^+_{T_2}(y_i)$ origin of a path $P\equiv a_2a_1...a_{\delta+4}$ in $T_2$ and thus $y_iy_{i+1}y_{i+2}vP^{-1}y_i\equiv a_3a_4...a_3=C$. So we can assume that $d^-_{T'_2}(y_i)\geq 2$ for all $y_i\in T_1$ and $s^-_{T'_2}(N^-_{T'_2}(T_1))\geq 3$ and $(T'_2;a_0a_{n-1}...a_{\delta+3})$ is one of the following exceptions:\\
Dual Exc $5$, $7$, $51$, $E'_{8}(|T'_2|),E'_9(|T'_2|)$ and $E'_{10}(|T'_2|)$: There exists $x\in N^+_{T_2}(y_i)$ origin of a path $P\equiv a_2a_1...a_{\delta+4}$ in $T_2$ so $y_iy_{i+1}y_{i+2}vP^{-1}y_i\equiv C$.\\
\textbf{Case $2$}: If $|B'|=2$. We distinguish two cases:\begin{itemize}
\item[(a)] If $a_2\rightarrow a_1$, suppose first that $(T_1;a_1...a_{\delta})$ is a Gr\"{u}nbaum's exception, every vertex $x_i\in T_1$ is an origin of a path $P_i\equiv a_0a_1...a_{\delta-1}$ and $d^-_{T_2}(x_i)\geq 2$, denote by $y_i$ the origin of $^*P_i$. If $x_i$ has an inneighbor origin of a path $P\equiv a_{n-1}a_{n-2}...a_{\delta+1}$ in $T_2$ then $P_ivP^{-1}x_i\equiv C$. Otherwise, $(T_2;a_{n-1}a_{n-2}...a_{\delta+1})$ is an exception, so by the building lemmas it is an biexception with respect to $x_i$ (Considering that $a_{n-1}a_{n-2}...a_{\delta+1}=-(1,Q,1)$ with odd number of blocks and that $d^-_{T_2}(y_i)\geq 2$), in any possible biexception except $Exc (0,1)$, there is an origin $w$ of a path $P\equiv a_{n-1}a_{n-2}...a_{\delta+1}$ in $T_2-S$ which is an outneighbor of $y_i$, let $x_j=y_i$, we have $P_jvP^{-1}y_i\equiv C$. And if it is the biexception $Exc (0,1)$ with respect to $x_i$, in this case, by induction on $i$, we may assume that $1\rightarrow T_1$, so $d^+(1)=4$, a contradiction. So we may assume that $(T_1;a_1...a_{\delta})$ is not a Gr\"{u}nbaum's exception, let $x_1...x_{\delta}\equiv a_1...a_{\delta}$ in $T_1$ such that $d^+_{T_2}(x_1)=max \,\{d^+_{T_2}(x);x\in Or(a_1...a_{\delta},T_1)\}$. If $x_1$ has an outneighbor $u$ origin of a path $P\equiv a_0a_{n-1}...a_{\delta+2}$ in $T_2$ then $ux_1...x_{\delta}vP^{-1}\equiv C$. Otherwise, if $d^+_{T_2}(x_1)\geq 2$ and $s^-_{T_2}(N^+_{T_2}(x_1))\geq 3$ then $(T_2;a_0a_{n-1}...a_{\delta+2})$ is an exception. If $d^+_{T_2}(x_2)=0$, then, by lemma \ref{2.8}, $\overline{T}$ contains $C$. So we may assume that $d^+_{T_2}(x_2)\geq 1$ and thus, by the building lemmas, we can suppose that $(T_2;a_0a_{n-1}...a_{\delta+2})$ is one of the following biexceptions with respect to $x_1$:\\
Dual Exc $(3,1)$: $T(N^-(1))\neq 4A$, then $T$ contains $C$ or it will be solved.\\
Dual Exc $(11,1)$: If $N^+(x_{\delta})\cap \{1,2,4,5\}$, one can find $C$, so we may assume that $\{1,2,4,5\}\rightarrow x_{\delta}$ and thus if $\delta >2$ then $T(N^+(2))\neq \overline{5E}$, so we may assume that $\delta =2$, therefore $13v2x_245x_1\equiv C$.\\
Dual Exc $(E_2(5),1)$: $T(N^-(x_2))\neq F_2(5)$, then $T$ contains $C$ or it will be solved.\\
Dual Exc $(E_8(6),1)$: $T(N^-(1))\neq F_8(6)$, then $T$ contains $C$ or it will be solved.\\
If $d^+_{T_2}(x_1)\geq 2$ and $S^-_{T_2}(N^+_{T_2}(x_1))=N^+_{T_2}(x_1)=\{u\rightarrow w\}$, if there exists $t\in N^+_{T_2\setminus \{u,w\}}(x_2)$, then $|T_2|\geq 4$, so there exists $Q\equiv a_{n-1}a_{n-2}...a_{\delta+2}$ in $T_2-t$ of origin distinct from $u,w$, we have $x_1tx_2...x_{\delta}vQ^{-1}x_1\equiv C$ and if $x_2\rightarrow u$, since $(T_2;a_{n-2}a_{n-3}...a_{\delta+2})$ is not a Gr\"{u}nbaum's exception, there exists $Q\equiv a_{n-3}...a_{\delta+2}$ in $T_2\setminus \{u,w\}$, we have $wux_2...x_{\delta}vQ^{-1}x_1w\equiv C$. So we may assume that $N^+_{T_2}(x_2)\subseteq \{w\}$, if $\delta\geq 3$ then $x_{\delta}\rightarrow u$, let $P\equiv a_{\delta+2}a_{\delta+3}...a_{n-2}$ be a path in $T_2\setminus \{u,w\}$, we have $wvx_2...x_{\delta}uPx_1w\equiv C$, and if $\delta=2$, set $s=min\,\{i\geq 3;d^+_C(a_i)=2\}$ and let $P\equiv a_2a_3...a_{s-1}$ be a path in $T_2$ such that $u,w\in V(P)$ and $T_2-V(P)$ contains a path $Q\equiv a_{s+1}a_{s+2}...a_{n-1}$, we have $x_1x_2PvQx_1\equiv C$.\\
Assume now that $d^+_{T_2}(x_1)=1$, set $N^+_{T_2}(x_1)=\{u\}$, by lemma \ref{2.4}, either $T_1=4B$ or $\delta=2$. In the first case, set $T_1=\{u_1,u_2,u_3\}\rightarrow u_4$, we can suppose that $\{u_1,u_2,u_3\}\rightarrow \{u\}$ and $d^-_{T_2}(u)=0$ since otherwise, by taking $u_1$ in the place of $v$ the result follows. If $a_1...a_4=-(2,1)$, there exist $w,t\in T_2-u$ such that $t$ is an origin of a path $Q\equiv a_{n-2}a_{n-3}...a_{\delta+2}$ in $T_2\setminus \{u,w\}$, we have $u_4uu_2wu_1vQ^{-1}u_3u_4\equiv C$, else, $a_1...a_4=-(1,1,1)$, if there exists $w\in N^-_{T_2}(u_4)$, then there exist a path $Q\equiv a_{n-1}a_{n-2}...a_{\delta+2}$ in $T_2-w$ of origin distinct from $u$, thus $u_3u_2wu_4u_1vQ^{-1}u_3\equiv C$, otherwise, let $w,t\in T_2-u$ such that $t$ is an origin of a path $Q\equiv \overline{a_{\delta+2}a_{\delta+3}...a_{n-2}}$ in $T_2\setminus \{u,w\}$, we have $u_4wu_2u_1uu_3Qvu_4\equiv \overline{C}$, then $T$ contains $C$. In the second case, suppose that there exists $w\in N^+_{T_2}(x_2)\cap N^-_{T_2}(x_1)$, if $T_2-w$ contains a path $P\equiv a_{n-1}a_{n-2}...a_{\delta+2}$ of origin an inneighbor of $x_1$ then $x_1wx_2vP^{-1}x_1\equiv C$, else, without loss of generality, we can suppose that $(T_2-w;a_{n-1}...a_{\delta+2})$ is a Gr\"{u}naum's exception, in this case, one can find a path $P\equiv a_{\delta+1}a_{\delta+2}...a_0a_1$ in $T_2\cup \{v\}$ with origin $w$ and end $u$, so $x_1x_2Px_2\equiv a_1a_2...a_1=C$. So we can assume that $(T_2-u)\rightarrow x_2$, if $u\rightarrow x_2$ then $T_2-u$ contains a path $P\equiv a_5a_6...a_0$ we have $vx_1ux_2Pv\equiv a_1a_2...a_1=C$. Now if $x_2\rightarrow u$, if there exists $w\in N^-_{T_2}(x_1)$ such that $T_2-w$ contains a path $P\equiv a_2a_3...a_{n-3}$ of origin distinct from $u$ so $x_1x_2Pvwx_1\equiv C$, otherwise, $|T_2|=4$ and $d^-_{T_2}(u)=0$, set $T_2=\{u,w_1,w_2,w_3\}$ we have $w_1w_2vw_3x_1x_2uw_1\equiv C$.
\item[(b)] If $a_1\rightarrow a_2$, suppose first that $(T_1;a_1...a_{\delta})$ is a Gr\"{u}nbaum's exception, let $a\in T_2$ such that $d^+_{T_1}(a)=min\,\{d^+_{T_1}(x);x\in T_2\}$ and set $T'_1=T_1\cup \{a\}$ and $T'_2=T_2-a$. Let $x_0...x_{\delta}\equiv a_0...a_{\delta}$ be a path in $T'_1$ such that $d^-_{T'_2}(x_0)=max\,\{d^-_{T'_2}(x);x\in Or(a_0...a_{\delta};T'_1)\}$, if $x_0$ has an inneighbor origin of a path $P\equiv a_{n-1}a_{n-2}...a_{\delta+2}$ in $T'_2$ then $x_0...x_{\delta}vP^{-1}x_0\equiv C$ or an inneighbor origin of a path $Q\equiv \overline{a_{\delta+1}a_{\delta+2}...a_{n-2}}$ then $x_{\delta}...x_0Qvx_{\delta}\equiv \overline{C}$. Otherwise, if $T'_2=3A$, there exists $x\in T_1$ such that $x\rightarrow a$, thus $vaO_{T_1}O_{T'_2}v\equiv a_1a_2...a_1=C$. Assume now that $T'_2\neq 3A$, if $d^-_{T'_2}(x_0)\geq 2$ then $(T'_2;a_{n-1}a_{n-2}...a_{\delta+2})$ is one of the following exceptions:\\
Dual Exc $2$: If $\delta=3$ by taking $C=a'_0a'_1...a'_0=a_4a_5...a_4$ it is solved in lemma \ref{3.3}. If $\delta=5$, since $a_0...a_{\delta}$ has at least $3$ origins in $T'_1$ then $d^-(1)\geq 6$, a contradiction.\\
Dual Exc $4$: If $\delta=3$, there exists $x\in T_1$ such that $d^-_{T_2}(x)\geq 3$, since otherwise, by the building lemmas, we may assume that $\{1,2,4,5\}\rightarrow x_1$, if $x_1=a$, by definition of $a$ we deduce that $\{3,4,5\}\rightarrow \{x_2,x_3\}$. And if $\delta=5$, obviously $d^-_{T_2}(x)\geq 3$ for all $x\in T_1$. So let $y_3...y_{\delta+2}\equiv a_3...a_{\delta+2}$ be a path in $T_1$ such that $d^-_{T_2}(y_{\delta+2})\geq 3$ and let $u\in N^+_{T_2}(y_3)$, $d^-_{T_2-u}(y_{\delta+2})\geq 2$ so $y_{\delta+2}$ has an inneighbor origin of a path $P=(1,1,2)$ in $T_2$, thus $vuy_3...y_{\delta+2}Pv\equiv a_1a_2...a_1=C$.\\
Dual Exc $8$: Since $d^-_{T'_2}(x_0)=2$ then $T_1=3A$, if $x_3\rightarrow 2$ then $x_0...x_321v345x_0\equiv C$, so $2\rightarrow x_3$. If $x_3\rightarrow 3$ then $x_0...x_33124v5x_0\equiv C$, so $3\rightarrow x_3$. If $x_3\rightarrow 4$ then $x_0...x_34351v2x_0\equiv C$, so $4\rightarrow x_3$. Thus $T'_1=4A$, if $x_2=a$ then $x_0x_1x_2v4x_33125x_0\equiv C$, so $x_1=a$ and thus $d^+_{T_1}(x)\geq 2$ for all $x\in T'_2$. So $4\rightarrow x_2$ and thus $x_0x_1x_243x_321v5x_0\equiv C$.\\
Dual Exc $9$: If $N^+(x_0)\cap \{4,5\}\neq \emptyset$ then by the building lemmas we may assume that $T'_2\rightarrow x_1$ and thus $x_1\neq a$. $x_1$ has an inneighbor origin of a path $P\equiv a_6a_7...a_0$ in $T_2-3$, thus $v3x_0x_3x_1Pv\equiv a_1a_2...a_1=C$. So we may suppose that $\{4,5\}\subseteq N^-(x_0)$, if $x_{\delta}\rightarrow 2$ then $x_0...x_{\delta}21v354x_0\equiv C$, so $2\rightarrow x_{\delta}$. If $x_{\delta}\rightarrow 3$ then $x_0...x_{\delta}32v514x_0\equiv C$, so $3\rightarrow x_{\delta}$. If $x_{\delta}\rightarrow 4$ then $x_0...x_{\delta}43v125x_0\equiv C$, so $\{4,5\}\rightarrow x_{\delta}$. If $\delta=5$ then $2\rightarrow x_0$ and $x_{\delta}\rightarrow 1$, thus $x_0...x_{\delta}1453v2x_0\equiv C$, so $\delta=3$ and $T'_1=4A$. If $x_2=a$ then $x_0x_1x_2v3x_32415x_0\equiv C$, so $x_1=a$. Thus $1\rightarrow \{x_2,x_3\}$, we have $x_0x_1x_214x_332v5x_0\equiv C$.\\
Dual Exc $15$: If $\delta=3$ then by taking $C=a_7a_8...a_7$ it is solved in lemma \ref{3.3}, so $\delta\geq 5$. So $T'_1$ contains at least $3$ origins of an antidirected inpath, thus $d^-(4)\geq 8$, a contradiction.\\
Dual Exc $33$: If $\delta=3$, let $P=-(1,1,1,2,1)$ be a path in $T'_2$ of end an inneighbor of $x_0$, we have $x_0...x_3vPx_0\equiv a_7a_8...a_7=C$, so $\delta\geq 5$. Let $y_{\delta-2}...y_{2\delta-2}\equiv a_{\delta-2}...a_{2\delta-2}$ be a path in $T_1$ and $u\in N^+_{T_2}(y_{\delta-2})$, $d^-_{T_2-u}(y_{2\delta-2})\geq 2$ so $T_2-u$ contains a path $P\equiv a_{2\delta-1}a_{2\delta}...a_{\delta-5}$ of origin an inneighbor of $y_{2\delta-2}$, thus $y_{\delta-2}...y_{2\delta-2}Pvuy_{\delta-2}\equiv a_{\delta-2}a_{\delta-1}...a_{\delta-2}=C$.\\
Dual Exc $37$: If $\delta=3$ then $5264731x_0...x_3v5\equiv C$, so $\delta\geq 5$.  Let $y_3...y_{\delta+2}\equiv a_3...a_{\delta+2}$ be a path in $T_1$ and $u\in N^+_{T_2}(y_3)$, $d^-_{T_2-u}(y_{\delta+2})\geq 2$ so $y_{\delta+2}$ has an outneighbor origin of a path $P\equiv a_{\delta+3}a_{\delta+4}...a_0$ in $T_2-u$, thus $vuy_3...y_{\delta+2}Pv\equiv a_1a_2...a_2=C$.\\
Dual Exc $E_1(4)$: If $\delta=3$ then by taking $C=a_4a_5...a_5$ it is solved in lemma \ref{3.3}, so $\delta=5$. Let $y_5...y_9\equiv a_5...a_9$ be a path in $T_1$, by theorem \ref{2.1} $y_5$ has an outneighbor origin of a path $P\equiv a_4a_3...a_0$ in $T_2$, thus $y_5...y_9vP^{-1}y_5\equiv a_5a_6...a_5=C$.\\
Dual Exc $E_8,E_9,E_{10}$: $d^-_{T'_2}(x_0)=2$ so $\delta=3$. By the building lemmas we may assume that $(T'_2-X)\rightarrow x_1$ so $T'_2=4A$. $d^-_{T_1}(w)=2$ so $x_1\neq a$. Let $x\rightarrow y$ be an arc in $T'_2-X$ and $P\equiv a_6a_7...a_{n-3}$ be a path in $T'_2-(X\cup \{x,y\})$, we have $yxx_1x_3x_2x_0uPvwy\equiv C$.\\
Suppose now that $d^-_{T'_2}(x_0)=1$, if $T'_2$ is a transitive triangle, we may assume that $T'_1=4A$, and thus $d^+_{T_1}(a)\geq 1$; set $V(T'_2)=\{u_1,u_2,u_3\}$ such that $u_1u_2u_3$ is a directed outpath and $V(T_1)=\{y_1,y_2,y_3\}$. Without loss of generality suppose that $y_1\rightarrow a$, if $a\rightarrow u_3$, $T_1\cup \{u_1\}$ contains a path $z_0...z_3\equiv a_0...a_3$ and thus $z_0...z_3vau_2u_3z_0\equiv C$ or $z_0...z_3vu_2au_3z_0\equiv C$, so $u_3\rightarrow a$. If $u_1\rightarrow y_3$ then $u_3vay_1y_2y_3u_1u_2u_3\equiv C$, so $y_3\rightarrow y_3$. There exist $i\in \{1,2,3\}$ and $\{j,k\}= \{1,2\}$ such that $y_i\rightarrow u_j$ and $u_k\rightarrow y_{i-1}$ then $avu_jy_iy_{i+1}y_{i-1}u_ku_3a\equiv C$. If $|T'_2|=4$ then by taking $C=a_4a_5...a_4$ it is solved in lemma \ref{3.3}, so we may assume that $|T'_2|\geq 5$. Set $\{u\}=N^-_{T'_2}(x_0)$, if $a=x_1$ then $(T'_2-u)\rightarrow \{x_2,x_3\}$, let $P\equiv a_6a_7...a_{n-3}$ be a path in $T'_2-u$ and $\{x,y\}=V(T'_2-u)-V(P)$, we have $x_0x_1x_2xyx_3Pvux_0\equiv C$, so $x_2=a$.
If $a_7\rightarrow a_6$, then every vertex $x$ of $T_1$ has an inneighbor origin of a path $P\equiv \overline{a_2a_1...a_7}$ in $T_2$ and thus $I_XvP^{-1}x\equiv \overline{a_3a_4...a_3}=\overline{C}$, so $a_6\rightarrow a_7$. Similarly we may assume that $a_{n-4}\rightarrow a_{n-3}$. Let $x,y\in N^-_{T_2}(x_1)$ such that $x\rightarrow y$, if $x_0$ has an outneighbor origin of a path $P\equiv a_{n-3}a_{n-4}...a_4$ in $T_2\setminus \{x,y\}$ then $x_1xyvP^{-1}x_0x_3x_1\equiv C$, otherwise, $(T_2\setminus \{x,y\};a_{n-3}a_{n-4}...a_4)$ is one of the following exceptions:\\
Dual Exc $1$: If $d^+_{T'_2-u}(x_3)\geq 2$ then $x_3$ has an outneighbor origin of a path $P\equiv a_4...a_7$ in $T'_2-u$ and thus $x_0...x_3Pvux_0\equiv C$ unless $(T'_2-u;a_4...a_7)=(4A;(1,1,1))$, in this case there exists $w\in N^+_{T'_2}(x_0)$ such that $x_3$ has an outneighbor origin of a path $P\equiv a_0a_10a_8a_7=-(2,1)$ in $T'_2-w$ and thus $x_3...x_0wvP^{-1}x_3\equiv a_1a_2...a_1=C$. So we may assume that $d^-_{T'_2}(x_3)\geq 3$, so $x_3$ has an inneighbor origin of a path $P\equiv a_6a_7...a_0=(1,1,2)$ in $T'_2$ and thus $vax_1x_0x_3Pv\equiv a_1a_2...a_1=C$.\\
Dual Exc $5,7$: If $d^+_{T'_2-u}(x_3)\geq 3$ then $x_3$ has an outneighbor origin of a path $P\equiv a_4a_5...a_{n-3}=(1,1,2)$ in $T_2-u$ and thus $x_0...x_3Pvux_0\equiv C$. Otherwise, let $w\in N^+_{T'_2}(x_0)$ such that $d^-_{T'_2-w}(x_3)\geq 3$, $x_3$ has an inneighbor origin of a path $P\equiv a_{n-5}a_{n-6}...a_2=(1,2,1)$ in $T'_2-w$ and thus $wvP^{-1}x_3...x_0w\equiv C$.\\
Dual Exc $24$: We have $541x_0x_3x_1xyv2365\equiv C$.\\
Assume now that $(T_1;a_1...a_{\delta})$ is not a Gr\"{u}nbaum's exception, let $x_1...x_{\delta}\equiv a_1...a_{\delta}$ be a path in $T_1$, if $x_1$ has an outneighbor $w$ origin of a path $P\equiv a_0a_{n-1}...a_{\delta+2}$ in $T_2$ then $wx_1...x_{\delta}vP^{-1}\equiv C$ and if there exists $w\in N^-_{T_2}(x_1)\cap N^-_{T_2}(x_2)$ such that $T_2-w$ contains a path $Q\equiv a_{n-1}a_{n-2}...a_{\delta+2}$ with origin an inneighbor of $x_1$ then $x_1wx_2...x_{\delta}vP^{-1}x_1\equiv C$. Otherwise, we distinguish $4$ cases:\begin{itemize}
\item[(i)] If $d^+_{T_2}(x_1)\geq 2$ and $s^-_{T_2}(N^+_{T_2}(x_1))\geq 3$ and so $(T_2;a_0a_{n-1}...a_{\delta+2})$ is an exception. By the building lemmas, we may assume that it is one of the following biexceptions with respect to $x_1$:\\
Dual Exc $(3,1)$: $d^-(1)\geq 5$, a contradiction.\\
Dual Exc $(11,1)$: $T(N^-(4))\neq 5E$ nor $F_2$, so $T$ contains $C$ or it will be solved.\\
Dual Exc $(E_2(5),1)$: $T(N^+(x_2))\neq \overline{F_2}$ nor $5\overline{E}$, so $T$ contains $C$ or it will be solved.\\
Dual Exc $(E_8(6),1)$: $d^-(1)\geq 7$, a contradiction.
\item[(ii)] If $d^+_{T_2}(x_1)=s^-_{T_2}(N^+_{T_2}(x_1))=2$, set $N^+_{T_2}(x_1)=\{u\rightarrow w\}$, clearly $N^-_{T_2}(x_2)=\{w\}$. If $|T_2|\geq 5$, since otherwise there is a vertex in $T_2$ of indegree greater than $4$. If $|T_2|=5$ then $T_2\setminus \{u,w\}=3A$, let $t_1t_2t_3=O_{T_2\setminus \{u,w\}}$, we have $t_3wx_2...x_{\delta}vux_1t_1t_2t_3\equiv C$. So $\delta\geq 6$, let $a\in T_2\setminus \{u,w\}$ such that $T_2\setminus \{u,w,a\}\neq 3A,5A,7A$, and let $y_1...y_{\delta}\equiv a_1...a_{\delta}$ be a path in $(T_1-x_2)\cup \{a\}$, $y_1$ has an outneighbor $t$ origin of a path $P\equiv a_0a_{n-1}...a_{\delta+2}$ in $(T_2-a)\cup \{x_2\}$ of end distinct from $x_2$, thus $ty_1...y_{\delta}vP^{-1}\equiv C$, unless $|T|=9$ and $a_4...a_7$ is directed, in this case $w$ is an origin of a path $P\equiv a_8a_7...a_3=-(1,3,1)$ in $T_2$ and thus $x_2x_1vP^{-1}x_2\equiv C$.
\item[(iii)] If $N^+_{T_2}(x_1)=\{u\}$, suppose that there exists $w\in N^-_{T_2}(x_1)\cap N^-_{T_2}(x_2)$ so $(T_2-w;a_{n-1}a_{n-2}...a_{\delta+2})$ is a Gr\"{u}nbaum's exception, suppose that $u\equiv 1$, if $T_2-w=5A$ or $7A$ then $N^+(w)\cap \{2,3\}\neq \emptyset$ and so $w$ is an origin of a path $P\equiv a_1a_2...a_{\delta+2}$ in $T_2\cup \{x_1\}$ with end distinct from $x_1$ and thus $x_2...x_{\delta}vP^{-1}x_2\equiv a_2a_3...a_2=C$. If $T_2-w=3A$, then $w\rightarrow 2$ and so $3\rightarrow w$. If $w\rightarrow u$ then $uwx_2...x_{\delta}v32x_1u\equiv C$, so $u\rightarrow w$. If $\delta=2$, by taking $C=u_0u_1...u_0\equiv a_3a_4...a_3$ it is solved in lemma \ref{3.3}. If $T_1$ is a transitive triangle then $x_1$ is an origin of an antidirected outpath $Q$ in $T_1$ and $x_1$ has an inneighbor origin of a path $P\equiv a_3a_2a_1a_0$, thus $QvP^{-1}x_1\equiv a_4a_5...a_4=C$. If $3$, which is an origin of a path $P\equiv a_0a_{n-1}...a_{\delta+2}$ in $T_2$, has an inneighbor in $T_1$ that is an origin of a path $a_1...a_{\delta}$ in $T_1$ then $T$ contains $C$. So we may assume that $|Or(a_1...a_{\delta};T_1)|leq 2$ and $3\rightarrow Or(a_1...a_{\delta};T_1)$. If $a_{\delta}\rightarrow w$ then $1x_1...x_{\delta}w2v31\equiv $, so we may assume that $N^+_{T_1}(w)=\{x_1,x_2,x_4\}$. We have $d^+_{T_1}(x_1)=2$ so $T_1=4A$. If $a_1...a_4$ is antidirected then $d^+_{T_2}(x_4)=2$, so $x_4$ has an outneighbor $x$ in $T_2$ such that $T_2-x$ contains a path $P=-(1,1)$ of origin an inneighbor of $x_1$ and thus $xvP^{-1}x_1...x_4x\equiv C$. If $a_1...a_4=(1,2)$ then $Or((1,2);T_1)=T_1-Or((2,1);T_1)$ and thus there exists $y_4...y_7\equiv a_4...a_7=(2,1)$ such that $y_4\rightarrow 3$ therefore $21w3y_4...y_7v2\equiv C$. If $a_1...a_4=(2,1)$ then $\{2,3\}\rightarrow \{x_1,x_3\}=Or((2,1);T_1)$, thus $x_33x_1x_4x_2vw13x_3\equiv C$. Assume now that $N^-_{T_2}(x_2)=\{u\}$, we may suppose that $d^-_{T_2}(u)=0$ since otherwise, by taking $x_2$ in the place of $v$ one can find $C$. If $a_{\delta+2}a_{\delta+3}...a_{n-1}$ has at least $4$ blocks, then $(T_2-u)\cup \{x_1\}$ contains a path $P\equiv a_0a_{n-1}...a_{\delta+2}$ with origin $w$ and end in $T_2$, thus $wux_2...x_{\delta}vP^{-1}\equiv C$, so $a_{\delta+2}...a_{n-1}$ has only two blocks. If $\delta=2$, let $x\rightarrow y$ be an arc in $T_2-u$, we have $yxx_1vI_{T_2\setminus \{u,x,y\}}ux_2y\equiv C$. If $\delta=3$, let $P\equiv a_3a_4...a_{n-3}$ be a path in $T_2-u$ and suppose that $x_1x_3x_2=-(1,1)$, we have $x_2x_3x_1Pvux_2\equiv C$. So $\delta=4$, if $a_1...a_4$ is antidirected, let $P\equiv a_5a_6...a_{n-1}$ be a path in $T_2$, thus $x_4...x_1vPx_4\equiv C$. So $a_1...a_4=(1,2)$, if $|T_2|=4$ then $u$ is an origin of a path $P=(2,1,1)$ in $T_2\cup \{x_1\}$ with end in $T_2$, thus $vx_4x_3x_2Pv\equiv a_3a_4...a_3=C$; and if $|T_2|=5$ then $x_1$ has an inneighbor origin of a path $P=(2,1,1)$ in $T_2$ and thus $vx_4...x_1Pv\equiv a_3a_4...a_3=C$.
\item[(iv)] If $d^+_{T_2}(x_1)=0$, so $x_1$ is the unique origin of the path $a_1...a_{\delta}$ in $T_1$. Let $w\in N^-_{T_2}(x_2)$, we have $(T_2-w;a_{n-1}a_{n-2}...a_{\delta+2})$ is a Gr\"{u}nbaum's exception. $x_{\delta}$ has an outneighbor in $T_2-w$, but any vertex in $T_2-w$ is an origin of a path $P\equiv a_{\delta+1}a_{\delta+2}...a_{n-3}$ in $T_2-w$, so let $P$ be a path of origin an outneighbor of $x_{\delta}$ and set $\{a\}=V(T_2-w)-V(P)$, we have $x_1wx_2...x_{\delta}Pvax_1\equiv C$.\end{itemize}
\end{itemize}
\textbf{Case $3$}: If $|B|=|B'|=3$. Since the first and the last arcs of the path $a_0...a_{\delta-1}$ have the same direction then $(T_1;a_0...a_{\delta-1})$ is not a Gr\"{u}nbaum's exception, so let $x_0...x_{\delta-1}\equiv a_0...a_{\delta-1}$ be a path in $T_1$ such that $d^-_{T_2}(x_0)=max\,\{d^-_{T_2}(x);x\in Or(a_0...a_{\delta-1};T_1)\}$, if $x_0$ has an inneighbor origin of a path $P\equiv a_{n-1}a_{n-2}...a_{\delta+1}$ in $T_2$ then $x_0...x_{\delta-1}vP^{-1}x_0\equiv C$, otherwise, we distinguish three cases:\begin{itemize}
\item[(i)] If $d^-_{T_2}(x_0)=2$ and $s^-_{T_2}(N^-_{T_2}(x_0))=2$, set $N^-_{T_2}(x_0)=\{u\rightarrow w\}$. Suppose that $x_1$ has an outneighbor $t\in T_2\setminus \{u,w\}$, if $T_2\setminus \{u,w,t\}$ contains a path $P\equiv a_{n-4}a_{n-5}...a_{\delta+1}$ then $tx_1...x_{\delta-1}vP^{-1}uwx_0t\equiv C$, else, if $t$ has an inneighbor $t'\in T_2\setminus \{u,w\}$ then $T_2\setminus \{u,w,t,t'\}$ contains a path $P\equiv a_{n-5}a_{n-6}...a_{\delta+1}$ and thus $tx_1...x_{\delta-1}vP^{-1}wux_0t't\equiv C$. Otherwise, $t\rightarrow (T_2\setminus \{u,w,t\}$, if $x_{\delta-1}$ has an outneighbor $x$ in $T_2\setminus \{u,w,t\}$ then $|T_2|\geq 8$ since otherwise $d^-(x)>\Delta^+(T)$, let $xyz$ be a direccted outpath in $T_2$ and $P\equiv a_{n-4}a_{n-5}...a_{\delta+4}$ be a path in $T_2\setminus \{u,w,t,x,y,z\}$, we have $tx_1...x_{\delta-1}xyzvP^{-1}uwx_0t\equiv C$, else, $(T_2\setminus \{u,w,t\})\rightarrow x_{\delta-1}$, if $t\rightarrow x_{\delta-2}$, let $x\rightarrow y$ be an arc in $T_2\setminus \{u,w,t\}$ and $P\equiv a_{n-4}a_{n-5}...a_{\delta+3}$ be a path in $T_2\setminus \{u,w,t,x,y\}$, we have $x_0...x_{\delta-2}txyx_{\delta-1}P^{-1}vuwx_0\equiv C$ and if there exists $x\in T_2\setminus \{u,w,t\}$ an inneighbor of $x_{\delta-2}$, let $P\equiv a_{n-4}a_{n-5}...a_{\delta+2}$ be a path in $T_2\setminus \{u,w,t,x\}$, we have $tx_1...x_{\delta-2}xx_{\delta-1}vP^{-1}uwx_0t\equiv C$. So we may assume that $N^-_{T_2}(x_{\delta-2})=\{w\}$, let $xy$ be an arc in $T_2\setminus \{u,w,t\}$ and $P\equiv a_{\delta+3}a_{\delta+4}...a_{n-4}$ in $T_2\setminus \{u,w,x,y\}$, we have $tx_1...x_{\delta-2}wxyx_{\delta-1}Pvux_0t\equiv C$. So we can suppose now that $(T_2\setminus \{u,w\})\rightarrow x_1$, we have to mention that $T_1\cup \{u\}$ contains a path $Q\equiv a_{n-1}a_0...a_{\delta-1}$ of origin $x_1$ and end distinct from $u$, in fact, if $a_2\rightarrow a_1$ then $x_1x_0ux_2...x_{\delta-1}\equiv Q$ and if $a_1\rightarrow a_2$, $T_1-x_1$ contains a path $Q'\equiv a_2...a_{\delta-1}$ of origin distinct from $x_0$ and thus $x_1uQ'\equiv Q$, let $P\equiv a_{\delta+1}a_{\delta+2}...a_{n-4}$ be a path in $T_2-w$, and set $\{t\}=V(T_2)-V(P)$, we have $QvPwtx_1\equiv a_{n-1}a_0...a_{n-1}=C$.
\item[(ii)] If $d^-_{T_2}(x_0)\geq 2$ and $s^-_{T_2}(N^-_{T_2}(x_))\geq 3$ then $(T_2;a_{n-1}a_{n-2}...a_{\delta+1})$ is an exception, so, by the building lemmas, we may assume that it is one of the following biexceptions with respect to $x_0$:\\
Dual Exc $(7,0)$: If $d^+_{T_2}(x_{\delta-1})\geq 1$ then there exist $xyz$ and $ab$ two directed outpaths such that $V(T_2)=\{a,b,x,y,z\}$, $x\in N^+(x_{\delta-1})$ and $b\in N^-(x_0)$, thus $x_0...x_{\delta-1}xyzvabx_0\equiv C$, so $T_2\rightarrow x_{\delta-1}$. Let $u\in N^-_{T_2}(x_{\delta-2})$, there exists $P\equiv a_{n-1}a_{n-2}...a_{\delta+2}$ in $T_2-u$ of origin an inneighbor of $x_0$, we have $x_0...x_{\delta-2}ux_{\delta-1}vP^{-1}x_0\equiv C$.\\
Dual Exc $(E'_8(6),1),(E'_9(8),1),(E'_{10}(10),1)$: If there exists $x\in N^+_{T_2-X}(x_{\delta-1})$ then $|T_2|\geq 8$, let $xyz$ be a directed outpath in $T_2-X$ and $P\equiv a_{n-1}a_{n-2}...a_{\delta+4}$ be a path in $T_2\setminus \{x,y,z\}$ of origin an inneighbor of $x_0$, we have $x_0...x_{\delta-1}xyzvP^{-1}x_0\equiv C$, so $(T_2-X)\rightarrow x_{\delta-1}$. If there exists $x\in N^-_{T_2-X}(x_{\delta-2})$, let $P\equiv a_{n-1}a_{n-2}...a_{\delta+2}$ be a path in $T_2-x$ of origin an inneighbor of $x_0$, we have $x_0...x_{\delta-2}xx_{\delta-1}vP^{-1}x_0\equiv C$, so $x_{\delta-2}\rightarrow (T_2-X)$ and thus $|T_2|\geq 8$. If $x_{\delta-2}$ has an inneighbor $x\in X$ such that there is an inneighbor $y\in X-x$ of $x_0$ of positive indegree in $T_2-x$, set $\{z\rightarrow y\}=X-x$ and let $a,b\in T_2-X$ and $P\equiv a_{\delta+3}a_{\delta+4}...a_{n-4}$ be a path in $T_2\setminus \{x,y,z,a,b\}$, we have $x_0...x_{\delta-2}xabx_{\delta-1}Pvzyx_0\equiv C$, so we may suppose that $T_2=\{u,w,t\}$ is transitive such that $u\rightarrow w\rightarrow t$, $N^-_{T_2}(x_0)=\{u,t\}$ and $N^-_{T_2}(x_{\delta-2})=\{t\}$. Let $ab$ be an arc in $T_2-X$ and $P\equiv a_{n-1}a_{n-2}...a_{\delta+4}$ be a path in $T_2\setminus \{a,b\}$ with origin an inneighbor of $x_0$, we have $x_0...x_{\delta-3}vx_{\delta-2}abx_{\delta-1}P^{-1}x_0\equiv C$.
\item[(iii)] If $d^-_{T_2}(x_0)=1$, set $N^-_{T_2}(x_0)=\{u\}$, we have $x_0$ is maximal in $T_1$. By lemma \ref{2.3}, $a_0...a_{\delta-1}$ is directed. If $T_1$ is transitive with $\delta=3$ or $4$, then $x_0$ is an origin of a path $Q\equiv a_1...a_{\delta}$ in $T_1$ and thus $x_0$ has an outneighbor $w$ origin of a path $P\equiv a_0a_{n-1}...a_{\delta+2}$ in $T_2$, we have $wQvP^{-1}\equiv C$. If $T_1=\overline{4B}$, we can suppose that, without loss of generality, $\{x_1,x_2,x_3\}\rightarrow u$, so $d^-_{T_2-u}(x_1)\geq 3$ and thus $x_1$ has either an inneighbor origin of a path $P\equiv a_{n-2}a_{n-3}...a_{\delta+1}$ in $T_2-u$ and thus $x_0ux_2x_3vP^{-1}x_1x_0\equiv C$, or $(T_2-u;a_{n-2}a_{n-3}...a_{\delta+2})$ is one of the dual of the exceptions $22$, $24$, $38$ or $47$, in this case $x_1$ has an outneighbor origin of a path $Q\equiv a_{\delta+4}a_{\delta+5}...a_1$ in $T_2-u$ and thus $vx_2x_3ux_0x_1x_0Pv\equiv a_2a_3...a_2=C$. Now if $|T_1|=2$, we can suppose that $a_5\rightarrow a_6$ since otherwise by lemma \ref{3.4} one can find $\overline{C}=a'_0a'_1...a'_0=\overline{a)a_5a_4...a_5}$, also we can suppose that $a_{n-5}\rightarrow a_{n-4}$ since otherwise one can find $C=a_{n-4}a_{n-3}...a_{n-4}$. The blocks containing the arcs $a_{n-5}\rightarrow a_{n-4}$ and $a_5\rightarrow a_6$ cannot be of length $2$ since otherwise it is solved in case $2$. Obviously $u$ is not maximal in $T_2$ so let $w\in N^+_{T_2}(u)$, if $d^-_{T_2\setminus \{u,w\}}(x_1)\geq 2$ then $x_1$ has an inneighbor origin of a path $P\equiv a_{n-3}a_{n-4}...a_3$ in $T_2\setminus \{u,w\}$ and thus $wux_0P^{-1}x_1vw\equiv C$ unless $(T_2\setminus \{u,w\},a_{n-3}a_{n-4}...a_3)$ is one the exceptions Exc $0,4,8,9,33,37,E_8(6),E_9(7)$ or $E_{10}(9)$, if it is the exception $0$ then if $d^+_{T_2}(x_1)\geq 2$ then $T_2$ contains a path $P\equiv a_0a_{n-1}...a_3$ of origin an outneighbor $x$ of $x_1$ and thus $xx_1x_0vP^{-1}\equiv C$ unless $s^-_{T_2}(N^+_{T_2}(x_1))\leq 3$, but since $T_2$ contains a circuit triangle with at least two inneighbors of $x_1$ then $T_2=\{u,w\}\rightarrow \{1,2,3\}$ and $N^+_{T_2}(x_1)=\{u,w\}$, thus $vx_0w12u3x_1v\equiv C$, so $d^-_{T_2}(x_1)\geq 4$, let $t\in N^+_{T_2}(x_0)$, $x_1$ has an inneighbor origin of a path $P=(1,2)$ in $T_2-t$, we have $x_1x_0tvP^{-1}x_1\equiv a_2a_3...a_2=C$. Otherwise, one can check that $(T_2\cup \{x_1\})\setminus \{u,w\}$ contains a path $P\equiv a_{n-1}a_0...a_{n-6}$ of origin and end in $T_2$, thus $vwux_0Pv\equiv a_{n-5}a_{n-4}...a_{n-5}=C$. So we may assume that $d^-_{T_2\setminus \{u,w\}}(x_1)\leq 1$, in other terms $d^+_{T_2}(x_1)\geq |T_2|-3$, if $s^-_{T_2}(N^+_{T_2}(x_1))\geq 4$ then $x_1$ has an ouitneighbor $x$ origin of a path $P\equiv a_0a_{n-1}...a_3$ in $T_2$ and thus $xx_1x_0vP^{-1}\equiv C$ unless $(T_2;a_0a_{n-1}...a_3)=E_8(7)$, in this case $x_1$ has an outneighbor origin of a path $P\equiv a_{n-4}a_{n-3}...a_0$ in $T_2$, then $vP^{-1}x_1x_0v\equiv a_{n-1}a_0...a_{n-1}=C$, otherwise, $s^-_{T_2}(N^+_{T_2}(x_1)\leq 3$, then $|T_2|=5$ and $N^+_{T_2}(x_1)\rightarrow N^-_{T_2}(x_1)$, in this case $u$ and $x_1$ have a common outneighbor $w'$ then $x_1$ has an inneighbor $x$ origin of a path $P=(1,1)$ in $T_2\setminus \{u,w'\}$, we have $w'ux_0P^{-1}x_1vw\equiv C$.
\end{itemize}
\textbf{Case $4$}: If $|B|>|B'|\geq 3$. Rotate $C$ until $a_{n-1}\rightarrow a_{n-2}$ for the first time. Suppose first that $(T_1;a_1...a_{\delta})$ is a Gr\"{u}nbaum's exception and for $x_i\in T_1$, let $P_i\equiv a_3a_4...a_{\delta+2}$ be a Hamiltonian path of $T_1$ of origin $x_i$. If an $x_i$ has an  outneighbor origin of a path $P\equiv a_2a_1...a_{\delta+4}$ in $T_2$ then $P_ivP^{-1}x_i\equiv a_3a_4...a_4=C$. Otherwise, if $(T_2;a_2a_1...a_{\delta+4})$ is an exception so it is one of the duals of the exceptions $(E_8;(\Delta-4,1,1,1))$, $(E_8;(\Delta,-6,1,1,1,1,1))$ or $(E_{10};(\Delta-8,1,1,1,1,1,1,1))$, then $|N^-_{T_2}(T_1)|\geq 3$; else, $(T_2;a_2a_1...a_{\delta+4})$ is not an exception, so $s^-_{T_2}(\displaystyle \bigcup_{x\in T_1} N^+_{T_2}(x))\leq 3$, then $s^-_{T_2}(N^-_{T_2}(T_1))=|T_2|$ and $|N^-_{T_2}(T_1)|\geq 3$ unless $|T_2|=5$ and $T_1=3A$. In the first case, $|N^-_{T_2}(T_1)|\geq 3$, let $u,w,t\in N^-_{T_2}(T_1)$, since $(T_2-u;a_1a_0a_{n-1}...a_{\delta+4})$ and $(T_2-u;a_0a_{n-1}...a_{\delta+3})$ cannot be two exceptions in the same time, then $w$ or $t$ is an origin of a path $Q\equiv a_1a_0...a_{\delta+4}$ or $R\equiv a_0a_{n-1}...a_{\delta+3}$ in $T_2-u$ and thus $x_iu^*P_ivQ^{-1}x_i\equiv a_2a_3...a_2=C$ or $P'_ivR^{-1}x_i\equiv a_1a_2...a_1=C$, where $P'_i\equiv a_1a_2...a_{\delta+1}$ is a path in $T_1\cup \{u\}$ with origin $x_i$. Now if $|T_2|=5$ and $T_1=3A$, we have $N^+_{T_2}(T_1)\rightarrow N^-_{T_2}(T_1)$ where $N^+_{T_2}(T_1)=\{u_1,u_2,u_3\}=3A$ and $N^-_{T_2}(T_1)=\{b_1,b_2\}$, thus $u_1u_2b_2b_1P_iu_3va_1\equiv C$.\\
Assume now that $(T_1;a_1...a_{\delta})$ is not a Gr\"{u}nbaum's exception and let $x_1...x_{\delta}\equiv a_1...a_{\delta}$ be a path in $T_1$ such that $d^-_{T_2}(x_1)=max\,\{d^-_{T_2}(x);x\in Or(a_1...a_{\delta},T_1)\}$. If $x_1$ has an inneighbor $u$ origin of a path $P\equiv a_0a_{n-1}...a_{\delta}$ in $T_2$ then $ux_1...x_{\delta}vP^{-1}\equiv C$. Otherwise, since $d^-_{T_2}(x_1)\geq 2$, $(T_2;a_0a_{n-1}...a_{\delta})$ is one of the following biexceptions with respect to $x_1$:\\
Dual Exc $(6,0),(31,0),(E_{11},0)$ and $(E_{13},0)$: In this case, $d^-_{T_2}(x_1)=2$, so by lemma \ref{2.4}, since $\delta \geq 3$ and $a_1...a_{\delta}$ is not directed, $T_1=4B$ and $a_1...a_{\delta}$ has $3$ origins, which implies that there exists a common inneighbor $z$ of the origins of $a_1...a_{\delta}$ in $T_2$ $d^+(z)>\Delta^+(T)$, a contradiction.\\
Dual Exc $(26,0)$: If $x_{\delta}\rightarrow 1$ then $4x_1...x_{\delta}16523v4\equiv C$. In a similar way, we can assume that $\{1,2,3,4,5\}\rightarrow x_{\delta}$, then $T(N^+(4)\neq 6H$, thus $T$ contains $C$.\\
\textbf{Case $5$}: If $|B'|=1$ and $a_1\rightarrow a_2$. Suppose first that $(T_1;a_1...a_{\delta})$ is a Gr\"{u}nbaum's exception. If $a_{\delta+4}\rightarrow a_{\delta+3}$, let $x_3...a_{\delta+2}\equiv \overline{a_3...a_{\delta+2}}$ be a path in $T_1$, since $d^-_{T_2}(x_3)\geq 2$ then $x_3$ has an inneighbor origin of a path $P\equiv \overline{a_2a_1...a_{\delta+4}}$ and thus $x_3...x_{\delta+2}vP^{-1}x_3\equiv \overline{a_3a_4...a_3}=\overline{C}$ unless $(T_2;\overline{a_2a_1...a_{\delta+4}})$ is an exception and by the building lemmas we may assume that it is a biexception with respect to $x_3$, but since any vertex of $T_1$ is an origin of a path $Q\equiv \overline{a_3...a_{\delta+2}}$ then if $N^-_{T_2}(x_3)\cup S\neq V(T_2)$, one can find $C$, thus we are in one of the following biexceptions with respect to $x_3$:\\
Exc $(22,1)$: We have $|T_1|=3$, set $V(T_1)=\{u_1,u_2,u_3\}$ we have $T_1\rightarrow \{4,5,6\}$ and $\{1,2,3\}\rightarrow T_1$, thus $u_11u_3u_2v6452u_1\equiv C$.\\
Exc $(24,1)$: $d^-(4)>\delta$ a contradiction.\\
So we may assume that $a_{\delta+3}\rightarrow a_{\delta+4}$. If there exists $x\in T_1$ such that $d^+_{T_2}(x)=1$ then $T_1=3A=\{x_1,x_2,x_3\}$, without loss of generality, suppose that $d^+_{T_2}(x_1)=1$. Let $u\in N^+_{T_2}(x_2)$, since $d^-_{T_2}(x_1)\geq |T_2-u|-1$ then $x_1$ has an inneighbor origin of a path $P\equiv a_{\delta+3}a_{\delta+4}...a_0$ in $T_2-u$ and thus $vux_2x_3x_1Pv\equiv a_1a_2...a_2=C$ unless four cases:\begin{itemize}
\item[(i)] If $T_2-u=\{u_1\rightarrow u_2\}$ with $x_1\rightarrow u_1$. If $x_3\rightarrow u$ then $u_2vux_3x_1x_2u_1u_2\equiv C$, so $u\rightarrow x_3$. If $u\rightarrow u_2$ then $u_2vu_1x_1x_2x_3uu_2\equiv C$, so $u_2\rightarrow u$ and without loss of generality $x_3\rightarrow u_2$. If $u_1\rightarrow u$ then $uvu_2x_3x_1x_2u_1u\equiv C$, so $u\rightarrow u_1$ and thus $T=7A$ which contains no $(2,1,1,1,1,1)$, $(T;C)=A_3$.
\item[(ii)] If $T_2-u=\{u_1,u_2,u_3\}$ is transitive such that $x_1\rightarrow u_1\rightarrow \{u_2,u_3\}$. We can suppose that $u_1\rightarrow x_2$ and thus $u\rightarrow x_3$, we have $x_1ux_3x_2vu_3u_1u_2x_1\equiv C$.
\item[(iii)] If $(T_2-u;a_6...a_0)=(4A;(1,1,1))$ and $N^+_{T_2}(x_1)=\{4\}$. If $1\rightarrow u$ then $21ux_2x_3x_134$\\$v2\equiv C$, if $3\rightarrow u$ then $43ux_2x_3x_112v4\equiv C$ and if $4\rightarrow u$ then $14ux_2x_3x_123v1\equiv C$; so $u\rightarrow \{1,3,4\}$. If $2\rightarrow x_3$ then $1v4x_1x_2x_323u1\equiv C$ and if $x_3\rightarrow 2$ then $x_12x_3x_2v1u43x_1\equiv C$.
\item[(iv)] If $(T_2-u;a_6...a_0)=(5C;(2,1,1))$ and $N^+_{T_2}(x_1)=\{5\}$. We can suppose that $5\rightarrow x_2$ since otherwise $x_1$ has an inneighbor origin of a path $P=(2,1,1)$ in $T_2-5$ and thus $v5x_2x_3x_1Pv\equiv a_1a_2...a_2=C$. If $N^-(u)\cap \{1,2,3\}\neq \emptyset$, say $1\rightarrow u$, then $21ux_2x_3x_1345v2\equiv C$, so $u\rightarrow \{1,2,3\}$. If $4\rightarrow u$ then $54ux_2x_3x_1123v5\equiv C$, so $u\rightarrow 4$. If $5\rightarrow u$ then $15ux_2x_3x_1234v1\equiv C$, so $u\rightarrow 5$, and thus $x_3\rightarrow u$. Therefore $3vux_3x_1x_251423\equiv C$.
\end{itemize}
Now we can suppose that $d^+_{T_2}(x)\geq 2$ for all $x\in T_1$. Let $a\in T_2$ such that $d^+_{T_2}(a)=\Delta^+(T_2)$ and set $T'_1=T_1\cup \{a\}$ and $T'_2=T_2-a$. Let $x_0...x_{\delta}\equiv a_0...a_{\delta}$ be a path in $T'_1$, if $x_0$ has an inneighbor origin of a path $P\equiv a_{n-1}a_{n-2}...a_{\delta+2}$ in $T'_2$ then $x_0...x_{\delta}vP^{-1}x_0\equiv C$. Otherwise, if $d^-_{T'_2}(x_0)\geq 2$ then $(T'_2;a_{n-1}a_{n-2}...a_{\delta+2})$ is an exception. By the building lemmas we may assume that it is one of the following biexceptions with respect to $x_0$:\\
Exc $(1,1)(1),(2),(4)$: By definition of $a$, $a$ cannot be $x_1$, so $\{2,4\}\rightarrow x_{\delta}$, since otherwise by taking $x_{\delta}...x_1\equiv a_1...a_{\delta}$ the result follows. If $x_{\delta}\rightarrow 3$ then $x_0...x_{\delta}34v21x_0\equiv C$, so $3\rightarrow x_{\delta}$. Since $d^+_{T_2}(x_{\delta})\geq 2$ then$x_{\delta}\rightarrow 1$, if $3\rightarrow x_1$ then $x_0...x_{\delta}12v43x_0\equiv C$, so $x_1\rightarrow 3$. If $\delta=5$, $T'_1$ has at least $4$ origins of an antidirected Hamiltonian inpath and since so $d^+(a)\geq 6$, a contradiction, in fact, $d^-_{T_2}(x)\geq 3$ for all $x\in T_1$ so we can suppose that any origin of an antidirected Hamiltonian inpath in $T'_1$ is an outneighbor of $a$. Thus $\delta=3$ and so $3vax_1x_0x_34213\equiv C$.\\
Exc $(1,1)(3)$: As above wen can assume that $\delta=3$. If $x_{\delta}\rightarrow 4$ then $x_0...x_{\delta}41v32x_0\equiv C$ and if $4\rightarrow x_{\delta}$ then $3vax_1x_0x_34213\equiv C$.\\
Exc $(5,1)$: By definition of $a$, $a$ cannot be $x_1$, so $x_{\delta}\rightarrow 5$ and thus $x_0...x_{\delta}51v432x_0\equiv C$.\\
Exc $(7,1)$: Similarly we can suppose that $\{1,2,3\}\rightarrow x_{\delta}$. If $x_{\delta}\rightarrow 4$ then $x_0...x_{\delta}45v321$\\$x_0\equiv C$ and if $x_{\delta}\rightarrow 5$ then $x_0...x_{\delta}51v432x_0\equiv C$.\\
Exc $(10,1)$: If $\delta=5$, since $T'_1$ contains at least $4$ origin of an antidirected Hamiltonian inpath and one of them is an outneighbor of $2$, so the result follows and thus $\delta=3$. Let $y_1y_2y_3$ be a directed outpath of $T_1$ of origin $x_0$, one can find an outneighbor $u\in T'_2$ of $x_0$ such that $y_3$ has an inneighbor in $T_2-u$ origin of a path $Q=(2,1,1)$ and thus $vuy_1y_2y_3Qv\equiv a_1a_2...a_1=C$.\\
Exc $(22,1)$: If $T_2=7A$ then there exists $a'\in T_2$ such that there is a path $x'_0...x'_{\delta}\equiv a_0...a_{\delta}$ in $T_1\cup \{a'\}$ with $x'_0$ has an inneighbor origin of a path $P\equiv a_{n-1}a_{n-2}...a_{\delta+2}$ in $T_2-a'$ and the result follows, so $d^+_{T_2}(a)\geq 4$. There exists $a'\in T_2$ such that $T_2-a'\neq 6F,6H$ and the result follows; in fact, the condition $d^+_{T_2}(a')\geq 1$ is sufficient instead of $d^+_{T_2}(a')=\Delta^+(T_2)$.\\
Exc $(24,1)$: By definition of $a$, $a$ cannot be $x_1$ so $\{1,2,3,6\}\rightarrow x_{\delta}$. So $N^+_{T'_2}(x_{\delta})\subseteq \{4,5\}$, say $x_{\delta}\rightarrow 4$, we have $x_0...x_{\delta}46v1253x_0\equiv C$.\\
Else, set $N^-_{T'_2}(x_0)=\{u\}$, since $d^-_{T_2}(x)\geq 2$ for all $x\in T_1$ then $a\equiv x_1$. Since $u$ is not an origin of a path $P\equiv a_{n-1}a_{n-2}...a_{\delta+2}$ then $N^+_{T'_2}(a)\cap N^+_{T'_2}(x_0)\neq \emptyset$, let $w\in N^+_{T'_2}(a)\cap N^+_{T'_2}(x_0)$. If $T'_2-w$ contains a path $P\equiv a_{n-2}a_{n-3}...a_{\delta+2}$ of origin distinct from $u$ then $wx_1...x_{\delta}vP^{-1}x_0w\equiv C$.  Otherwise, two cases may arise:\begin{itemize}
\item[(i)] If $(T'_2-w;a_{n-2}...a_{\delta+2})=(4A;(1,1,1))$, there exists $w'\in N^+_{T'_2}(a)$ such that $T'_2-w'$ contains a path $P$ and the result follows. In fact, otherwise, $u$ is an origin of a path $Q=(2,1,1)$ in $T'_2$ and thus $x_0...x_{\delta}vQ^{-1}x_0\equiv C$.
\item[(ii)] If $(T'_2-w;a_{n-2}...a_{\delta+2})$ is a Gr\"{u}nbaum's exception, if there exists $w'\in T'_2-w$ such that $T'_2-w'\neq 3A,5A,7A$, the result follows, so we can assume that $N^+_{T'_2-u}(a)\subseteq \{x\in T'_2-u/ T'_2-x\neq 3A,5A,7A\}$. By completing $T_2$ under the conditions $|N^+_{T'_2-u}(a)|\leq 2$ and $d^+_{T_2}(a)=\Delta^+(T_2)$, $u$ is an origin of a path $P\equiv a_{n-1}a_{n-2}...a_{\delta+2}$ in $T'_2$ which is a contradiction. 
\end{itemize}
Assume now that $(T_1;a_1...a_{\delta})$ is not a Gr\"{u}nbaum's exception and let $x_1...x_{\delta}$ be a path in $T_1$ such that $d^+_{T_2}(x_1)=max\,\{d^+_{T_2}(x);x\in Or(a_1...a_{\delta},T_1)\}$. If $x_1$ has an outneighbor $u$ origin of a path $P\equiv a_0a_{n-1}...a_{\delta+2}$ in $T_2$ then $ux_1...x_{\delta}vP^{-1}\equiv C$. Otherwise, if $d^+_{T_2}(x_1)\geq 2$ then $(T_2;a_0a_{n-1}...a_{\delta+2})$ is an exception. By the building lemmas we may assume that it is one of the following biexceptions with respect to $x_1$:\\
Dual Exc $(0,1)$: $T(N^-(x_2))=T(x_1,1,2)\neq 3A$ then $T$ contains $C$.\\
Dual Exc $(2,1)$: $T(N^-(1))\neq 4A,F_1(4)$ so $T$ contains $C$.\\
Dual Exc $(4,1)$, $(9,1)(2),(3),(4)$, $(15,1)$, $(16,1)(2)$, $(33,1)$ and $(E_5,1)(2)$: Since $T(N^+(x_2))$ has $v$ as minimal vertex then $(T(N^+(x_2));a_0a_{n-1}...a_{\delta+2})$ is not an exception, unless it is $E_8(5)$, this case will treated below.\\
Dual Exc $(8,1)$: If $d^+_{T_2}(x)=4$ for an $x\in T_1$ then $T(N^+(x))\neq 5A,5D,5E$, so we may assume that $d^-_{T_2}(x)\geq 2$ for all $x\in T_1$. If $x_{\delta}\rightarrow 2$ then $5x_1...x_{\delta}v21v345\equiv C$, if $x_{\delta}\rightarrow 3$ then $5x_1...x_{\delta}3124v5\equiv C$ and if $x_{\delta}\rightarrow 4$ then $5x_1...x_{\delta}43v125\equiv C$, so $\{2,3,4\}\rightarrow x_{\delta}$ and thus $x_2\neq x_{\delta}$ i.e. $\delta\geq 2$. If $\delta=3$ then $35x_2x_1v12x_343\equiv C$. If $\delta=4$ then $(T_2;\overline{a_0a_1...a_4})$ is not an exception, so let $y_1...y_4\equiv \overline{a_6...a_9}$, since $d^-_{T_2}(y_4)\geq 2$ then $y_4$ has an inneighbor $u$ origin of a path $P\equiv \overline{a_0...a_4}$ in $T_2$ and thus $Pvy_1...y_4u\equiv C$. So $\delta=5$ and $x_{\delta}...x_1\equiv x_1...x_{\delta}$ since otherwise there is a block of length $3$. Since $d^-_{T_2}(x_4)\geq 2$, one can find an inneighbor of $x_4$ origin of a path $P\equiv a_1a_0...a_7$ in $T_2\cup \{x_5\}$  of end distinct from $x_5$ and thus $x_4x_3...x_1vP^{-1}x_4\equiv a_2a_3...a_2=C$.\\
Dual Exc $(9,1)(1)$: If $N^+(x_{\delta})\cap \{2,3,4,5\}\neq \emptyset$, one can find $w,t\in T_2$ and $P=-(1,1)$ of origin an outneighbor $u$ of $x_1$ in $T_2\setminus \{w,t\}$ such that $x_{\delta}\rightarrow w\rightarrow t$ and thus $ux_1...x_{\delta}wtvP^{-1}\equiv C$, and if $x_{\delta}\rightarrow 1$ then $2x_1...x_{\delta}1453v2\equiv C$.\\
Dual Exc $(16,1)$: If $d^+_{T_2}(x_{\delta})\geq 1$, one can check that $x_{\delta}$ has an outneighbor origin of a path $P\equiv a_{\delta+1}a_{\delta+2}...a_0$ in $T_2\cup \{v\}$ of end $4$ and so $4x_1...x_{\delta}P\equiv C$, so $a_{\delta}\rightarrow a_{\delta-1}$. Let $u\in N^-_{T_2}(x_{\delta-1})$, $u$ is an origin of a path $P=(2,1)$ in $T_2\setminus \{w,t\}$ where $w\in N^+_{T_2}(x_1)$ and $t\in T_2$, thus $wx_1...x_{\delta-1}Px_{\delta}tvw\equiv C$.\\
Dual Exc $(32,1)$: If $d^+_{T_2}(x_{\delta})\geq 1$, one can check that $x_{\delta}$ has an outneighbor origin of a path $P\equiv a_{\delta+1}a_{\delta+2}...a_0$ in $T_2\cup \{v\}$ of end $5$ or $6$ and thus $T$ contains $C$. Otherwise, let $u\in N^-_{T_2}(x_{\delta-1})$, as above one can find $C$.\\
Dual Exc $(E_1,1)$: Let $a_s$ be the origin of $B$, $T$ contains $\overline{C}=a'_0a'_1...a'_0=\overline{a_{s+\delta+2}a_{s+\delta+1}...}$\\$\overline{a_{s+\delta+2}}$.\\
Dual Exc $(E_5,1)(1)$: If $x_{\delta}\rightarrow 1$ then $wx_1...x_{\delta}12v3uw\equiv C$, so we may assume that $\{1,2,3\}\rightarrow x_{\delta}$. If $a_{\delta-1}\rightarrow a_{\delta}$ then $wx_1...x_{\delta-1}v1x_{\delta}23uw\equiv C$, so $a_{\delta}\rightarrow a_{\delta-1}$. Then $\delta=2$ and thus $uw3vx_2x_112u\equiv \overline{C}$.\\
Dual Exc $(E_8(5),1)$: If $N^+(x_{\delta})\cap \{1,2,3\}\neq \emptyset$, say $1\in N^+(x_{\delta})$, we have $ux_1...x_{\delta}12v3w$\\$u\equiv C$, so $\{1,2,3\}\rightarrow x_{\delta}$ and thus $x_2\neq x_{\delta}$ i.e. $\delta\geq 3$. If $T_1=3A,5A$ or $7A$, $x_2$ is an origin of a path $y_1...y_{\delta}\equiv a_1...a_{\delta}$ and so $1y_1...y_{\delta}v3u2w1\equiv C$, so $T_1\neq 3A,5A,7A$. If $a_{\delta-1}\rightarrow a_{\delta}$, let $u_0...u_{\delta-1}\equiv a_0...a_{\delta-1}$ be a path in $T_1$, any vertex of $T_2$ is an origin of a path $P\equiv a_{n-1}...a_{\delta+1}$ so $u_0...u_{\delta-1}vP^{-1}u_0\equiv C$, so $a_{\delta}\rightarrow a_{\delta-1}$. If $N^-(x_{\delta-1})\cap \{u,w\}\neq \emptyset$, say $u\in N^-(x_{\delta-1})$, we have $wx_1...x_{\delta-1}u1x_{\delta}23vw\equiv C$ so $x_{\delta-1}\rightarrow \{u,w\}$ and thus $x_{\delta-1}\neq x_2$ i.e. $\delta\geq 4$. If $a_1...a_{\delta}$ is symmetric then $(T_2;a_0a_{n-1}...a_{\delta+2})$ is not a biexception with respect to $x_{\delta}$, so by taking $x_{\delta}...x_1\equiv a_1...a_{\delta}$ one can find $C$. So we may assume that $\delta=4$, let $u_1...u_4\equiv \overline{a_6...a_9}$, every vertex of $T_2$ is an origin of a path $P\equiv \overline{a_0a_1...a_4}$, thus $u_4Pvu_1...u_4\equiv \overline{a_9a_0...a_9}=C$.\\
Else, $d^+_{T_2}(x_1)\leq 1$, if there exist $w\in N^-_{T_2}(x_1)\cap N^-_{T_2}(x_2)$ and $P\equiv a_{n-1}a_{n-2}...a_{\delta+2}$ in $T_2-w$ of origin an inneighbor of $x_1$ then $x_1wx_2...a_{\delta}vP^{-1}x_1\equiv C$. Otherwise, since $(T_2-w;a_{n-1}a_{n-2}...a_{\delta+2})$ is not a Gr\"{u}nbaum's exception, $d^+_{T_2}(x_1)=1$, set $N^+_{T_2}(x_1)=\{u\}$. Two cases may arise:\begin{itemize}
\item[(i)] If $d^-_{T_2}(x_2)\geq 2$ then either $(T_2-w;a_{n-1}a_{n-2}...a_{\delta+2})=(4A;(1,1,1))$ or $(5C;(2,1,1))$ or $u$ is minimal in $T_2$ and $a_{n-1}a_{n-2}...a_{\delta+2}$ is directed. If $(T_2-w;a_{n-1}...a_{\delta+2})=(4A;(1,1,1))$, $u\equiv 4$, if $4\rightarrow x_2$ then, without loss of generality suppose that $3\rightarrow w$, $24x_2...x_{\delta}vw3x_112\equiv C$; so $x_2\rightarrow 4$. Now we suppose that for any $w'\in N^-_{T_2}(x_2)$ $T_2-w'=4A$ with $u\equiv 4$. So if $3\rightarrow x_2$ then $u$ is an origin of a path $P\equiv a_0a_{n-1}...a_{\delta+2}$, a contradiction. So $4\rightarrow w$ and $w\rightarrow 3$, if there is a path $P=(1,1)$ in $T(1,2,w)$ of origin an inneighbor of $x_2$ then $x_2...x_{\delta}v43x_1P^{-1}x_2\equiv a_2a_3...a_2=C$, so $T(1,2,w)=3A$. If $\delta=2$ then $x_2x_1v4213wx_2\equiv C$, if $\delta=3$ and $T_1=3A$ then $3x_2x_3x_1v241w3\equiv C$, if $\delta=3$ and $T_1$ is transitive, if $a_1a_2a_3$ is antidirected then $2v4x_1x_3x_2w312\equiv C$ and if $a_1a_2a_3$ is directed then $x_2x_1x_3v4213wx_2\equiv C$, finally if $\delta=4$, since $d^-(3)=5$ then $3\rightarrow x_4$, and if $N^+(x_4)\cap \{1,2,w\}\neq \emptyset$, say $x_4\rightarrow 1$, then $4x_1...x_4132wv4\equiv C$ so $\{1,2,w\}\rightarrow x_4$ and thus $T(N^+(1))=T(x_1,x_2,x_4,w,3)\neq T_2$ and the result follows. Now if $(T_2-w;a_{n-1}a_{n-2}...a_{\delta+2})=(5C;(2,1,1))$, if $N^+(w)\cap \{1,2,3\}\neq \emptyset$, say $w\rightarrow 1$, then $1wx_2...x_{\delta}v43x_1251\equiv C$ so $\{1,2,3\}\rightarrow w$. If $5\rightarrow x_2$ then $15x_2...a_{\delta}vw243x_11\equiv C$ so $x_2\rightarrow 5$. If $N^-(x_2)\cap \{1,2,3\}\neq \emptyset$, say $1\rightarrow x_2$, then $T_2-1\neq 5C$ and the result follows, so $x_2\rightarrow \{1,2,3\}$ and thus $4\rightarrow x_2$ and $w\rightarrow 5$. If $N^+(x_{\delta})\cap \{1,2,3\}\neq \emptyset$, say $x_{\delta}\rightarrow 1$, then $5x_1...x_{\delta}142w3v5\equiv C$, so $\{1,2,3\}\rightarrow x_{\delta}$ and thus $x_2\neq x_{\delta}$ i.e. $\delta\geq 3$. If $a_{\delta-1}\rightarrow a_{\delta}$ then $x_1wx_2...x_{\delta-1}v1x_{\delta}2354x_1\equiv C$, so $a_{\delta}\rightarrow a_{\delta-1}$. If $x_{\delta}\rightarrow 5$ then $x_1wx_2...a_{\delta}51v432x_1\equiv C$, so $5\rightarrow x_{\delta}$. If $(T_1;a_1...a_{\delta})=(4B;(1,2))$ then $x_1x_3x_4x_2\equiv a_1...a_4$ with $x_2\rightarrow 1$ so $T$ contains $C$. Otherwise, since $(T_1;a_1...a_{\delta})\neq (4A;(1,1,1))$ nor $(5C;(2,1,1))$ then there exists $y_1...y_{\delta}\equiv a_1...a_{\delta}$ such that $y_1\neq x_1$, and if $y_2$ is distinct from $x_1$ and $x_2$ then $d^-(5)\geq 7$, a contradiction. If $y_2=x_1$ then $1\rightarrow y_2$ and so $T$ contains $C$, so $y_2=x_2$. But $T(N^-(5))\neq \overline{T_2}$, and the result follows. Assume now that $u$ is minimal in $T_2-w$ and $a_{n-1}a_{n-2}...a_{\delta+2}$ is directed. If $|T_2|\geq 4$, let $a_s$ be the origin of $B$, then $T$ contains $\overline{C}=a'_0a'_1...a'_0=\overline{a_{s+\delta+2}a_{s+\delta+1}...a_{s+\delta+2}}$. Otherwise, if $\delta=2$ and $T_2$ is transitive, let $P=(1,1)$ be a path in $T_2$ of origin an inneighbor of $x_2$, we have $x_2x_1vP^{-1}x_2\equiv C$, if $\delta=2$ and $T_2=3A$, then $vx_1O_{T_2}x_2v\equiv C$, and if $\delta=3$ then $T_1=3A$ since otherwise $d^+(x_1)\geq 4$, and thus $T=7A$ where $12457631\equiv C$.
\item[(ii)] If $d^-_{T_2}(x_2)=1$ then we may assume that $u$ is minimal in $T_2$ and $a_{n-1}a_{n-2}...a_{\delta+2}$ is directed, this is possible by taking $x_2$ in the place of $v$. If $|T_2|\geq 4$, let $a_s$ be the origin of $B$, then $T$ contains $\overline{C}=a'_0a'_1...a'_0=\overline{a_{s+\delta+2}a_{s+\delta+1}...a_{s+\delta+2}}$. Otherwise, set $T_2-u=\{u_1\rightarrow u_2\}$, if $|T_1|=2$ then $u_2\rightarrow x_2$ and so $x_2\rightarrow u_1$ thus $x_1x_2u_1u_2vux_1\equiv C$, and if $|T_1|=3$ then $T_1$ is transitive and $T_2\rightarrow x_3$, we have $x_3u_2x_2x_1vu_1ux_3\equiv C$.
\end{itemize}
\textbf{Case $6$}: If $|B'|=1$ and $a_2\rightarrow a_1$. Suppose first that $(T_1;a_1...a_{\delta})$ is a Gr\"{u}nbaum's exception and for $x_i\in T_1$ let $y_0^i...x_{\delta-1}^i\equiv a_0...a_{\delta-1}$ be a path of origin $x_i$ in $T_1$. If an $x_i$ has an inneighbor origin of a path $P\equiv a_{n-1}a_{n-2}...a_{\delta+1}$ in $T_2$ then $y_0^i...y_{\delta-1}^ivP^{-1}y_0^i\equiv C$. Otherwise, $(T_2;a_{n-1}a_{n-2}...a_{\delta+1})$ is an exception, so by the building lemmas it is one of the biexceptions $(0,1),(3,1),(4,1)$, $(8,1),(9,1),(15,1),(16,1),(32,1),(33,1),(E_8(5),1)$\\ and $(E_8(6),1)$ with respect to $x_i$. For $Exc (0,1)$, $T=7A$ and thus $13465721\equiv C$, for $Exc (4,1)$ and $Exc (33,1)$, $(T_2;a_{n-1}a_{n-2}...a_{\delta+1})$ is not a biexception with respect to $y_1^i$. And for the remaining biexceptions, there is an outneighbor of $x_i$ origin of a path $P\equiv a_{\delta+2}a_{\delta+3}...a_0$ and therefore $vy_{\delta-1}^iy_{\delta-2}^i...y_0^iPv\equiv a_1a_2...a_1=C$.\\
Assume now that $(T_1;a_1...a_{\delta})$ is not a Gr\"{u}nbaum's exception and let $x_1...x_{\delta}\equiv a_1...a_{\delta}$ be a path in $T_1$ such that $d^+_{T_2}(x_1)=max\,\{d^+_{T_2}(x);x\in Or(a_1...a_{\delta};T_1)\}$. If $x_1$ has an outneighbor $u$ origin of a path $P\equiv a_0a_{n-1}...a_{\delta+2}$ in $T_2$ then $ux_1...x_{\delta+2}vP^{-1}\equiv C$. Otherwise, if $d^+_{T_2}(x_1)\geq 2$ and $d^+_{T_2}(x_2)=0$ then, by lemma \ref{2.8}, $(T_2;a_0a_{n-1}...a_{\delta+2})$ is one of the following exceptions:\\ Dual Exc $0,4,19,33$: If $T_1-x_2$ contains a path $y_2...y_{\delta}\equiv a_2...a_{\delta}$ then every vertex of $T_2$ is an origin of a path $P\equiv a_1a_0...a_{\delta+2}$ in $T_2\cup \{x_2\}$ of end distinct from $x_2$ and thus $y_2...y_{\delta}vP^{-1}y_2\equiv a_2a_3...a_2=C$; and if $(T_1;a_2...a_{\delta})=(3A;(1,1))$, let $a\in T_2$ and $y_1...y_4\equiv a_1...a_4$ be a path in $(T_1-x_2)\cup \{a\}$, $a\neq y_4$, every vertex of $T_2-a$ is an origin of a path $P\equiv a_0a_{n-1}...a_{\delta+2}$ in $(T_2-a)\cup \{x_2\}$, thus $y_1...y_4vP^{-1}y_1\equiv a_1a_2...a_1=C$.\\
Dual Exc $E_1(4)$: If $\delta=2$ then $x_14321vx_2x_1\equiv C$. If $\delta=3$ then $x_3x_2O_{T_2}vx_1x_3\equiv C$. If $\delta=4$, let $y_1...y_4=(2,1)$ be a path in $T_1$ such that $y_1\neq x_2$, $y_1$ has an outneighbor origin of a path $P\equiv \overline{a_5a_4a_3a_2}$ in $T_2$ and thus $y_4vP^{-1}y_1...y_4\equiv \overline{C}$.\\
Dual Exc $E_1(6)$: Every vertex of $T_2$ is an origin of a path $P\equiv \overline{a_2a_1...a_{n-3}}=-(2,1,2)$, so let $y_3...y_{\delta+2}\equiv \overline{a_3...a_{\delta+2}}$ be a path in $T_1$, $y_3\neq x_2$, thus $y_3...y_{\delta+2}vP^{-1}y_3\equiv \overline{a_3a_4...a_3}=\overline{C}$.\\
And if $d^+_{T_2}(x_1)\geq 2$ and $d^+_{T_2}(x_2)\geq 1$ then $(T_2;a_0a_{n-1}...a_{\delta+2})$ is an exception, so by the building lemmas we may assume that it is one of the following biexceptions with respect to $x_1$:\\ (When $|T_2|=5$ and $a_{\delta}\rightarrow a_{\delta-1}$ we can suppose that $\delta=2$. In fact, if $\delta=4$ or $5$ take $\overline{C}=a'_0a'_1...a'_0=\overline{a_0a_{n-1}...a_0}$ and it is solved in an another case, and if $\delta=3$ then $a_0a_1...a_3$ is a block greater than $B$.)\\
Dual Exc $(0,1)$: If $\delta=2$ then $x_2v123x_1x_2\equiv C$. If $\delta=3$ then $x_2\rightarrow 3$ and thus $x_3v123x_2x_1x_3\equiv C$.\\
Dual Exc $(2,1)$: If $\delta=2$ then $1vx_2x_14231\equiv C$. Else, we have $x_{\delta}\rightarrow 2$ so $4x_1x_2213v4\equiv C$.\\
Dual Exc $(4,1)$: By the building lemmas, if $N^+(x_{\delta})\cap \{1,2,5\}\neq \emptyset$ then $x_{\delta}$ has an outneighbor origin of a path $P\equiv a_{\delta+1}a_{\delta+2}...a_0$ in $T_2\cup \{x_1\}$ and thus $vx_2...x_{\delta}Pv\equiv a_1a_2...a_1=C$, so we may assume that $\{1,2,5\}\rightarrow x_{\delta}$. If $a_{\delta-1}\rightarrow a_{\delta}$ then $1x_1...x_{\delta-1}v435$\\$x_{\delta}21\equiv C$, so $a_{\delta}\rightarrow a_{\delta-1}$. Then $\delta=2$ and thus $421vx_2x_1534\equiv \overline{C}$.\\
Dual Exc $(8,1)$: By the building lemmas, if $N^+(x_{\delta})\cap \{2,3,4\}\neq \emptyset$ then $x_{\delta}$ has an outneighbor origin of a path $P\equiv a_{\delta+1}a_{\delta+2}...a_0$ in $T_2\cup \{x_1\}$ and thus $vx_2...x_{\delta}Pv\equiv a_1a_2...a_1=C$, so we may assume that $\{2,3,4\}\rightarrow x_{\delta}$. If $a_{\delta-1}\rightarrow a_{\delta}$ then $5x_1...x_{\delta-1}v213$\\$x_{\delta}45\equiv C$, so $a_{\delta}\rightarrow a_{\delta-1}$. Then $\delta=2$ and thus $253vx_2x_1142\equiv \overline{C}$.\\
Dual Exc $(9,1)(1),(3),(4)$: By the building lemmas, if $N^+(x_{\delta})\cap \{3,5\}$ then $x_{\delta}$ has an outneighbor origin of a path $P\equiv a_{\delta+1}a_{\delta+2}...a_0$ in $T_2\cup \{x_1\}$ and thus $vx_2...x_{\delta}Pv\equiv a_1a_2...a_1=C$, so we may suppose that $\{3,5\}\rightarrow x_{\delta}$. If $a_{\delta-1}\rightarrow a_{\delta}$ then $4x_1...x_{\delta-1}v5x_{\delta}31$\\$24\equiv C$, so $a_{\delta}\rightarrow a_{\delta-1}$. Then $\delta=2$ and thus $253vx_2x_1142\equiv \overline{C}$.\\
Dual Exc $(9,1)(2)$: By the building lemmas, if $N^+(x_{\delta})\cap \{3,4\}$ then $x_{\delta}$ has an outneighbor origin of a path $P\equiv a_{\delta+1}a_{\delta+2}...a_0$ in $T_2\cup \{x_1\}$ and thus $vx_2...x_{\delta}Pv\equiv a_1a_2...a_1=C$, so we may suppose that $\{3,4\}\rightarrow x_{\delta}$. If $a_{\delta-1}\rightarrow a_{\delta}$ then $5x_1...x_{\delta-1}v4x_{\delta}31$\\$25\equiv C$, so $a_{\delta}\rightarrow a_{\delta-1}$. Then $\delta=2$ and thus $253vx_2x_1142\equiv \overline{C}$.\\
Dual Exc $(15,1)$: If $d^+_{T_2}(x_{\delta})\geq 1$, one can find an outneighbor of $x_{\delta}$ origin of a path $P\equiv a_{\delta+1}a_{\delta+2}...a_1$ in $T_2\cup \{v\}$ of end and outneighbor $u$ of $x_1$ and thus $ux_1...x_{\delta}vP\equiv C$. So we may suppose that $d^+_{T_2}(x_{\delta})=0$ and so $a_{\delta}\rightarrow a_{\delta-1}$. If $\delta=5$ or $6$ take $\overline{C}=a'_0a'_1...a'_0=\overline{a_0a_{n-1}...a_0}$ and it is solved in an another case; if $\delta=4$ then $x_4v435261x_1x_2x_3x_4\equiv C$; if $\delta=3$ then $a_0a_1...a_3$ is a block greater than $B$ and if $\delta=2$ then $4vx_2x_1162534\equiv C$.\\
Dual Exc $(16,1)$: If $d^+_{T_2}(x_{\delta})\geq 1$, one can check that $x_{\delta}$ has an outneighbor origin of a path $P\equiv a_{\delta+1}a_{\delta+2}...a_0$ in $T_2\cup \{v\}$ of end $4$ and so $4x_1...x_{\delta}P\equiv C$, so $a_{\delta}\rightarrow a_{\delta-1}$. Let $u\in N^-_{T_2}(x_{\delta-1})$, $u$ is an origin of a path $P=(2,1)$ in $T_2\setminus \{w,t\}$ where $w\in N^+_{T_2}(x_1)$ and $t\in T_2$, thus $wx_1...x_{\delta-1}Px_{\delta}tvw\equiv C$.\\
Dual Exc $(32,1)$: If $d^+_{T_2}(x_{\delta})\geq 1$, one can check that $x_{\delta}$ has an outneighbor origin of a path $P\equiv a_{\delta+1}a_{\delta+2}...a_0$ in $T_2\cup \{v\}$ of end $5$ or $6$ and thus $T$ contains $C$. Otherwise, let $u\in N^-_{T_2}(x_{\delta-1})$, as above one can find $C$.\\
Dual Exc $(33,1)$: By the building lemmas, $N^+_{T_2}(x_1)=\{1,2,3,6\}$ and $\{1,2,4,5,6,7\}\subseteq N^-_{T_2}(x_{\delta})$ (since $a_{\delta+1}a_{\delta+1}...a_0\equiv a_0a_{n-1}...a_{\delta+1}$). For any $u\in T_2$, one can find a path $P\equiv a_{\delta}a_{\delta+1}...a_0$ in $T_2\cup \{x_{\delta},v\}$ ending by an outneighbor of $x_1$, so in both cases, $a_{\delta-1}\rightarrow a_{\delta}$ and $a_{\delta}\rightarrow a_{\delta-1}$, $x_{\delta-1}$ is an origin of a path $Q\equiv a_{\delta-1}a_{\delta}...a_0$ ending by an outneighbor of $x_1$ and thus $x_1x_2...x_{\delta-2}Px_1\equiv a_1a_2...a_1=C$ unless $T_2\rightarrow x_{\delta-1}$, in this case $T(N^-(x_{\delta})\neq 7A$.\\
Dual Exc $(E_8(5),1)$: If $x_{\delta}\rightarrow 1$ then $wx_1...x_{\delta}12v3uw\equiv C$, so we may assume that $\{1,2,3\}\rightarrow x_{\delta}$. If $a_{\delta-1}\rightarrow a_{\delta}$ then $wx_1...x_{\delta-1}v1x_{\delta}23uw\equiv C$, so $a_{\delta}\rightarrow a_{\delta-1}$. Then $\delta=2$ and thus $uw3vx_2x_112u\equiv \overline{C}$.\\
Else, by lemma \ref{2.4}, either $T_1=\overline{4B}$ or $\delta=2$. If $T_1=\overline{4B}$, set $T_1=\{u_1,u_2,u_3\}\rightarrow u_4$, we have $N^+_{T_2}(u_i)=\{u\}$ for all $i=1,2,3$ and $u\rightarrow u_4$. Suppose that $a_1...a_4=-(2,1)$, if $a_{n-2}\rightarrow a_{n-3}$, let $P\equiv a_3a_4...a_{n-3}$ be a path in $T_2$ of origin distinct from $u$, we have $u_4u_2u_1Pvu_3u_4\equiv C$; else, let $w\in N^+_{T_2}(u)$, $t\in T_2\setminus \{u,w\}$ and $P\equiv a_{n-3}a_{n-4}...a_{\delta+4}$ be a path in $T_2\setminus \{u,w,t\}$ (this is possible since $d^+_{T_2}(u)\geq 3$), we have $u_2u_1wuu_4vP^{-1}u_3tu_2\equiv C$. Assume now that $a_1...a_4=-(1,1,1)$, if $|T_2|=4$, let $T_2-u=\{y_1,y_2,y_3\}$ we have $y_3vu_1u_4uy_1u_3u_2y_2y_3\equiv C$; else, let $w\in N^+_{T_2}(u)$, $t\in T_2\setminus \{u,w\}$ and $P\equiv a_{n-3}a_{n-4}...a_{\delta+2}$ be a path in $T_2\setminus \{u,w,t\}$, we have $u_2wuu_4u_1vP^{-1}u_3tu_2\equiv C$. If $\delta=2$ and $N^+_{T_2}(x_1)=\{u\}$, if $|T_2|>3$, let $P\equiv \overline{a_{n-2}a_{n-3}...a_2}$ be a path in $T_2$ of origin distinct from $u$, we have $x_2vP^{-1}x_1x_2\equiv \overline{C}$. Otherwise, $|T_2|=3$, set $T_2-u=\{w\rightarrow t\}$. If $T_2=3A$, i.e. $t\rightarrow u \rightarrow w$, then $u\rightarrow x_2$, if $w\rightarrow x_2$ then $x_2ux_1vtwx_2\equiv C$ and if $x_2\rightarrow w$ then $x_1wx_2vutx_1\equiv C$. So we may assume that $T_2$ is a transitive triangle, also, for every vertex $x\in T$, $T(N^+(x))$ and $T^(N^-(x))$ are transitive, since, otherwise, $T$ or $\overline{T}$ contains $C$. Thus $\{w,t\}\rightarrow x_2$, if $u\rightarrow x_2$ then $x_2ux_1vtwx_2\equiv C$, and if $x_2\rightarrow u$ then $(T;C)=A_{14}$.    $\square$
\end{proof}
\begin{lemma}\label{3.6}
Let $T$ be a tournament of order $n$ with $\delta^-(T)\geq 2$ and $C$ a non-directed and non-antidirected cycle of order $n$. Let $B$ be a block of $C$ of maximal length and set $C=a_0a_1...a_{n-1}a_0$ such that $B$ is forward and it ends at $a_{\delta^-(T)+2}$. If $a_{n-1}\leftarrow a_0\rightarrow a_1$, then $T$ contains $C$ if and only if $(T;C)$ is not one of the exceptions $A_2$ and $A_{14}$.
\end{lemma}
\begin{proof}
By theorem \ref{3.1}, we may suppose that $\delta^-(T)\leq \delta^+(T)$. Set $\delta=\delta^-(T)$. Let $v\in T$ be a vertex such that $d^-(v)=\delta$ and denote by $T_1=T(N^-(v))$ and $T_2=T(N^+(v))$. Denote by $B'$ the block of $C$ containing the arc $a_0\rightarrow a_{n-1}$. We distinguish seven cases:\\
\textbf{Case $1$}: If $|B|=|B'|\geq 4$. Rotate $C$ until $a_{\delta}\rightarrow a_{\delta-1}$. $(T_1;a_1...a_{\delta})$ is not a Gr\"{u}nbaum's exception since its first block has the same direction of its last one. Let $x_1...x_{\delta}\equiv a_1...a_{\delta}$ in $T_1$ such that $d^+_{T_2}(x_1)=max\,\{d^+_{T_2}(x);x\in Or(a_1...a_{\delta};T_1)\}$. If there is an outneighbor $u$ of $x_1$ origin of $P\equiv a_0a_{n-1}...a_{\delta+2}$ in $T_2$, then $ux_1...x_{\delta}vP^{-1}\equiv C$. Otherwise, if $d^+_{T_2}(x_1)\geq 2$ then $(T_2;a_0a_{n-1}...a_{\delta+2})$ is an exception, and since $a_1a_0...a_{\delta+2}=(3,...,2)$ then it is not a biexception with respect to $x_1$, so by the building lemmas $T$ contains $C$. Else, $d^+_{T_2}(x_1)=1$, then by lemma \ref{2.4}, $(T_1;a_1...a_{\delta})=(\overline{4B};-(1,1,1))$ and $N^+_{T_2}(x)=\{u\}$ for all $x\in Or(a_1...a_{\delta};T_1)$ and so $d^+_{T_2}(u)\geq 3$, since for $x\in Or(a_1...a_{\delta};T_1)$ and $d^+_{T_2}(x)\geq 2$ we have $s^+_{T_2}(N^+_{T_2}(x))\geq 3$ and $(T_2;a_0a_{n-1}...a_{\delta+2})$ is not an exception. We have $a_0a_{n-1}...a_{\delta+3}=+(2,...,1)$ with an even number of blocks; so there is an outneighbor $w$ of $u$ in $T_2$ origin of a path $P\equiv a_0a_{n-1}...a_{\delta+3}$ and thus $wux_1x_4x_2x_3vP^{-1}\equiv C$.\\
\textbf{Case $2$}: If $|B|=|B'|=3$. Rotate $C$ until $a_{\delta}\rightarrow a_{\delta-1}$. Suppose first that $(T_1;a_1...a_{\delta})$ is a Gr\"{u}nbaum's exception. If $d^+_{T_2}(x)\geq 2$ for all $x\in T_1$, let $a\in T_2$ such that $d^+_{T_2}(a)=\Delta^+(T_2)$ and set $T'_1=T_1\cup \{a\}$ and $T'_2=T_2-a$. Let $x_0x_1...x_{\delta}\equiv a_0a_1...a_{\delta}$ in $T'_1$ such that $d^+_{T'_2}(x_0)\geq 2$, this is possible by lemma \ref{2.4}. If $T'_2$ contains a path $P\equiv a_{n-1}a_{n-2}...a_{\delta+2}$ of origin an outneighbor of $x_0$ then $x_0...x_{\delta}vP^{-1}x_0\equiv C$. Otherwise $(T'_2;a_{n-1}...a_{\delta+2})$ is an exception, so by the building lemmas we may assume that it is one of the following biexceptions with respect to $x_0$:\\
Exc $(0,0),(4,0)$ and $(33,0)$: $x_{\delta}$ has an outneighbor origin of a path $Q\equiv a_{\delta+1}a_{\delta+2}...a_{n-2}$ in $T_2$ so $x_0...x_{\delta}Qvx_0\equiv C$.\\
Exc $(9,0)$: $a\neq x_1$ so $x_{\delta}\leftarrow \{1,3,4\}$. In fact, if not, the path $x_1...x_{\delta}\equiv x_{\delta}...x_1$. If $x_{\delta}\rightarrow 2$ then $x_0...x_{\delta}23541vx_0\equiv C$. So we may assume that $N^-_{T_2}=\{1,2,3,4\}$, thus $x_0...x_{\delta}513v24x_0\equiv C$.\\
Exc $(15,0)$: $a\neq x_1$ so $\{2,4,5,6\}\subseteq N^-(x_{\delta})$. Then $x_{\delta}$ has an outneighbor in $\{1,3\}$ which is an origin of a path $Q\equiv a_{\delta+1}...a_{n-2}$, we have $x_0...x_{\delta}Qvx_0\equiv C$.\\
Exc $(16,0)$: $x_{\delta}\rightarrow 6$ so $x_0...x_{\delta}651423vx_0\equiv C$.\\
Suppose now that there exists $x\in T_1$ such that $d^+_{T_2}(x)=1$ so $T_1=3A$. Let $u_3$ be a such vertex and $u_1u_2u_3$ be a directed outpath of $T_1$. There exist $w,t\in T_2$ such that $w\rightarrow \{u_1,t\}$ and $u_3$ has an inneighbor origin of a path $Q\equiv a_7a_8...a_0$ in $T_2\setminus \{w,t\}$ and thus $vtwu_1u_2u_3Qv\equiv a_1a_2...a_2=C$ unless two case:\begin{itemize}
\item[(i)] $(T_2\setminus \{w,t\};a_7a_8...a_0)=(4B;(1,2))$, where $u_3\rightarrow \{u\}\rightarrow \{y_1,y_2,y_3\}$, in this case we can suppose that $\{u,y_1,y_2,y_3\}\rightarrow \{w,t\}$ since otherwise one can find $C$, so there exist $w',t'$ such that $w'\rightarrow \{t',u_3\}$ and $Q\equiv a_7a_8...a_0$ in $T_2\setminus \{w',t'\}$ with origin an inneighbor of $u_2$, thus $vt'w'u_3u_1u_2Qv\equiv C$.
\item[(ii)] $T_2\setminus \{w,t\}= \{y_1\rightarrow y_2\}$ where $u_3\rightarrow y_1$, if $w$ or $t\rightarrow y_1$ one can find $C$ so we can suppose that $d^+_{T_2}(y_1)=3$. Then if there exists $u_i\leftarrow y_1$ so $T$ contains $C$. Assume now that $T_1\rightarrow y_1$, if $T_2-y_1=3A$ then there is an outneighbor of a $u_i$ origin of a directed outpath in $T_2-y_1$ thus $u_{i+2}y_1u_{i+1}u_iO_{T_2-y_1}vu_{i+2}\equiv C$, else, there is an outneighbor of a $u_i$ origin of an antidirected outpath $R$ in $T_2-y_1$ thus $u_iy_1u_{i+2}u_{i+1}vR^{-1}u_i\equiv C$.
\end{itemize}
So we may suppose that $(T_1;a_1...a_{\delta})$ is not a Gr\"{u}nbaum's exception. Let $x_1...x_{\delta}\equiv a_1...a_{\delta}$ in $T_1$ such that $d^+_{T_2}(x_1)=max\,\{d^+_{T_2}(x);x\in Or(a_1...a_{\delta};T_1)\}$. If $d^+_{T_2}(x_1)\geq 2$ then $s^+_{T_2}(N^+_{T_2}(x_1))\geq 3$ and so $(T_2;a_0a_{n-1}...a_{\delta+2})$ is an exception, by the building lemmas we may suppose that it is one of the following biexceptions with respect to $x_1$: Exc $(7,0),(E_8(|T_2|),0),(E_9(|T_2|),0)$ or $(E_{10}(|T_2|),0)$. In these exceptions, we have respectively $T(N^+(x_2))\neq 5C$, $N^-(u)\neq \overline{F_8}, \overline{F_9}$ or $\overline{F_{10}}$, where $u\in X$; thus $T$ contains $C$. Without loss of generality, we can suppose that $d^+_{T_2}(x_1)\leq 1$ and $d^+_{T_2}(x_{\delta})\leq 1$, so assume that $N^+_{T_2}(x_1)=\{u\}$, if $d^+_{T_2}(u)\geq 2$ then $(T_2-u;a_{n-1}a_{n-2}...a_{\delta+2})$ is one of the exceptions: $0,4,8,9,15,16,32,33,37,E_8,E_9$ or $E_{10}$. It is easy to check that there is an inneighbor $w$ of $x_2$ origin of a path $P\equiv a_1a_0...a_{\delta+2}$ in $T_2\cup \{x_1\}$ and thus $x_2...x_{\delta}vP^{-1}x_2\equiv a_2a_3...a_2=C$. Otherwise, $d^+_{T_2}(u)=1$, take $u$ in the place of $v$ in $\overline{T}$ and set $T'_1=\overline{T}(N^-(u))$ and $T'_2=\overline{T}(N^+(u))$, $x_1$ and $v$ are origins of two paths $P\equiv a_0a_{n-1}...a_{\delta+2}$ and $P'\equiv a_{\delta+1}a_{\delta+2}...a_{n-1}$ in $T'_2$. Let $y_1...y_{\delta}\equiv a_1...a_{\delta}$ in $T'_1$, $y_1$ or $y_{\delta}$ is an inneighbor of $x$ or $v$, thus $y_1...y_{\delta}uPy_1\equiv a_1a_2...a_1=C$ or $uy_1...y_{\delta}P'u\equiv C$. Therefore $\overline{T}$ contains $C$ and so $T$ contains $C$.\\
\textbf{Case $3$}: If $|B|>|B'|\geq 3$. Rotate $C$ until $a_{n-2}\rightarrow a_{n-1}$ for the first time. Obviously we have $a_2\rightarrow a_1$. Assume first that $(T_1;a_1...a_{\delta})$ is a Gr\"{u}nbaum's exception. If there exists $x_3x_4...x_{\delta+2}\equiv a_3a_4...a_{\delta+2}$ such that $x_3$ has an inneighbor $u\in T_2$ origin of a path $P\equiv a_2a_1...a_{\delta+4}$ then $x_3x_4...x_{\delta+2}vP^{-1}x_3\equiv a_3a_4...a_3=C$. Otherwise, since there is no a such exception, we have $s^+_{T_2}(N^-_{T_2}(T_1))\leq 3$ so $s^+_{T_2}(N^+_{T_2}(T_1))=|T_2|$, therefore, for $y_2y_3...y_{\delta+1}\equiv a_2a_3...a_{\delta+1}$, there is an outneighbor $u\in T_2$ of $y_2$ origin of a path $P\equiv a_1a_0...a_{\delta+3}$ and thus $y_2...y_{\delta+1}vP^{-1}y_2\equiv a_2a_3...a_2=C$. So we can suppose that $(T_1;a_1...a_{\delta})$ is not a Gr\"{u}nbaum's exception. Let $x_1...x_{\delta}\equiv a_1...a_{\delta}$ in $T_1$ such that $d^+_{T_2}(x_1)=max\,\{d^+_{T_2}(x);x\in Or(a_1...a_{\delta};T_1)\}$. If $x_1$ has an outneighbor $u\in T_2$ origin of a path $P\equiv a_0a_{n-1}...a_{\delta+2}$ then $ux_1...x_{\delta}vP^{-1}\equiv C$. Otherwise, if $d^+_{T_2}(x_1)\geq 2$ and $d^+_{T_2}(x_2)=0$ then, since $a_0a_{n-1}...a_{\delta+2}$ has an even number of blocks with first block of length $1$ and last block of length $2$, by lemma \ref{2.8}, $(T_2;a_0a_{n-1}...a_{\delta+2})$ is one of the exceptions $19$ or $E_1(4)$. If $a_{\delta-1}\rightarrow a_{\delta}$, let $y_0...y_{\delta-1}\equiv a_0...a_{\delta-1}$ be a path in $T_1$ such that $d^+_{T_2}(y_0)\geq 2$, $T_2$ has at least $|T_2|-1$ origin of a path $P\equiv a_{n-1}a_{n-2}...a_{\delta+1}$, so there is an outneighbor of $y_0$ origin of a path $P$ in $T_2$ and thus $y_0...y_{\delta-1}vP^{-1}y_0\equiv C$. Otherwise, $a_{\delta}\rightarrow a_{\delta-1}$, then $T_1$ has at least $3$ origins of $a_1...a_{\delta}$, and so $Or(a_1...a_{\delta};T_1)\cup \{x_2\}\subseteq N^+(x)$ for all $x\in Or(a_0a_{n-1}...a_{\delta+2};T_2)$, then $d^+(x)>\Delta^+(T)$ for all $x\in Or(a_0a_{n-1}...a_{\delta+2};T_2)$, a contradiction. Then $d^+_{T_2}(x_2)\geq 1$ and $(T_2;a_0a_{n-1}...a_{\delta+2})$ is an exception, by the building lemmas we may assume that it is one of the following biexceptions with respect to $x_1$:\\
Exc $(2,0)$: If $x_{\delta}$ has an outneighbor in $T_2$ then it is an origin of a directed outpath, we have $vx_1...x_{\delta}O_{T_2}v\equiv C$. Else, $\delta=2$, we have $1vx_2x_13421\equiv C$.\\
Exc $(E_1(|T_2|),0)$: $\overline{T}(N^-(2))\neq F_1,F_3,F_5$, so $T$ contains $C$.\\
Exc $(E_3(|T_2|),0)$: $\overline{T}(N^-(1))\neq F_1,F_3,F_5$, so $T$ contains $C$.\\
Exc $(E_5(|T_2|),0)$: $\overline{T}(N^-(1))\neq F_1,F_3,F_5$, so $T$ contains $C$.\\
Exc $(E_{11}(|T_2|),0)$: $T(N^+(3))\neq F_{11},F_{13}$, so $T$ contains $C$.\\
Exc $(E_{13}(|T_2|),0)$: $\overline{T}(N^-(1))\neq F_{13}$, so $T$ contains $C$.\\
Assume now that $d^+_{T_2}(x_1)=1$, so by lemma \ref{2.3}, $(T_1;a_1...a_{\delta})=(\overline{4B};-(1,1,1))$, $(\overline{4B};-(2,1))$ or $(3A;-(2))$ or $|T_1|=2$. First, if $T_1=\overline{4B}$, we can assume that $N^+_{T_2}(x_1,x_2,x_4)=\{u\}$ since otherwise there is a vertex $x$ of minimal indegree with $\overline{T}(N^-(x))\neq \overline{4B}$. We have $d^+_{T_2}(u)\geq 3$, so if $a_1...a_{\delta}=-(1,1,1)$ then $(T_2-u,a_{n-1}a_{n-2}...$\\$a_{\delta+2})=(\overline{6D};-(1,2,2))$, we have $62ux_1x_4x_2x_3v43516\equiv C$. And if $a_1...a_{\delta}=-(2,1)$ then $(T_2-u;a_{n-1}a_{n-2}...a_{\delta+})=(\overline{6H};-(1,1,3),(\overline{7B};-(2,1,3))$ or $(\overline{7C};-(2,1,3))$. In this case $d^+_{T_2-u}(x_4)\geq 4$ so it has an outneighbor $w\in T_2-u$ origin of a path $P\equiv a_{n-1}a_{n-2}...a_{\delta+1}$ in $T_2$. We have $x_4x_1x_2x_3vP^{-1}x_4\equiv C$. Now if $T_1=3A$, we can suppose that $N^+_{T_2}(x_i)$ are pairwaisly different, which is a contradiction. So assume that $|T_1|=2$, let $u\in N^+_{T_2}(x_2)$ and $P\equiv a_{n-1}a_{n-2}...a_3$ in $T_2-u$ of origin and inneighbor of $x_1$, we have $x_2uvP^{-1}x_1x_2\equiv C$.\\
\textbf{Case $4$}: If $|B'|=2$ and $a_1\rightarrow a_2$. Without loss of generality, we can suppose that $(T_1;a_1...a_{\delta})$ is not a Gr\"{u}nbaum's exception. In fact, if it is, set $C=a'_0a'_1...a'_{n-1}a'_0=a_{\delta}a_{\delta-1}...a_{\delta+1}a_{\delta}$. So let $x_1...x_{\delta}\equiv a_1...a_{\delta}$ in $T_1$ such that $d^-_{T_2}(x_1)=max\,\{d^-_{T_2}(x);x\in Or(a_1...a_{\delta},T_1)\}$. We can suppose that if $s^+_{T_2}(N^-_{T_2}(x_1))=2$ then $T(N^+_{T_2}(x_1)\neq 3A,5A,7A$ since otherwise, set $N^-_{T_2}(x_1)=\{u\rightarrow w\}$, in $\overline{T}$, taking $w$ in the place of $v$, $\overline{T}(N^+(w))$ has $v$ as a maximal vertex and $u$ as a minimal vertex. So let $Q\equiv a_{n-2}a_{n-3}...a_{\delta+2}$ be a path in $T(N^+_{T_2}(x_1))$, we have $uwx_2...x_{\delta}Q^{-1}x_1u\equiv C$. Otherwise, if $x_1$ has an inneighbor $u$ origin of a path $P\equiv a_0a_{n-1}...a_{\delta+2}$ in $T_2$ then $ux_1...x_{\delta}vP^{-1}\equiv C$. Else, $(T_2;a_0a_{n-1}...a_{\delta+2})$ is an exception, so by the building lemmas we may assume that it is one of the following biexceptions with respect to $x_1$:\\
Exc $(5,1)$: $T(N^+(x_2))\neq 5B,5C$.\\
Exc $(7,1)$: We have $\{1,2\}\subseteq N^-(x_1)$, if $d^+_{T_2}(x_{\delta})\geq 1$ one can find $Q\equiv a_{\delta+1}a_{\delta+2}...a_0$ in $T_2\cup \{v\}$ of origin an outneighbor of $x_{\delta}$ and end an inneighbor $u$ of $x_1$, thus $ux_1...x_{\delta}Q\equiv C$. So we may assume that $T_2\rightarrow x_{\delta}$ and $a_{\delta}\rightarrow a_{\delta-1}$. If $a_1...a_{\delta}$ is antidirected then $5x_{\delta}...x_1v21435\equiv C$. So $\delta=4$ and $a_1...a_4=(1,2)$, we have $x_4v32145x_1x_2x_3x_4\equiv C$.\\
\textbf{Case $5$}: If $|B'|=2$ and $a_2\rightarrow a_1$. Suppose first that $(T_1;a_1...a_{\delta})$ is a Gr\"{u}nbaum's exception. If $d^+_{T_2}(x)\geq 2$ for all $x\in T_1$, let $a\in T_2$ such that $d^+_{T_2}(a)=\Delta^+(T_2)$ . Set $T'_1=T_1\cup \{a\}$ and $T'_2=T_2-a$ and let $x_0...x_{\delta}\equiv a_0...a_{\delta}$ be a path in $T'_1$ such that $d^+_{T'_2}(x_{\delta})=max\,\{d^+_{T'_2}(x);x\in Or(a_{\delta}...a_0,T'_1)\}$, if $x_{\delta}$ has an outneighbor in $T'_2$ origin of a path $P\equiv a_{\delta+1}a_{\delta+2}...a_{n-2}$ then $x_0...x_{\delta}Pvx_0\equiv C$. Otherwise, by lemma \ref{2.4}, $d^+_{T_2}(x_{\delta})\geq 2$ then $(T'_2;a_{\delta+1}...a_{n-2})$ is an exception. By the building lemmas we may suppose that it is one of the following biexceptions with respect to $x_{\delta}$:\\
Exc $(1,0)$: If $\delta=5$, there exists $y_0y_{n-1}...y_{\delta+2}\equiv \overline{a_0a_{n-1}...a_{\delta+2}}$ such that $y_0$ has an inneighbor origin of a path $Q\equiv \overline{a_1a_2...a_{\delta}}$ in $T_2$, we have $y_0Qvy_{\delta+2}...y_0\equiv \overline{C}$, thus $T$ contains $C$. Assume now that $\delta=3$, let $u_1u_2u_3$ be a directed outpath of $T_1$ such that $u_1\rightarrow a$. If $u_3\leftarrow 1$ or $2$, say $1$, then $u_1av4231u_3u_2u_1\equiv C$. So assume that $\{1,2\}\rightarrow u_3$, if $d^-_{T_2}(a)\geq 1$, say $a\rightarrow 1$, then $1au_1u_2u_32v431\equiv C$. So we may suppose that $d^+_{T_2}(a)=4$, if $u_1$ has an outneighbor in $T_2-a$, say $1$, we have $v1u_1u_2u_32a43v\equiv C$. So $T_2-a\rightarrow u_1$, we have $v21u_1u_2u_334av\equiv C$.\\
Exc $(27,0)$: Taking the block $a_{\delta+3}a_{\delta+4}a_{\delta+5}a_{\delta+6}$ of maximal length, $T$ contains $C$.\\
Exc $(48,0),(E'_8,0),(E'_9,0)$ and $(E'_{10},0)$: There is an outneighbor of $x_{\delta}$ of outdegree greater than $\Delta^+(T)$ which is a contradiction.\\
Suppose now that there exists $x\in T_1$ such that $d^+_{T_2}(x)= 1$, let $u_1$ a such vertex and $u_1u_2u_3$ be a directed outpath of $T_1$. There exist $u\in N^+_{T_2}(u_1)$ and an inneighbor of $u_3$ origin of a path $Q\equiv a_6a_7...a_{n-1}$ in $T_2-u$ and thus $vuu_1u_2u_3Qv\equiv C$, unless $T$ is the Paley tournament which contains $C$ or $(T_2-u;a_6a_7...a_{n-1})=(4B;(1,2))$, set $T_2-u=y_1\rightarrow \{y_2,y_3,y_4\}$, by supposing that $u_i\rightarrow y_j$ $j=2,3$ or $4$ one can find $C$ so we can assume that $\{y_2,y_3,y_4\}\rightarrow T_1$ and $u\rightarrow \{y_2,y_3,y_4\}$; we have $vy_3y_2u_1u_2u_3uy_4y_1v\equiv C$.\\
Assume now that $(T_1;a_1...a_{\delta})$ is not a Gr\"{u}nbaum's exception and let $x_1...x_{\delta}\equiv a_1...a_{\delta}$ be a path in $T_1$ such that $d^-_{T_2}(x_1)=max\,\{d^-_{T_2}(x);x\in Or(a_1...a_{\delta},T_1)\}$. If $x_1$ has an inneighbor $u$ origin of a path $P\equiv a_0a_{n-1}...a_{\delta}$ in $T_2$ then $ux_1...x_{\delta}vP^{-1}\equiv C$. Otherwise, suppose that $d^-_{T_2}(x_1)\geq 2$, if $s^+_{T_2}(N^-_{T_2}(x_1))=2$ set $N^-_{T_2}(x_1)=\{u\rightarrow w\}$, taking $w$ in $\overline{T}$ in the place of $v$, then $\overline{T}(N^+(w))$ has $u$ as a minimal vertex and $v$ as a maximal vertex. So we may suppose that $T_2$ has a minimal vertex $x$. If $x_2\rightarrow x$, if $|T_2|\geq 4$ let $Q\equiv a_{n-1}a_{n-2}...a_{\delta+2}$ in $T_2-x$ of origin distinct from $u$ and $w$, we have $x_1xx_2...x_{\delta+2}vQ^{-1}x_1\equiv C$ and if $|T_2|=3$, we may assume that $|T|=6$, since, otherwise, an another origin of a path $y_1y_2y_3\equiv a_1a_2a_3$ in $T_1$ is an outneighbor of $x$ which is an origin of a path $R\equiv a_0a_6a_5$ in $T_2$, and so $xy_1y_2y_3vR^{-1}\equiv C$; thus $(T;C)=A_{14}$. So $x\rightarrow x_2$, taking $x$ in the place of $w$, the maximal component has at least $3$ vertices, $u,w$ and $x_2$. So we can assume that $s^+_{T_2}(N^-_{T_2}(x_1))\geq 3$ and so if $d^+_{T_2}(x_2)=0$ then, by lemma \ref{2.8} $(T_2;a_0a_{n-1}...a_{\delta+2})$ is the exception $7$, if $T_1-x_2$ contains a path $y_2...y_{\delta}\equiv a_2...a_{\delta}$ then every vertex of $T_2$ is an origin of a path $P\equiv a_1a_0...a_{\delta+2}$ in $T_2\cup \{x_2\}$ and so $y_2...y_{\delta}vP^{-1}y_2\equiv a_2a_3...a_2=C$, otherwise, without loss of generality we may assume that $x_4\rightarrow 5$ and $1\rightarrow x_1$, thus $31x_1...x_{\delta}52v23\equiv C$. And if $d^+_{T_2}(x_2)\geq 1$ then $(T_2;a_0a_{n-1}...a_{\delta+2})$ is an exception. By the building lemmas we may suppose that it is one of the following biexceptions with respect to $x_1$:\\
Exc $(5,0),(7,0)$: If $d^+_{T_2}(x_{\delta})\geq 1$ one can find $Q\equiv a_{\delta+1}a_{\delta+2}...a_0$ in $T_2\cup \{v\}$ of origin an outneighbor of $x_{\delta}$ and end an inneighbor $u$ of $x_1$, thus $ux_1...x_{\delta}Q\equiv C$. So we may assume that $T_2\rightarrow x_{\delta}$ and $|T_1|=2$ since otherwise $T(N^+(1))=T(\{x_1,x_2,x_{\delta},2,4\})\neq 5B$. Thus $v25x_2x_1143v\equiv C$.\\
So we may assume that $d^-_{T_2}(x_1)=1$, set $N^-_{T_2}(x_1)=\{u\}$, and, by lemma \ref{2.3}, either $(T_1;a_1...a_{\delta})=(\overline{4B};-(1,2))$ or $|T_1|=2$. In the first case we can suppose that $(T_1-x_1)\rightarrow u$ since otherwise $T(N^-(x_1))\neq \overline{4B}$, then $d^+_{T_2}(u)\geq 3$ and so $(T_2;a_{n-1}a_{n-2}...a_{\delta+2})$ is one of the exceptions $18,24,26,34,38,47$, and in these ones, $u$ has an outneighbor origin of a path $Q\equiv a_{n-2}a_{n-3}...a_{\delta+1}$, thus $x_2x_1x_3x_4vQ^{-1}x_2\equiv C$. Now, if $|T_1|=2$. If $x_2\rightarrow u$ then, taking $x_1$ in the place of $v$ we can suppose that $d^-_{T_2}(u)=0$ which is a contradiction, so assume that $u\rightarrow x_2$. Also taking $x_1$ in the place of $v$ we can suppose that $d^-_{T_2}(u)=1$, thus either $(T_2-u;a_{n-1}a_{n-2}...a_{\delta+2})=(4A;(1,1,1))$ and this is a contradiction since $d^+(3)=1$ or $T_2=\{u,w,t\}$ which is similar to the above case ($d^-_{T_2}(x_1)=2$).\\
\textbf{Case $6$}: If $|B'|=1$ and $a_1\rightarrow a_2$. Suppose first that $(T_1;a_1...a_{\delta})$ is a Gr\"{u}nbaum's exception, if $a_{\delta+4}\rightarrow a_{\delta+3}$, $T$ contains $\overline{C}=u_0u_1...u_{n-1}u_0=\overline{a_2a_3...a_2}$ by case $5$, so $a_{\delta+3}\rightarrow a_{\delta+4}$. Let $a\in T_2$ such that $d^+_{T_1}(a)< |T_1|$ and $d^-_{T_2}(a)=max\,\{d^-_{T_2}(x);x\in T_2, d^+_{T_1}(x)< |T_1|\}$, set $T'_1=T_1\cup \{a\}$ and $T'_2=T_2-a$. Since $T'_1\neq 4B$, there exists $x_0...x_{\delta}\equiv a_0...a_{\delta}$ such that $a\neq x_{\delta}$, if $x_0$ has an outneighbor origin of a path $P\equiv a_{n-1}a_{n-2}...a_{\delta+2}$ in $T'_2$ then $x_0...x_{\delta}vP^{-1}x_0\equiv C$. Otherwise, if $d^+_{T'_2}(x_0)\geq 2$ then $(T'_2;a_{n-1}a_{n-2}...a_{\delta+2})$ is an exception, by the building lemmas we may assume that it is one of the following biexceptions with respect to $x_0$:\\
Dual Exc $(0,1)$: Since $d^-_{T_2}(x_0)\geq 2$ then $a=x_2$ and thus $x_3\rightarrow 3$, therefore $x_0...x_332v1x_0\equiv C$.\\
Dual Exc $(3,1)$: If $x_{\delta}\rightarrow 2$ then $x_0...x_{\delta}24v31x_0\equiv C$ so $2\rightarrow x_{\delta}$ and then if $x_1\neq a$ $x_{\delta}...x_1\equiv x_1...a_{\delta}$, so $T$ contains $C$, thus $x_1=a$. One can check that $x_0\rightarrow x_{\delta}$ and so $x_0x_{\delta}x...x_1\equiv a_0a_1...a_{\delta}$, but $x_1\rightarrow 2$, so $T$ contains $C$.\\
Dual Exc $(4,1)$: If $x_1\neq a$ then $N^+_{T'_2}(x_1)\subseteq N^+_{T'_2}(x_{\delta})$, thus $x_0...x_{\delta}54v132x_0\equiv C$. Otherwise, $x_0x_{\delta}...x_154v132x_0\equiv C$.\\
Dual Exc $(8,1)$: If $x_1\neq a$ then $N^+_{T'_2}(x_1)\subseteq N^+_{T'_2}(x_{\delta})$, thus $x_0...x_{\delta}21v345x_0\equiv C$. Otherwise, $x_0x_{\delta}...x_121v345x_0\equiv C$.\\
Dual Exc $(9,1)$: If $x_{\delta}\rightarrow 3$ then $x_0$ has an outneighbor origin of an antidirected inpath $P$ in $T(2,4,5)$ and thus $x_0...x_{\delta}31vP^{-1}x_0\equiv C$, so $3\rightarrow x_{\delta}$ and $x_1=a$, therefore $x_0x_{\delta}...x_131vP^{-1}x_0\equiv C$.\\
Dual Exc $(15,1)$: If $x_1\neq a$ then $N^+_{T'_2}(x_1)\subseteq N^+_{T'_2}(x_{\delta})$, thus $x_0...x_{\delta}54v3261x_0\equiv C$. Otherwise, $x_0x_{\delta}...x_154v3261x_0\equiv C$.\\
Dual Exc $(16,1)$: If $x_1\neq a$ then $N^+_{T'_2}(x_1)\subseteq N^+_{T'_2}(x_{\delta})$, thus if $x_0\rightarrow 5$, $x_0...x_{\delta}21v6435x_0$\\$\equiv C$ and if $x_0\rightarrow 4$ $x_0...x_{\delta}21v6534x_0\equiv C$. Otherwise, if $x_0\rightarrow 5$, $x_0x_{\delta}...x_121v6435x_0\equiv C$ and if $x_0\rightarrow 4$ $x_0x_{\delta}...x_121v6534x_0\equiv C$\\
Dual Exc $(32,1)$: If $x_1\neq a$ then $N^+_{T'_2}(x_1)\subseteq N^+_{T'_2}(x_{\delta})$, thus $x_0...x_{\delta}21v6435x_0\equiv C$. Otherwise, $x_0x_{\delta}...x_121v6435x_0\equiv C$.\\
Dual Exc $(33,1)$: If $x_1\neq a$ then $N^+_{T'_2}(x_1)\subseteq N^+_{T'_2}(x_{\delta})$, thus $x_0...x_{\delta}43v67521x_0\equiv C$. Otherwise, $x_0x_{\delta}...x_143v67521x_0\equiv C$.\\
Dual Exc $(E_8(5),1)$: If $x_1\neq a$ then $N^+_{T'_2}(x_1)\subseteq N^+_{T'_2}(x_{\delta})$, thus $x_0...x_{\delta}21v3wux_0\equiv C$. Otherwise, $x_0x_{\delta}...x_121v3wux_0\equiv C$.\\
Dual Exc $(E_8(6),1)$: If $x_{\delta}\rightarrow 1$ then $x_0...x_{\delta}13v2O_Xx_0\equiv C$ so $1\rightarrow x_{\delta}$ and then if $x_1\neq a$ $x_{\delta}...x_1\equiv x_1...a_{\delta}$, so $T$ contains $C$, thus $x_1=a$. One can check that $x_0\rightarrow x_{\delta}$ and so $x_0x_{\delta}...x_1\equiv a_0a_1...a_{\delta}$, but $x_1\rightarrow 1$, so $T$ contains $C$.\\
Otherwise, $d^+_{T_2}(x_0)\leq 1$, if there exists $w\in N^-_{T'_2}(x_1)\cap N^-_{T'_2}(x_0)$ then if $x_0$ has an inneighbor origin of a path $P\equiv a_{n-2}a_{n-3}...a_5$ in $T'_2-w$ then $wx_1x_2x_3vP^{-1}x_0w\equiv C$ and else $(T'-w;a_{n-3}a_{n-4}...a_5)$ is one of the following exceptions:\\
Exc $0$: $\{1,2\}\rightarrow x_0$, if $w\rightarrow 3$ then $wx_1x_2x_3v21x_03w\equiv C$, so we can suppose that $3\rightarrow w$. If $x_3\rightarrow 2$ then $wx_1x_2x_323v1x_0w\equiv C$ and if $x_3\rightarrow 3$ then $wx_1x_2x_331v2x_0w\equiv C$, so $\{2,3\}\rightarrow x_3$ and thus $x_1=a$. If $1\rightarrow x_2$ then $vwax_0x_3x_2123v\equiv C$, so $x_2\rightarrow \{1,2,3\}$ and thus $v12x_2x_0x_33wav\equiv C$.\\
Exc $1$: $N^+_{T'_2-w}(x_0)=\{4\}$, if $w\rightarrow 2$ then $wx_1...x_{\delta}v43x_012w\equiv C$ and if $w\rightarrow 4$ then $wx_1...x_{\delta}v21x_034w\equiv C$, so we can suppose that $\{2,4\}\rightarrow w$ since otherwise one can find $C$. If $d^+_{T'_2-w}(x_{delta})\geq 1$, set $V(T'_2-w)=\{x,y,x',y'\}$ such that $x\rightarrow \{x_0,y\}$, $x'\rightarrow y'$ and $x_{\delta}\rightarrow x'$, then $wx_1...x_{\delta}x'y'vyxx_0w\equiv C$. Otherwise, if $3\rightarrow x_1$ then $3x_1...x_{\delta}v14w2x_03\equiv C$ so $x_1\rightarrow 3$ and $x_1=a$ and thus $1x_{\delta}...x_134vw2x_01\equiv C$.\\
Exc $7$:  $N^+_{T'_2-w}(x_0)=\{5\}$, in a similar way we can suppose that $\{1,2,3\}\rightarrow w$. If $d^+_{T'_2-w}(x_{\delta})\geq 1$, set $V(T'_2-w)=\{x,y,z,x',y'\}$ such that $xyz$ and $x'y'$ are two directed outpaths, $x\rightarrow x_0$ and $x_{\delta}\rightarrow x'$, then $wx_1...x_{\delta}x'y'vzyxx_0w\equiv C$. Otherwise, if there exists $w'\in N^-(x_1)\cap \{1,2,3\}$, say $w'=3$, then $w'x_1...x_{\delta}vw1254x_0w'\equiv C$, so $x_1\rightarrow \{1,2,3\}$ and $x_1=a$ and thus $x_0x_{\delta}...x_11w243v5x_0\equiv C$.\\
Exc $4$ and $33$: In a similar way we can suppose that $(T'_2-w)\rightarrow w$. If $N^+(x_{\delta})\cap (T'_2-w)\neq \emptyset$ then there exist $x\in N^+_{T'_2-w}(x_{\delta})$, $y\in N^+_{T'_2-w}(x)$ and $P\equiv a_{n-3}a_{n-4}...a_{\delta+4}$ in $T'_2\setminus \{w,x,y\}$ of origin an inneighbor of $x_0$, thus $wx_1...x_{\delta}xyvP^{-1}x_0w\equiv C$. So we may assume that $(T'_2-w)\rightarrow x_{\delta}$ and so $x_1=a$. Let $x'\in N^+_{T'_2}(x_1)$, $w'\in N^-_{T'_2}(x_0)\cap N^-_{T'_2}(x_{\delta})$ and $P\equiv a_{n-3}a_{n-4}...a_{\delta+4}$ be a path in $T'_2\setminus \{w,w',x'\}$ of origin an inneighbor of $x_0$, we have $w'x_{\delta}...x_1x'wvP^{-1}x_0w'\equiv C$.\\
Otherwise, $d^-_{T'_2}(x_1)\leq 1$, then $x_1=a$, in fact, suppose that $x_1\neq a$, if $|T_1|=3$ then $x_1\in T_1$ with $d^-_{T_2}(x_1)\leq 1$, a contradiction, and if $|T_1|=5$ then $x_1\in T_1$ with $d^-_{T_2}(x_1)\leq 2$, a contradiction. Suppose first that $T_1=3A$, if $x_3\rightarrow x_1$ then $y_0...y_4=x_2x_0x_1x_3\equiv a_0...a_4$ with $d^-_{T'_2}(y_1)\geq 2$ then $T$ contains $C$, so $x_1\rightarrow x_3$ and thus $d^-_{T'_2}(x_1)=1$, set $N^-_{T'_2}(x_1)=\{u\}$. We also have $N^+_{T'_2}(x_0)=\{u\}$. If $a_{n-2}a_{n-3}a_{n-4}$ is directed then by $T$ contains $C=u_0u_1...u_0=a_{n-2}a_{n-1}...a_{n-2}$ by lemmas \ref{3.3} and \ref{3.4}, so $a_{n-2}a_{n-3}a_{n-4}$ is not directed. If $a_{n-2}a_{n-3}a_{n-4}=(1,1)$, since the number of blocks of $a_4a_5...a_{n-3}$ is odd then $T'_2-u$ contains a path $P\equiv a_4a_5...a_{n-3}$ and so $x_0x_3x_2x_1Pvux_0\equiv C$, and if $a_{n-2}a_{n-3}a_{n-4}=-(1,1)$, if $C=(2,1,2,1,1,1,2,1)$, by taking $C=u_0u_1...u_0\equiv a_5a_6...a_5$ it is solved in lemma \ref{3.5}, otherwise, let $u'\in N^-_{T'_2}(x_3)$ and $P\equiv a_3a_4...a_{n-4}$ be a path in $T'_2-u'$ of origin in $N^-(x_0)$, we have $x_1x_2x_0Pvu'x_3x_1\equiv C$. Assume now that $T_1=5A$, if $a_{n-3}\rightarrow a_{n-2}$ then $T$ contains $\overline{C}=u_0u_1...u_0\equiv \overline{a_4a_5...a_4}$ by lemma \ref{3.5}, so we may suppose that $a_{n-2}\rightarrow a_{n-3}$. Let $y_0...y_4\equiv a_0...a_4$ be a path in $T_1$ with end $x_0$ and let $w\in N^+_{T_2}(y_0)$, since $d^-_{T_2}(x_0)\geq |T_2-w|-2$ then by theorem \ref{2.1} $x_0$ has an inneighbor origin of a path $P\equiv a_5a_6...a_{n-3}$ in $T_2-w$ and so $y_0...y_4Pvwy_0\equiv C$, unless $(T_2-w;a_5...a_{n-3})=(5C;(2,1,1))$, this case is solved in lemma \ref{3.5} by taking $C=u_0u_1...u_{n-1}u_0\equiv a_8a_9...a_8$.

Assume now that $(T_1;a_1...a_{\delta})$ is not a Gr\"{u}nbaum's exception and let $x_1...x_{\delta}\equiv a_1...a_{\delta}$ be a path in $T_1$ such that $d^-_{T_2}(x_1)=max\,\{d^-_{T_2}(x);x\in Or(a_1...a_{\delta},T_1)\}$, if $x_1$ has an inneighbor $u$ origin of a path $P\equiv a_0a_{n-1}...a_{\delta+2}$ then $ux_1...x_{\delta}vP^{-1}\equiv C$. Otherwise, $(T_2;a_0a_{n-1}...a_{\delta+2})$ is one of the following biexceptions with respect to $x_1$:\\
Exc $(1,1)$: If $d^-_{T_2}(x_1)=2$ then, by lemma \ref{2.4}, $(T_1;a_1...a_{\delta})=(4B;(2,1)),(4B;(1,1,$\\$1)),(3A;(2))$ or $|T_1|=2$. If it is $(4B;(2,1))$ then the block $a_0...a_3$ is greater than $B$, a contradiction. If it is $(4B;(1,1,1))$ or $(3A;(2))$ then any vertex of $\{x_1,x_2,x_{\delta}\}$ may play the role of $x_1$, so if $N^-_{T_2}(x_1)=\{1,3\}$ then $d^-(x_2)<\delta$, a contradiction and if $N^-_{T_2}(x_1)=\{1,2\}$ then $N^-_{T_2}(x_2)=\{1,3\}$, a contradiction and if $N^-_{T_2}(x_1)=N^-_{T_2}(x_2)=\{2,3\}$ then $d^-(x_1)=\delta$ with $T(N^-(x_1))\neq 3A$ nor $4B$, the result follows. So we may assume that $\delta=2$. If $x_2\rightarrow 1$ then $x_1x_2123v4x_1\equiv C$, else, $x_1v4231x_2x_1\equiv C$. So we may assume that $N^+_{T_2}(x_1)=\{1,2,3\}$ and $\{2,4\}\subseteq N^+(x_2)$. Since $d^-(2)=4$ then $2\rightarrow x_{\delta}$. By supposing that $d^+_{T_2}(x_{\delta})\geq 1$ one can find $C$ so we may assume that $T_2\rightarrow x_{\delta}$ and so $a_{\delta}\rightarrow a_{\delta-1}$ and $d^+_{T_1}(x_{\delta}=\delta-1$. $x_{\delta}$ is an origin of a path $Q\equiv a_1...a_{\delta}$ so $4Qv3124\equiv C$.\\
Exc $(6,1)$: Since $d^-_{T_2}(x_1)=2$ then $\delta=2$ and thus $5x_1x_2142v35\equiv C$.\\
Exc $(10,1)$: $d^-(2)\geq 6$ a contradiction.\\
Exc $(12,1)$: Since $d^-_{T_2}(x_1)=2$ then $\delta=2$ and thus $1x_1x_2354v21\equiv C$.\\
Exc $(22,1)$: If $x_{\delta}\rightarrow 4$ then $1x_1...x_{\delta}46v3521\equiv C$, so we may assume that $\{4,5,6\}\rightarrow x_{\delta}$. If $a_{\delta-1}\rightarrow a_{\delta}$ then $2x_1...x_{\delta-1}v5x_{\delta}41632\equiv C$, so $a_{\delta}\rightarrow a_{\delta-1}$. If $1\rightarrow x_{\delta-1}$ then $2x_1...x_{\delta-1}15x_{\delta}46v32\equiv C$, so $x_{\delta-1}\rightarrow \{1,2,3\}$. Then $\{4,5,6\}\cap N^-(x_{\delta-1})\neq \emptyset$, without loss of generality we may suppose that $4\rightarrow x_{\delta-1}$ then $3x_1...x_{\delta-1}46x_{\delta}52v13\equiv C$.\\
Exc $(24,1)$: By supposing that $d^+_{T_2}(x_{\delta})\geq 1$ one can find $C$, otherwise, $x_{\delta}$ is an origin of a path $Q\equiv a_1...a_{\delta}$ in $T_1$ and so $6Qv524316\equiv C$.\\
Exc $(26,1)$: $T(N^-(4))\neq \overline{6H}$ so the result follows.\\
\textbf{Case $7$}: If $|B'|=1$ and $a_2\rightarrow a_1$. Suppose first that $(T_1;a_1...a_{\delta})$ is a Gr\"{u}nbaum's exception, if $a_{\delta+4}\rightarrow a_{\delta+3}$ (in this case $|T_2|\geq 5$), for $x_i\in T_1$ let $P_i=y_3^iy_4^i...y_{\delta+2}^i\equiv \overline{a_3a_4...a_{\delta+2}}$ be a path of origin $x_i$ in $T_1$. If there exists $x_i\in T_1$ having an outneighbor origin of a path $P\equiv \overline{a_2a_1...a_{\delta+4}}$ in $T_2$ then $P_ivP^{-1}x_i\equiv \overline{a_3a_4...a_4}=\overline{C}$. Otherwise, if $\delta=5$ or $7$, then $d^-_{T_2}(x_i)\geq 2$ and $b_i(P_i)=1$ and so $(T_2;a_2a_1...a_{\delta+4})$ is an exception. By the building lemmas we may suppose that it is a biexception with respect to $x_i$; the sole biexception verifying that $\overline{a_3a_2...a_{\delta+4}}=+(1,1,1,1,Q)$ and $T_2-S\subseteq N^+(x_i)$ is Dual Exc $(22,1)$: $d^-(1)=8$, a contradiction. So we may suppose that $\delta=3$, if there exists $x_i\in T_1$ such that $d^+_{T_2}(x_i)\geq 2$ then $(T_2;\overline{a_2a_1...a_{\delta+4}})$ is an exception, by the building lemmas we may suppose that it is one of the following biexceptions with respect to an $x_i$:\\
Dual Exc $(4,1),(8,1),(9,1),(24,1),(33,1)$: $(T_2;\overline{a_6a_7...a_3})$ is not a biexception with respect to $y_4^i$, thus $y_5^i$ is an origin of a path $Q\equiv \overline{a_5a_6...a_3}$ in $T_2\cup \{y_4^i,y_5^i\}$ with end in $T_2$, thus $vy_3^iQv\equiv \overline{a_4a_5...a_4}=\overline{C}$.\\
Dual Exc $(22,1),(E_8(5),1)$: $(T_2;\overline{a_8a_9...a_3})$ is not a biexception with respect to $x_i$, thus $y_4^i$ is an origin of a path $Q\equiv \overline{a_6a_7...a_3}$ in $T_2\cup \{y_4^i,x_i\}$ with end in $T_2$, thus $vy_4^iQv\equiv \overline{a_4a_5...a_4}=\overline{C}$.\\
Otherwise, $d^+_{T_2}(x_i)=1$ for all $x_i\in T_1$, if $a_7\rightarrow a_8$, there exists $u\in N^-_{T_2}(y^i_5)$ such that $T_2-u\neq 4A$ with $y^i_3\rightarrow 4$, so $(T_2-u;a_2a_1...a_8)\neq (4A;(1,1,1))$ with $y^i_3\rightarrow 4$ and thus $y^i_3$ has an inneighbor origin of a path $P\equiv a_2a_1...a_8$ in $T_2$, we have $y^i_3y^i_4y^i_5uvP^{-1}y^i_3\equiv a_3a_4...a_3=C$. Else, $a_8\rightarrow a_7$, let $u\in N^-_{T_2}(y^i_3)$ and $w\in N^+_{T_2}(u)$, $T_2\setminus \{u,w\}$ has at least $2$ origins of a path $a_6a_7...a_{n-1}$ so $y^i_5$ has an inneighbor origin of a path $P\equiv a_6a_7...a_{n-1}$ in $T_2\setminus \{u,w\}$, we have $vwuy^i_3y^i_4y^i_5Pv\equiv C$.\\
So we may assume that $a_{\delta+3}\rightarrow a_{\delta+4}$, if there exists $x\in T_1$ such that $d^+_{T_2}(x)=1$ then $T_1=3A$, set $N^+_{T_2}(x)=\{u\}$ and let $y_3y_4y_5\equiv a_3a_4a_5$ be a path in $T_1$ of end $x$ and $w,t\in N^-(y_3)$ such that $w\rightarrow t$. $T_2\setminus \{w,t\}$ contains a path $Q\equiv a_6a_7...a_{n-1}$ of origin an inneighbor of $x$ and then $vtwy_3y_4y_5Qv\equiv C$ unless three cases:\begin{itemize}
\item[(i)] If $|T_2|=4$, set $V(T_2)=\{u,w,t,y\}$, if $y\leftarrow w$ or $t$, say $w$, then $ty_5y_4vuywy_3t\equiv C$, so $y\rightarrow \{w,t\}$. Without loss of generality we may suppose that $u\rightarrow T_1$, so $y_4\leftarrow w$ or $y$, say $w$, we have $y_3uy_5vytwy_4y_3\equiv C$.
\item[(ii)] If $(T_2\setminus \{w,t\};a_6a_7...a_{n-1})=(4A;(1,1,1))$, there exists $w',t'\in T_2\setminus \{w,t\}$ such that $T_2\setminus \{w',t'\}\neq 4A$, thus $y_3$ has an inneighbor origin of a path $P=(1,1,1)$ in $T_2\setminus \{w',t'\}$, therefore $y_3y_4y_5w't'vP^{-1}y_3\equiv a_3a_4...a_3=C$.
\item[(iii)] If $(T_2\setminus \{w,t\};a_6a_7...a_{n-1})=(5C;(2,1,1))$ or $T_2\setminus \{w,t\}$ is a transitive triangle and $a_6a_7a_8$ is directed, both cases are solved in lemma \ref{3.3} by taking $\overline{C}=a'_0a'_1...a'_0=\overline{a_7a_6...a_7}$.
\end{itemize}
Now $d^+_{T_2}(x)\geq 2$ for all $x\in T_1$, let $a\in T_2$ such that $d^+_{T_2}(a)=\Delta^+(T_2)$ and set $T'_1=T_1\cup \{a\}$ and $T'_2=T_2-a$. Let $x_{n-1}x_0...x_{\delta-1}\equiv a_{n-1}a_0...a_{\delta-1}$ be a path in $T'_1$ such that $d^-_{T'_2}(x_{n-1})=max\,\{d^-_{T'_2}(x);x\in Or(a_{n-1}a_0...a_{\delta-1})\}$. If $x_{n-1}$ has an inneighbor origin of a path $P\equiv a_{n-2}a_{n-3}...a_{\delta+1}$ in $T'_2$ then $x_0...x_{\delta-1}vP^{-1}x_{n-1}x_0\equiv C$. Otherwise, if $d^-_{T'_2}(x_{n-1})\geq 2$ then $(T'_2;a_{n-2}a_{n-3}...a_{\delta+1})$ is an exception. By the building lemmas we may assume that it is one of the following biexceptions with respect to $x_{n-1}$:\\
Dual Exc $(1,0)$: We have $x_{\delta-2}x_{\delta-3}...x_{n-1}3421vx_{\delta-1}x_{\delta-2}\equiv \overline{C}$.\\
Dual Exc $(7,0)$: If $\delta=3$ it is solved in lemma \ref{3.3} by taking $C=a'_0a'_1...a'_0=a_4a_5...a_4$ and if $\delta=5$, let $Q\equiv a_5a_6...a_9$ be a path in $T_1$ of origin $x_0$, since $T_2\neq 6H$ and $d^-_{T_2}(x_0)=4$ then $x_0$ has an inneighbor origin of a path $R\equiv a_4a_3...a_{n-1}$ in $T_2$, thus $QvR^{-1}x_0\equiv a_5a_6...a_5=C$.\\
Dual Exc $(48,0)$: Since $d^+_{T_2}(x_0)=2$ then $\delta \neq 7$. If $\delta=3$ then it is solved in lemma \ref{3.3} by taking $\overline{C}=a'_0a'_1...a'_0=\overline{a_7a_6...a_7}$ and if $\delta=5$ then it is solved in lemma \ref{3.5} by taking $\overline{C}=a'_0a'_1...a'_0=\overline{a_9a_8...a_9}$.\\
Dual Exc $(E'_9(7),0),(E'_{10}(9),0)$: Let $Q\equiv a_{n-4}a_{n-3}...a_{\delta-5}$ be a path of origin $x_{n-1}$ in $T_1$, since $d^-_{T_2}(x_{n-1})\geq 5$ and $(T_2;a_{\delta-5}a_{\delta-4}...a_{n-6})$ is not an exception then $x_{n-1}$ has an inneighbor origin of a path $R\equiv a_{\delta-5}a_{\delta-4}...a_{n-6}$ in $T_2$ and thus $vQRv\equiv a_{n-6}a_{n-5}...a_{n-6}=C$.\\
Dual Exc $(E'_9(8),0),(E'_{10}(10),0)$: Let $Q\equiv a_{n-5}a_{n-4}...a_{\delta-6}$ be a path of origin $x_{n-1}$ in $T_1$, since $d^+_{T_2}(x_{n-1})\geq 3$ and $(T_2;a_{\delta-6}a_{\delta-5}...a_{n-7})$ is not an exception then $x_{n-1}$ has an outneighbor origin of a path $R\equiv a_{\delta-6}a_{\delta-5}...a_{n-7}$ in $T_2$ and thus $vQRv\equiv a_{n-7}a_{n-6}...a_{n-7}=C$.\\
Exc $(0,1)$: We have $d^+_{T_2}(x_{n-1})=2$ then $x_{n-1}\rightarrow a$ and thus $T_1\rightarrow a\rightarrow T_2$, $x_0\rightarrow 1$ or $2$, say $1$, we have $a3x_{n-1}x_2x_01v2a\equiv C$.\\
Exc $(4,1)$: $d^-_{T_2}(x_0)=4$ then $x_0\neq a$ and so $d^-_{T_2}(x_{\delta-1})=4$, by theorem \ref{2.1} $x_0$ or $x_{\delta-1}$ is an origin of a path $y_{n-1}y_0...a_{\delta-1}\equiv a_{n-1}a_0...a_{\delta-1}$ in $T'_1$ with $d^-_{T'_2}(y_{n-1})=4$ then $T$ contains $C$.\\

Exc $(8,1)$: $d^-_{T_2}(x_0)\geq 3$ then $x_0\neq a$ and thus $\{2,3,4\}\rightarrow x_{\delta-1}$, if $x_{\delta-1}\rightarrow 1$ then $x_0...x_{\delta-1}134v25x_{n-1}x_0\equiv C$ and if $x_{\delta-1}\rightarrow 5$ then $x_0...x_{\delta-1}541v32x_{n-1}x_0\equiv C$.\\
Exc $(9,1)$: For all $u\in N^+_{T'_2}(x_{\delta-1})$, there exist $u_1,u_2,w_1,w_2\in T'_2$ such that $uu_1u_2$ and $w_1w_2$ are two directed outpaths, $w_1\in N^-_{T'_2}(x_{n-1})$ and $V(T'_2)=\{u,u_1,u_2,w_1,w_2\}$, thus $x_0...x_{\delta-1}uu_1u_2vw_2w_1x_{n-1}x_0\equiv C$.\\
Exc $(16,1)$: If $\delta=3$ then it is solved in lemma \ref{3.3} by taking $\overline{C}=a'_0a'_1...a'_0=\overline{a_7a_6...a_7}$. If $\delta=5$ or $7$, there exists $x\in T_1$ such that $6\rightarrow x$ ($x$ can be $x_{n-1}$ or $x_0$), if $T_2\neq \overline{7C}$ then, by theorem \ref{2.1}, $x$ has an inneighbor origin of an antidirected Hamiltonian inpath $P$ in $T_2$, also if $T_2=\overline{7C}$, $6$ is an origin of $P$, let $Q\equiv a_4a_5...a_{\delta+3}$ be a path in $T_1$ of origin $x$, we have $QvP^{-1}x\equiv a_4a_5...a_4=C$.\\
Exc $(32,1)$: Similar to Exc $(16,1)$.\\
Exc $(33,1)$: $d^-_{T_2}(x_0)=6$ then $x_0\neq a$ and so $d^-_{T_2}(x_{\delta-1})=6$, by theorem \ref{2.1} $x_0$ or $x_{\delta-1}$ is an origin of a path $y_{n-1}y_0...a_{\delta-1}\equiv a_{n-1}a_0...a_{\delta-1}$ in $T'_1$ with $d^-_{T'_2}(y_{n-1})=6$ then $T$ contains $C$.\\
Exc $(E_8(5),1),(E_8(6),1)$: $d^-_{T_2}(x_0)\geq 3$ then $x_0\neq a$ and thus $\{1,2,3\}\rightarrow x_{\delta-1}$. So $d^+_X(x_{\delta-1})\geq 1$ and then $x_{\delta-1}$ has an outneighbor origin of a path $P\equiv a_{n-1}a_{n-2}...a_{\delta+2}$ in $T'_2$, therefore $x_{\delta-1}x_{\delta-2}...x_{n-1}vP^{-1}x_{\delta-1}\equiv C$.\\
So we may assume now that $d^-_{T'_2}(x_{n-1})=1$, and thus $T_1=3A$, if $x_7\rightarrow x_8$ it is solved in lemma \ref{3.5} by taking $\overline{C}=a'_0a'_1...a'_0=\overline{a_7a_6...a_7}$, so $a_8\rightarrow a_7$. Let $y_1y_2y_3$ be a directed outpath in $T_1$ such that $y_1=x_{n-1}$, and let $u\in N^+_{T_2}(y_3)$. Set $P=a_2a_1a_0...a_7$, we have $P=+(1,1,1,Q,1)$ with an odd number of blocks ($P$ can be $+(1,1,1)$) and $d^+_{T_2-u}(y_1)\geq |T_2-u|-2$, so by theorem \ref{2.1} $y_1$ has an outneighbor origin of a path $^*P$ in $T_2-u$ and thus $y_1y_2y_3wv(^*P)^{-1}y_1\equiv a_2a_3...a_2=C$ unless $(T_2;^*P)$ is an exception. By the building lemmas we may assume that $(T_2-u;^*P)$ is a biexception with respect to $y_1$, and the sole possible biexceptions are:\\
Dual Exc $(4,1)(2)$: $a\equiv 4$ or $5$ and $a\rightarrow u$ (since $a$ has a maximal outdegree in $T_2$), if $a\equiv 4$ then $21y_1y_2y_3u43v52\equiv C$ and if $a\equiv 5$ then $32y_1y_2y_3u54v13\equiv C$.\\
Dual Exc $(4,1)(3)$: $a\equiv 3$ or $5$, if $a\equiv 3$ then $14y_1y_2y_3u32v51\equiv C$ and if $a\equiv 5$ then $32y_1y_2y_3u54v13\equiv C$.\\
Dual Exc $(9,1)(1)$: $a\equiv 1$ or $3$, if $a\equiv 1$ then $v2y_1y_2y_3u1453v\equiv C$ and if $a\equiv 3$ then $v5y_1y_2y_3u3124v\equiv C$.\\
Dual Exc $(15,1)(1)$: $a\equiv 4$ or $5$, say $a\equiv 4$, then $21y_1y_2y_3u46v352\equiv C$.\\
Assume now that $(T_1;a_1...a_{\delta})$ is not a Gr\"{u}nbaum's exception and let $x_1...x_{\delta}\equiv a_1...a_{\delta}$ be a path in $T_1$ such that $d^-_{T_2}(x_1)=max\,\{d^-_{T_2}(x);x\in Or(a_1...a_{\delta},T_1)\}$, if $x_1$ has an inneighbor $u$ origin of a path $P\equiv a_0a_{n-1}...a_{\delta+2}$ in $T_2$ then $ux_1...x_{\delta}vP^{-1}\equiv C$. Otherwise, if $d^-_{T_2}(x_1)\geq 2$ and $d^+_{T_2}(x_2)=0$ then, by lemma \ref{2.8} $(T_2;a_0a_{n-1}...a_{\delta+2})$ is one of the exceptions $1,18$ or $22$. If $T_1-x_2$ contains a path $y_2...y_{\delta}\equiv a_2...a_{\delta}$ then every vertex of $T_2$ is an origin of a path $P\equiv a_1a_0...a_{\delta+2}$ in $T_2\cup \{x_2\}$ with end distinct from $x_2$ so $y_2...y_{\delta}vP^{-1}y_2\equiv a_2a_3...a_2=C$; and if $(T_1-x_2;a_2a_3a_4)=(3A;(1,1))$, then without loss of generality we may suppose that $Or(a_0a_{n-1}...a_{\delta+2};T_2)\subseteq N^+(x)$ for all $x\in T_1-x_2$, so if $(T_2;a_0a_{n-1}...a_{\delta+2})$ is the exception Exc $18$ then $14365x_2x_3x_4x_1v21\equiv C$ and if it is the exception Exc $22$ then one can check that there exist $w\in N^+_{T_2}(x_4)$ and $w'\in N^+_{T_2}(w)$ such that $T_2\setminus \{w,w'\}$ contains a path $P\equiv a_0a_{n-1}...a_{\delta+4}$ of origin an inneighbor $x$ of $x_1$ then $xx_1...x_4ww'vP^{-1}\equiv C$. If $d^-_{T_2}(x_1)\geq 1$ and $d^+_{T_2}(x_2)\geq 1$ then $(T_2;a_0a_{n-1}...a_{\delta+2})$ is an exception. By the building lemmas we may assume that it is one of the following biexceptions with respect to $x_1$:\\
Exc $(1,1)(1)$: If $x_{\delta}\rightarrow 3$ then $1x_1...x_{\delta}34v21\equiv C$ and if $x_{\delta}\rightarrow 4$ then $2x_1...x_{\delta}41v32\equiv C$, so $\{3,4\}\rightarrow x_{\delta}$. If $a_{\delta-1}\rightarrow a_{\delta}$ then $1x_1...x_{\delta-1}v3x_{\delta}421\equiv C$ so $a_{\delta}\rightarrow a_{\delta-1}$. If $\delta=2$ then $v43x_2x_112v\equiv C$. If $\delta=3$, we can assume that $T_1$ is transitive since otherwise $x_{\delta}$ is an origin of a directed inpath which has $4$ as inneighbor, $2\leftarrow x_2$ or $x_3$ since otherwise $d^+(2)=5$, without loss of generality suppose that $x_2\rightarrow 2$, we have $x_1x_3v3241x_2x_1\equiv \overline{C}$. If $\delta=4$ then $a_1...a_4$ is antidirected and $T_1=4A$ and thus $x_2x_4x_3x_1=-(1,2)$. We have $v14x_2x_4x_3x_123v\equiv C$.\\
Exc $(1,1)(2)$: By supposing that $N^+(x_{\delta})\cap \{1,3,4\}\neq \emptyset$ one can find $C$, and since $d^-(2)=4$ then $T_2\rightarrow x_{\delta}$. If $\delta=2$ then $v43x_2x_112v\equiv C$, else, $N^-(x_1)\neq \overline{4A}$, so $T$ contains $C$.\\
Exc $(1,1)(3)$: By supposing that $N^+(x_{\delta})\cap \{1,4\}\neq \emptyset$ one can find $C$ so $\{1,4\}\rightarrow x_{\delta}$. If $a_{\delta-1}\rightarrow a_{\delta}$ then $2x_1...x_{\delta-1}v4x_{\delta}132\equiv C$ so $a_{\delta}\rightarrow a_{\delta-1}$. If $\delta=2$ then $v14x_2x_123v\equiv C$. If $\delta=3$ then $T_1$ is transitive, if $2\rightarrow x_2$ then $x_2x_3v143x_12x_2\equiv \overline{C}$ and if $x_2\rightarrow 2$ then $x_1x_3v3142x_2x_2\equiv \overline{C}$. If $\delta=4$ then $a_1...a_4$ is transitive and $T_1=4A$. We have $4\rightarrow x_2$ or $x_3$, say $x_2$, then $v14x_2x_4x_3x_123v\equiv C$.\\
Exc $(1,1)(4)$: By supposing that $N^+(x_{\delta})\cap \{1,3,4\}\neq \emptyset$ one can find $C$ so $\{1,3,4\}\rightarrow x_{\delta}$. Similarly we can assume that $a_{\delta}\rightarrow a_{\delta-1}$. If $\delta=2$ or $3$ it is solved as Exc $(1,1)(1)$. If $\delta=4$, since $d^-_{T_1}(x_1)=1$ then $T_1\neq 4A$ and so $x_2$ or $x_4$ is an origin of an antidirected inpath in $T_1$ which has $4$ as inneighbor and the result follows.\\
Exc $(6,1)$: If $x_{\delta}\rightarrow 1$ then $5x_1...x_{\delta}142v35\equiv C$, so $\{1,2,3\}\rightarrow x_{\delta}$. If $a_{\delta-1}\rightarrow a_{\delta}$, if $x_{\delta}\rightarrow 4$ then $5x_1...x_{\delta-1}v14x_{\delta}235\equiv C$ and if $4\rightarrow x_{\delta}$ then $5x_1...x_{\delta-1}v1x_{\delta}4235\equiv C$, so $a_{\\delta}\rightarrow a_{\delta-1}$. If $\delta=2$ then $241x_2x_15v32\equiv C$. If $\delta=3$ then $34vx_1x_2x_32513\equiv C$. If $\delta=4$ then $a_1...a_4$ is antidirected and $T_1=4A$, we have $N^+(x_4)\cap \{4,5\}\neq \emptyset$ so let $x\in N^+_{T_2}(x_4)$, $x_1$ has an outneighbor in $T_2$ origin of a path $Q=-(2,1)$ in $T_2-x$ then $vxx_4...x_1Qv\equiv C$. If $\delta=5$, if $a_1...a_5=-(2,1,1)$ then $53vx_5...x_14125\equiv C$ and if $a_1...a_5=-(1,1,2)$ then by theorem \ref{2.1} $x_2$ or $x_5$ is an origin of a path $Q=-(1,1,2)$ and thus $1Qv43521\equiv C$.\\
Exc $(10,1)$: By supposing that $N^+(x_{\delta})\cap \{1,2,3\}\neq \emptyset$ one can find $C$ so $\{1,2,3\}\rightarrow x_{\delta}$. If $a_{\delta-1}\rightarrow a_{\delta}$ then $5x_1...x_{\delta-1}v1x_{\delta}3425\equiv C$, so $a_{\delta}\rightarrow a_{\delta-1}$. If $\delta=2$ then $4213x_1x_2v54\equiv C$. If $\delta=3$ then $5243x_3x_1x_2v15\equiv C$. If $\delta=4$ then $x_2x_135v421x_4x_3x_2$\\$\equiv C$. If $\delta=5$, if $a_1...a_5=-(2,1,1)$ then $41vx_5...x_13524\equiv C$ and if $a_1...a_5=-(1,1,2)$ then by theorem \ref{2.1} $x_2$ or $x_5$ is an origin of a path $Q=-(1,1,2)$ and thus $1Qv43521\equiv C$.\\
Exc $(12,1)$: By supposing that $N^+(x_{\delta})\cap \{3,4,5\}\neq \emptyset$ one can find $C$ so $\{3,4,5\}\rightarrow x_{\delta}$. If $a_{\delta-1}\rightarrow a_{\delta}$ then $1x_1...x_{\delta-1}v3x_{\delta}5421\equiv C$, so $a_{\delta}\rightarrow a_{\delta-1}$. If $\delta=2$ then $v54x_2x_1213v\equiv C$. If $\delta=3$ then $x_3x_11v4532x_2x_3\equiv C$. If $\delta=4$ then $T_1=4A$, we have $453x_2x_4x_3x_11v24\equiv C$. If $\delta=5$, if $a_1...a_5=-(2,1,1)$ then by taking $C=a'_0a'_1...a'_{n-1}a'_0=a_6a_5...a_7a_6$ it is solved in lemma \ref{3.3}, and if $a_1...a_5=-(1,1,2)$ then by theorem \ref{2.1} $x_2$ or $x_5$ is an origin of a path $Q=-(1,1,2)$ and thus $1Qv43521\equiv C$.\\
Exc $(22,1)$: If $x_{\delta}\rightarrow 4$ then $1x_1...x_{\delta}46v3521\equiv C$, so $\{4,5,6\}\rightarrow x_{\delta}$. If $a_{\delta-1}\rightarrow a_{\delta}$ then $3x_1...x_{\delta-1}v6x_{\delta}52413\equiv C$, so $a_{\delta}\rightarrow a_{\delta-1}$. If $1\rightarrow x_{\delta-1}$ then $2x_1...x_{\delta-1}15x_{\delta}46v32\equiv C$ and if $4\rightarrow x_{\delta-1}$ then $1x_1...x_{\delta-1}46x_{\delta}53v21\equiv C$.\\
Exc $(24,1)$: By supposing that $d^+_{T_2}(x_{\delta})\geq 2$ one can find $C$ so $T_2\rightarrow x_{\delta}$ and thus $a_{\delta}\rightarrow a_{\delta-1}$. If $1\rightarrow x_{\delta-1}$ then $3x_1...x_{\delta-1}15x_{\delta}62v43\equiv C$. If $4\rightarrow x_{\delta-1}$ then $1x_1...x_{\delta-1}46x_{\delta}35v21\equiv C$. If $6\rightarrow x_{\delta-1}$ then $2x_1...x_{\delta-1}61x_{\delta}54v32\equiv C$.\\
Exc $(26,1)(1),(2)$: By supposing that $N^+(x_{\delta})\cap \{1,2,3,6\}\neq \emptyset$ one can find $C$, so $\{1,2,3,6\}\rightarrow x_{\delta}$ and then $5x_1...x_{\delta-1}v6x_{\delta}32145\equiv C$.\\
Exc $(26,1)(3)$: By supposing that $N^+(x_{\delta})\cap \{1,2,3,5\}\neq \emptyset$ one can find $C$, so $\{1,2,3,5\}\rightarrow x_{\delta}$. If $x_{\delta}\rightarrow 4$ then $T(N^-(4))\neq \overline{6H}$ and so $T$ contains $C$, thus $4\rightarrow x_{\delta}$ and $6x_1...x_{\delta-1}v1x_{\delta}45236\equiv C$.\\
So we may suppose that $d^-_{T_2}(x_1)=1$, set $\{u\}=N^-_{T_2}(x_1)$,so by lemma \ref{2.3} three cases may arise:\begin{itemize}
\item[(i)] If $T_1=\overline{4B}$ and $a_1...a_4=-(1,2)$. We can suppose that $u$ is minimal in $T_2$ and $T_2-u=3A,5A$ or $7A$, since otherwise by taking $x_1$ in the place of $v$, $v$ is minimal in $T(N^+(x_1))$. If $T_2\neq 4B$ then $x_2$ has an inneighbor origin of a path $P\equiv a_{n-2}a_{n-3}...a_{\delta}$ in $T_2$ and thus $x_4x_1x_3vP^{-1}x_2x_4\equiv C$. Otherwise, set $V(T_2)=\{u,u_1,u_2,u_3\}$, we have $x_2x_1uu_1x_4vu_3u_2x_3x_2\equiv C$.
\item[(ii)] If $T_1$ is a transitive triangle and $a_1a_2a_3$ is directed, if $a_{\delta+4}\rightarrow a_{\delta+3}$ it is solved in lemma \ref{3.3} by taking $C=a'_0a'_1...a'_0=a_{\delta+3}a_{\delta+2}...a_{\delta+3}$, so we may assume that $a_{\delta+3}\rightarrow a_{\delta+4}$. Similarly we may assume that $a_{n-2}\rightarrow a_{n-3}$ since otherwise $T$ contains $\overline{C}=\overline{a_{n-1}a_0...a_{n-1}}$. Suppose that there exists $w\in N^+_{T_2-u}(x_2)$, if $x_1$ has an outneighbor origin of a path $P\equiv a_{n-1}a_{n-2}...a_{\delta+2}$ in $T_2-w$ then $x_1wx_2x_3vP^{-1}x_1\equiv C$, otherwise $(T_2-w;a_{n-1}...a_{\delta+2})$ is a Gr\"{u}nbaum's exception and then $w$ is an origin of a path $P\equiv \overline{a_{n-2}a_{n-3}...a_{\delta}}$ in $T_2$ and thus $x_1x_3vP^{-1}x_2x_1\equiv \overline{C}$. So we may assume that $N^+_{T_2}(x_2)=\{u\}$, and thus $x_3\rightarrow u$ (to obtain $T(N^-(x_1))$ a transitive triangle). We have $|T_2|\geq 5$ since otherwise $T(N^-(x_3))=3A$, then $a_{n-1}a_{n-2}...a_5$ has at least four blocks, let $s=max\,\{i;d^+_C(a_i)=0\}$, let $w\in N^+_{T_2}(u)$, $Q\equiv \a_{n-1}a_{n-2}...a_{s+1}$ and $Q'\equiv a_{s-1}a_{s-2}...a_5$ be two paths such that $V(Q)\cup V(Q')=V(T_2)\setminus \{u,w\}$, we have $x_1wux_3vQ'^{-1}x_2Q^{-1}x_1\equiv C$.
\item[(iii)] If $\delta=2$, if $|T_2|=2$ then $(T;C)=A_2$ which is an exception, so we may assume that $|T_2|\geq 4$. If there exists $w\in N^+_{T_2-u}(x_2)$, since $(T_2-w;a_4a_5...a_{n-1})\neq (\overline{4B};-(1,2))$ (in fact, otherwise, it is solved in lemma \ref{3.3} by taking $C=a'_0a'_1...a'_0=a_3a_4...a_3$), $x_1$ has an outneighbor origin of a path $P\equiv a_4a_5...a_{n-1}$ in $T_2-w$ and thus $vwx_2x_1Pv\equiv C$. So we can suppose that $(T_2-u)\rightarrow x_2$, there exist $w_1,w_2\in T_2-u$ and $P\equiv a_5a_6...a_{n-1}$ a path in $T_2\setminus \{w_1,w_2\}$ of origin an inneighbor of $x_2$ and thus $vw_1x_1w_2x_2Pv\equiv C$ unless $T_2=4B$, in this case set $V(T_2)=\{u,u_1,u_2,u_3\}$, we have $vu_2u_1x_2x_1uu_3v\equiv C$.    $\square$
\end{itemize}
\end{proof}
\begin{theorem}\label{3.7}
Let $T$ be a tournament of order $n=2k\geq 4$, $T$ contains an antidirected Hamiltonian cycle $C$ if and only if $(T;C)$ is not one of the exceptions.
\end{theorem}
\begin{proof}
By theorem \ref{3.1}, we may suppose that $\delta^-(T)\leq \delta^+(T)$. Set $\delta=\delta^-(T)$. Let $v\in T$ be a vertex such that $d^-(v)=\delta$ and denote by $T_1=T(N^-(v))$ and $T_2=T(N^+(v))$.\\
\textbf{Case $1$}: If $\delta=0$, $T$ contains $C$ if and only if $T_2$ contains an antidirected inpath $P$, in this case $vPv\equiv C$. So the sole exceptions are $A_4,A_5$ and $A_6$.\\
\textbf{Case $2$}: If $\delta=1$, set $T_1=\{u\}$. $3$ cases may arise:\begin{itemize}
\item[(a)] If $d^+_{T_2}(u)\geq 3$, let $w\in N^+_{T_2}(u)$ such that $d^-_{T_2}(w)\geq 1$. If $u$ has an outneighbor origin of an antidirected inpath $P$ in $T_2-w$ then $uPvwu\equiv C$, otherwise, $(T_2-w;-(1,1,...,1))$ is one of the following biexceptions with respect to $u$: (We may suppose that $w\rightarrow (T_2\setminus \{S\cup \{w\}\})$ since, otherwise, let $R=(1,1,...,1)$ be a path in $(T_2-w)\cup \{u\}$ of origin in $T_2\setminus \{S\cup \{w\}\}$ and end in $T_2$, we have $wRvw\equiv C$)\\
Dual Exc $(0,1)$: If $1\rightarrow w$ then $uw13v2u\equiv C$ and if $w\rightarrow w$ then $u1w3v2u\equiv C$.\\
Dual Exc $(4,1)$: We have $u1w5v342u\equiv C$.\\
Dual Exc $(8,1)$: If $w\rightarrow 1$ then $u2w1v345u\equiv C$ and if $w\rightarrow 5$ then $u5w4v132u\equiv C$. So we may assume that $\{1,5\}\rightarrow w$, in this case $T-v=7A$ and thus $(T;C)=A_6$.\\
Dual Exc $(9,1)(1),(3)$: We have $u415w3v2u\equiv C$.\\
Dual Exc $(9,1)(2)$: We have $u514w3v2u\equiv C$.\\
Dual Exc $(9,1)(4)$: We have $u5w231v4u\equiv C$.\\
Dual Exc $(33,1)$: We have $u1w62435v7u\equiv C$.\\
Dual Exc $(E_8(5),1)$: We have $(T;C)=A_{10}$.
\item[(b)] If $d^-_{T_2}(u)\geq 4$, let $w\in N^-_{T_2}(u)$ such that $d^+_{T_2}(w)\geq 2$ and $t\in N^+_{T_2}(w)$. If $u$ has an inneighbor origin of an antidirected outpath $P$ in $T_2\setminus \{w,t\}$ then $uPvtwu\equiv C$. Otherwise, $(T_2\setminus \{w,t\};(1,1,...,1))$ is one of the following biexceptions with respect to $u$:\\
Exc $(1,1)(4)$: We can suppose that $T(w,t,1,2)=4A$ since otherwise it contains $P=(1,1,1)$ with origin an inneighbor of $u$, so $u\rightarrow t$ and $t\rightarrow \{1,2\}$. If $t\rightarrow 4$ then $u34t1vw2u\equiv C$, else, $u124tv3wu\equiv C$ if $w\rightarrow 3$ and $u124tvw3u\equiv C$ if $3\rightarrow w$.\\
Exc $(1,1)(1),(2),(3)$: There exist $w',t'\in T_2$ such that $w'\rightarrow \{t',u\}$ and $d^-_{T_2\setminus \{w',t'\}}(u)\geq 3$, so $T$ contains $C$.\\
Exc $(24,1)$: By taking $T_2\setminus \{4,6\}$, we can suppose that $T(w,t,1,2,3,5)=6H$ and $u\rightarrow t$, so $\{1,2,3\}\rightarrow w$ and thus $u6541w23vtu\equiv C$.\\
Exc $(22,1)$: Taking $T_2\setminus \{1,5\}$, we obtain $d^-_{T_2\setminus \{1,5\}}(u)\geq 4$, and the result follows.
\item[(c)] If $d^+_{T_2}(u)\leq 2$ and $d^-_{T_2}(u)\leq 3$, then $|T_2|\leq 4$. If $|T_2|=2$ then $T=3A\cup \{x\}$ which is the exception $A_4$. So we may assume that $|T_2|=4$, set $V(T_2)=\{w_1,w_2,w_3,w_4\}$, if $N^+_{T_2}(u)=\{w_1\}$ then $T$ contains $C$ if and only if $d^-_{T_2}(w_1)\geq 1$, in fact, without loss of generality, suppose that $w_2\in N^-(w_1)$ we have $uw_2w_1vw_3w_4u\equiv C$; so $C$ is not contained in $T$ if and only if $(T;C)$ is the exception $A_{11}$. If $N^+_{T_2}(u)=\{w_1,w_2\}$, say $w_3\rightarrow w_4$, then $T$ contains $C$ if and only if $w_1$ or $w_2\rightarrow w_3$, so $C$ is not contained in $T$ if and only if $(T;C)$ is the exception $A_{12}$.
\end{itemize}
\textbf{Case $3$}: If $\delta$ is a positive even number. Let $x_1...x_{\delta}$ be a Hamiltonian antidirected inpath of $T_1$ such that $d^-_{T_2}(x_1)=max\,\{d^-_{T_2}(x);x\in Or(-(1,1,...,1);T_1)\}$ and $d^+_{T_2}(x_{\delta})\geq 1$ (This is possible unless $\delta=2$). Let $u\in N^+_{T_2}(x_{\delta})$, if $x_1$ has an inneighbor origin of an antidirected outpath $P$ in $T_2-u$ then $x_1...x_{\delta}uvP^{-1}x_1\equiv C$. Otherwise, suppose first that $d^-_{T_2-u}(x_1)\geq 2$ then $(T_2-u;(1,1,...,1))$ is one of the following exceptions:\\
Exc $1$: We distinguish four cases:\begin{itemize}
\item[(i)] If $N^-_{T_2}(x_1)=\{1,2,3\}$ then by the building lemmas we may suppose that $\{2,4\}\subseteq N^-(x_2)$. By supposing that there exists $i\in \{1,3,4\}\cap N^-(u)$ one can find $C$, so we may assume that $u\rightarrow \{1,3,4\}$. 
If $T_1=4A$ then $x_1x_3vx_4x_241u32x_1\equiv C$, so $T_1\neq 4A$. If there exists $u'\in N^+_{T_2-u}(x_{\delta})$ the result follows so $N^+_{T_2}(x_{\delta})=\{u\}$. If $\delta=4$ then $x_{\delta}$ is an origin of an antidirected inpath in $T_1$ with $d^-_{T_2}(x_{\delta})>d^-_{T_2}(x_1)$, a contradiction. If $\delta=2$, if $x_1\rightarrow u$ then $ux_143x_212vu\equiv C$, and if $u\rightarrow x_1$ then $4ux_11x_223v4\equiv C$.
\item[(ii)] If $N^-_{T_2}(x_1)=\{1,2\}$ then by the building lemmas we may suppose that $\{2,4\}\subseteq N^-(x_2)$. By supposing that there exists $i\in \{3,4\}\cap N^-(u)$ one can find $C$, so we may assume that $u\rightarrow \{3,4\}$ and thus $3\rightarrow \{x_2,...,x_{\delta}\}$. If there exists $u'\in N^+_{T_2-u}(x_{\delta})$ the result follows so $N^+_{T_2}(x_{\delta})=\{u\}$. If $\delta=4$, if $u\rightarrow 2$ then $x_1vx_4x_3x_234u21x_1\equiv C$ and if $2\rightarrow u$ then $43x_3x_2x_11x_42uv4\equiv C$. If $\delta=2$, if $x_1\rightarrow u$ then $ux_143x_212vu\equiv C$ and if $u\rightarrow x_1$ then $4ux_11x_223v4\equiv C$.
\item[(iii)] If $N^-_{T_2}(x_1)=\{1,3\}$ then by the building lemmas we may suppose that $\{2,3,4\}\subseteq N^-(x_2)$. By supposing that there exists $i\in \{1,3,4\}\cap N^-(u)$ one can find $C$, so we may assume that $u\rightarrow \{1,3,4\}$. If $x_{\delta}\rightarrow 2$ then $x_1...x_{\delta}2v4u31x_1\equiv C$, if $x_{\delta}\rightarrow 3$ then $x_1vx_2...x_{\delta}3u142x_1\equiv C$ and if $x_{\delta}\rightarrow 4$ then $x_1vx_2...x_{\delta}4u312x_1\equiv C$, so $\{2,3,4\}\rightarrow x_{\delta}$. If $\delta=4$, suppose that $x_1\rightarrow x_4$, since $2$ and $4$ are origins of an antidirected outpath $P$ with end in $\{2,4,u\}$ then if $x_3$ has an inneighbor in $\{2,4,u\}$ we obtain $x_4x_1vx_2x_3P^{-1}x_4\equiv C$, so $x_3\rightarrow \{2,4,u\}$ and thus $2\rightarrow u$. If $1\rightarrow x_2$ then $4x_1x_43x_212x_3uv4\equiv C$ and if $x_2\rightarrow 1$ then $1x_2x_13x_42ux_34v1\equiv C$. So $x_4\rightarrow x_1$, and we can suppose that $x_3\rightarrow 4$. $x_2$ or $x_4$ is an origin of an antidirected inpath $y_1...y_4$ in $T_1$ such that $y_4\in \{x_1,x_3\}$ and thus $y_1...y_44v1u32y_1\equiv C$. If $\delta=2$, if $x_2\rightarrow 1$, if $u\rightarrow 2$ then $x_1x_21v2u43x_1\equiv C$, so $2\rightarrow u$; if $u\rightarrow x_1$ then $1ux_13x_242v1\equiv C$ and if $x_1\rightarrow u$ then $T-v=7A$ and $(T;C)$ is the exception $A_6$. So we may assume that  $1\rightarrow x_2$, if $x_1\rightarrow u$ then $ux_121x_234vu\equiv C$ and if $u\rightarrow x_1$ then $24x_21x_1u3v2\equiv C$.
\item[(iv)] If $N^-_{T_2}(x_1)=\{2,3\}$. By supposing that there exists $i\in \{1,4\}\cap N^-(u)$ one can find $C$, so we may assume that $u\rightarrow \{1,4\}$. If $x_{\delta}\rightarrow 2$ then $x_1...x_{\delta}2v1u43x_1\equiv C$, so $2\rightarrow x_{\delta}$.\\
Suppose first that $u\rightarrow 3$, by supposing that there exists $i\in \{1,4\}\cap N^+(x_{\delta})$ one can find $C$, so we may assume that $x_{\delta}\rightarrow \{1,4\}$. Suppose that $x_{\delta}\rightarrow 3$ then $3\rightarrow \{x_2,x_3\}$ if it exist. If $u\rightarrow x_1$ then $x_1...x_{\delta}3v241ux_1\equiv C$, so $x_1\rightarrow u$. If $u\rightarrow 2$ then $x_1vx_2...x_{\delta}312u4x_1\equiv C$, so $2\rightarrow u$. If $\delta=4$ then if $x_2\rightarrow 2$  and so $x_2x_3x_4vx_14u312x_2\equiv C$ and if $2\rightarrow x_2$, without loss of generality suppose that $x_4\rightarrow x_2$, so $1x_143x_3x_4x_22uv1\equiv C$. If $\delta=2$ then $31x_22ux_14v3\equiv C$. Otherwise, $3\rightarrow x_{\delta}$, suppose that $\delta=4$, if $u\rightarrow x_3$ then $3ux_3x_2x_12x_441v3\equiv C$, so $x_3\rightarrow u$. If $u\rightarrow x_1$, let $a\in N^-_{T_2}(x_3)$, $b\in N^+_{T_2-a}(u)$ and $\{c\rightarrow d\}=T_2\setminus \{a,b,u\}$, we have $bux_1x_2x_3ax_4cdvb\equiv C$, so $x_1\rightarrow u$. Suppose that $x_1\rightarrow x_3$, if $x_2\rightarrow x_4$ then $ux_1x_3x_2x_44213vu\equiv C$, so $x_4\rightarrow x_2$ and thus, similarly $x_2\rightarrow 4$. So $4\rightarrow x_3$ and since $d^+(x_1)=5$ then $x_4\rightarrow x_1$ and so $ux_4x_1x_2x_34213vu\equiv C$. So we may assume that $x_3\rightarrow x_1$, if $x_1\rightarrow x_4$ then if $x_3$ has an inneighbor origin of an antidirected outpath $P$ in $T_2$ of end distinct from $u$ then $x_4x_1vx_2x_3Px_4\equiv C$, otherwise, $2\rightarrow u$, $N^-_{T_2}(x_3)=\{1,3\}$ and $u\rightarrow x_2$. If $x_2\rightarrow 3$ then $x_2vx_3ux_1x_44213x_2\equiv C$, otherwise, $x_23x_4x_1vx_3241ux_2\equiv C$. So $x_4\rightarrow x_1$, if $4\rightarrow x_2$ then $x_3x_1x_4x_24213vux_3\equiv C$ and if $x_2\rightarrow 4$ then $4\rightarrow x_3$ and thus $x_4x_1x_2x_34213vux_4\equiv C$. Now if $\delta=2$ then $ux_143x_212vu\equiv C$. So we may assume that $3\rightarrow u$, if $\delta=2$ then if $4\rightarrow x_2$ then $u3x_12x_241vu\equiv C$, so $x_2\rightarrow 4$ and thus $u\rightarrow 2$ and $x_1\rightarrow u$ since otherwise we are in the case $(i)$ or $(ii)$. If $x_2\rightarrow 3$ then $x_1\rightarrow u$ and thus $3x_2ux_1142v3\equiv C$ and of $3\rightarrow x_2$ then $ux_143x_212vu\equiv C$. So $\delta=4$, suppose that $x_4\rightarrow x_1$, similarly we may assume that $\{1,2\}\rightarrow x_2$, also $x_3\rightarrow 4$ and thus $4\rightarrow \{x_2,x_4\}$. If $x_3\rightarrow x_1$ then $x_1\rightarrow u$ and $x_2x_4x_1x_3$ or $x_4x_2x_1x_3$ is an antidirected inpath with $d^-_{T_2}(x_2)\geq d^-_{T_2}(x_1)$ and $d^-_{T_2}(x_4)\geq d^-_{T_2}(x_1)$, a contradiction, so $x_1\rightarrow x_3$. Similarly we may suppose that $x_4\rightarrow x_2$. We have $x_1x_3x_4x_212vu34x_1\equiv C$, so $x_1\rightarrow x_4$. If $u\rightarrow x_3$ then $31x_42x_1x_2x_3uv3\equiv C$, so $x_3\rightarrow u$. Suppose that $u\rightarrow x_1$. If there exists $u'\in N^+_{T_2}(x_4)$ such that $x_1$ has an inneighbor in $T_2-u'$ of outdegree at least $2$ in $T_2-u'$ then $T$ contains $C$, so we may assume that $\{1,2,3,4\}\rightarrow x_4$. Similarly, by taking $x_4$ in the place of $v$, we may assume that $x_2\rightarrow 4$, so $4\rightarrow x_3$ and thus $1ux_1x_2x_34x_423v1\equiv C$. So $x_1\rightarrow u$ and then $u\rightarrow \{2,x_2\}$ and $x_3\rightarrow x_1$. And thus $T_1$ cannot contains a vertex of outdegree $0$ in $T_1$, so without loss of generality, if there exists $x\in T$ such that $d^-(x)=4$ and $T(N^-(x))$ contains a vertex of outdegree $0$ then $T$ contains $C$, by duality, if $d^+(x)=4$ and $T(N^+(x))$ contains a vertex of indegree $0$ then $T$ contains $C$, so we may assume that $x_2\rightarrow 4$ and thus $4\rightarrow \{x_3,x_4\}$. We have $u3x_1x_2x_34x_412vu\equiv C$.
\end{itemize}
Exc $20$: One can check that for any $i\in \{3,4,5,6\}$, $1$ or $2$ is an origin of an antidirected outpath $P$ in $T_2\setminus \{i,i+1,u\}$, so if there exists $i\in N^-(u)\cap \{3,4,5,6\}$ then $x_1...x_{\delta}ui(i+1)vP^{-1}x_1\equiv C$, thus $u\rightarrow \{3,4,5,6\}$. If there exists $u'\in N^+_{T_2-u}(x_{\delta})$, $1$ or $2$ is an origin of an antidirected ouptath $P$ in $T_2-u'$ and thus $x_1...x_{\delta}u'vP^{-1}x_1\equiv C$, so $N^+_{T_2}(x_{\delta})=\{u\}$. There exist $w\in N^-_{T_2}(x_{\delta-1})$, $w'\in N^+_{T_2}(w)$ and $a\in \{1,2\}$ pairwaisly distinct such that $x_{\delta}$ has an inneighbor origin of an antidirected outpath $P$ in $T_2\setminus \{a,w,w'\}$, thus $x_{\delta}ax_1...x_{\delta-1}ww'vP^{-1}x_{\delta}\equiv C$.\\
Exc $22$: Without loss of generality we can suppose that $\{1,2\}\subseteq N^-(x_1)$. One can check that for any $i\in \{4,5,6\}$, there exists an outneighbor $j$ of $i$ such that $1$ or $2$ is an origin of an antidirected outpath $P$ in $T_2\setminus \{i,j,u\}$, so if there exists $i\in N^-(u)\cap \{1,2,3,4,5\}$ then $x_1...x_{\delta}uijvP^{-1}x_1\equiv C$, thus $u\rightarrow \{4,5,6\}$. If there exists $u'\in N^+_{T_2-u}(x_{\delta})$ distinct from $4$, $1$ or $2$ is an origin of an antidirected ouptath $P$ in $T_2-u'$ and thus $x_1...x_{\delta}u'vP^{-1}x_1\equiv C$, and if $4\in N^+(x_{\delta})$ then $x_1...x_{\delta}4536v21x_1\equiv C$, so $N^+_{T_2}(x_{\delta})=\{u\}$. There exist $w\in N^-_{T_2}(x_{\delta-1})$, $w'\in N^+_{T_2}(w)$ and $a\in \{4,6\}$ pairwaisly distinct such that $x_{\delta}$ has an inneighbor origin of an antidirected outpath $P$ in $T_2\setminus \{a,w,w'\}$, thus $x_{\delta}ax_1...x_{\delta-1}ww'vP^{-1}x_{\delta}\equiv C$.\\
Exc $23$: One can check that for any $i\in \{1,2,3,4,5\}$, there exists an outneighbor $j$ of $i$ such that $4$ or $6$ is an origin of an antidirected outpath $P$ in $T_2\setminus \{i,j,u\}$, so if there exists $i\in N^-(u)\cap \{1,2,3,4,5\}$ then $x_1...x_{\delta}uijvP^{-1}x_1\equiv C$, thus $u\rightarrow \{1,2,3,4,5\}$. If there exists $u'\in N^+_{T_2-u}(x_{\delta})$, $4$ or $6$ is an origin of an antidirected ouptath $P$ in $T_2-u'$ and thus $x_1...x_{\delta}u'vP^{-1}x_1\equiv C$, so $N^+_{T_2}(x_{\delta})=\{u\}$. There exist $w\in N^-_{T_2}(x_{\delta-1})$, $w'\in N^+_{T_2}(w)$ and $a\in \{4,6\}$ pairwaisly distinct such that $x_{\delta}$ has an inneighbor origin of an antidirected outpath $P$ in $T_2\setminus \{a,w,w'\}$, thus $x_{\delta}ax_1...x_{\delta-1}ww'vP^{-1}x_{\delta}\equiv C$.\\
Exc $24$: Without loss of generality we may suppose that $1\in N^-(x_1)$. One can check that for any $i\in \{2,4,5,6\}$, there exists an outneighbor $j$ of $i$ such that $1$ is an origin of an antidirected outpath $P$ in $T_2\setminus \{i,j,u\}$, so if there exists $i\in N^-(u)\cap \{2,4,5,6\}$ then $x_1...x_{\delta}uijvP^{-1}x_1\equiv C$, thus $u\rightarrow \{2,4,5,6\}$. If $3\rightarrow u$ then $x_1...x_{\delta}u3456v21x_1\equiv C$, so $u\rightarrow 3$. If there exists $u'\in N^+_{T_2-u}(x_{\delta})$, $1$ is an origin of an antidirected ouptath $P$ in $T_2-u'$ and thus $x_1...x_{\delta}u'vP^{-1}x_1\equiv C$, so $N^+_{T_2}(x_{\delta})=\{u\}$. There exist $w\in N^-_{T_2}(x_{\delta-1})$, $w'\in N^+_{T_2}(w)$ and $a\in N^-_{T_2}(x_1)\cap N^-_{T_2}(x_{\delta})$ pairwaisly distinct such that $x_{\delta}$ has an inneighbor origin of an antidirected outpath $P$ in $T_2\setminus \{a,w,w'\}$, thus $x_{\delta}ax_1...x_{\delta-1}ww'vP^{-1}x_{\delta}\equiv C$.\\
Exc $29$ :One can check that for any $i\in \{2,3,5\}$, there exists an outneighbor $j$ of $i$ such that $4$ or $6$ is an origin of an antidirected outpath $P$ in $T_2\setminus \{i,j,u\}$, so if there exists $i\in N^-(u)\cap \{2,3,5\}$ then $x_1...x_{\delta}uijvP^{-1}x_1\equiv C$, thus $u\rightarrow \{2,3,5\}$. Similarly if $1\in N^-(u)$ then $x_1...x_{\delta}u1453v26x_1\equiv C$, thus $u\rightarrow 1$. If there exists $u'\in N^+_{T_2-u}(x_{\delta})$, $4$ or $6$ is an origin of an antidirected ouptath $P$ in $T_2-u'$ and thus $x_1...x_{\delta}u'vP^{-1}x_1\equiv C$, so $N^+_{T_2}(x_{\delta})=\{u\}$. There exist $w\in N^-_{T_2}(x_{\delta-1})$, $w'\in N^+_{T_2}(w)$ and $a\in \{4,6\}$ pairwaisly distinct such that $x_{\delta}$ has an inneighbor origin of an antidirected outpath $P$ in $T_2\setminus \{a,w,w'\}$, thus $x_{\delta}ax_1...x_{\delta-1}ww'vP^{-1}x_{\delta}\equiv C$.\\
Exc $30$: One can check that for any $i\in \{3,4,5\}$, there exists an outneighbor $j$ of $i$ such that $4$ or $6$ is an origin of an antidirected outpath $P$ in $T_2\setminus \{i,j,u\}$, so if there exists $i\in N^-(u)\cap \{3,4,5\}$ then $x_1...x_{\delta}uijvP^{-1}x_1\equiv C$, thus $u\rightarrow \{3,4,5\}$. Similarly if there exists $i\in N^-(u)\cap \{1,2,6\}$, set $\{u_1\rightarrow u_2\}=T(1,2,6)-i$, we have $x_1...x_{\delta}i5u_1u_2v43x_1\equiv C$, thus $u\rightarrow \{1,2,6\}$. If there exists $u'\in N^+_{T_2-u}(x_{\delta})$, $4$ or $6$ is an origin of an antidirected ouptath $P$ in $T_2-u'$ and thus $x_1...x_{\delta}u'vP^{-1}x_1\equiv C$, so $N^+_{T_2}(x_{\delta})=\{u\}$. There exist $w\in N^-_{T_2}(x_{\delta-1})$, $w'\in N^+_{T_2}(w)$ and $a\in \{4,6\}$ pairwaisly distinct such that $x_{\delta}$ has an inneighbor origin of an antidirected outpath $P$ in $T_2\setminus \{a,w,w'\}$, thus $x_{\delta}ax_1...x_{\delta-1}ww'vP^{-1}x_{\delta}\equiv C$.\\
Otherwise, $d^-_{T_2-u}(x_1)\leq 1$. If $\delta \geq 4$ then, by lemma \ref{2.6}, $T_1\in \{\overline{4A},4B\}$. If $T_1=4B$, taking $x_3x_2x_4x_1$ as antidirected inpath in $T_1$ with $u'\in N^+_{T_2}(x_1)$ such that $d^+_{T_2-u'}(x_3)\geq 2$ the result follows. So we may suppose that $T_1=\overline{4A}\equiv 4A$, $N^-_{T_2}(x_1)=\{u,w\}$ and $\{u\}\subseteq N^+_{T_2}(x_{\delta})\subseteq \{u,w\}$. If there exists $x\in N^+_{T_2}(x_1)\cap N^+_{T_2}(x_2)$ then $x_1$ has an outneighbor origin of an antidirected inpath $P$ in $T_2\setminus \{u,x\}$ and thus $x_1xx_2x_3x_4uvP^{-1}x_1\equiv C$ unless $T_2\setminus \{u,x\}=3A$, set $T_2\setminus \{u,x\}=\{w,w_1,w_2\}$, if $w\rightarrow x$ then $xx_2x_3x_4uvw_2x_1w_1wx\equiv C$ so $x\rightarrow w$. If $x_2\rightarrow w_1$, without loss of generality suppose that $w_2\rightarrow x$, then $w_1x_2x_3x_4uvww_2xx_1w_1\equiv C$, so $w_1\rightarrow x_2$. If $w\rightarrow u$ then $ux_4x_3x_2xvw_2x_1w_1wu\equiv C$, so $u\rightarrow w$. If $w_2\rightarrow x$ then $x_2x_3x_4ux_1w_1vww_2xx_2\equiv C$ so $x\rightarrow w_2$ and thus $w_1\rightarrow x$ (Since otherwise $w$ is an origin of an antidirected outpath in $T_2-u$). If $x_2\rightarrow w_2$ then $x_2x_3x_4ux_1w_1vwxw_2xx_2\equiv C$ so $w_2\rightarrow x_2$ and if $x_2\rightarrow w$ then $x_2x_3x_4ux_1xw_1w_2vwx_2\equiv C$ so $w\rightarrow x_2$. If $x_4\rightarrow w$ then $T(N^-(x_4))=T(x_1,w_1,w_2,x)\neq 4A$ so $T$ contains $C$, thus $w\rightarrow x_4$. If $x_2\rightarrow u$ then $T(N^+(x_4))\neq 4A$ so $u\rightarrow x_2$. If $x_3\rightarrow x$ then $T(N^+(x_2))\neq 4A$ so $x\rightarrow x_3$. If $u\rightarrow w_2$ then $x_1x_2x_3xx_4ww_1vw_2ux_1\equiv C$ and if $w_2\rightarrow u$ then $x_1x_2x_3xx_4w_2uvw_1wx_1\equiv C$. So we may assume that $N^+_{T_2}(x_1)\rightarrow x_2$. If $x_2\rightarrow w$, since $w$ is not an origin of an antidirected Hamiltonian outpath in $T_2-u$ then if $w$ is an origin of an antidirected inpath $P$ of order $3$ in $T_2-u$, let $Q$ be an antidirected Hamiltonian inpath of $T_2\setminus (\{u\}\cup V(P))$, we have $x_2x_3x_4uvQ^{-1}x_1P^{-1}x_2\equiv C$, and if not then $w\rightarrow (T_2\setminus \{u,w\})$ and $T_2\setminus =3A,5A$ or $7A$, since $\{x_3,u\}\subseteq N^+(x_4)\cap N^-(x_1)$ then there exists $u_1u_2u_3u_4\equiv x_1x_2x_3x_4$ in $T(x_1,x_3,x_4,u)$ such that $u_1=x_1$ and $u_4=x_4$, $x_2$ has an inneighbor origin of an antidirected outpath $P$ in $T_2\setminus \{u,w\}$ of order $|T_2|-3$, let $\{w'\}=V(T_2)-(V(P)\cup \{u,w\})$, thus $u_1u_2u_3u_4x_2Pvw'wx_1\equiv C$, so $w\rightarrow x_2$ and thus $N^+_{T_2}(x_2)=\{u\}$. So without loss of generality we may assume that $N^+_{T_2}(x_4)=\{u\}$ since otherwise we can take $x_4x_3x_2\equiv x_2x_3x_4$ and we find $C$. If $d^-_{T_2}(u)=0$ then $N^-(u)=\{x_2,x_3,x_4,v\}$ with $T(N^-(u))\neq 4A$ and thus $T$ contains $C$ so we may assume that $d^-_{T_2}(u)\geq 1$. If $d^+_{T_2-u}(x_3)\geq 2$, let $u'\in N^-_{T_2}(u)$, $w'\in N^+_{T_2-u'}(x_3)$ and $P$ be an antidirected outpath in $T_2\setminus \{u',w'\}$, we have $ux_2x_1x_3w'vP^{-1}u'u\equiv C$. Otherwise, $x_3\rightarrow u$ so $d^-_{T_2-u}(x_3)\geq 2$, let $u'\in N^+_{T_2}(u)$, $w'\in N^-_{T_2-u'}(x_3)$ and $P$ be an antidirected outpath in $T_2\setminus \{u',w'\}$, we have $x_3x_2x_1uu'vP^{-1}w'x_3\equiv C$.\\
Now if $\delta=2$, if $d^-_{T_2-u}(x_1)=1$ set $\{w\}=N^-_{T_2-u}(x_1)$. Suppose that $d^+_{T_2-w}(x_2)\geq 2$, then $|T_2|\geq 5$ since otherwise $d^-(x_2)\leq 1$, let $u'\in N^+_{T_2-w}(x_2)$ distinct from $u$, we may assume that $x_1\rightarrow u$. If $x_1$ has an outneighbor origin of an antidirected inpath $P$ in $T_2\setminus \{u,u'\}$ then $x_1u'x_2uvP^{-1}x_1\equiv C$, otherwise, $(T_2\setminus \{u,u'\};-(1,1,...,1))$ is a Gr\"{u}nbaum's exception and $d^+_{T_2}(x_2)=2$. Since $w$ is not an origin of an antidirected outpath in $T_2-u$ then there exists $w'\in N^-_{T_2-u}(u')$, let $P$ be an antidirected outpath in $T_2\setminus \{u,u',w'\}$, we have $x_2w'u'x_1uvP^{-1}x_2\equiv C$. So we may assume that $N^+_{T_2-w}(x_2)=\{u\}$. If $|T_2|\geq 7$, one can find $u_1,u_2,u_3\in T_2\setminus \{u,w\}$ such that $T_2\setminus \{w,u,u_1,u_2,u_3\}$ is not $4A$ nor $6H$, thus $x_2$ has an inneighbor origin of an antidirected outpath $P$ in $T_2\setminus \{w,u,u_1,u_2,u_3\}$, suppose that $u_3\rightarrow u_2$, therefore $x_1u_2u_3x_2Pvu_1x_1\equiv C$. So we may suppose that $|T_2|\in \{3,5\}$ and we have $3$ cases:\begin{itemize}
\item[(i)] If $T_2\setminus \{u,w\}$ is a transitive triangle, set $T'_2=T_2\setminus \{u,w\}=\{u_1,u_2,u_3\}$ such that $u_i\rightarrow u_{i+1}$ for $i\in \{1,2\}$. Obviously $\{u_2,u_3\}\rightarrow w$. If $x_2\rightarrow w$ then $x_2wu_2u_3x_1u_1vux_2\equiv C$, so $w\rightarrow x_2$. If $w\rightarrow u_1$, if $x_1\rightarrow u$ then $x_2wu_1vux_1u_3u_2x_2\equiv C$ so $u\rightarrow x_1$, if $u$ has an outneighbor $u_i$ in $T'_2$ then $x_1uu_ivu_{i+2}u_{i+1}x_2wx_1\equiv C$ so $T'_2\rightarrow u$ and thus $x_1x_2uu_2u_3vu_1wx_1\equiv C$; so $u_1\rightarrow w$. If $u$ has an inneighbor $u_i$ in $T'_2$ then $x_1vx_2uu_iwu_{i+1}u_{i+2}x_1\equiv C$, so $u\rightarrow T'_2$. If $x_1\rightarrow u$ then $ux_1u_2u_1x_2u_3wvu\equiv C$, so $u\rightarrow x_1$ and thus $u_1ux_1wx_2u_2u_3vu_1\equiv C$.
\item[(ii)] If $T_2\setminus \{u,w\}$ is a circuit triangle, set $T'_2=T_2\setminus \{u,w\}=\{u_1,u_2,u_3\}$. If $x_2\rightarrow w$, if $w$ has an inneighbor $u_i$ in $T'_2$ then $x_2wu_iu_{i+1}x_1u_{i+2}vux_2\equiv C$, and if $w\rightarrow T'_2$, if $d^-_{T'_2}(u)\geq 1$, say $u_1\rightarrow u$, then $x_1x_2uu_1u_2vu_3wx_1\equiv C$ and if $d^-_{T'_2}(u)=0$ then $(T;C)=A_{10}$, so we may assume that $w\rightarrow x_2$. If $x_1\rightarrow u$, if $T'_2\rightarrow w$ then $ux_1u_2u_1x_2u_3wvu\equiv C$ and otherwise, by taking $x_1$ in the place of $v$ one can find $C$; so $u\rightarrow x_1$. If $u$ has an outneighbor $u_i$ in $T'_2$ then $u_iux_1wx_2u_{i+1}u_{i+2}vu_i\equiv C$, so $T'_2\rightarrow u$. If $w$ has an outneighbor $u_i$ in $T'_2$ then $u_iwx_1x_2uu_{i+1}u_{i+2}vu_i\equiv C$, so $T'_2\rightarrow w$. We have $x_1vx_2uu_1wu_2u_3x_1\equiv C$.
\item[(iii)] If $T_2\setminus \{u,w\}=\{z\}$. We have $T-u=5A$ then $(T;C)=A_5$.
\end{itemize}
If $d^-_{T_2-u}(x_1)=0$ i.e. $N^-_{T_2}(x_1)=N^+_{T_2}(x_2)=\{u\}$, if there exists $u_1,u_2,u_3,u_4\in T_2-u$ such that $u_4$ is an origin of an antidirected outpath $P$ in $T_2\setminus \{u_1,u_2,u_3\}$, suppose that $u_3\rightarrow u_2$, then $u_1x_1u_2u_3x_2Pvu_1\equiv C$. Otherwise, if $|T_2|=5$, then $d^-_{T_2}(u)=0$ and thus $(T;C)=A_{10}$. If $|T_2|=3$ then $T-v=5A$ and thus $(T;C)=A_5$.\\
So we may assume now that $T_2\rightarrow x_2$. If $|T|=6$ then there exists $u\in T_2$ such that $T-u=5A$ and thus $(T;C)=A_5$, so we may suppose that $|T|\geq 8$. If $d^+_{T_2}(x_1)\geq 3$, let $\{u_1,u_2,u_3\}\subseteq N^+_{T_2}(x_1)$ such that $u_2\rightarrow u_3$ and let $P$ be an antidirected outpath in $T_2\setminus \{u_1,u_2,u_3\}$, we have $u_1x_1u_3u_2x_2Pvu_1\equiv C$, so $d^+_{T_2}(x_1)\leq 2$. We have $d^-_{T_2}(x_1)\geq 3$, let $w_1,w_2,w_3\in N^-_{T_2}(x_1)$ such that$w_2\rightarrow w_3$ and let $P$ be an antidirected outpath in $T_2\setminus \{w_1,w_2,w_3\}$, we have $w_3w_2x_1w_1x_2Pvw_3\equiv C$.\\
\textbf{Case $4$}: If $\delta\geq 3$, $\delta$ is odd and, without loss of generality, $T_1\neq 3A,5A,7A$. Let $x_1...x_{\delta}$ be a Hamiltonian antidirected outpath of $T_1$ such that $x_{\delta}\rightarrow x_1$ and $d^-_{T_2}(x_2)=max\,\{d^-_{T_2}(y_2);y_1...y_{\delta}$ is an antidirected outpath in $T_1$ with $y_{\delta}\rightarrow y_1\}$. We distinguish two cases:\begin{itemize}
\item[(a)] If $d^+_{T_2}(x_1)\geq 3$. Suppose first that $d^+_{T_2}(x_{\delta})\geq 1$ and let $u\in N^+_{T_2}(x_{\delta})$. Obviously $d^+_{T_2-u}(x_1)\geq 2$ so there is an outneighbor of $x_1$ origin of an antidirected inpath $P$ in $T_2-u$ and thus $x_1...x_{\delta}uvP^{-1}x_1\equiv C$; unless $(T_2-u;-(1,1,...,1))$ is one of the following biexceptions with respect to $x_1$:\\
Dual Exc $(4,1)$: If there exists $u'\in N^+_{T_2-u}(x_1)$ such that $u$ is an origin of an antidirected inpath $P$ in $T_2-u'$ then $u'x_1...x_{\delta}Pvu'\equiv C$, so we may assume that $u\rightarrow \{1,5\}$. If $3\rightarrow x_2$, $3$ is an origin of an antidirected outpath $P$ in $T_2\setminus \{4,5\}$ then $45x_1x_{\delta}...x_2Pv4\equiv C$. Otherwise, $4\rightarrow x_2$ and we may assume that $u\rightarrow 3$, $3$ is an origin of an antidirected outpath $P$ in $T_2\setminus \{2,4\}$ so $24x_2...x_{\delta}x_1Pv2\equiv C$.\\
Dual Exc $(8,1)$: If $2\rightarrow u$ then $5x_1...x_{\delta}u21v345\equiv C$ so $u\rightarrow 2$. If $3\rightarrow u$ then $5x_1...x_{\delta}u3124v5\equiv C$ and if $4\rightarrow u$ then $2x_1...x_{\delta}u4351v2\equiv C$, so $u\rightarrow \{3,4\}$. If $1\rightarrow x_2$, if $u\rightarrow 5$ then $23x_1x_{\delta}...x_214u5v2\equiv C$ and if $5\rightarrow u$ then $41x_2...x_{\delta}u5x_132v4\equiv C$. And if $5\rightarrow x_2$ then $15x_2...x_{\delta}x_132u4v1\equiv C$.\\
Dual Exc $(9,1)(1)$: If $1\rightarrow u$ then $2x_1...x_{\delta}u2453v2\equiv C$, if $2\rightarrow u$ then $4x_1...x_{\delta}u21v354\equiv C$, if $3\rightarrow u$ then $4x_1...x_{\delta}u31v524\equiv C$, if $4\rightarrow u$ then $5x_1...x_{\delta}u43v125\equiv C$ and if $5\rightarrow u$ then $4x_1...x_{\delta}u53v124\equiv C$, so $u\rightarrow \{1,2,3,4,5\}$. If $1$ or $2\rightarrow x_2$, say $1$, then $23x_1x_{\delta}...x_214u5v2\equiv C$ and if $4\rightarrow x_2$ then $51x_1x_{\delta}...x_243u2v5\equiv C$.\\
Dual Exc $(9,1)(2)$: If $4\rightarrow u$ then $5x_1...x_{\delta}u43v125\equiv C$ and if $5\rightarrow u$ then $2x_1...x_{\delta}u54v132\equiv C$, so $u\rightarrow \{4,5\}$ and thus $23x_1x_{\delta}...x_214u5v2\equiv C$.\\
Dual Exc $(9,1)(3),(4)$: If $2\rightarrow u$ then $4x_1...x_{\delta}u21v354\equiv C$ and if $5\rightarrow u$ then $4x_1...x_{\delta}u53v124\equiv C$, so $u\rightarrow \{2,5\}$ and thus $34x_2...x_{\delta}x_115u2v3\equiv C$.\\
Dual Exc $(33,1)$: In $\overline{7A}$, any vertex $i$ is an origin of a path $P_i=-(1,1,1,1,2)$ of end $i+1$, so if there exists $i\in \{1,4,6,7\}$ such that $i\rightarrow u$ then $(i+1)x_1...x_{\delta}uP_i^*v(i+1)\equiv C$, so $u\rightarrow \{1,4,6,7\}$. Thus $62x_2...x_{\delta}x_14537v1u6\equiv C$.\\
Dual Exc $(E_8(5),1)$: Set $X=\{u_1,u_2\}$. If $1\rightarrow u$ then $u_1x_1...x_{\delta}u13u_22vu_1\equiv C$ so $u\rightarrow \{1,2,3\}$. Without loss of generality suppose that $u_1\rightarrow x_2$ and $u_2\rightarrow u$, we have $31x_1x_{\delta}...x_2u_12u_2uv3\equiv C$.\\
Or $(T_2-u;(1,1,..,1))$ is an exception with $N^-_{T_2}(x_2)=N^+_{T_2}(x_{\delta})=\{u\}$, so $T_1$ is a transitive triangle (easy to prove) and $N^+_{T_2}(x_{\delta})=\{u\}$. If $u$ is an origin of an antidirected inpath $P$ in $T_2-z$ for a $z\in N^+(x_1)$ then $zx_1x_2x_3Pvz\equiv C$, otherwise, let $u_1,u_2\in N^+_{T_2-u}(x_1)$ such that $u_1\rightarrow u_2$, if $T_2\setminus \{u_1,u_2\}$ contains an antidirected outpath $P$ of origin and end in $T_2-u$ then $u_2x_1vx_2P^{-1}x_3u_1u_2\equiv C$, otherwise, $T_2\setminus \{u_1,u_2\}=4A$ and $y_1...y_4$ is the sole antidirected outpath in $T_2\setminus \{u_1,u_2\}$. If $u=y_4$, since $u$ is not an origin of an antidirected inpath in $T_2-u_1$ then $u_1\rightarrow y_1$ and thus $u_2x_1x_2x_3y_4y_3y_2vy_1u_1u_2\equiv C$. Else, $u=y_1$, since $u_2u_1y_4y_3y_2$ is not antidirected $y_4\rightarrow u_1$. If $y_3\rightarrow u_2$ then $u_2x_1x_2x_3y_1y_4u_1vy_2y_3u_2\equiv C$ so $u_2\rightarrow y_3$ and thus $u_1x_1vx_2y_2y_1y_3u_2x_3y_4u_1\equiv C$.\\
Assume now that $N^-(x_{\delta})=T_2$, let $a\in N^+_{T_2}(x_1)$, $b\in N^-_{T_2}(x_{\delta-1})$ and $P$ be an antidirected oupath in $T_2\setminus \{a,b\}$, we have $ax_1...x_{\delta-1}bx_{\delta}Pva\equiv C$.

\item[(b)] If $d^+_{T_2}(x_1)\leq 2$ i.e. $d^-_{T_2}(x_1)\geq |T_2|-2$. Let $u\in N^-_{T_2}(x_2)$, if $d^+_{T_2}(u)=0$ and $N^-_{T_2}(x_2)=\{u\}$, let $a\in N^-_{T_2-u}(x_1)$ and $P$ be an antidirected inpath in $T_2\setminus \{a,u\}$, we have $x_1x_{\delta}...x_3vx_2Puax_1\equiv C$. So we may suppose that there exists $w\in N^+_{T_2}(u)$, if $x_1$ has an inneighbor $y$ origin of an antidirected outpath $P$ in $T_2\setminus \{u,w\}$ then $yx_1x_{\delta}...x_2uwvP^{-1}\equiv C$, otherwise, either $(T_2\setminus \{u,w\};(1,1,...,1))$ is one of the following exceptions:\\
Exc $1$: If $V(T_2\setminus \{u,w\})=\{v_1,v_2,v_3,v_4\}$ such that $v_2\in N^-(x_1)\cap N^-(v_1)$ and $v_4\in N^-(w)\cap N^-(v_3)$ then $x_1x_{\delta}...x_2uwv_3v_4vv_1v_2x_1\equiv C$. So we may assume that $w\rightarrow 4$. Suppose first that $\{1,2\}\subseteq N^-(x_1)$, so $w\rightarrow 3$. If $u\rightarrow 1$ or $2$, say $1$, then $x_1x_{\delta}...x_2u1v4w32x_1\equiv C$, so $\{1,2\}\rightarrow u$. If $x_2\rightarrow 4$ then $x_1x_{\delta}...x_3vx_24w32u1x_1\equiv C$ so $4\rightarrow x_2$. Without loss of generality suppose that $u\rightarrow 4$, since $1$ is an origin of an antidirected outpath $P$ in $T_2\setminus \{u,4\}$, then $x_1x_{\delta}..x_2u4vP^{-1}x_1\equiv C$. So we may assume that $|N^-(x_1)\cap \{1,2\}|=1$, and thus $\{u,w,3\}\subseteq N^-(x_1)$. Suppose that $1\in N^-(x_1)$, so $w\rightarrow \{1,3\}$, if $u\rightarrow 2$ then $x_1x_{\delta}...x_2u2v4w31x_1\equiv C$, so $2\rightarrow u$. If $u\rightarrow 3$ then $x_1x_{\delta}...x_2u3v241wx_1\equiv C$, so $3\rightarrow u$. If $u\rightarrow 4$ then $x_1x_{\delta}...x_2u4v213wx_1\equiv C$, so $4\rightarrow u$. If $1\rightarrow x_2$ then $x_1x_{\delta}...x_212vu34wx_1\equiv C$ and if $x_2\rightarrow 1$ then $x_1x_{\delta}...x_3vx_214u23wx_1\equiv C$. So we may assume that $N^+_{T_2}(x_1)=\{1,4\}$ and thus $4\rightarrow \{1,4\}$. If $u\rightarrow 2$, without loss of generality suppose that $w\rightarrow 3$, then $x_1x_{\delta}...x_2u241v3wx_1\equiv C$, so $2\rightarrow u$. If $u\rightarrow 3$ then $x_1x_{\delta}...x_2u3v241wx_1\equiv C$, so $3\rightarrow u$. If $1\rightarrow x_2$, since $d^-_{T_2\setminus \{1,2\}}(4)\geq 2$ then $T_2\setminus \{1,2\}$ contains an antidirected outpath of origin distinct from $4$, thus $x_1x_{\delta}...x_212vP^{-1}x_1\equiv C$ and if $x_2\rightarrow 1$ then $x_1x_{\delta}...x_3vx_21w43u2x_1\equiv C$.\\
Exc $24$: If $d^-_{T_2}(w)\geq 2$, one can find and inneighbor $u_1$ of $w$ in $T_2-u$ having an outneighbor $u_2$ in $T_2\setminus \{u,w\}$ such that $T_2\setminus \{u,w,u_1,u_2\}$ contains an antidirected outpath $P$ of origin an inneighbor of $x_1$ and thus $x_1x_{\delta}...x_2uwu_1u_2vP^{-1}x_1\equiv C$, so we may assume that $d^-_{T_2}(w)=1$. If $d^+_{T_2}(u)\geq 2$, let $w'\in N^+_{T_2-w}(u)$, $T_2\setminus \{u,w'\}$ contains an antidirected outpath $P$ of origin an inneighbor of $x_1$ and thus $x_1x_{\delta}...x_2uw'vP^{-1}x_1\equiv C$, so $d^+_{T_2}(u)=1$. We may suppose that $x_2\rightarrow (T_2-u)$, let $a\in N^-_{T_2\setminus \{u,w\}}(x_1)$, $b\in N^+_{T_2\setminus \{u,w\}}(x_2)$ and $P$ be an antidirected inpath in $T_2\setminus \{u,w,a,b\}$, we have $x_1x_{\delta}...x_2bwPuax_1\equiv C$.\\
Or $T_2\setminus \{u,w\}=\{u_1\rightarrow u_2\}$ such that $u_1\in N^+(x_1)$. Suppose first that $u\rightarrow u_1$, then $w\rightarrow u_2$ and $x_1\rightarrow w$. And since $d^+(u)=4$ then $x_3\rightarrow u$ and thus, without loss of generality suppose that $w\rightarrow u_1$, $x_1x_2x_3uvu_2wu_1x_1\equiv C$. So we may assume that $u_1\rightarrow u$. If $u_1\rightarrow x_2$, if $u\rightarrow x_1$ then $x_1x_3x_2u_1u_2vwux_1\equiv C$ and if $w\rightarrow x_1$ then, without loss of generality suppose that $w\rightarrow u_2$, we have $x_1x_3x_2u_1uvu_2wx_1\equiv C$, so we can assume that $x_2\rightarrow u_1$. If $T_2\rightarrow x_3$, without loss of generality we may assume that $u_2\rightarrow w$, thus $x_1x_2ux_3u_2wvu_1x_1\equiv C$. So we may assume that $d^+_{T_2}(x_3)=1$. Suppose that $x_3\rightarrow w$, if $u\rightarrow u_2$ then $x_1x_2x_3wuu_2vu_1x_1\equiv C$, so $u_2\rightarrow u$. If $x_1\rightarrow u$ then $x_1x_2x_3wvu_2u_1ux_1\equiv C$, so $u\rightarrow x_1$. If $u_2\rightarrow w$ then $x_1x_2x_3wu_2uvu_1x_1\equiv C$, so $w\rightarrow u_2$. If $x_1\rightarrow w$ then $u_2u_1x_3ux_2x_1wvu_2\equiv C$, so $w\rightarrow x_1$. If $u_1\rightarrow w$ then there exists an outneighbor of $x_2$ in $\{u_2,w\}$ origin of an antidirected inpath $P$ of end in $\{u,w,u_2\}$ and thus $x_1x_3vx_2Px_1\equiv C$, so $w\rightarrow u_1$ and thus $x_2\rightarrow \{u_2,w\}$. Therefore $uu_2x_1x_3wx_2u_1vu\equiv C$, so we can assume that $w\rightarrow x_3$. Suppose that $x_3\rightarrow u_1$, if $u_2\rightarrow x_2$ then $wux_2u_2x_1x_3u_1vw\equiv C$, so $x_2\rightarrow u_2$. If $u\rightarrow x_1$ then $x_2u_2u_1wux_1x_3vx_2\equiv C$, so $x_1\rightarrow u$ and thus $w\rightarrow x_1$. If $x_2\rightarrow w$ then $w\rightarrow u_2$ and $u_2\rightarrow u$, so $x_1x_3vx_2wu_1uu_2x_1\equiv C$, so $w\rightarrow x_2$. If $w\rightarrow u_2$ then $x_1vx_3x_2wu_2u_1ux_1\equiv C$, so $u_2\rightarrow w$. If $u_2\rightarrow u$ then $x_1x_2x_3u_1vwu_2ux_1\equiv C$, so $u\rightarrow u_2$. Therefore $u_2ux_2wx_1x_3u_1vu_2\equiv C$. Now suppose that $x_3\rightarrow u_2$, if $x_1\rightarrow w$ then $x_1x_2x_3u_2u_2uvwx_1\equiv C$, so $w\rightarrow x_1$. If $u_1\rightarrow w$ then $x_2x_1vx_3u_2u_1wux_3\equiv C$, so $w\rightarrow u_1$ then $x_2\rightarrow \{w,u_2\}$ and thus $u_2\rightarrow w$; if $u\rightarrow x_1$ then $wux_1x_3u_2x_2u_1vw\equiv C$, so $x_1\rightarrow u$ and then $T-v=7A$, therefore $(T;C)=A_6$. Suppose now that $x_3\rightarrow u$. We have $2$ cases:\begin{itemize}
\item[(i)] If $u_1\rightarrow w$. If $x_1\rightarrow w$ then $x_1x_2x_3uvu_2u_1wx_1\equiv C$, so $w\rightarrow x_1$. Suppose that $u_2\rightarrow u$ then $u\rightarrow x_1$; if $u_2\rightarrow w$ then $x_1x_2x_3uu_2wvu_1x_1\equiv C$, so $w\rightarrow u_2$. If $x_2\rightarrow u_2$ then $x_1x_3vx_2u_2u_1wux_1\equiv C$, so $u_2\rightarrow x_2$ and thus $x_2\rightarrow w$, therefore $x_2vx_3x_1u_2uu_1wx_2\equiv C$. So we may assume that $u\rightarrow u_2$. If $w\rightarrow u_2$, then $u_2\rightarrow x_2$ and thus $x_2\rightarrow w$. If $u\rightarrow x_1$ then $x_2vx_3x_1uu_2u_1wx_3\equiv C$, and if $x_1\rightarrow u$ then $u_2x_3u_1wvux_1x_2u_2\equiv C$. So we can suppose now that $u_2\rightarrow w$. Then $w\rightarrow x_2$ and thus $x_2\rightarrow u_2$, if $u\rightarrow x_1$ then $x_2vx_3x_1uwu_1u_2x_2\equiv C$ and if $x_1\rightarrow u$ then $wx_3u_1u_2vux_1x_2w\equiv C$.
\item[(ii)] If $w\rightarrow u_1$. If $w\rightarrow u_2$ then $x_1x_2x_3uvu_2wu_1x_1\equiv C$, so $u_2\rightarrow w$. If $u_2\rightarrow u$ then $x_1x_2x_3uu_2wvu_1x_1\equiv C$, so $u\rightarrow u_2$. And thus $x_1\rightarrow w$, therefore $x_1x_2x_3uu_1u_2vwx_1\equiv C$.
\end{itemize}
\end{itemize}
\textbf{Case $5$}: If $\delta=3$ and for all $x\in T$ such that $d^-(x)=3$ we have $T(N^-(x))=3A$. Set $V(T_1)=\{y_1,y_2,y_3\}$. If $d^+_{T_1}(x)\geq 2$ for any $x\in T_2$.  If $|T_2|=4$ then $d^+_{T_2}(x)\leq 2$ for all $x\in T_2$, and so $T_2=4A$. Also, $d^+_{T_1}(x)=1$ for all $x\in T_2$. Without loss of generality we may suppose that $y_1\rightarrow 3$. If $4\rightarrow y_3$ then $y_1vy_2y_34213y_1\equiv C$ and if $y_3\rightarrow 4$ then $4\rightarrow \{y_1,y_2\}$. Since $T(N^+(2))=3A$ then $N^+(2)=\{3,y_1,y_3\}$. Thus $T-1=7A$ and $(T;C)=A_6$. Now if $|T_2|\geq 6$, without loss of generality suppose that $d^+_{T_2}(y_1)=max\{d^+_{T_2}(y_i),\,i\in \{1,2,3\}\}$. Let $a\in N^-_{T_2}(y_2)\cap N^-_{T_2}(y_3)$ and set $T'_1=T_1\cup \{a\}$ and $T'_2=T_2-a$. We have $x_1x_2x_3x_4=y_3ay_2y_1=-(1,1,1)$. If $d^+_{T'_2}(x_4)=0$ then there exists $a'\in T'_2$ such that $d^+_{T_1}(a')=3$, by taking $x_1a'x_3x_4\equiv x_1...x_4$ we obtain $d^+_{T_2-a'}(x_4)=1$. So we may suppose that there exists $u\in N^+_{T'_2}(x_4)$ with $d^-_{T'_2-u}(x_1)\geq 2$. In fact, otherwise, $d^+_{T_2}(x_1)\geq |T_2|-3$ with $N^+_{T_2}(x_1)\subseteq N^-_{T_2}(x_4)$, thus $d^-_{T_2}(x_4)\geq |T_2|-3$ and so $d^+_{T_2}(x_4)\leq |T_2|-(|T_2|-3)$, then $|T_2|-3\leq 3$, but $|T_2|\geq 6$, so $|T_2|=6$; let $a',u'\in N^+_{T_2}(x_1)$, we have $x_3a'x_4x_1=-(1,1,1)$ with $d^-_{T_2\setminus \{a',u'\}}(x_3)\geq 2$. So if $x_1$ has an inneighbor origin of an antidirected outpath $P$ in $T'_2-u$ then $x_1...x_4uvP^{-1}x_1\equiv C$. Otherwise, $(T'_2-u;(1,1,...,1))$ is an exception, and it will be treated below. Otherwise, there exist $y\in T_2$ such that $d^+_{T_1}(y)\leq 1$, suppose that, without loss of generality, $\{y_1,y_3\}\rightarrow y$. Set $T'_1=T_1\cup \{y\}$ and $T'_2=T_2-y$. We have $x_1...x_4=y_2y_1yy_3=-(1,1,1)$. Since $T(x_2,x_3,v)\neq 3A$ then $d^+_{T'_2}(x_4)\geq 1$, so let $u\in N^+_{T'_2}(x_4)$. If $x_1$ has an inneighbor origin of an antidirected outpath $P$ in $T'_2-u$ then $x_1...x_4uvP^{-1}x_1\equiv C$. Otherwise, if $d^-_{T'_2-u}(x_1)\leq 1$, if $|T_2|=4$, set $T_2\setminus \{y,u\}=\{w_1\rightarrow w_2\}$. We have $T(N^-(y_3))=T(y_2,w_1,w_2)=3A$. If $y_1\rightarrow w_1$ then $T(y_3,w_2,u)=3A$ and so $u\rightarrow w_2$. Since $T(w_1,u,v)\neq 3A$ then $y\rightarrow w_2$. Similarly, we find that $T-v=7A$ and $(T;C)=A_6$. So we may suppose that $w_1\rightarrow y_1$. Similarly, we find that $T-v=7A$ and $(T;C)=A_6$. Now if $|T_2|\geq 6$, two cases may arise:\begin{itemize}
\item[(i)] If $N^-_{T'_2-u}(x_1)=\{w\}$. Suppose first that $x_1\rightarrow u$, then $N^-(x_1)=\{w,x_2,x_3\}$ and thus $w\in N^-(x_2)$. We have $x_2x_4ux_1=-(1,1,1)$, if $d^-_{T_2-u}(x_2)\geq 2$ then there exists $u'\in N^+_{T_2-u}(x_1)$ such that $d^-_{T_2\setminus \{u,u'\}}(x_2)\geq 2$ and the result follows, so $N^-_{T_2}(x_2)=\{u,w\}$ and thus $w\in N^-(x_4)$. We have $T'_2\setminus \{u,w\}\,\subseteq N^+(x_1)\cap N^+(x_2)$ with $|T'_2\setminus \{u,w\}|\geq 3$ and $d^-_{T'_2\setminus \{u,w\}}(x_4)\geq 1$ so we may find $u"\in T'_2\setminus \{u,w\}$ such that $d^-_{T_2\setminus \{u",x_3\}}(x_4)\geq 2$ with $x_4x_1u"x_2=-(1,1,1)$ and the result follows. So we may assume that $u\rightarrow x_1$, and so $N^+_{T'_2}(x_1)\subseteq N^-(x_4)$. Set $T^+=T(N^+_{T'_2}(x_1))$. If $x_2$ has an outneighbor $a\in T^+$ then $x_4x_1ax_2=-(1,1,1)$ with $x_2\rightarrow x_3$ and $d^-_{T_2\setminus \{a,x_3\}}(x_4)\geq 2$, so the result follows. Thus $T^+\subseteq N^-(x_2)$. If there exists an antidirected outpath $P$ in $T^+$ such that $x_1$ has an outneighbor origin of an antidirecrted inpath $P$ of order $|T^+|-1$ in $T'_2V(P)$ then $x_3x_4x_2u_1u_2x_1Pvx_3\equiv C$. Otherwise, suppose that $u$ has an outneighbor $u_1\in T^+$, then $T(u,u_1,w)=3A$ and thus $T^+\rightarrow w$ and $w\rightarrow u$. If $u_1$ has an inneighbor $u_2\in T^+$, let $P$ be an antidirected outpath in $T'_2\setminus \{u,u_1,u_2\}$, then $x_3x_2x_1uu_1u_2x_4Pvx_3\equiv C$. Else, $N^+_{T^+}(u)=\{u_1\}$, if $|T^+|=5$ or $7$, let $P$ be an antidirected outpath in $T^+-u_1$ of order $|T^+|-3$, there exists and antidirecte outpath $Q$ in $T'_2-V(P)$ of origin in $T^+$, thus $x_3x_4x_2Px_1Qvx_3\equiv C$. If $|T^+|=3$, set $\{u_2\rightarrow u_3\}=T^+-u_1$, if $u\rightarrow x_2$ then $x_3x_4x_2uu_1x_1u_3u_2wvx_3\equiv C$, and if $x_2\rightarrow u$, $x_3$ has an outneighbor $u_i$ distinct from $u_1$, in fact, $d^+_{T'_2}(x_3)\geq 2$, so if $\{u_1,w\}\subseteq N^+(x_3)$ then $ux_4x_2u_2u_3x_1u_1x_3wvu\equiv C$ and if $N^+_{T'_2}(x_3)=\{u,w\}$ or $N^+_{T'_2}(x_3)=\{u,u_1\}$ then $T(N^+(x_3))\neq 3A$, so without loss of generality we may suppose that $x_3\rightarrow u_2$ and thus $ux_2x_1x_3u_2u_1x_4u_3wvu\equiv C$. So we can suppose that $T^+\rightarrow u$, if $w$ has an outneighbor $u_1\in T^+$, let $P$ be an antidirected outpath in $T^+-u_1$, we have $x_1...x_4uPvu_1wx_1\equiv C$, so $T^+\rightarrow w$. Since $T(x_1,x_2,x_3)\neq 3A$ then $u\rightarrow w$. If $|T^+|=5$ or $7$, let $P$ be an antidirected outpath in $T^+$ of order $|T^+|-3$, there exists an antidirected outpath $Q$ in $T'_2-V(P)$ of origin in $T^+$, thus $x_3x_4x_2Px_1Qvx_3\equiv C$. Now if $|T^+|=3$, set $V(T^+)=\{u_1,u_2,u_3\}$ such that $u_1\rightarrow u_2$, since $T(w,x_1,x_3)\neq 3A$ then $x_2\in N^+(u)$, thus $x_3x_4x_2uwu_1u_2x_1u_3vx_3\equiv C$.
\item[(ii)] If $N^-(x_1)=\{x_2,x_3,u\}$, then $x_3\rightarrow u$ and $u\rightarrow x_2$. We may assume that $d^+_{T'_2}(x_4)=1$ and thus $d^+_{T'_2}(x_2)=0$. There exist $a,b\in T'_2-u$ such that $x_2$ has an inneighbor origin of an antidirected outpath $P$ in $T_2\setminus \{a,b\}$, thus $ax_1x_4bx_2Pva\equiv C$.
\end{itemize}
So we may assume that $d^-_{T'_2-u}(x_1)\geq 2$ and then $(T'_2-u;(1,1,...,1))$ is an exception. For the both cases, if $d^-_{T'_2}(x_1)\geq 2$ and $(T'_2-u;(1,1,...,1))$ is an exception then it is one of the following ones:\\
Exc $1$: If $d^+_{T_1}(x)\geq 2$ for all $x\in T_2$, $3$ cases may arise:\begin{itemize}
\item[(i)] If $\{1,2\}\subseteq N^-(x_1)$, by the building lemmas we can suppose that $\{2,4\}\rightarrow x_2$. If $N^-(u)\cap \{3,4\}\neq \emptyset$ one can find $C$, so $u\rightarrow \{3,4\}$. So by supposing that $N^+(x_4)\cap \{1,2,3,4\}\neq \emptyset$ we can find $C$, so $\{1,2,3,4\}\rightarrow x_4$. Thus $d^-_{T'_2}(x_3)\geq 3$. If $x_4\rightarrow x_2$ then $x_1ux_3x_4x_242v31x_1\equiv C$ and if $x_2\rightarrow x_4$, since $T'_2-4\neq 4A$ and $d^-_{T'_2-4}(x_3)\geq 2$ then $x_3$ has an inneighbor origin of an antidirected outpath $P$ in $T'_2-4$ and thus $x_3x_2x_4x_14vP^{-1}x_3\equiv C$.
\item[(ii)] If $N^-_{T'_2}(x_1)=\{1,3\}$, by the building lemmas we can suppose that $\{2,3,4\}\rightarrow x_2$. By supposing that there exists $i\in \{1,3,4\}\cap N^-(u)$ one can find $C$, so we may suppose that $u\rightarrow \{1,3,4\}$. If $x_4\rightarrow 2$ then $x_1...x_42v4u31x_1\equiv C$, if $x_4\rightarrow 3$ then $x_1vx_2x_3x_43u142x_1\equiv C$ and if $x_4\rightarrow 4$ then $x_1vx_2x_3x_44u312x_1\equiv C$, so $\{2,3,4\}\rightarrow x_4$. If $x_4\rightarrow x_2$ then $x_1ux_3x_4x_242v31x_1\equiv C$ and if $x_2\rightarrow x_4$, since $T'_2-4\neq 4A$ and $d^-_{T'_2-4}(x_3)\geq 2$ then $x_3$ has an inneighbor origin of an antidirected outpath $P$ in $T'_2-4$ and thus $x_3x_2x_4x_14vP^{-1}x_3\equiv C$.
\item[(iii)] If $N^-_{T'_2}(x_1)=\{2,3\}$. By supposing that there exists $i\in \{1,4\}\cap N^-(u)$ one can find $C$, so we may assume that $u\rightarrow \{1,4\}$. If there exists $u'\in \{1,4\}\cap N^+(x_4)$ then $d^-_{T'_2-u'}(x_1)\geq 3$ and the result follows, and if there exists $u'\in \{2,3\}\cap N^+(x_4)$ then $d^-_{T'_2-u'}(x_1)\geq 2$ with an inneighbor of outdegree at least $2$ in $T'_2-u'$, so the result follows. Thus $\{1,2,3,4\}\rightarrow x_4$. We have $y_1...y_4=x_31x_4x_1=-(1,1,1)$ with $d^-_{T_2\setminus \{1,4\}}(y_1)\geq 2$, so the result follows unless $T(x_2,u,2,3)=4A$ with $N^+_{T_2\setminus \{1,4\}}(y_1)\cap N^-(4)\neq \emptyset$, so the result follows.
\end{itemize}
Otherwise, we distinguish $3$ cases:\begin{itemize}
\item[(i)] If $\{1,2\}\subseteq N^-(x_1)$. If $3\rightarrow u$ then $x_1...x_4u34v21x_1\equiv C$ and if $4\rightarrow u$ then $x_1...x_4u41v32x_1\equiv C$, so $u\rightarrow \{3,4\}$. If $x_4\rightarrow 1$ or $2$, say $1$, then $x_1...x_41v4u32x_1\equiv C$, so $\{1,2\}\rightarrow x_4$ and thus $x_3x_2x_11x_423u4vx_3\equiv C$.
\item[(ii)] If $N^-_{T'_2-u}(x_1)=\{1,3\}$, by the building lemmas we can suppose that $\{2,3,4\}\subseteq N^-(x_2)$. By supposing that there exists $i\in \{1,3,4\}\cap N^-(u)$ one can find $C$, so we may suppose that $u\rightarrow \{1,3,4\}$. If $x_4\rightarrow 2$ then $x_1...x_42v4u31x_1\equiv C$, if $x_4\rightarrow 3$ then $x_1vx_2x_3x_43u142x_1\equiv C$ and if $x_4\rightarrow 4$ then $x_1vx_2x_3x_44u312x_1\equiv C$, so $\{2,3,4\}\rightarrow x_4$. If $x_1\rightarrow x_3$, since $T(x_1,x_3,v)\neq 3A$ then $x_2$ has an outneighbor $x\in \{1,u\}$, but $x_4$ has an inneighbor origin of an antidirected outpath $P$ in $T'_2-x$ so $x_4x_1x_3x_2xvP^{-1}x_4\equiv C$, so $x_3\rightarrow x_1$. If $4\rightarrow x_3$, if $2\rightarrow u$ then $14x_3x_2x_13x_42uv1\equiv C$ and if $u\rightarrow 2$ then $x_3x_2x_13x_441u2vx_3\equiv C$, thus $x_3\rightarrow 4$ and so $ux_4x_234x_3x_113vu\equiv C$.
\item[(iii)] If $N^-_{T'_2-u}(x_1)=\{2,3\}$. By supposing that there exists $i\in \{1,4\}\cap N^-(u)$ one can find $C$, so we may assume that $u\rightarrow \{1,4\}$. If $x_4\rightarrow 2$ then $x_1...x_42v1u43x_1\equiv C$, so $2\rightarrow x_4$. If $3\rightarrow x_4$ then $x_3x_2x_12x_434u1vx_3\equiv C$, so $x_4\rightarrow 3$. If $u\rightarrow x_1$ then $241ux_1...x_43v2\equiv C$, so $x_1\rightarrow u$. If $2\rightarrow x_2$ then $x_3x_4x_22x_134u1vx_3\equiv C$ and if $x_2\rightarrow 2$ then $312x_2x_3x_4ux_14v3\equiv C$.
\end{itemize}
Exc $20$: One can check that for any $i\in \{3,4,5,6\}$, $1$ or $2$ is an origin of an antidirected outpath $P$ in $T'_2\setminus \{i,i+1,u\}$, so if there exists $i\in N^-(u)\cap \{3,4,5,6\}$ then $x_1...x_4ui(i+1)vP^{-1}x_1\equiv C$, thus $u\rightarrow \{3,4,5,6\}$. If there exists $u'\in N^+_{T'_2-u}(x_4)$, $1$ or $2$ is an origin of an antidirected ouptath $P$ in $T'_2-u'$ and thus $x_1...x_4u'vP^{-1}x_1\equiv C$, so $N^+_{T'_2}(x_4)=\{u\}$. If $d^-_{T'_2}(x_3)\geq 1$ then there exist $w\in N^-_{T'_2}(x_3)$, $w'\in N^+_{T'_2}(w)$ and $a\in \{1,2\}$ pairwaisly distinct such that $x_4$ has an inneighbor origin of an antidirected outpath $P$ in $T'_2\setminus \{a,w,w'\}$, thus $x_4ax_1...x_3ww'vP^{-1}x_4\equiv C$. Otherwise, since $T(x_1,x_3,v)\neq 3A$ then $d^+_{T'_2}(x_2)\geq 1$, thus there exist $x\in N^+_{T'_2}(x_2)$ and $P$ an antidirected outpath $P$ in $T'_2-x$ of origin distinct from $u$, therefore $x_4x_1x_3x_2xvP^{-1}x_4\equiv C$.\\
Exc $22$: Without loss of generality we can suppose that $\{1,2\}\subseteq N^-(x_1)$. One can check that for any $i\in \{4,5,6\}$, there exists an outneighbor $j$ of $i$ such that $1$ or $2$ is an origin of an antidirected outpath $P$ in $T'_2\setminus \{i,j,u\}$, so if there exists $i\in N^-(u)\cap \{1,2,3,4,5\}$ then $x_1...x_4uijvP^{-1}x_1\equiv C$, thus $u\rightarrow \{4,5,6\}$. If there exists $u'\in N^+_{T'_2-u}(x_4)$ distinct from $4$, $1$ or $2$ is an origin of an antidirected ouptath $P$ in $T'_2-u'$ and thus $x_1...x_4u'vP^{-1}x_1\equiv C$, and if $4\in N^+(x_4)$ then $x_1...x_44536v21x_1\equiv C$, so $N^+_{T'_2}(x_4)=\{u\}$. If $d^-_{T'_2}(x_3)\geq 1$, then there exist $w\in N^-_{T'_2}(x_3)$, $w'\in N^+_{T'_2}(w)$ and $a\in \{4,6\}$ pairwaisly distinct such that $x_4$ has an inneighbor origin of an antidirected outpath $P$ in $T_2\setminus \{a,w,w'\}$, thus $x_4ax_1...x_3ww'vP^{-1}x_3\equiv C$. Otherwise, since $T(x_1,x_3,v)\neq 3A$ then $d^+_{T'_2}(x_2)\geq 1$, thus there exist $x\in N^+_{T'_2}(x_2)$ and $P$ an antidirected outpath $P$ in $T'_2-x$ of origin distinct from $u$, therefore $x_4x_1x_3x_2xvP^{-1}x_4\equiv C$.\\
Exc $23$: One can check that for any $i\in \{1,2,3,4,5\}$, there exists an outneighbor $j$ of $i$ such that $4$ or $6$ is an origin of an antidirected outpath $P$ in $T'_2\setminus \{i,j,u\}$, so if there exists $i\in N^-(u)\cap \{1,2,3,4,5\}$ then $x_1...x_4uijvP^{-1}x_1\equiv C$, thus $u\rightarrow \{1,2,3,4,5\}$. If there exists $u'\in N^+_{T'_2-u}(x_4)$, $4$ or $6$ is an origin of an antidirected ouptath $P$ in $T'_2-u'$ and thus $x_1...x_4u'vP^{-1}x_1\equiv C$, so $N^+_{T'_2}(x_4)=\{u\}$. If $d^-_{T'_2}(x_3)\geq 1$, then there exist $w\in N^-_{T'_2}(x_3)$, $w'\in N^+_{T'_2}(w)$ and $a\in \{4,6\}$ pairwaisly distinct such that $x_4$ has an inneighbor origin of an antidirected outpath $P$ in $T_2\setminus \{a,w,w'\}$, thus $x_4ax_1...x_3ww'vP^{-1}x_4\equiv C$. Otherwise, since $T(x_1,x_3,v)\neq 3A$ then $d^+_{T'_2}(x_2)\geq 1$, thus there exist $x\in N^+_{T'_2}(x_2)$ and $P$ an antidirected outpath $P$ in $T'_2-x$ of origin distinct from $u$, therefore $x_4x_1x_3x_2xvP^{-1}x_4\equiv C$.\\
Exc $24$: Without loss of generality we may suppose that $1\in N^-(x_1)$. One can check that for any $i\in \{2,4,5,6\}$, there exists an outneighbor $j$ of $i$ such that $1$ is an origin of an antidirected outpath $P$ in $T'_2\setminus \{i,j,u\}$, so if there exists $i\in N^-(u)\cap \{2,4,5,6\}$ then $x_1...x_4uijvP^{-1}x_1\equiv C$, thus $u\rightarrow \{2,4,5,6\}$. If $3\rightarrow u$ then $x_1...x_4u3456v21x_1\equiv C$, so $u\rightarrow 3$. If there exists $u'\in N^+_{T'_2-u}(x_4)$, $1$ is an origin of an antidirected ouptath $P$ in $T'_2-u'$ and thus $x_1...x_4u'vP^{-1}x_1\equiv C$, so $N^+_{T'_2}(x_4)=\{u\}$. If $d^-_{T'_2}(x_3)\geq 1$, then there exist $w\in N^-_{T'_2}(x_3)$, $w'\in N^+_{T'_2}(w)$ and $a\in N^-_{T'_2}(x_1)\cap N^-_{T'_2}(x_4)$ pairwaisly distinct such that $x_4$ has an inneighbor origin of an antidirected outpath $P$ in $T'_2\setminus \{a,w,w'\}$, thus $x_4ax_1...x_3ww'vP^{-1}x_4\equiv C$.\\
Exc $29$ :One can check that for any $i\in \{2,3,5\}$, there exists an outneighbor $j$ of $i$ such that $4$ or $6$ is an origin of an antidirected outpath $P$ in $T'_2\setminus \{i,j,u\}$, so if there exists $i\in N^-(u)\cap \{2,3,5\}$ then $x_1...x_4uijvP^{-1}x_1\equiv C$, thus $u\rightarrow \{2,3,5\}$. Similarly if $1\in N^-(u)$ then $x_1...x_4u1453v26x_1\equiv C$, thus $u\rightarrow 1$. If there exists $u'\in N^+_{T'_2-u}(x_4)$, $4$ or $6$ is an origin of an antidirected ouptath $P$ in $T'_2-u'$ and thus $x_1...x_4u'vP^{-1}x_1\equiv C$, so $N^+_{T'_2}(x_4)=\{u\}$. If $d^-_{T'_2}(x_3)\geq 1$, the there exist $w\in N^-_{T'_2}(x_3)$, $w'\in N^+_{T'_2}(w)$ and $a\in \{4,6\}$ pairwaisly distinct such that $x_4$ has an inneighbor origin of an antidirected outpath $P$ in $T'_2\setminus \{a,w,w'\}$, thus $x_4ax_1...x_3ww'vP^{-1}x_4\equiv C$. Otherwise, since $T(x_1,x_3,v)\neq 3A$ then $d^+_{T'_2}(x_2)\geq 1$, thus there exist $x\in N^+_{T'_2}(x_2)$ and $P$ an antidirected outpath $P$ in $T'_2-x$ of origin distinct from $u$, therefore $x_4x_1x_3x_2xvP^{-1}x_4\equiv C$.\\
Exc $30$: One can check that for any $i\in \{3,4,5,\}$, there exists an outneighbor $j$ of $i$ such that $4$ or $6$ is an origin of an antidirected outpath $P$ in $T'_2\setminus \{i,j,u\}$, so if there exists $i\in N^-(u)\cap \{3,4,5\}$ then $x_1...x_4uijvP^{-1}x_1\equiv C$, thus $u\rightarrow \{3,4,5\}$. Similarly if there exists $i\in N^-(u)\cap \{1,2,6\}$, set $\{u_1\rightarrow u_2\}=T(1,2,6)-i$, we have $x_1...x_4i5u_1u_2v43x_1\equiv C$, thus $u\rightarrow \{1,2,6\}$. If there exists $u'\in N^+_{T'_2-u}(x_4)$, $4$ or $6$ is an origin of an antidirected ouptath $P$ in $T'_2-u'$ and thus $x_1...x_4u'vP^{-1}x_1\equiv C$, so $N^+_{T'_2}(x_4)=\{u\}$. If $d^-_{T'_2}9x_3)\geq 1$, then there exist $w\in N^-_{T'_2}(x_3)$, $w'\in N^+_{T'_2}(w)$ and $a\in \{4,6\}$ pairwaisly distinct such that $x_4$ has an inneighbor origin of an antidirected outpath $P$ in $T'_2\setminus \{a,w,w'\}$, thus $x_4ax_1...x_3ww'vP^{-1}x_4\equiv C$. Otherwise, since $T(x_1,x_3,v)\neq 3A$ then $d^+_{T'_2}(x_2)\geq 1$, thus there exist $x\in N^+_{T'_2}(x_2)$ and $P$ an antidirected outpath $P$ in $T'_2-x$ of origin distinct from $u$, therefore $x_4x_1x_3x_2xvP^{-1}x_4\equiv C$.\\
\textbf{Case $6$}: If $\delta=5$ and for all $x\in T$ such that $d^-(x)=5$ we have $T(N^-(x))=5A$. Let $a\in T_2$ such that $d^-_{T_2}(a)\geq 1$ and set $T'_1=T_1\cup \{a\}$ and $T'_2=T_2-a$. Let $x_1...x_6=-(1,1,...,1)$ be a path in $T'_1$ and let $u\in N^+_{T'_2}(x_6)$. If $x_1$ has an inneighbor origin of an antidirected outpath $P$ in $T'_2-u$ then $x_1...x_6uvP^{-1}x_1\equiv C$. Otherwise, if $d^-_{T'_2-u}(x_1)\geq 2$ then $(T'_2-u;(1,1,...,1))$ is one of the following exceptions:\\
Exc $1$: Set $V(T_1)=\{y_1,y_2,...,y_5\}$ such that $y_i\equiv i$ in $5A$, we distinguish two cases:\begin{itemize}
\item[(i)] If there exists $y\in T_2$ and $i\in \{1,2,3,4,5\}$ such that $\{y_i,y_{i+1}\}\rightarrow y$, without loss of generality suppose that $i=1$. Let $u'\in N^+_{T_2-y}(y_1)$ and suppose that $T_2\setminus \{u',y\}=4A$. We have $P=y_5y_3y_4y_2yy_1=-(1,1,...,1)$. Since $T(N^+(x))=5A$ for all $x\in T$ such that $d^+(x)=5$ then $d^+_{T_2-y}(y_1)\geq 2$. If $\{1,2\}\subseteq N^-(y_5)$. If $3\rightarrow u'$ then $Pu'34v21y_5\equiv C$ and if $4\rightarrow u'$ then $Pu'41v32y_5\equiv C$, so $u'\rightarrow \{3,4\}$. If $N^+(y_1)\cap \{3,4\}\neq \emptyset$ then $d^-(y_1)=5$ with $T(N^-(y_1))\neq 5A$, a contradiction, so $y_1\rightarrow 1$ or $2$, say $1$, then $P1v4u'32y_5\equiv C$. If $N^-_{T_2\setminus \{u',y\}}(y_5)=\{1,3\}$, by the building lemmas we can suppose that $\{2,3,4\}\rightarrow y_3$. If $3\rightarrow u'$ then $Pu'34v21y_5\equiv C$, if $4\rightarrow u'$ then $Pu'42v31y_5\equiv C$ and if $1\rightarrow u'$ then $Pu'12v43y_5\equiv C$, so $u'\rightarrow \{1,3,4\}$. If $y_1\rightarrow 2$ then $P2v4u'31y_5\equiv C$, if $y_1\rightarrow 3$ or $4$, say $3$, then $P3u'4v21y_5\equiv C$, so $N^-(y_1)=\{y_4,y_5,2,3,4\}$ and thus $\{2,4\}\rightarrow y_4$ and $y_4\rightarrow 3$. We have $21y_53y_34y_4y_2yy_5u'v2\equiv C$. Now if $N^-_{T_2\setminus \{u',y\}}(y_5)=\{2,3\}$, if $1\rightarrow u'$ then $Pu'12v43y_5\equiv C$ and if $4\rightarrow u'$ then $Pu'41v32y_5\equiv C$, so $u'\rightarrow \{1,4\}$. If $y_1\rightarrow 2$ then $P2v1u'43y_5\equiv C$, so $2\rightarrow y_1$.  If $3\rightarrow y_1$ then $T(N^-(y_1))\neq 5A$ so $y_1\rightarrow 3$ and thus $\{1,4\}\rightarrow y_1$. But $T(N^-(y_1))=5A$ so $y_4\rightarrow 4$ and $\{1,2\}\rightarrow y_4$. In any case $T(N^-(y_2))\neq 5A$ and $T(N^+(y_2))\neq 5A$, a contradiction.
\item[(ii)] If for all $x\in T_2$ and for all $i$, $N^+(x)\cap \{y_iy_{i+1}\}\neq \emptyset$. Let $y\in T_2$ such that $d^+_{T_2}(y)=\Delta^+(T_2)$, without loss of generality suppose that $y\rightarrow \{y_1,y_5\}$. We have $P=y_1yy_5y_3y_4y_2=-(1,1,...,1)$. Let $u'\in N^+_{T_2-y}(y_2)$, we may assume that $T_2\setminus \{y,u'\}=4A$. If $\{1,2\}\subseteq N^-(y_1)$, by the building lemmas we can suppose that $\{2,4\}\rightarrow y$, and since $y_1\rightarrow 4$ then $4\rightarrow y_5$ and thus $31y_12y4y_5y_3y_4y_2u'v3\equiv C$. If $N^-_{T_2\setminus \{y,u'\}}(y_1)=\{1,3\}$, by the building lemmas, we may assume that $\{2,3,4\}\rightarrow y$, a contradiction. So we can suppose that $N^-_{T_2\setminus \{y,u'\}}(y_1)=\{2,3\}$. If $1\rightarrow u'$ then $Pu'12v43y_1\equiv C$ and if $4\rightarrow u'$ then $Pu'41v32y_1\equiv C$, so $u'\rightarrow \{1,4\}$. Suppose that $y\rightarrow y_4$ then $P'=y_5yy_4y_2y_3y_1=-(1,1,...,1)$. We have $\{1,4\}\rightarrow y_5$, if $2\rightarrow u'$ then $P'4vu'231y_5\equiv C$, so $u'\rightarrow 2$ and thus $3\rightarrow u'$ (since $\Delta^(T_2)=3$), we have $P'43u'v21y_5\equiv C$. So we may assume that $y_4\rightarrow y$ and so $y\rightarrow y_3$, thus $y_2\rightarrow y$. So $y_5y_4y_1yy_3y_2=-(1,1,...,1)$, but $4\rightarrow y_5$, so $y_5y_4y_1yy_3y_2u'v3124y_5\equiv C$.
\end{itemize}
Exc $20$: One can check that for any $i\in \{3,4,5,6\}$, $1$ or $2$ is an origin of an antidirected outpath $P$ in $T'_2\setminus \{i,i+1,u\}$, so if there exists $i\in N^-(u)\cap \{3,4,5,6\}$ then $x_1...x_6ui(i+1)vP^{-1}x_1\equiv C$, thus $u\rightarrow \{3,4,5,6\}$. If there exists $u'\in N^+_{T'_2-u}(x_6)$, $1$ or $2$ is an origin of an antidirected ouptath $P$ in $T'_2-u'$ and thus $x_1...x_6u'vP^{-1}x_1\equiv C$, so $N^+_{T'_2}(x_6)=\{u\}$. There exist $w\in N^-_{T'_2}(x_5)$, $w'\in N^+_{T'_2}(w)$ and $a\in \{1,2\}$ pairwaisly distinct such that $x_6$ has an inneighbor origin of an antidirected outpath $P$ in $T'_2\setminus \{a,w,w'\}$, thus $x_6ax_1...x_5ww'vP^{-1}x_6\equiv C$.\\
Exc $22$: Without loss of generality we can suppose that $\{1,2\}\subseteq N^-(x_1)$. One can check that for any $i\in \{4,5,6\}$, there exists an outneighbor $j$ of $i$ such that $1$ or $2$ is an origin of an antidirected outpath $P$ in $T'_2\setminus \{i,j,u\}$, so if there exists $i\in N^-(u)\cap \{1,2,3,4,5\}$ then $x_1...x_6uijvP^{-1}x_1\equiv C$, thus $u\rightarrow \{4,5,6\}$. If there exists $u'\in N^+_{T'_2-u}(x_6)$ distinct from $4$, $1$ or $2$ is an origin of an antidirected ouptath $P$ in $T'_2-u'$ and thus $x_1...x_6u'vP^{-1}x_1\equiv C$, and if $4\in N^+(x_6)$ then $x_1...x_64536v21x_1\equiv C$, so $N^+_{T'_2}(x_6)=\{u\}$. There exist $w\in N^-_{T'_2}(x_5)$, $w'\in N^+_{T'_2}(w)$ and $a\in \{4,6\}$ pairwaisly distinct such that $x_6$ has an inneighbor origin of an antidirected outpath $P$ in $T'_2\setminus \{a,w,w'\}$, thus $x_6ax_1...x_5ww'vP^{-1}x_6\equiv C$.\\
Exc $23$: One can check that for any $i\in \{1,2,3,4,5\}$, there exists an outneighbor $j$ of $i$ such that $4$ or $6$ is an origin of an antidirected outpath $P$ in $T'_2\setminus \{i,j,u\}$, so if there exists $i\in N^-(u)\cap \{1,2,3,4,5\}$ then $x_1...x_6uijvP^{-1}x_1\equiv C$, thus $u\rightarrow \{1,2,3,4,5\}$. If there exists $u'\in N^+_{T'_2-u}(x_6)$, $4$ or $6$ is an origin of an antidirected ouptath $P$ in $T'_2-u'$ and thus $x_1...x_6u'vP^{-1}x_1\equiv C$, so $N^+_{T'_2}(x_6)=\{u\}$. There exist $w\in N^-_{T_2}(x_5)$, $w'\in N^+_{T'_2}(w)$ and $a\in \{4,6\}$ pairwaisly distinct such that $x_6$ has an inneighbor origin of an antidirected outpath $P$ in $T'_2\setminus \{a,w,w'\}$, thus $x_6ax_1...x_5ww'vP^{-1}x_6\equiv C$.\\
Exc $24$: Without loss of generality we may suppose that $1\in N^-(x_1)$. One can check that for any $i\in \{2,4,5,6\}$, there exists an outneighbor $j$ of $i$ such that $1$ is an origin of an antidirected outpath $P$ in $T_2\setminus \{i,j,u\}$, so if there exists $i\in N^-(u)\cap \{2,4,5,6\}$ then $x_1...x_6uijvP^{-1}x_1\equiv C$, thus $u\rightarrow \{2,4,5,6\}$. If $3\rightarrow u$ then $x_1...x_6u3456v21x_1\equiv C$, so $u\rightarrow 3$. If there exists $u'\in N^+_{T'_2-u}(x_6)$, $1$ is an origin of an antidirected ouptath $P$ in $T'_2-u'$ and thus $x_1...x_6u'vP^{-1}x_1\equiv C$, so $N^+_{T'_2}(x_6)=\{u\}$. There exist $w\in N^-_{T'_2}(x_5)$, $w'\in N^+_{T'_2}(w)$ and $a\in N^-_{T'_2}(x_1)\cap N^-_{T'_2}(x_6)$ pairwaisly distinct such that $x_6$ has an inneighbor origin of an antidirected outpath $P$ in $T'_2\setminus \{a,w,w'\}$, thus $x_6ax_1...x_5ww'vP^{-1}x_6\equiv C$.\\
Exc $29$ :One can check that for any $i\in \{2,3,5\}$, there exists an outneighbor $j$ of $i$ such that $4$ or $6$ is an origin of an antidirected outpath $P$ in $T'_2\setminus \{i,j,u\}$, so if there exists $i\in N^-(u)\cap \{2,3,5\}$ then $x_1...x_6uijvP^{-1}x_1\equiv C$, thus $u\rightarrow \{2,3,5\}$. Similarly if $1\in N^-(u)$ then $x_1...x_6u1453v26x_1\equiv C$, thus $u\rightarrow 1$. If there exists $u'\in N^+_{T'_2-u}(x_6)$, $4$ or $6$ is an origin of an antidirected ouptath $P$ in $T'_2-u'$ and thus $x_1...x_6u'vP^{-1}x_1\equiv C$, so $N^+_{T'_2}(x_6)=\{u\}$. There exist $w\in N^-_{T'_2}(x_5)$, $w'\in N^+_{T'_2}(w)$ and $a\in \{4,6\}$ pairwaisly distinct such that $x_6$ has an inneighbor origin of an antidirected outpath $P$ in $T'_2\setminus \{a,w,w'\}$, thus $x_6ax_1...x_5ww'vP^{-1}x_6\equiv C$.\\
Exc $30$: One can check that for any $i\in \{3,4,5,\}$, there exists an outneighbor $j$ of $i$ such that $4$ or $6$ is an origin of an antidirected outpath $P$ in $T'_2\setminus \{i,j,u\}$, so if there exists $i\in N^-(u)\cap \{3,4,5\}$ then $x_1...x_6uijvP^{-1}x_1\equiv C$, thus $u\rightarrow \{3,4,5\}$. Similarly if there exists $i\in N^-(u)\cap \{1,2,6\}$, set $\{u_1\rightarrow u_2\}=T(1,2,6)-i$, we have $x_1...x_6i5u_1u_2v43x_1\equiv C$, thus $u\rightarrow \{1,2,6\}$. If there exists $u'\in N^+_{T'_2-u}(x_6)$, $4$ or $6$ is an origin of an antidirected ouptath $P$ in $T'_2-u'$ and thus $x_1...x_6u'vP^{-1}x_1\equiv C$, so $N^+_{T'_2}(x_6)=\{u\}$. There exist $w\in N^-_{T'_2}(x_5)$, $w'\in N^+_{T'_2}(w)$ and $a\in \{4,6\}$ pairwaisly distinct such that $x_6$ has an inneighbor origin of an antidirected outpath $P$ in $T'_2\setminus \{a,w,w'\}$, thus $x_6ax_1...x_5ww'vP^{-1}x_6\equiv C$.\\
So we may assume that $d^+_{T'_2}(x_1)\geq 3$ and by Theorem \ref{2.1}, $x_1$ is always an end of an antidirected inpath $y_1...y_6$ in $T'_1$ and thus we may find $u'\in N^+_{T'_2}(y_6)$ such that $d^-_{T'_2-u'}(y_1)\geq 2$, therefore $T$ contains $C$.\\
\textbf{Case $7$}: If $\delta=7$ and for all $x\in T$ such that $d^-(x)=7$ we have $T(N^-(x))=7A$, let $a\in T_2$ such that $d^+_{T_2}(a)=\Delta^+(T_2)$ and set $T'_1=T_1\cup \{a\}$ and $T'_2=T_2-a$, there exists $x_1...x_8=-(1,1,...,1)$ such that $d^-_{T'_2}(x_1)\geq 3$ and $d^+_{T'_2}(x_8)\geq 2$. Let $u\in N^+_{T'_2}(x_8)$, $T'_2-u$ has an antidirected Hamiltonian outpath $P$ of origin an inneighbor of $x_1$ and thus $x_1...x_8uvP^{-1}x_1\equiv C$ unless $(T'_2-u;(1,1,...,1))$ is one of the following biexceptions with respect to $x_1$:\\
Exc $(22,1)$: If $N^-(u)\cap \{4,5,6\}\neq \emptyset$, one can find $C$. Otherwise, let $u'\in N^+_{T'_2}(x_8)$, $T'_2-u'$ is not an exception and the result follows.\\
Exc $(24,1)$: If $d^-_{T'_2}(u)\geq 1$ one can find $C$, otherwise, let $u'\in N^+_{T'_2}(x_8)$ distinct from $u$, $T'_2-u'$ is not an exception, so the result follows.    $\square$
\end{proof}
\\The following corollary discuss the existence of any non-directed cycle in tournaments.
\begin{corollary}
Let $T$ be a tournament on $n\geq 3$ vertices and $C$ be a non-directed cycle of order $m$, where $3\leq m \leq n$. $T$ contains $C$ if and only if $(T;C)$ is not one of the exceptions $A_i$, $1\leq i\leq 18$.
\end{corollary}
\begin{proof}
If $m=n$, theorems \ref{3.2} and \ref{3.7} and lemmas \ref{3.3}, \ref{3.4}, \ref{3.5} and \ref{3.6} give the result. Suppose now that $m=n-1$. If there exists $x\in T$ such that $(T-x;C)$ is not an exception, the result follows. Otherwise, we discuss according to the values of $n$:\begin{itemize}
\item[(i)] If $n=5$, if $T=5A$ then $(T;C)=A_{15}$. Otherwise, let $x\in T$ such that $d^+(x)=3$ and let $y$ be its inneighbor. it is easy to check that $d^-(y)\geq 2$ with $T\setminus \{x,y\}=3A$ and so $(T;C)=A_{16}$.
\item[(ii)] If $n=7$, if $T-x=6F$ for all $x\in T$ then $(T;C)=A_{18}$. If $T-x=6M$ for all $x\in T$ then $(T;C)=A_{17}$. It is easy to check that if there exists $x\in T$ such that $T-x=6M$ then $T-y\notin \{6O,6P\}$ for all $y\in T$. Also, one can check that $|\{x\in T; T-x=6O\}|\leq 2$ and $|\{x\in T; T-x=6P\}|\leq 3$, so if $T\neq 7K$ then there exists $x\in T$ such that $T-x\notin \{6O,6P\}$ and the result follows.
\end{itemize}
Finally, it is easy to check that there exists $x\in T$ such that $T-x\notin \{5A,5F,7A,7K\}$, so for all $m\leq n-2$, $T$ contains $C$.    $\square$
\end{proof}
\renewcommand{\abstractname}{Acknowledgments}
\begin{abstract}
I would like to express my gratitude to Professor Amine El Sahili and Doctor Maydoun Mortada for their insightful comments.
\end{abstract}

\end{document}